\documentclass{ociamthesis}

\usepackage{mathpazo} 
\usepackage[scaled]{helvet} 
\usepackage{courier} 
\normalfont
\usepackage[T1]{fontenc}

\usepackage{amsmath}
\usepackage{amssymb}
\usepackage{stmaryrd}
\usepackage{bbold}
\usepackage{amsthm}
\usepackage{xspace}
\usepackage{hyperref}
\usepackage{graphicx}

\newtheorem{theorem}{Theorem}[section]
\newtheorem*{theorem*}{Theorem}
\newtheorem{proposition}[theorem]{Proposition}
\newtheorem{lemma}[theorem]{Lemma}
\newtheorem{corollary}[theorem]{Corollary}
\theoremstyle{definition}\newtheorem{example}[theorem]{Example}
\theoremstyle{definition}\newtheorem{examples}[theorem]{Examples}
\theoremstyle{definition}\newtheorem{definition}[theorem]{Definition}
\theoremstyle{definition}\newtheorem{definitions}[theorem]{Definitions}
\theoremstyle{definition}\newtheorem{remark}[theorem]{Remark}
\theoremstyle{definition}\newtheorem{conjecture}[theorem]{Conjecture}
\theoremstyle{definition}
\theoremstyle{definition}
\theoremstyle{definition}
\theoremstyle{definition}\newtheorem{notation}[theorem]{Notation}

\newcommand{\Tr}{\ensuremath{\textrm{Tr}\xspace}}
\newcommand{\Int}{\ensuremath{\textrm{\rm Int}\xspace}}
\newcommand{\op}{{\textrm{op}}}

\newcommand{\catSet}{\ensuremath{\textrm{\bf Set}\xspace}}
\newcommand{\catRel}{\ensuremath{\textrm{\bf Rel}\xspace}}

\newcommand{\catVect}{\ensuremath{\textrm{\bf Vect}\xspace}}
\newcommand{\catFVect}{\ensuremath{\textrm{\bf FVect}\xspace}}
\newcommand{\catFHilb}{\ensuremath{\textrm{\bf FHilb}\xspace}}
\newcommand{\catHilb}{\ensuremath{\textrm{\bf Hilb}\xspace}}
\newcommand{\catSuperHilb}{\ensuremath{\textrm{\bf SuperHilb}\xspace}}
\newcommand{\catAb}{\ensuremath{\textrm{\bf Ab}\xspace}}

\newcommand{\catCHaus}{\ensuremath{\textrm{\bf CHaus}\xspace}}
\newcommand{\catHaus}{\ensuremath{\textrm{\bf Haus}\xspace}}
\newcommand{\catGraph}{\ensuremath{\textrm{\bf Graph}}\xspace}
\newcommand{\catMat}{\ensuremath{\textrm{\bf Mat}}\xspace}
\newcommand{\catGr}{\ensuremath{\textrm{\bf Gr}}\xspace}
\newcommand{\catSpek}{\ensuremath{\textrm{\bf Spek}}\xspace}

\newcommand{\catMonCat}{\ensuremath{\textrm{\bf MonCat}}\xspace}
\newcommand{\catSymMonCat}{\ensuremath{\textrm{\bf SymMonCat}}\xspace}
\newcommand{\catSymTraceCat}{\ensuremath{\textrm{\bf SymTraceCat}}\xspace}
\newcommand{\catCCCat}{\ensuremath{\textrm{\bf CCCat}}\xspace}

\newcommand{\catMonSig}{\ensuremath{\textrm{\bf MonSig}}\xspace}
\newcommand{\catSGraph}{\ensuremath{\textrm{\bf SGraph}}\xspace}

\newcommand{\catSGraphTG}{\ensuremath{\catSGraph_{T}}}

\newcommand{\dom}{\ensuremath{\textrm{dom}}}
\newcommand{\cod}{\ensuremath{\textrm{cod}}}
\newcommand{\In}{\textrm{\rm In}}
\newcommand{\Out}{\textrm{\rm Out}}
\newcommand{\Bound}{\textrm{\rm Bound}}
\newcommand{\FCsp}{\textrm{\rm FCsp}}
\newcommand{\Csp}{\textrm{\rm Csp}}
\newcommand{\csp}{\textrm{\rm csp}}

\newcommand{\sgn}{\textrm{sgn}}

\newcommand{\tikznode}[2]{\tikz{\node[shape=rectangle,inner sep=0,#1]{#2};}}

\newcommand{\normalised}{\ensuremath{\!\downarrow\,}}
\newcommand{\cmdrewritesto}{\tikz[baseline=-0.25em] { \draw [-open triangle 45, line width=0.2pt] (0,0) -- (0.5,0); }\,}
\newcommand{\cmdrewriteequiv}{\tikz[baseline=-0.25em] { \draw [open triangle 45-open triangle 45, line width=0.2pt] (0,0) -- node [auto,yshift=-1.2mm] {$*$} (0.7,0); }\,}
\newcommand{\cmdrewritetrans}{\tikz[baseline=-0.25em] { \draw [-open triangle 45, line width=0.2pt] (0,0) -- node [auto,pos=0.3,yshift=-1.2mm] {$*$} (0.5,0); }\,}

\newcommand{\cmdleftrewritetrans}{\tikz[baseline=-0.25em] { \draw [open triangle 45-, line width=0.2pt] (0,0) -- node [auto,pos=0.7,yshift=-1.2mm] {$*$} (0.5,0); }\,}

\DeclareMathOperator{\rewritesto}{\cmdrewritesto}
\DeclareMathOperator{\rewriteequiv}{\cmdrewriteequiv}
\DeclareMathOperator{\rewritetrans}{\cmdrewritetrans}

\DeclareMathOperator{\leftrewritetrans}{\cmdleftrewritetrans}

\newcommand{\bra}[1]{\ensuremath{\left\langle #1 \right|}}
\newcommand{\ket}[1]{\ensuremath{\left|  #1 \right\rangle}}
\newcommand{\braket}[2]{\ensuremath{\langle#1|#2\rangle}}
\newcommand{\ketbra}[2]{\ensuremath{\ket{#1}\!\bra{#2}}}

\newcommand{\ketGHZ} {\ket{\textit{GHZ}\,}}
\newcommand{\ketW}   {\ket{\textit{W\,}}}
\newcommand{\ketBell}{\ket{\textit{Bell\,}}}
\newcommand{\ketEPR} {\ket{\textit{EPR\,}}}
\newcommand{\ketGHZD}{\ket{\textrm{GHZ}^{(D)}}}
\newcommand{\ketWD}  {\ket{\textrm{W}^{(D)}}}

\newcommand{\braBell}{\bra{\textit{Bell\,}}}
\newcommand{\braEPR} {\bra{\textit{EPR\,}}}

\newcommand{\numket}[1]{\ket{\underbar{$#1$}}}

\newcommand{\updot}{\ensuremath{\blacktriangle}}
\newcommand{\downdot}{\ensuremath{\blacktriangledown}}
\newcommand{\cpline}{\ensuremath{\textrm{\bf CP}^1}}

\usepackage{tikz}
\usepackage{pgfplots}
\usetikzlibrary{trees}
\usetikzlibrary{topaths}
\usetikzlibrary{decorations.pathmorphing}
\usetikzlibrary{decorations.markings}
\usetikzlibrary{matrix,backgrounds,folding}
\usetikzlibrary{chains,scopes,positioning,fit}
\usetikzlibrary{arrows,shadows}
\usetikzlibrary{calc} 
\usetikzlibrary{chains}
\usetikzlibrary{shapes,shapes.geometric,shapes.misc}

\usepackage{tikzfig}
\tikzstyle{every picture}=[baseline=-0.25em]
\tikzstyle{dotpic}=[scale=0.6]
\tikzstyle{diredges}=[every to/.style={diredge}]
\tikzstyle{dot graph}=[shorten <=-0.1mm,shorten >=-0.1mm,scale=0.6]
\tikzstyle{digraph}=[-latex]
\tikzstyle{plot point}=[circle,fill=black,minimum width=2mm,inner sep=0]
\tikzstyle{string graph}=[scale=0.6]
\tikzstyle{sg diredge}=[-stealth]
\tikzstyle{rewrite edge}=[-open triangle 45]
\tikzstyle{sg bold diredge}=[-stealth,thick,shorten >=-1pt]
\tikzstyle{sg vertex}=[circle,minimum width=2.2mm,fill=white,draw=black,inner sep=0mm]
\tikzstyle{labelled sg vertex}=[circle,minimum width=7mm,fill=white,draw=black,inner sep=0mm]
\tikzstyle{sg grey vertex}=[sg vertex,fill=gray!30!white]
\tikzstyle{sg black vertex}=[sg vertex,fill=black]
\tikzstyle{sg bold vertex}=[circle,minimum width=2.2mm,fill=white,draw=black,very thick,inner sep=0mm]
\tikzstyle{sg wire vertex}=[circle,minimum width=1mm,fill=black,inner sep=0mm]

\tikzstyle{control qubit}=[circle,minimum width=1mm,fill=black,inner sep=0mm,yshift=-0.6mm]

\newcommand{\phantombox}[1]{\tikz[baseline=(current bounding box).east]{\path [use as bounding box] (0,0) rectangle #1;}}
\tikzstyle{braceedge}=[decorate,decoration={brace,amplitude=2mm,raise=-1mm}]
\tikzstyle{small braceedge}=[decorate,decoration={brace,amplitude=1mm,raise=-1mm}]
\tikzstyle{left hook arrow}=[left hook-latex]
\tikzstyle{right hook arrow}=[right hook-latex]

\newcommand{\measurement}{\tikz[scale=0.6]{ \draw [use as bounding box,draw=none] (0,-0.1) rectangle (1,0.7); \draw [fill=white] (1,0) arc (0:180:5mm); \draw (0,0) -- (1,0) (0.5,0) -- +(60:7mm);}}

\tikzstyle{dot}=[inner sep=0.7mm,minimum width=0pt,minimum height=0pt,fill=black,draw=black,shape=circle]
\tikzstyle{white dot}=[dot,fill=white]
\tikzstyle{alt white dot}=[white dot,label={[xshift=2.9mm,yshift=-0.1mm]left:$\cdot$}]
\tikzstyle{gray dot}=[dot,fill=gray!50]
\tikzstyle{box vertex}=[draw=black,rectangle]
\tikzstyle{whitebg}=[fill=white,inner sep=2pt]
\tikzstyle{graph state vertex}=[sg vertex,fill=black]

\tikzstyle{wide point}=[fill=white,draw=black,shape=isosceles triangle,shape border rotate=90,isosceles triangle stretches=true,inner sep=1pt,minimum width=1.5cm,minimum height=5mm]
\tikzstyle{wide copoint}=[fill=white,draw=black,shape=isosceles triangle,shape border rotate=-90,isosceles triangle stretches=true,inner sep=1pt,minimum width=1.5cm,minimum height=5mm,yshift=-0.7mm]
\tikzstyle{symm}=[ultra thick,shorten <=-1mm,shorten >=-1mm]

\tikzstyle{square box}=[rectangle,fill=white,draw=black,minimum height=6mm,minimum width=6mm,yshift=-0.7mm]
\tikzstyle{square gray box}=[rectangle,fill=gray!30,draw=black,minimum height=6mm,minimum width=6mm]
\tikzstyle{point}=[regular polygon,regular polygon sides=3,draw=black,scale=0.75,inner sep=-0.5pt,minimum width=7mm,fill=white]
\tikzstyle{copoint}=[point,regular polygon rotate=180,fill=white]
\tikzstyle{gray point}=[point,fill=gray!40!white]
\tikzstyle{gray copoint}=[copoint,fill=gray!40!white]

\tikzstyle{open graph}=[baseline=-0.25em]
\tikzstyle{greybg}=[background rectangle/.style={fill=black!5,draw=black!30,rounded corners=1ex}, show background rectangle]

\tikzstyle{edge point}=[circle,minimum width=1mm,fill=black,inner sep=0mm]
\tikzstyle{vertex point}=[circle,minimum width=2.2mm,fill=white,draw=black,inner sep=0mm]
\tikzstyle{gray vertex point}=[circle,minimum width=2.2mm,fill=gray!30!white,draw=black,inner sep=0mm]
\tikzstyle{edge label}=[inner sep=2pt, font=\small]
\tikzstyle{on edge label}=[fill=white, font=\footnotesize, inner sep=1 pt]

\newcommand{\edgearrow}{{\arrow[black]{>}}}
\newcommand{\edgetick}{{\arrow[black,scale=0.7,very thick]{|}}}
\tikzstyle{short diredge}=[->]

\tikzstyle{diredge}=[short diredge]
\tikzstyle{medium diredge}=[short diredge]

\tikzstyle{halfedge}=[-)]
\tikzstyle{other halfedge}=[(-]
\tikzstyle{freeedge}=[(-)]
\tikzstyle{white edge}=[line width=5pt,white]
\tikzstyle{tick}=[postaction=decorate,decoration={markings, mark=at position 0.5 with \edgetick}]
\tikzstyle{small map edge}=[|-latex, gray!60!blue, shorten <=0.9mm, shorten >=0.5mm]
\tikzstyle{thick dashed edge}=[very thick,dashed,gray!40]
\tikzstyle{map edge}=[|-latex,very thick, gray!40, shorten <=1mm, shorten >=0.5mm]
\tikzstyle{tickedge}=[postaction=decorate,
  decoration={markings, mark=at position 0.5 with \edgetick}]
  
\tikzstyle{dirtickedge}=[->,postaction=decorate,
  decoration={markings, mark=at position 0.5 with \edgetick}]
\tikzstyle{dirdoubletickedge}=[->,postaction=decorate,
  decoration={markings, mark=at position 0.4 with \edgetick},
  decoration={markings, mark=at position 0.6 with \edgetick}]

\tikzstyle{arrs}=[-latex,font=\small,auto]
\tikzstyle{arrow plain}=[arrs]
\tikzstyle{arrow dashed}=[dashed,arrs]
\tikzstyle{arrow bold}=[very thick,arrs]
\tikzstyle{arrow hide}=[draw=white!0,-]
\tikzstyle{arrow reverse}=[latex-]
\tikzstyle{cdnode}=[]

\tikzstyle{cnot}=[fill=white,shape=circle,inner sep=-1.4pt]
\tikzstyle{bang box}=[draw=black,dashed,minimum height=12mm,minimum width=12mm,fill=gray!20]
\tikzstyle{wire label}=[font=\footnotesize, auto]

\newcommand{\emptybox}[2]{\tikz{\node[rectangle,inner sep=0,minimum width=#1,minimum height=#2] {};}}

\newcommand{\dotunit}[1]{%
\,\begin{tikzpicture}[dotpic,yshift=-1mm]
\node [#1] (a) at (0,0.25) {}; 
\draw [medium diredge] (a)--(0,-0.2);
\end{tikzpicture}\,}
\newcommand{\dotcounit}[1]{%
\,\begin{tikzpicture}[dotpic,yshift=1.5mm]
\node [#1] (a) at (0,-0.25) {}; 
\draw [medium diredge] (0,0.2)--(a);
\end{tikzpicture}\,}
\newcommand{\dotmult}[1]{%
\,\begin{tikzpicture}[dotpic,yshift=0.5mm]
	\node [#1] (a) {};
	\draw [medium diredge] (a) -- (-90:0.55);
	\draw [medium diredge] (45:0.6) -- (a);
	\draw [medium diredge] (135:0.6) -- (a);
\end{tikzpicture}\,}

\newcommand{\dotcomult}[1]{%
\,\begin{tikzpicture}[dotpic]
	\node [#1] (a) {};
	\draw [medium diredge] (90:0.55) -- (a);
	\draw [medium diredge] (a) -- (-45:0.6);
	\draw [medium diredge] (a) -- (-135:0.6);
\end{tikzpicture}\,}
\newcommand{\dottickunit}[1]{%
\begin{tikzpicture}[dotpic,yshift=-1mm]
\node [#1] (a) at (0,0.35) {}; 
\draw [postaction=decorate,
       decoration={markings, mark=at position 0.3 with \edgetick},
       decoration={markings, mark=at position 0.85 with \edgearrow}] (a)--(0,-0.25);
\end{tikzpicture}}
\newcommand{\dottickcounit}[1]{%
\begin{tikzpicture}[dotpic,yshift=1mm]
\node [#1] (a) at (0,-0.35) {}; 
\draw [postaction=decorate,
       decoration={markings, mark=at position 0.8 with \edgetick},
       decoration={markings, mark=at position 0.45 with \edgearrow}] (0,0.25) -- (a);
\end{tikzpicture}}
\newcommand{\dotonly}[1]{%
\,\begin{tikzpicture}[dotpic,yshift=-0.3mm]
\node [#1] (a) at (0,0) {};
\end{tikzpicture}\,}

\newcommand{\dotcap}[1]{%
\,\begin{tikzpicture}[dotpic,yshift=2.5mm]
	\node [#1] (a) at (0,0) {};
	\draw [bend right,medium diredge] (a) to (-0.4,-0.6);
	\draw [bend left,medium diredge] (a) to (0.4,-0.6);
\end{tikzpicture}\,}

\newcommand{\tick}{%
\,\,\begin{tikzpicture}[dotpic]
	\node [style=none] (a) at (0,0.35) {};
	\node [style=none] (b) at (0,-0.35) {};
	\draw [dirtickedge] (a) -- (b);
\end{tikzpicture}\,\,}

\newcommand{\lolli}{%
\,\begin{tikzpicture}[dotpic,yshift=-1mm]
	\path [use as bounding box] (-0.25,-0.25) rectangle (0.25,0.5);
	\node [style=dot] (a) at (0, 0.15) {};
	\node [style=none] (b) at (0, -0.25) {};
	\draw [medium diredge] (a) to (b.center);
	\draw [diredge, out=45, looseness=1.00, in=135, loop] (a) to ();
\end{tikzpicture}\,}

\newcommand{\cololli}{%
\,\begin{tikzpicture}[dotpic]
	\path [use as bounding box] (-0.25,-0.5) rectangle (0.25,0.5);
	\node [style=none] (a) at (0, 0.5) {};
	\node [style=dot] (b) at (0, 0) {};
	\draw [diredge, in=-45, looseness=2.00, out=-135, loop] (b) to ();
	\draw [medium diredge] (a.center) to (b);
\end{tikzpicture}\,}

\newcommand{\blackdot}{\dotonly{dot}}
\newcommand{\unit}{\dotunit{dot}}
\newcommand{\counit}{\dotcounit{dot}}
\newcommand{\mult}{\dotmult{dot}}
\newcommand{\comult}{\dotcomult{dot}}
\newcommand{\tickunit}{\dottickunit{dot}}
\newcommand{\tickcounit}{\dottickcounit{dot}}

\newcommand{\whitedot}{\dotonly{white dot}}
\newcommand{\whiteunit}{\dotunit{white dot}}
\newcommand{\whitecounit}{\dotcounit{white dot}}
\newcommand{\whitemult}{\dotmult{white dot}}
\newcommand{\whitecomult}{\dotcomult{white dot}}

\newcommand{\whitecap}{\dotcap{white dot}}

\newcommand{\altwhiteunit}{\dotunit{alt white dot}}
\newcommand{\altwhitecounit}{\dotcounit{alt white dot}}
\newcommand{\altwhitemult}{\dotmult{alt white dot}}
\newcommand{\altwhitecomult}{\dotcomult{alt white dot}}

\newcommand{\graydot}{\dotonly{gray dot}}
\newcommand{\grayunit}{\dotunit{gray dot}}
\newcommand{\graycounit}{\dotcounit{gray dot}}
\newcommand{\graymult}{\dotmult{gray dot}}
\newcommand{\graycomult}{\dotcomult{gray dot}}

\newcommand{\blacktranspose}{\ensuremath{{\,\blackdot\!\textrm{\rm\,T}}}}
\newcommand{\whitetranspose}{\ensuremath{{\,\whitedot\!\textrm{\rm\,T}}}}
\newcommand{\graytranspose}{\ensuremath{{\,\graydot\!\textrm{\rm\,T}}}}

\newcommand{\circl}{\begin{tikzpicture}[dotpic]
		\node [style=none] (a) at (-0.25, 0.25) {};
		\node [style=none] (b) at (0.25, 0.25) {};
		\node [style=none] (c) at (-0.25, -0.25) {};
		\node [style=none] (d) at (0.25, -0.25) {};
		\draw [in=45, out=135] (b.center) to (a.center);
		\draw [in=135, out=225] (a.center) to (c.center);
		\draw [in=225, out=-45] (c.center) to (d.center);
		\draw [style=diredge, in=-45, out=45] (d.center) to (b.center);
\end{tikzpicture}}

\tikzstyle{cdiag}=[matrix of math nodes, row sep=3em, column sep=3em, text height=1.5ex, text depth=0.25ex,inner sep=0.5em]
\tikzstyle{arrow above}=[transform canvas={yshift=0.5ex}]
\tikzstyle{arrow below}=[transform canvas={yshift=-0.5ex}]

\newcommand{\csquare}[8]{
\begin{tikzpicture}
    \matrix(m)[cdiag,ampersand replacement=\&]{
    #1 \& #2 \\
    #3 \& #4  \\};
    \path [arrs] (m-1-1) edge node {$#5$} (m-1-2)
                 (m-2-1) edge node {$#6$} (m-2-2)
                 (m-1-1) edge node [swap] {$#7$} (m-2-1)
                 (m-1-2) edge node {$#8$} (m-2-2);
\end{tikzpicture}
}

\newcommand{\NWbracket}[1]{%
\draw #1 +(-0.25,0.5) -- +(-0.5,0.5) -- +(-0.5,0.25);}
\newcommand{\NEbracket}[1]{%
\draw #1 +(0.25,0.5) -- +(0.5,0.5) -- +(0.5,0.25);}
\newcommand{\SWbracket}[1]{%
\draw #1 +(-0.25,-0.5) -- +(-0.5,-0.5) -- +(-0.5,-0.25);}
\newcommand{\SEbracket}[1]{%
\draw #1 +(0.25,-0.5) -- +(0.5,-0.5) -- +(0.5,-0.25);}
\newcommand{\THETAbracket}[2]{%
\draw [rotate=#1] #2 +(0.25,0.5) -- +(0.5,0.5) -- +(0.5,0.25);}

\newcommand{\posquare}[8]{
\begin{tikzpicture}
    \matrix(m)[cdiag,ampersand replacement=\&]{
    #1 \& #2 \\
    #3 \& #4  \\};
    \path [arrs] (m-1-1) edge node {$#5$} (m-1-2)
                 (m-2-1) edge node [swap] {$#6$} (m-2-2)
                 (m-1-1) edge node [swap] {$#7$} (m-2-1)
                 (m-1-2) edge node {$#8$} (m-2-2);
	\NWbracket{(m-2-2)}
\end{tikzpicture}
}

\newcommand{\pbsquare}[8]{
\begin{tikzpicture}
    \matrix(m)[cdiag,ampersand replacement=\&]{
    #1 \& #2 \\
    #3 \& #4  \\};
    \path [arrs] (m-1-1) edge node {$#5$} (m-1-2)
                 (m-2-1) edge node [swap] {$#6$} (m-2-2)
                 (m-1-1) edge node [swap] {$#7$} (m-2-1)
                 (m-1-2) edge node {$#8$} (m-2-2);
  \SEbracket{(m-1-1)}
\end{tikzpicture}
}

\newcommand{\ctri}[6]{
	\begin{tikzpicture}[-latex]
		\matrix (m) [cdiag,ampersand replacement=\&] { #1 \& #2 \\ #3 \& \\ };
		\path [arrs] (m-1-1) edge node {$#4$} (m-1-2)
		      (m-1-1) edge node [swap] {$#5$} (m-2-1)
		      (m-2-1) edge node [swap] {$#6$} (m-1-2);
	\end{tikzpicture}
}

\newcommand{\crun}[5]{
\ensuremath{#1 \overset{#2}{\longrightarrow} #3 \overset{#4}{\longrightarrow} #5}
}

\newcommand{\cspan}[5]{
\ensuremath{#1 \overset{#2}{\longleftarrow} #3
               \overset{#4}{\longrightarrow} #5}
}

\newcommand{\ccospan}[5]{
\ensuremath{#1 \overset{#2}{\longrightarrow} #3
               \overset{#4}{\longleftarrow} #5}
}

\newcommand{\cpair}[4]{
\begin{tikzpicture}
    \matrix(m)[cdiag,ampersand replacement=\&]{
    #1 \& #2 \\};
    \path [arrs] (m-1-1.20) edge node {$#3$} (m-1-2.160)
                 (m-1-1.-20) edge node [swap] {$#4$} (m-1-2.-160);
\end{tikzpicture}
}

\newcommand{\csquareslant}[9]{
\begin{tikzpicture}[-latex]
    \matrix(m)[cdiag,ampersand replacement=\&]{
    #1 \& #2 \\
    #3 \& #4  \\};
    \path [arrs] (m-1-1) edge node {$#5$} (m-1-2)
                 (m-2-1) edge node {$#6$} (m-2-2)
                 (m-1-1) edge node [swap] {$#7$} (m-2-1)
                 (m-1-2) edge node {$#8$} (m-2-2)
				 (m-1-2) edge node [swap] {$#9$} (m-2-1);
\end{tikzpicture}
}

\newcommand{\dpoobjects}[6]{%
\def\dpooA{#1}
\def\dpooB{#2}
\def\dpooC{#3}
\def\dpooD{#4}
\def\dpooE{#5}
\def\dpooF{#6}}

\newcommand{\dpoarrows}[7]{%
\def\dpoaA{#1}
\def\dpoaB{#2}
\def\dpoaC{#3}
\def\dpoaD{#4}
\def\dpoaE{#5}
\def\dpoaF{#6}
\def\dpoaG{#7}}

\newcommand{\dposquares}[2]{%
  #1
  #2
	\begin{tikzpicture}
		\matrix (m) [cdiag,ampersand replacement=\&] {
        \dpooA \& \dpooB \& \dpooC \\
        \dpooD \& \dpooE \& \dpooF \\
		};
		\path [arrs]
	  	  	(m-1-2) edge node [swap] {$\dpoaA$} (m-1-1)
	  	  	(m-1-2) edge node {$\dpoaB$} (m-1-3)

	  	  	(m-1-1) edge node [swap] {$\dpoaC$} (m-2-1)
	  	  	(m-1-2) edge node {$\dpoaD$} (m-2-2)
	  	  	(m-1-3) edge node {$\dpoaE$} (m-2-3)

	  	  	(m-2-2) edge node {$\dpoaF$} (m-2-1)
	  	  	(m-2-2) edge node [swap] {$\dpoaG$} (m-2-3);
	  \NEbracket{(m-2-1)};
    \NWbracket{(m-2-3)};
	\end{tikzpicture} 
}

\newcommand{\vkcleararrows}{%
\tikzstyle{vkarrow1}=[arrow plain]
\tikzstyle{vkarrow2}=[arrow plain]
\tikzstyle{vkarrow3}=[arrow plain]
\tikzstyle{vkarrow4}=[arrow plain]
\tikzstyle{vkarrow5}=[arrow plain]
\tikzstyle{vkarrow6}=[arrow plain]
\tikzstyle{vkarrow7}=[arrow plain]
\tikzstyle{vkarrow8}=[arrow plain]
\tikzstyle{vkarrow9}=[arrow plain]
\tikzstyle{vkarrow10}=[arrow plain]
\tikzstyle{vkarrow11}=[arrow plain]
\tikzstyle{vkarrow12}=[arrow plain]}
\vkcleararrows

\newcommand{\topobjects}[4]{%
\def\topoA{#1}
\def\topoB{#2}
\def\topoC{#3}
\def\topoD{#4}}

\newcommand{\bottomobjects}[4]{%
\def\botoA{#1}
\def\botoB{#2}
\def\botoC{#3}
\def\botoD{#4}}

\newcommand{\toparrows}[4]{%
\def\topaA{#1}
\def\topaB{#2}
\def\topaC{#3}
\def\topaD{#4}}

\newcommand{\bottomarrows}[4]{%
\def\botaA{#1}
\def\botaB{#2}
\def\botaC{#3}
\def\botaD{#4}}

\newcommand{\sidearrows}[4]{%
\def\sideaA{#1}
\def\sideaB{#2}
\def\sideaC{#3}
\def\sideaD{#4}}

\newcommand{\vksquare}[5]{%
#1
#2
#3
#4
#5
\begin{tikzpicture}
\matrix (m) [matrix, text height=1.5ex, text depth=0.25ex,inner sep=0.5em,ampersand replacement=\&,row sep=0.5cm,column sep=0.5cm] {
  \node(m-1-1){$\topoA$}; \& \& \& \& \node(m-1-5){$\topoB$}; \\
  \& \node[scale=0.8] (m-2-2){$\topoC$}; \& \& \node[scale=0.8](m-2-4){$\topoD$}; \& \\
  \& \& \& \\
  \& \node[scale=0.8](m-4-2){$\botoA$}; \& \& \node[scale=0.8](m-4-4){$\botoB$}; \& \\
  \node(m-5-1){$\botoC$}; \& \& \& \& \node(m-5-5){$\botoD$}; \\
};
\path [arrs]
  (m-1-1) edge [vkarrow1] node {$\topaA$} (m-1-5)
  (m-2-2) edge [vkarrow2] node [swap,pos=0.35,scale=0.9] {$\topaB$} (m-1-1)
  (m-2-4) edge [vkarrow3] node [pos=0.35,scale=0.9] {$\topaC$} (m-1-5)
  (m-2-2) edge [vkarrow4] node [scale=0.8] {$\topaD$} (m-2-4)
                     
  (m-4-2) edge [vkarrow9] node [swap,scale=0.8] {$\botaA$} (m-4-4)
  (m-4-2) edge [vkarrow10] node [swap,pos=0.35,scale=0.9] {$\botaB$} (m-5-1)
  (m-4-4) edge [vkarrow11] node [pos=0.35,scale=0.9] {$\botaC$} (m-5-5)
  (m-5-1) edge [vkarrow12] node [swap] {$\botaD$} (m-5-5)
                     
  (m-1-1) edge [vkarrow5] node [swap] {$\sideaA$} (m-5-1)
  (m-2-2) edge [vkarrow6] node [swap,scale=0.8] {$\sideaB$} (m-4-2)
  (m-2-4) edge [vkarrow7] node [scale=0.8] {$\sideaC$} (m-4-4)
  (m-1-5) edge [vkarrow8] node {$\sideaD$} (m-5-5);
\end{tikzpicture}}

\title{Pictures of Processes\\\bigskip {\Large Automated Graph Rewriting for Monoidal Categories and Applications to Quantum Computing}}
\author{Aleks Kissinger}
\college{St. Catherine's College}

\degree{Doctor of Philosophy}
\degreedate{Michaelmas 2011}

\begin{document}
	\maketitle
	\baselineskip=18pt plus1pt
	
	\newpage
	\tableofcontents
	\newpage
	
	
	\begin{abstract}
    This work is about diagrammatic languages, how they can be represented, and what they in turn can be used to represent. More specifically, it focuses on representations and applications of string diagrams. String diagrams are used to represent a collection of processes, depicted as ``boxes'' with multiple (typed) inputs and outputs, depicted as ``wires''. If we allow plugging input and output wires together, we can intuitively represent complex compositions of processes, formalised as morphisms in a monoidal category.
    
    While string diagrams are very intuitive, existing methods for defining them rigorously rely on topological notions that do not extend naturally to automated computation. The first major contribution of this dissertation is the introduction of a discretised version of a string diagram called a \textit{string graph}. String graphs form a partial adhesive category, so they can be manipulated using double-pushout graph rewriting. Furthermore, we show how string graphs modulo a rewrite system can be used to construct free symmetric traced and compact closed categories on a monoidal signature.
    
    The second contribution is in the application of graphical languages to quantum information theory. We use a mixture of diagrammatic and algebraic techniques to prove a new classification result for strongly complementary observables. Namely, maximal sets of strongly complementary observables of dimension $D$ must be of size no larger than 2, and are in 1-to-1 correspondence with the Abelian groups of order $D$. We also introduce a graphical language for multipartite entanglement and illustrate a simple graphical axiom that distinguishes the two maximally-entangled tripartite qubit states: GHZ and W. Notably, we illustrate how the algebraic structures induced by these operations correspond to the (partial) arithmetic operations of addition and multiplication on the complex projective line.
    
    The third contribution is a description of two software tools developed in part by the author to implement much of the theoretical content described here. The first tool is Quantomatic, a desktop application for building string graphs and graphical theories, as well as performing automated graph rewriting visually. The second is QuantoCoSy, which performs fully automated, model-driven theory creation using a procedure called conjecture synthesis.
	\end{abstract}
	
	\newpage
	



\chapter{Introduction}\label{ch:intro}

Quantum information theory is the study of how data can be encoded and manipulated using microscopic systems subject to quantum effects. Over the past two decades, it has grown into a large and diverse field, with applications in security, where quantum effects are used to design ``unlistenable'' data channels, foundations of physics, where fundamental principles of information are used to derive physical theories, and perhaps most notably quantum computing, where classically intractable computations such as factorisation of huge numbers can happen in the blink of an eye. Virtually all of these applications use quantum theory exactly as John von Neumann described it in 1932. However, amidst the increasing scale the problems considered, it becomes clear that this is analogous to writing complex computer programs using circuit diagrams. As in the case with software development, abstracting away from the low-level is crucial to progress.

In this dissertation, we seek out this abstraction by identifying and exploiting the behaviour of graphical representations of quantum systems. We develop a tool set for graphical reasoning by drawing a connection between categorical algebra and graph rewriting. We then show how these this tool set can be applied to the description of quantum phenomena using the language of \textit{string diagrams}.

String diagrams consist of boxes, which represent processes (physical, logical, algebraic, ...) that have some inputs and some outputs. Some of those inputs and outputs can be connected together using wires.
\ctikzfig{tensor_diagram}

The only real requirement we impose on string diagrams is that their ``value'' (typically as some sort of map, relation, or process) is unaffected by topological deformations. Due to the strongly physical and spatial qualities of string diagrams, it should come as no surprise that they were originally formulated by a physicist. String diagrams originated with Roger Penrose in 1971~\cite{Penrose1971} as an alternative notation for contractions of what he called \textit{abstract tensors}, which are essentially just morphisms with some named inputs and outputs. Furthermore, the idea of representing spatially and temporally composed processes using these types of diagrams dates back at least to the 1948 advent of Feynman diagrams~\cite{Kaiser2005}.

String diagrams make sense for any mathematical structure that has a well-behaved notion of horizontal (i.e. spatial) and vertical (i.e. temporal) composition. A very general way to formalise such structures is to use monoidal categories, which were introduced by Mac Lane~\cite{MacLane} to describe a wide variety of categories admitting associative, product-like structures (e.g. cartesian products, direct sums, tensor products).

A connection between the notions of string diagrams and monoidal categories was inevitable. Twenty years after the introduction of string graphs, Joyal and Street~\cite{JS} formalised this idea by using string diagrams (considered as topological graphs with extra structure) to build \textit{free monoidal categories}. Intuitively, a ``free X'' is an object for which the axioms of an ``X'' are true, but nothing else. So, a free monoidal category is a monoidal category where two morphisms are equal \textit{if and only if} they are equal by the axioms of a monoidal category. In other words, string diagrams, compared up to topological deformations (of a particular kind) exactly represent morphisms compared up to the axioms of a monoidal category.

A monoidal category is a category equipped with a bifunctor $\otimes : \mathcal V \times \mathcal V \rightarrow \mathcal V$ that is associative, up to isomorphism and has a left and right unit $I \in \textrm{ob}\mathcal V$. There are many notions of monoidal categories with additional structure that have an extremely wide variety of applications in areas such as the study of braids and knots, linear algebra and representation theory, quantum field theory, higher-dimensional algebra, enriched and internal category theory, homotopy theory, linear logic, and programming language semantics. We introduce a few of these extended notions of monoidal category in chapter \ref{ch:monoidal-categories}, namely strict and non-strict (planar) monoidal categories, braided and symmetric monoidal categories, symmetric traced categories, left- and right-autonomous categories, compact closed categories, $\dagger$-monoidal (pronounced dagger-monoidal) categories, and $\dagger$-compact closed categories. We offer a summary of what is known about the relationships between these kinds of categories, coherence results, and most importantly, graphical language theorems. A much more comprehensive collection of graphical language definitions, as well as the state of the art in what is and is not known about these languages is available in Selinger's excellent survey paper~\cite{selinger2009survey}.

We also review the notion of abstract tensor networks, roughly in the form it was introduced by Penrose in \cite{Penrose1971} and relate its formulation to monoidal categories and the topological graphical languages introduced by Joyal and Street.

One of the most useful aspects of a monoidal category is that it allows one to define algebraic structures \textit{internal} to a monoidal category. That is, an algebraic structure can be defined as a collection of morphisms in some monoidal category satisfying some axioms. Since such a definition only relies on the structure of a monoidal category, it makes sense in \textit{any} monoidal category. For instance, one can define a monoid in $\mathcal V$ as a triple $(A, \mu, \eta)$ where $A$ is an object and $\mu : A \otimes A \rightarrow A$ and $\eta : I \rightarrow A$ are two morphisms satisfying some equations (namely, associativity and unit laws). A monoid in the category of sets and functions is just the usual notion of a monoid, i.e. a unital semigroup. In the category of vector spaces and linear maps, it is a unital, associative algebra. In the opposite category, it is a counital, coassociative coalgebra. In the category of categories and functors, it is a (strict) monoidal category, justifying the intuition that a monoidal category is just a ``categorified'' monoid. In chapter \ref{ch:monoidal-algebra}, we define algebraic structures internal to a monoidal category and give various examples that will be used throughout this dissertation. These include monoids, commutative monoids, comonoids, Frobenius algebras, bialgebras, and Hopf algebras. We also provide many concrete examples of these algebraic structures, as they occur in familiar (and some less-familiar) categories.

Building on this background, the bulk of the thesis is organised into two roughly independent parts. The first part is about applying techniques from the theory of graph rewriting to string diagrams. The second part is about applying monoidal category theory and graphical languages to the study of quantum mechanics. In particular, diagrams are used to study quantum computing and quantum entanglement theory. A third, shorter part focuses on implementing the theoretical work from the previous two parts in a program called Quantomatic.

Part \ref{part:rewriting} opens with an introduction to rewrite systems in chapter \ref{ch:rewrite-systems}. Rewrite systems provide a very general means of reasoning systematically about algebraic theories. In fact, this reasoning is so systematic that it can be done by a computer. Rewriting lives at the heart of most computer algebra systems (CAS), automated reasoning tools, and proof assistants. The idea behind rewriting is very simple. Rather than considering \textit{equations} $(s = t)$ between terms, as one typically does in (universal) algebra, one considers directed reductions called \textit{rewrite rules} $(s \rewritesto t)$. The application of rewrite systems from an algebraists point of view is that they can help solve \textit{word problems}.

A word problem is a question of the form, ``Is term $s$ equivalent to term $t$ by the axioms of an equational theory $E$?'' It is well known that word problems are not decidable in general. However, given a suitably nice algebraic theory and some elbow grease, it often \textit{is} possible to solve a word problem by turning $E$ into a rewrite system $R$. If we end up with a nice enough rewrite system, we can solve the word problem by rewriting $s$ repeatedly until no rule from $R$ applies (called normalising $s$), doing the same to $t$, and comparing the two results to see if they are equal.

If we just randomly pick directions for each of the equations in $E$, this technique is very unlikely to work. However, if we can find a nice rewrite system (i.e. one that is terminating and confluent), normal forms always exist and are unique, and there is an evident solution to the word problem for all terms. Thus a large portion of the rewriting literature is about how to go about turning sets of equations into nice rewrite systems, turning ill-behaved rewrite systems into better ones, and coping with ill-behaved systems using more sophisticated strategies than ``normalise and compare''.

Rewrite systems are not just restricted to terms. Just as term rewriting can be thought of us replacing certain subtrees (corresponding to subterms) with other trees, we can consider replacing certain subgraphs with another graph. This is called \textit{graph rewriting}. In 1973, Ehrig, Pfender, and Schneider introduced the double pushout (DPO) approach to graph rewriting~\cite{Ehrig1973}. We explain this technique in detail and with examples in section \ref{sec:graph-rewriting}. While DPO rewriting can be formulated in many categories (including the category of sets, the category of graphs, and any topos), DPO rewriting is not well-defined in all categories with pushouts. In 1979, Ehrig and Kreowski identified certain abstract properties of a category with pushouts that allow one to do double-pushout rewriting~\cite{Ehrig1979}. One abstract formulation of categories in which DPO rewriting makes sense are adhesive categories, introduced by Lack and Soboci\'nski in 2003~\cite{LackAdh2003}. In these categories, pushouts involving monomorphisms behave like coverings in the sense that they form so-called \textit{van Kampen Squares}. All toposes are adhesive categories, and in a recent result~\cite{Lack2011}, Lack showed that any adhesive category embeds fully and faithfully in a topos, and that embedding preserves all the adhesive structure. So, another way to think of adhesive categories is ``categories where pushouts of monomorphisms behave as they do in toposes''.

In section \ref{sec:partial-adhesive}, we generalise adhesive categories to \textit{partial adhesive categories}. These are categories $\mathcal C$ that embed fully and faithfully in an adhesive category $S : \mathcal C \rightarrow \mathcal A$, such that $S$ preserves monomorphisms. Intuitively, these are categories whose objects are the objects of an adhesive category (e.g. directed graphs) that satisfy certain additional axioms (e.g. simple graphs: at most one edge connecting any vertex to another). We then illustrate that the adhesive-like properties of pushouts in $\mathcal A$ are inherited by the $S$-pushouts in $\mathcal C$ (i.e. the pushouts in $\mathcal C$ that exist and are preserved by $S$). As a result, DPO rewriting is well-defined for partial adhesive categories as long as one restricts to certain matching morphisms called $S$-matchings.

Graph rewriting can be applied to string diagrams, but not directly. This is for the simple reason that string diagrams are not graphs in a strict sense. The wires in string diagrams need not be connected to boxes at both ends. They can even be connected to themselves to form circles. Wires that are not connected to a box at their source serve as inputs for string diagrams, and wires that are not connected at their target serve as outputs. Wires that are not connected to a box at either end are called \textit{free wires}, and represent identity maps. A directed graph $G$ consists of a set of vertices $V_G$ and a set of edges $E_G$, as well as \textit{total} functions $s,t : E_G \rightarrow V_G$. Therefore we cannot represent string diagrams as digraphs. Even if we relax the requirement that $s$ and $t$ be total functions, there is no way to distinguish circles from free wires.

We solve this problem by defining \textit{string graphs}. These are typed graphs whose vertices fall into two categories: wire-vertices and box-vertices. The wires in string diagrams are replaced by chains of wire-vertices.
\ctikzfig{wire}

Representing boxes as box-vertices, we can translate string diagrams into string graphs.
\ctikzfig{string_diagram_to_graph}

In chapter \ref{ch:string-graphs}, we define the partial adhesive category of string graphs and string graph homomorphisms. We also define special pushouts called \textit{pluggings}, which are used to plug the outputs of one string graph into the inputs of another string graph to form the composed graph. We also define string graph rewrite rules in such a way that any monomorphism is an $S$-matching. Therefore double-pushout rewriting is always well-defined.

The vigilant reader will notice that the correspondence between string diagrams and string graphs is \textit{nearly} 1-to-1. The only obstacle is that wires in string diagrams can be converted to chains of wire-vertices of any length. To eliminate this redundancy, we consider two string graphs to be equivalent if the only difference between the two is the length of the wires.
\ctikzfig{wire_homeo_example}

This corresponds to the wires of the associated string diagrams being homeomorphic, when considered as subspaces of topological graphs. For that reason, this equivalence relation is called \textit{wire-homeomorphism}. We define wire-homeomorphism using a confluent, terminating rewrite system on string graphs in section \ref{sec:homeomorphism}.

Like their topological counterparts, string graphs can be used to construct free monoidal categories. We do this by defining a framed cospan construction over the category of string graphs. Recall that for any category $\mathcal C$ with pushouts, the bicategory of cospans $\Csp(\mathcal C)$ has as objects the objects of $\mathcal C$, $1$-morphisms cospans $\ccospan{X}{f}{F}{g}{Y}$, and $2$-morphisms cospan homomorphisms. Composition is performed by pushout, and identities are cospans of identity maps. We form the category of framed cospans of string graphs by restricting the objects in the cospan construction to discrete graphs consisting of wire-vertices and the cospans $\ccospan{X}{}{G}{}{Y}$ to maps covering the inputs and outputs of $G$. Composition by pushout then reduces to the intuitive notion of plugging together string graphs.
\ctikzfig{cospan_plugging}

In section \ref{sec:sg-free-monoidal-categories}, we show that the free symmetric traced category and the free compact closed category on a monoidal signature can be constructed as a category of framed cospans of string graphs. By shifting from a topological representation to a combinatoric one, these morphisms can be represented straightforwardly on a computer, and they can be manipulated using automated graph rewriting techniques. This is explored in Part \ref{part:automation}, following the introduction of graphical theories for quantum computing.

Part \ref{part:entanglement} describes in detail two graphical theories that are of particular interest for quantum computing. These theories were formulated in the context of Categorical Quantum Mechanics (CQM), a program initiated by Abramsky and Coecke in 2004~\cite{AC2003,AC2004} whose purpose was to study quantum phenomena from the point of view of monoidal category theory. More than any one particular result, CQM represents a set of principles and an approach to the study of quantum theory. In this approach, compositionality is at the forefront. CQM asserts that all of the interesting and important aspects of quantum theory can be witnessed by studying systems and processes and the ways in which they compose. It emphasises the role of compound systems, information flow, and diagrammatic reasoning while de-emphasising the role of Hilbert spaces as a crucial component to the understanding of quantum phenomena.

This part opens with a brief introduction for the non-physicist to quantum mechanics, quantum computing, and quantum information theory in chapter \ref{ch:quantum}. Chapter \ref{ch:categorical-quantum-mechanics} introduces categorical quantum mechanics and illustrates the role of monoidal categories in quantum teleportation and the study of complementary observables. The latter was explored in detail by Coecke and Duncan in \cite{Coecke2008}. In quantum mechanics, an observable comes with a basis of orthonormal eigenstates corresponding to measurement outcomes. Two non-degenerate observables $O$ and $O'$ are called complementary if their associated bases of eigenstates are mutually unbiased. That is, for bases $\{ \ket{u_i} \}$, $\{ \ket{v_j} \}$ and for all $i,j$ (where $D$ is the dimension of the space):
\[ |\braket{u_i}{v_j}|^2 = \frac{1}{D} \]

Intuitively, if we measure a quantum state in an eigenstate of $O$ with respect to the $O'$ observable, we are equally likely to get any outcome. That is, maximal knowledge of a state with respect to $O$ implies minimal knowledge with respect to $O'$. A familiar example of complementary observables is position and momentum.

Mutually unbiased bases can be understood algebraically using particular kinds of interacting Frobenius algebras. Frobenius algebras in a monoidal category consist of a monoid $(A, \mult, \unit)$ and a comonoid $(A, \comult,\counit)$ satisfying the \textit{Frobenius identity}.
\ctikzfig{frobenius_comonoid_def}

$\dagger$-Frobenius algebras are Frobenius algebras whose comonoid structure is the adjoint of the monoid structure, i.e. $\comult = \left(\mult\right)^\dagger$ and $\counit = \left(\unit\right)^\dagger$.

Commutative $\dagger$-Frobenius algebras have attracted attention in recent years for exhibiting precisely the identities of 2-dimensional cobordisms, whose understanding is a crucial stepping stone to the formulation of topological quantum field theories. For more details, see e.g.~\cite{Atiyah1989,Kock2003}.

Special Frobenius algebras satisfy an additional identity on the loop map.
\ctikzfig{special}

As the name suggests, $\dagger$-special commutative Frobenius algebras ($\dagger$-SCFAs) are commutative $\dagger$-Frobenius algebras that are special. Coecke, Pavlovic, and Vicary showed that orthonormal bases over finite-dimensional complex Hilbert spaces are in 1-to-1 correspondence with $\dagger$-SCFAs. So, rather than studying mutually unbiased bases themselves, we can study their associated $\dagger$-SCFAs. From this point of view, the mutually unbiased basis condition can be summed up in a simple graphical identity, where $(\whitemult, \whiteunit, \whitecomult, \whitecounit)$ is the $\dagger$-SCFA induced by an orthonormal basis and $(\graymult, \grayunit, \graycomult, \graycounit)$ is the $\dagger$-SCFA induced by another, mutually unbiased basis.
\begin{equation}\label{eq:intro-compl}
\beginpgfgraphicnamed{op_dir_hopf}
\InputIfFileExists{op_dir_hopf.tikz}{}{\input{./figures/op_dir_hopf.tikz}}
\endpgfgraphicnamed
\end{equation}

In \cite{Coecke2008}, the Coecke and Duncan introduced several stronger forms of complementarity. One example is the case where the induced algebras of the two bases extend to a bialgebra. That is, the following equations are satisfied.
\begin{equation}\label{eq:intro-strong-compl}
\beginpgfgraphicnamed{zx_bialg}
\InputIfFileExists{zx_bialg.tikz}{}{\input{./figures/zx_bialg.tikz}}
\endpgfgraphicnamed
\end{equation}

In this dissertation, we refer to a pair of observables whose bases satisfy this condition as \textit{strongly complementary} observables. Coecke and Duncan showed that under certain additional assumptions, the equations in (\ref{eq:intro-strong-compl}) imply (\ref{eq:intro-compl}) in an arbitrary compact closed category. In section \ref{sec:complementary-obs}, we simplify this result in the case of the category of finite-dimensional Hilbert spaces by providing a new proof that (\ref{eq:intro-strong-compl}) \textit{always} implies (\ref{eq:intro-compl}). We also provide a new classification result for pairs of strongly complementary observables.

\begin{theorem*}
  Strongly complementary pairs of observables in a Hilbert space of dimension $D$ are in 1-to-1 correspondence with the finite Abelian groups of order $D$.
\end{theorem*}

Furthermore, we show that it is impossible for three distinct observables to be pairwise strongly complementary. This then classifies maximal sets of strongly complementary observables for all dimensions.

In chapter \ref{ch:monoidal-entanglement}, we turn to the application of diagrammatic techniques in the study of multipartite entanglement. The classification, computation, and manipulation of complex, many-body entangled quantum states is one of the most difficult problems facing quantum physicists and quantum information theorists. Any na\"ive approach to the problem of classifying multipartite entanglement is doomed to fail, and brute-force calculations involving many entangled quantum systems are untenable on today's computers. This suggests the need for more sophisticated techniques that capture and exploit as many symmetries and fundamental structure within a quantum system as possible. Rather than studying a multipartite state as a single, monolithic entity, we study it in terms of its components and explore how those components interact. We call this the \emph{compositional approach} to multipartite entanglement.

By way of the Choi-Jamio\l{}kowski isomorphism, we can consider quantum states and processes on the same footing. In that sense, a bipartite quantum state in $H \otimes H$ can be thought of as a quantum channel from $H$ to $H$. Similarly, we can treat a tripartite state as a map from $H \otimes H$ to $H$, i.e. a \textit{binary operation} on quantum states in $H$. Nearly all algebraic objects of interest (e.g. groups, rings, vector spaces) are sets equipped with one or more binary operations satisfying certain axioms. For that reason, we adopt a motto: ``Just as binary operations have a special status in the study of algebra, so too should tripartite states in the study of multipartite entanglement.''

To justify this assertion, we develop a methodology for representing and studying arbitrary qubit states using tripartite states as building blocks. It is a well-known result from quantum entanglement theory that there exist two canonical, genuinely-entangled tripartite states over qubits, up to equivalence by stochastic local operations and classical communication~\cite{DVC}. These states are the Greenberger-Horne-Zeilinger (GHZ) state and the W state.\footnote{W states were also first introduced by Greenberger, Horn, and Zeilinger in 1991~\cite{Zeilinger1991}, but they were not named until D\"ur, Vidal, and Cirac highlighted their significance in 2000~\cite{DVC}. It is generally believed that the W is for Wolfgang (D\"ur)~\cite{Cabello2002}.}
\[ \ketGHZ = \frac{1}{\sqrt{2}} \left( \ket{000} + \ket{111} \right) \qquad\qquad
   \ketW   = \frac{1}{\sqrt{3}} \left( \ket{100} + \ket{010} + \ket{001} \right) \]

In section \ref{sec:strong-symmetry-strong-max}, we identify two properties shared by GHZ and W states, which we call \textit{strong symmetry} and \textit{strong SLOCC-maximality}. Strongly symmetric states are symmetric states that extend naturally to larger symmetric states on any number of systems. For instance, the $N$-partite versions of GHZ and W are defined as:
\begin{align*}
  \ket{\textit{GHZ}_N} & := \ket{00\ldots0} + \ket{11\ldots1}                            \\
  \ket{\textit{W}_N}   & := \ket{10\ldots0} + \ket{010\ldots0} + \ldots + \ket{0\ldots01}
\end{align*}

SLOCC-maximal states are states that are maximal with respect to conversion by stochastic local operations and classical communication. That is, $\ket\Psi$ is SLOCC-maximal precisely when any state $\ket{\Psi'}$ that can be converted into $\ket\Psi$ by way of a SLOCC protocol must already be SLOCC-equivalent to $\ket\Psi$. Strongly SLOCC-maximal states are multipartite states that are \textit{inductively} SLOCC-maximal. That is, $\ket\Psi$ is a strongly SLOCC-maximal $N$ partite state if it is SLOCC-maximal and it is possible to obtain a strongly SLOCC-maximal $(N-1)$ partite state from $\ket\Psi$ by projecting out any of the subsystems $H$.

States that satisfy both of these conditions are called \textit{Frobenius states}. In section \ref{sec:frobenius-states}, we show that any Frobenius state extends to commutative Frobenius algebra. Conversely, commutative Frobenius algebras can be used to construct a Frobenius state. The commutative Frobenius algebra $\mathcal G$ associated with the Frobenius state $\ketGHZ$ is:
\begin{equation*}
  \begin{split}
    \whitemult & = \ket{0}\bra{00} + \ket{1}\bra{11} \qquad\qquad
    \whiteunit = \sqrt{2}\, \ket{+} = \ket{0}+\ket{1} \\
    \whitecomult & = \ket{00}\bra{0} + \ket{11}\bra{1} \qquad\qquad
    \whitecounit = \sqrt{2} \bra{+} = \bra{0}+\bra{1}
  \end{split}\vspace{-1.5mm}
\end{equation*}

For the Frobenius state $\ketW$, the associated Frobenius algebra $\mathcal W$ is:
\begin{equation*}
  \begin{split}
    \mult & = \ket{1}\bra{11} + \ket{0}\bra{01} + \ket{0}\bra{10} 
    \qquad\qquad\qquad
    \unit = \ket 1\qquad\qquad \\
    \comult & = \ket{00}\bra{0} + \ket{01}\bra{1} + \ket{10}\bra{1} 
    \qquad\qquad\qquad
    \counit = \bra 0\qquad\qquad
  \end{split}
\end{equation*}  

In \cite{CoeckeKissinger2010}, Coecke and Kissinger produced a unique characterisation of GHZ and W states in terms of properties of their associated Frobenius algebras, which is summarised in section \ref{sec:classification}. We highlight two types of commutative Frobenius algebras based on the value of the ``loop map''.
\ctikzfig{scfa_acfa_recap}

We prove that commutative Frobenius algebras are \textit{special} if and only if their associated tripartite states are SLOCC-equivalent to GHZ. Similarly, commutative Frobenius algebras are \textit{anti-special} if and only if their associated tripartite states are SLOCC-equivalent to W. Thus, the two canonical tripartite qubit states can be distinguished by two simple graphical identities.

Taking inspiration from the interaction properties of GHZ and W states, we define a \textit{GW-pair} as a pair of commutative Frobenius algebras (one special, one anti-special) satisfying certain graphical identities.
\ctikzfig{ghzw_bialg}
\ctikzfig{antiunit_cp}

We prove that the axioms of a GW-pair subsume the axiomatisation for GHZ and W states given in~\cite{CoeckeKissinger2010}. Again inspired by the example of the GHZ and W states, we introduce the notion of a \textit{distributive GW-pair}. Distributive GW-pairs behave similarly to rings, in that the ``multiplication'' induced by $\whitemult$ distributes, up to a scalar factor, over the ``addition'' defined by $\mult$.
\begin{equation}\label{eq:dist-propto}
\beginpgfgraphicnamed{distributive_propto}
\InputIfFileExists{distributive_propto.tikz}{}{\input{./figures/distributive_propto.tikz}}
\endpgfgraphicnamed
\end{equation}

Using a GW-pair, we can construct the abstract analogue of a CNOT gate, and verify graphically that it behaves as a CNOT. For the specific GW-pair defined by $\ketGHZ$ and $\ketW$, this is actually a CNOT gate. Using this fact, we prove that the generators of the pair $(\mathcal G, \mathcal W)$ are universal for quantum computation. By Choi-Jamio\l{}kowski, this means they can also be used to construct arbitrary multipartite qubit states.

The precise sense in which the GHZ algebra behaves like multiplication and the W algebra behaves like addition is explained in section \ref{sec:arithmetic}. Qubits defined on the Bloch sphere can equivalently be considered as points on the complex projective line $\cpline$. We can define (partial) addition and subtraction operations on the points in $\cpline$, considered as the set $\mathbb C$ with an additional point at infinity.
\begin{equation*}
\begin{aligned}
  k \cdot \infty      & = \infty & \qquad\qquad\qquad & k + \infty      & = \infty \\
  0 \cdot \infty      & = \bot   &                    & 0 + \infty      & = \infty \\
  \infty \cdot \infty & = \infty &                    & \infty + \infty & = \bot  
\end{aligned}
\end{equation*}

The GHZ algebra corresponds to the multiplication operation on $\cpline$, and the W algebra corresponds to addition. For finitary elements, multiplication distributes over addition. The failure of distributivity for $\infty$ in $\cpline$ (i.e. distributivity up to a non-zero scalar) is reflected in equation (\ref{eq:dist-propto}) by the fact that $k = 0$ when $\ket{a} = \numket{\infty}$.

There is still much to be learned about the algebras induced by quantum states, but already the compositional approach has yielded insights about the GHZ and W states that would not have been possible otherwise. One could picture quantum algorithms or protocols that leverage the behavioural qualities identified in this dissertation. For instance, treating inputs to graphs consisting of $\mult$ and $\whitemult$ as variables, we can think of such graphs as encoding polynomials in quantum states. Work is in progress to apply this insight to the development of quantum algorithms for hard problems such as finding the roots of diophantine polynomials.

In Part \ref{part:automation}, we introduce Quantomatic and QuantoCoSy, which are software tools for working with string graphs. Quantomatic allows a user to create and modify string graphs, graphical theories, and string graph rewrite systems. It also lets one selectively apply rewrite rules and normalise graphs with respect to a rewrite system. QuantoCoSy is a tool for synthesising \textit{new} graphical theories from concrete models using a technique called \textit{conjecture synthesis}. This procedure, introduced by Johansson, Dixon, and Bundy in 2010~\cite{Johansson:2010tk}, is a procedure for enumerating and checking equality for all terms of a certain size in an algebraic theory. The thing that makes this technique so effective is it builds a rewrite system dynamically during the enumeration procedure and actively avoids checking for redundant equalities, i.e. those that are already derivable using the rules it has discovered previously. It does this by only enumerating terms that are \textit{irreducible} with respect to a rewrite system. QuantoCoSy adapts this technique from term rewrite systems to string graph rewrite systems, and is showing potential to be a valuable tool in the generation of graphical theories from concrete, linear algebraic models.

In chapter \ref{ch:conclusion}, we review the major results of the dissertation and discuss future work, particularly in the area of automation.

\chapter{Monoidal Categories}\label{ch:monoidal-categories}

It is often useful to reason in a very general sense about processes and how they compose. Category theory provides the tool to do this. A category consists of a collection of objects $A, B, C, \ldots$, a collection morphisms $f, g, \cdots$, an associative operation $\circ$ for (vertical) composition, and for every object $A$ an identity morphism $1_A$. Objects can be thought of as types. They dictate which morphisms can be composed together. We shall primarily be interested in categories that have not only a vertical composition operation, but a horizontal composition as well. Such categories are called monoidal categories.

\begin{definition}\label{def:monoidal-category}
  A \emph{monoidal category} consists of a category $\mathcal V$, an object $I \in \mathcal V$ called the monoidal unit, a bifunctor $\otimes : \mathcal V \times \mathcal V \rightarrow \mathcal V$ called the monoidal product, and natural isomorphisms $\alpha_{A,B,C} : A \otimes (B \otimes C) \rightarrow (A \otimes B) \otimes C$, $\lambda_A : I \otimes A \rightarrow A$, and $\rho_A : A \otimes I \rightarrow A$, such that $\lambda_I = \rho_I$ and the following diagrams commute:
  \begin{equation}\label{eq:monidal-assoc}
    \begin{tikzpicture}
      \matrix (m) [cdiag] {
        A \otimes (B \otimes (C \otimes D)) &
        (A \otimes B) \otimes (C \otimes D) &
        ((A \otimes B) \otimes C) \otimes D \\
        A \otimes ((B \otimes C) \otimes D) & &
        (A \otimes (B \otimes C)) \otimes D \\
      };
      \path [arrs]
        (m-1-1) edge node {$\alpha$} (m-1-2)
        (m-1-2) edge node {$\alpha$} (m-1-3)
        (m-1-1) edge node {$A \otimes \alpha$} (m-2-1)
        (m-2-1) edge node {$\alpha$} (m-2-3)
        (m-2-3) edge node {$\alpha \otimes D$} (m-1-3);
    \end{tikzpicture}
  \end{equation}
  
  \begin{equation}\label{eq:monoidal-unit}
    \begin{tikzpicture}
      \matrix (m) [cdiag] {
        A \otimes (I \otimes B) & & (A \otimes I) \otimes B \\
        & A \otimes B & \\
      };
      \path [arrs]
        (m-1-1) edge node {$\alpha$} (m-1-3)
        (m-1-1) edge node [swap] {$A \otimes \lambda$} (m-2-2)
        (m-1-3) edge node {$\rho \otimes B$} (m-2-2);
    \end{tikzpicture}
  \end{equation}
\end{definition}

We shall refer to $(\otimes, \alpha, \lambda, \rho)$ as the \emph{monoidal structure} of $\mathcal V$. We often drop $\alpha$, $\lambda$, and $\rho$ when they are clear from the context. Monoidal categories where all three natural isomorphisms are actually equalities are called \emph{strict} monoidal categories.

\begin{examples}\label{exs:monoidal-categories}
  The condition of being a monoidal category is very weak. Most categories of interest admit at least one monoidal structure, and many admit several. Some examples:
  \begin{itemize}
    \item $(\catSet, \times, 1)$: the category of sets and total functions with the cartesian product $\times$ and the one-element set $1$ make $\catSet$ into a monoidal category.
    \item Disjoint union $+$ and the empty set form another monoidal structure on $\catSet$.
    \item More generally, any category with finite products or coproducts is monoidal.
    \item $(\catVect_K, \otimes, K)$: The category of $K$-vector spaces and $K$-linear maps is monoidal, with monoidal product taken as tensor product of vector spaces and tensor unit $K$, the 1-dimensional space.
    \item $(\catFVect_K, \otimes, K)$: The same as above, but restricted to finite-dimensional vector spaces.
    \item $(\catMat(K), \otimes, 1)$: The category whose objects are natural numbers and whose arrows $M : m \rightarrow n$ are $n \times m$ matrices taking values in $K$. Composition is matrix multiplication, the monoidal product is multiplication of natural numbers (on objects) and the Kronecker product of matrices (on arrows). This category is essentially $\catFVect_K$, with a chosen basis for all of its objects.
    \item $(\catRel, \times, 1)$: the category of sets and relations. Note that the cartesian product $\times$ is a monoidal product, but not a product in the categorical sense.
    \item $(\catRel, \oplus, \{ \})$: where $\oplus$ is the disjoint union of sets (on objects) and the disjoint union of relations (on arrows). We write the disjoint union using the $\oplus$ symbol to highlight the fact that it is actually a \textit{biproduct} in $\catRel$. As such, it is automatically a monoidal product.
  \end{itemize}
\end{examples}

In any monoidal category, $\alpha$, $\lambda$, and $\rho$ can be used to construct a natural isomorphism from some object to any other bracketing of that object, with or without monoidal units. E.g.
\[ (A \otimes I) \otimes (B \otimes (I \otimes C)) \cong (A \otimes (B \otimes (C \otimes I))) \]

It was shown by Mac Lane that the equations in Definition \ref{def:monoidal-category} suffice to show that \emph{any} such natural isomorphism is equal to any other one \cite{MacLane}. Such a theorem is known as a \emph{coherence theorem}, and it was the first of many concerning monoidal categories. By a minor abuse of notation, we shall often treat monoidal categories as if they were strict. That is, we often omit brackets, $\alpha$, $\lambda$, and $\rho$, simply assuming they are included where necessary. Coherence assures us that we can omit these details without ambiguity.

Though we shall occasionally use normal, algebraic notation for morphisms in monoidal categories, it is often vastly more convenient to use a graphical notation. In fact, the majority of this dissertation concerns formalising and exploiting graphical notation. For now, we shall treat this notation informally and fill in the details later. We represent objects as labelled wires:
\ctikzfig{cat_wires}

Morphisms can be thought of as processes. A morphism takes something of type $A$ and produces something of type $B$. For that reason, we'll draw morphisms as \emph{boxes} with a wire coming in labelled with a morphism's input type and a wire going out labelled with a morphism's output type.
\ctikzfig{cat_boxes}

Identity morphisms are special ``do nothing'' processes, which take something of type $A$ and return the thing itself. We represent these as empty wires.
\ctikzfig{cat_wires}

Morphisms are composed by plugging an output wire into an input wire.
\ctikzfig{cat_compose}

Implicit in this box-and-wire notation is the assumption that composition is associative.
\ctikzfig{cat_assoc}

\noindent ...and unital.
\ctikzfig{cat_unital}

We express the monoidal product of two objects as juxtaposition of wires.
\ctikzfig{cat_tensor_wires}

\noindent ...and the monoidal product of morphisms as the juxtaposition of boxes.
\ctikzfig{cat_tensor_boxes}

The monoidal product is also associative and unital, but possibly only up to isomorphism. We denote the monoidal unit $I$ as the empty graph. Note that the bifunctoriality of the tensor product is implicit in this notation.

\begin{center}
  $(g \otimes g') \circ (f \otimes f') = (g \circ f) \otimes (g' \circ f')$
  
  \bigskip

\beginpgfgraphicnamed{cat_bifunc}
\InputIfFileExists{cat_bifunc.tikz}{}{\input{./figures/cat_bifunc.tikz}}
\endpgfgraphicnamed
\end{center}

The following proposition is a simple consequence of binfunctoriality:

\begin{proposition}\label{prop:slide-boxes}
  For any morphisms $f : A \rightarrow B$ and $g : A' \rightarrow B'$ in a monoidal category,
  \[ (B \otimes g) \circ (f \otimes A') = (f \otimes B') \circ (A \otimes g) \]
\end{proposition}

Proposition \ref{prop:slide-boxes} can be interpreted graphically by ``sliding boxes'' past each other:

\begin{equation}\label{eq:slide-boxes}
\beginpgfgraphicnamed{cat_slide_boxes}
\InputIfFileExists{cat_slide_boxes.tikz}{}{\input{./figures/cat_slide_boxes.tikz}}
\endpgfgraphicnamed
\end{equation}

It was proved by Joyal and Street that planar box-and-wire diagrams can unambiguously represent morphisms in a monoidal category. They showed furthermore than this representation is sound and complete with respect to the algebraic definition of a monoidal category. This and similar results will be discussed at length in section \ref{sec:js-construction}.

\section{Types of Monoidal Categories}\label{sec:types-of-monoidal-categories}

So far, we have introduced monoidal categories. These are sometimes referred to as \emph{planar} monoidal categories, as diagrams of morphisms are always planar. However, many monoidal categories come with a notion of ``crossing'' wires. The weakest such category is called a braided monoidal category.

\begin{definition}\label{def:braided-monoidal-category}
  A \emph{braided monoidal category} is a monoidal category $(\mathcal V, \otimes, I, \alpha, \lambda, \rho)$ with an additional natural isomorphism $\gamma_{A,B} : A \otimes B \rightarrow B \otimes A$ called a \emph{braiding}, such that $\rho_A := \lambda_A \circ \gamma_{A, I}$ and the follow diagrams commute:
  \begin{equation}\label{eq:symmetry-hex}
    \begin{tikzpicture}
      \matrix (m) [cdiag,scale=0.4] {
        & (B \otimes A) \otimes C & B \otimes (A \otimes C) & \\
        (A \otimes B) \otimes C & & & B \otimes (C \otimes A) \\
        & A \otimes (B \otimes C) & (B \otimes C) \otimes A & \\
      };
      \path [arrs]
        (m-2-1) edge node {$\gamma \otimes C$} (m-1-2)
        (m-1-2) edge node {$\alpha$} (m-1-3)
        (m-1-3) edge node {$B \otimes \gamma$} (m-2-4)
        (m-2-1) edge node {$\alpha$} (m-3-2)
        (m-3-2) edge node {$\gamma$} (m-3-3)
        (m-3-3) edge node {$\alpha$} (m-2-4);
    \end{tikzpicture}
  \end{equation}
    
  \begin{equation}
    \begin{tikzpicture}
      \matrix (m) [cdiag,scale=0.4] {
        & (B \otimes A) \otimes C & B \otimes (A \otimes C) & \\
        (A \otimes B) \otimes C & & & B \otimes (C \otimes A) \\
        & A \otimes (B \otimes C) & (B \otimes C) \otimes A & \\
      };
      \path [arrs]
        (m-2-1) edge node {$\gamma^{-1} \otimes C$} (m-1-2)
        (m-1-2) edge node {$\alpha$} (m-1-3)
        (m-1-3) edge node {$B \otimes \gamma^{-1}$} (m-2-4)
        (m-2-1) edge node {$\alpha$} (m-3-2)
        (m-3-2) edge node {$\gamma^{-1}$} (m-3-3)
        (m-3-3) edge node {$\alpha$} (m-2-4);
    \end{tikzpicture}
  \end{equation}
\end{definition}

The braiding $\gamma$ and its inverse $\gamma^{-1}$ are drawn as wire crossings. Note how one wire is explicitly drawn over the top of the other.

\ctikzfig{cat_braid_crossings}

This is to emphasise that $\gamma$ may not be equal to $\gamma^{-1}$. That is, we cannot simply pass wires through each other in the graphical language. However, the naturality of $\gamma$ and the diagrams from Definition \ref{def:braided-monoidal-category} suffice to prove any equation about morphisms in a braided monoidal category that we can prove geometrically with braid diagrams. In particular, $\gamma$ satisfies the Yang-Baxter equation:

\ctikzfig{cat_yang_baxter}

We shall primarily be interested in a special case of a braided monoidal category called a symmetric monoidal category.

\begin{definition}\label{def:smc}
  A \textit{symmetric monoidal category} (SMC) is a monoidal category with a braiding $\sigma$ such that $\sigma_{A,B} = \sigma_{B,A}^{-1}$.
\end{definition}

When a braided monoidal category is symmetric, we refer to the braiding $\sigma$ as the \textit{symmetry map}. To emphasise the fact that $\sigma_{A,B} = \sigma_{B,A}^{-1}$, we do not distinguish over- and under-crossings:

\ctikzfig{cat_sym_crossing}

\begin{examples}
  All of the categories from Examples \ref{exs:monoidal-categories} are symmetric monoidal categories.
  \begin{itemize}
    \item $(Set, \times, \{*\})$, with $\sigma_{A,B} : A \times B \rightarrow B \times A$ defined as the canonical swap map $(a,b)\mapsto (b,a)$.
    \item $(Set, +, \{\})$, with $\sigma_{A,B} : A + B \rightarrow B + A$ the map that interchanges the two components of the disjoint union.
    \item For any category with finite products, the projection maps induce a canonical symmetry map:
    \begin{center}
      \begin{tikzpicture}
        \matrix (m) [cdiag] {
           & A \times B &   \\
         A &            & B \\
           & B \times A &   \\
        };
        \path [arrs]
          (m-1-2) edge node [swap] {$\pi_1$} (m-2-1)
          (m-1-2) edge node {$\pi_2$} (m-2-3)
          (m-3-2) edge node {$\pi_2$} (m-2-1)
          (m-3-2) edge node [swap] {$\pi_1$} (m-2-3)
          (m-1-2) edge [dashed] node {$\sigma_{A,B}$} (m-3-2);
      \end{tikzpicture}
    \end{center}
    \item The swap map is induced similarly in any category with finite coproducts.
    \item $(\catVect_k, \otimes, k)$ is an SMC, with $\sigma$ the tensor swap map. I.e. it is the linear extension of $\sigma_{V,W}(v \otimes w) = w \otimes v$.
    \item $(\catRel, \times, \{*\})$ is an SMC, with $\sigma$ defined as the swap map of the cartesian product.
  \end{itemize}
\end{examples}

We can interpret any progressive diagram (i.e. a diagram with no feedback loops) as a morphism in a symmetric monoidal category. Like in the case of planar monoidal categories, the axioms of a symmetric monoidal category ensure that there can be no ambiguity.
\ctikzfig{cat_smc_diagram}

The natural next question to ask would be, ``Is there a meaningful way to interpret diagrams \textit{with} feedback loops?'' The answer to this question is yes. There are actually two meaningful ways to do this. The first is a \textit{traced category}, and the second, which subsumes the first, is a \textit{compact closed category}.

\begin{definition}\label{def:traced-category}
  A symmetric traced category $\mathcal V$ is a symmetric monoidal category with a function
  \[ \Tr^X : \hom_{\mathcal V}(A \otimes X, B \otimes X) \rightarrow \hom_{\mathcal V}(A, B) \]
  defined for all objects $A, B, X$, satisfying the following five axioms:
  \begin{enumerate}
    \item $\Tr^X((g \otimes X) \circ f \circ (h \otimes X)) = g \circ \Tr^X(f) \circ h$
    \item $\Tr^Y(f \circ (A \otimes g)) = \Tr^X((B \otimes g) \circ f)$
    \item $\Tr^I(f) = f$ and $\Tr^{X \otimes Y}(f) = \Tr^X(\Tr^Y(f))$
    \item $\Tr^X(g \otimes f) = g \otimes \Tr^X(f)$
    \item $\Tr^X(\sigma_{X,X}) = 1_X$
  \end{enumerate}
  
  We refer to $\Tr^X$ is the \textit{trace operation} and $X$ as the object being \textit{traced out}.
\end{definition}

We depict this graphically by connecting the $X$-output of a map $f : A \otimes X \rightarrow B \otimes X$ to the $X$-input.
\ctikzfig{cat_trace}

Axioms 1 and 4 are implicit in this notation, since we do not draw the bounds of the trace operation. Axiom 2 is a ``box-sliding'' identity:
\ctikzfig{cat_trace_ax2}

A special case ($A = B = I$) of this axiom is the familiar property of matrix traces in linear algebra: $\Tr(M N) = \Tr (N M)$, i.e. the value of the trace is not affected by cyclic permutations.

Axiom 3 makes the trace operation respect the monoidal product on objects:
\ctikzfig{cat_trace_ax3}

Axiom 5 allows us to pull out loops in the diagram:
\ctikzfig{cat_trace_ax5}

\begin{examples}
  There are at least two ways in which the category $\catRel$ can be made into a symmetric traced category: one for each of the monoidal products defined from Examples \ref{exs:monoidal-categories}.
  \begin{itemize}
    \item $(\catRel, \otimes)$: For a relation $R : A \times X \rightarrow B \times X$, we define a new relation $\Tr^X(R) : A \rightarrow B$ as follows:
    \[ a \left[ \Tr^X(R) \right] b \Leftrightarrow \forall x \in X . (a,x) R (b,x) \]
    Thinking of a relation as a matrix over the booleans, this is analogous to the usual partial trace of a matrix.
    \item The trace for $(\catRel, \oplus)$ was defined in \cite{JSV}. Let $R : A \oplus X \rightarrow B \oplus X$ be a relation. We define the trace as:
    \[ a \left[ \Tr^X(R) \right] b \Leftrightarrow
         \left(a R b\right) \vee
         \left(\exists x, x' \in X\ .\ a R x \wedge x' R x \wedge x' R b \right) \]
    The new relation incorporates ``feedback'' from the $X$-output of $R$ to the $X$-input of $R$ via term ``$x' R x$'' on the RHS.
  \end{itemize}
\end{examples}

\begin{example}\label{ex:partial-trace-vect}
  The category $\catFVect_K$ of finite-dimensional vector spaces and linear maps is a traced monoidal category, with $\Tr$ given by the partial trace operation on a linear map. Suppose $f : A \otimes X \rightarrow B \otimes X$ is a linear map. Then, by fixing bases $x^i \in X$, $a^i \in A$ and $b^i \in B$, then $f$ is uniquely determined by an indexed collection $f^{i,j}_{k,l} \in K$ called a \textit{tensor}.
  \[ f(a^i \otimes x^j) = \sum_{k,l} f^{i,j}_{k,l} b^k \otimes x^l \]
  
  We can then define a new tensor by summing together the lower $X$-index with the upper $X$-index.
  \[ \tilde f = \sum_{k} f^{i,k}_{j,k} \]
  
  This new tensor defines a linear map from $A$ to $B$:
  \[ \Tr^X(f)(a^i) = \sum_j \tilde f^i_j b^j \]
  
  Equivalently, for $x_i \in X^*$ the corresponding basis of the dual space of $X$, we can define the partial trace as:
  \[ \Tr^X(f) := \sum_{i} (B \otimes x_i) \circ f \circ (A \otimes x^i) \]
  
  When $A = B = K$, this is just the usual trace of a matrix. With this is mind, we can see why $\catFVect_K$ is an example of a traced monoidal category, but $\catVect_K$ is not. For an infinite-dimensional vector space $V$, $\Tr^V(1_V)$ is also infinite. In particular, it is not an element of $\hom(K,K) \cong K$.
\end{example}

Note how the dual space plays a role in the definition of the trace. The dual space is actually a special case of a general categorical notion called a \textit{dual}. 

\begin{definition}\label{def:dual}
  Let $A$ and $A^*$ objects in a monoidal category. $A^*$ is called the \textit{right dual} of $A$ (equivalently, $A$ is called the \textit{left dual} of $A^*$) if there exist maps $d_A : A \otimes A^* \rightarrow I$, $e_A : I \rightarrow A^* \otimes A$ satisfying the ``line-yank'' identities:
  \begin{center}
    \begin{tikzpicture}
      \matrix (m) [cdiag] {
        A & \\
        A \otimes A^* \otimes A & A \\
      };
      \path [arrs]
        (m-1-1) edge node {$1_A$} (m-2-2)
        (m-1-1) edge node [swap] {$A \otimes e$} (m-2-1)
        (m-2-1) edge node [swap] {$d_A \otimes A$} (m-2-2);
    \end{tikzpicture}
    \quad
    \begin{tikzpicture}
      \matrix (m) [cdiag] {
        A^* & A^* \otimes A \otimes A^* \\
            & A^* \\
      };
      \path [arrs]
        (m-1-1) edge node [swap] {$1_A$} (m-2-2)
        (m-1-1) edge node {$e_A \otimes A^*$} (m-1-2)
        (m-1-2) edge node {$A^* \otimes d$} (m-2-2);
    \end{tikzpicture}
  \end{center}
  
  $e_A$ is called the \textit{cap}, and $d_A$ is called the \textit{cup}, of the compact structure.
\end{definition}

In the graphical notation, the object $A^*$ is represented as a wire labelled $A$, but directed upward instead of downward.
\ctikzfig{cat_dual_wire}

We represent $e_A$ and $d_A$ as half-turn of wire, forming a cup or a cap. 
\ctikzfig{cat_cap_cup}

The diagrams from Definition \ref{def:dual} are called the ``line-yank'' identities, because their graphical representations literally look like pulling a wire straight.
\ctikzfig{cat_line_yank}

\begin{definition}
  A monoidal category where every object has a right (resp. left) dual is called a \textit{right} (resp. left) \textit{autonomous category}.
\end{definition}

\begin{definition}\label{def:compact-closed-category}
  A \textit{compact closed category} is a category that is right autonomous and symmetric. For an object $A$ in a compact closed category, the dual maps $d_A$ and $e_A$ are called a \textit{compact structure} for $A$.
\end{definition}

Note that compact closed categories are automatically left autonomous. Any right dual $A^*$ of $A$ can be made into a left dual by choosing maps $d' := d \circ \sigma$ and $e' := \sigma \circ e$. The left line-yank identity can then be derived from symmetry and the right line-yank.
\ctikzfig{cat_reverse_link_yank}

Also note that in a compact category, and map $f : A \rightarrow B$ can also be considered as a map $f^* : A^* \rightarrow B^*$ by using caps and cups to ``bend the wires'' around.
\ctikzfig{bend_wires}

\begin{example}\label{ex:vect-compact}
  The category $\catFVect_K$ is compact closed. For any finite-dimensional vector space $A$, fix bases $a^i \in A$ and $a_i \in A^*$ such that $a_j \circ a^i = \delta_i^j$. Then, define a compact structure for $A$ as follows:
  \[ d_A :: a_i \otimes a^j \mapsto (a_i \circ a^j) \ \ \ \ \  e_A :: 1 \mapsto \sum_i a_i \otimes a^i \]
\end{example}

Compact closed categories are automatically symmetric traced categories. A cap and a cup can be used to construct a ``feedback loop'' that acts as a trace operation.
\begin{equation}\label{eq:compact-partial-trace}
  \Tr^X(f) := (B \otimes d_X) \circ (f \otimes X^*) \circ (A \otimes e'_X)
\end{equation}

The axioms of a compact structure then suffice to prove the five trace axioms given in Definition \ref{def:traced-category}. Therefore, compact closed categories subsume symmetric traced categories. Often these categories are easier to work with than their traced counterparts, especially using the graphical language. It would be convenient to use the compact structure axioms when proving identities in symmetric traced categories. This is possible due to a result by Joyal, Street, and Verity.

\begin{theorem}[\cite{JSV}]\label{thm:traced-full-emb}
  Any symmetric traced category can be fully and faithfully embedded in a compact category.
\end{theorem}

They prove this result by defining the free compact closure of a symmetric traced category, using a technique called the ``Int construction''. For a symmetric traced category $\mathcal V$, they build a compact closed category $\Int(\mathcal V)$ into which $\mathcal V$ embeds fully and faithfully. Thus $f = g$ in $\mathcal V$ if and only if $f = g$ in $\Int(\mathcal V)$. This construction is also free over $\mathcal V$, i.e. for a compact closed category $\mathcal V'$, a traced symmetric functor $F : \mathcal V \rightarrow \mathcal V'$ extends uniquely to a compact closed functor $\widetilde F : \Int(\mathcal V) \rightarrow \mathcal V'$.

We shall introduce one more type of monoidal category, introduced by Abramsky and Coecke, for the sake of reasoning about quantum information~\cite{AC2004}.

\begin{definition}
  A category $\mathcal C$ is called a \textit{$\dagger$-category} if there exists an identity-on-objects functor $(-)^\dagger : \mathcal C \rightarrow \mathcal C^\op$ such that $((-)^\dagger)^\dagger = 1_{\mathcal C}$.
\end{definition}

In particular, $\dagger$-categories are always isomorphic with their opposite category. As the notation might suggest, the $\dagger$ functor is an abstract version of the conjugate-transpose of a complex linear map. We can use it to define an abstract notion of unitarity.

\begin{definition}
  A morphism $U$ in a $\dagger$-category is called \textit{unitary} if it is an isomorphism and $U^\dagger = U^{-1}$.
\end{definition}

$\dagger$-monoidal categories are simply monoidal $\dagger$-categories, where all of the structural natural isomorphisms are actually unitary isomorphisms.

\begin{definition}
  A $\dagger$-symmetric monoidal category is a symmetric monoidal category where all of the components of $\alpha$, $\lambda$, and $\sigma$ are unitary.
\end{definition}

$\dagger$-compact closed categories have the additional property that the $\dagger$ reflects caps and cups vertically.

\begin{definition}
  A $\dagger$-compact closed category is a compact closed category where all of the components of $\alpha$, $\lambda$, and $\sigma$ are unitary and for all $A$, $d_A^\dagger = e_A'$ and $e_A^\dagger = d_A'$.
\end{definition}

\begin{examples}
  Several examples of $\dagger$-compact closed categories are:
  \begin{itemize}
    \item $(\catRel, \times)$, with $\dagger$ given by relational converse: $a R^\dagger b \Leftrightarrow b R a$
    \item $(\catMat(K), \otimes)$, with $\dagger$ given by transposition of matrices
    \item $(\catFHilb, \otimes)$, the category of finite-dimensional complex Hilbert spaces. The compact structure is the same as for $\catFVect_{\mathbb C}$, and $\dagger$ is given by the adjoint of linear operators with respect to $\langle - | - \rangle.$
  \end{itemize}
\end{examples}

The category $\catHilb$ is a $\dagger$-symmetric monoidal category, but it is not traced (and hence not compact closed).

For each of the categories we have introduced, we can form the category of (symmetric, traced, compact, ...) monoidal categories and suitably structure-preserving functors.

\begin{definition}\label{def:monoidal-functor}
  For monoidal categories $\mathcal C$, $\mathcal D$, a strong monoidal functor consists of functor $F : \mathcal C \rightarrow \mathcal D$ an isomorphism $\phi: I \rightarrow F(I)$ and a natural isomorphism $\psi_{A,B} : FA \otimes FB \rightarrow F(A \otimes B)$ such that following diagrams commute:
  \begin{center}
    \begin{tikzpicture}
      \matrix (m) [cdiag] {
        (FA \otimes FB) \otimes FC & F(A \otimes B) \otimes FC & F((A \otimes B) \otimes C) \\
        FA \otimes (FB \otimes FC) & FA \otimes F(B \otimes C) & F(A \otimes (B \otimes C)) \\
      };
      \path [arrs]
        (m-1-1) edge node {$\psi \otimes FC$} (m-1-2)
        (m-1-2) edge node {$\psi$} (m-1-3)
        (m-2-1) edge node {$FA \otimes \psi$} (m-2-2)
        (m-2-2) edge node {$\psi$} (m-2-3)
        (m-1-1) edge node [swap] {$\alpha$} (m-2-1)
        (m-1-3) edge node {$F(\alpha)$} (m-2-3);
    \end{tikzpicture}
    
    \begin{tikzpicture}
      \matrix (m) [cdiag] {
        FA \otimes I & FA \\
        FA \otimes FI & F(A \otimes I) \\
      };
      \path [arrs]
        (m-1-1) edge node {$\rho$} (m-1-2)
        (m-2-1) edge node {$\psi$} (m-2-2)
        (m-1-1) edge node [swap] {$FA \otimes \phi$} (m-2-1)
        (m-1-2) edge node {$F(\rho^{-1})$} (m-2-2);
    \end{tikzpicture}
    \qquad
    \begin{tikzpicture}
      \matrix (m) [cdiag] {
        I \otimes FA & FA \\
        FI \otimes FA & F(I \otimes A) \\
      };
      \path [arrs]
        (m-1-1) edge node {$\lambda$} (m-1-2)
        (m-2-1) edge node {$\psi$} (m-2-2)
        (m-1-1) edge node [swap] {$\phi \otimes FA$} (m-2-1)
        (m-1-2) edge node {$F(\lambda^{-1})$} (m-2-2);
    \end{tikzpicture}
  \end{center}
\end{definition}

The adjective \textit{strong} is used to distinguish from a \textit{lax monoidal functor}, where the isomorphisms are replaced by arbitrary maps. In the case of strict monoidal categories, this definition simplifies, as coherence diagrams need not be considered.

\begin{definition}\label{def:strict-monoidal-functor}
  For strict monoidal categories $\mathcal C$, $\mathcal D$, a \textit{strict monoidal functor} is a functor $F : \mathcal C \rightarrow \mathcal D$ where $I = F(I)$ and $F(A \otimes B) = FA \otimes FB$.
\end{definition}

For symmetric, traced, compact closed, and $\dagger$-compact closed categories, we can make similar definitions, where the functor application must commute with all of the structure in sight. There are also strict versions of all of these categories, where the associativity and unit natural isomorphisms are all identities.

\begin{definitions}\label{defs:moncats}
  The following are all categories of (small) monoidal categories:
  \begin{itemize}
    \item $\catMonCat$: the category of monoidal categories and monoidal functors
    \item $\catSymMonCat$: the category of symmetric monoidal categories and symmetric monoidal functors
    \item $\catSymTraceCat$: the category of symmetric traced categories and symmetric traced functors
    \item $\catCCCat$: the category of compact closed categories and compact functors
    \item $\catMonCat_s$, $\catSymMonCat_s$, $\catSymTraceCat_s$, $\catCCCat_s$: the strict versions
  \end{itemize}
\end{definitions}

\section{Free Monoidal Categories}

An important question for monoidal categories is as follows.

\begin{center}
  Given a suitable description for the generators of a (symmetric, traced, compact closed, ...) monoidal category, can we generate the \textit{free} such category?
\end{center}

This question is important because two arrows $f$ and $g$ are equal in a free monoidal category \textit{if and only if} their equality can be established only using the axioms of that category. We can make this precise by expounding on the usual universal property satisfied by free objects, but first, we define the notion of a \textit{monoidal signature}, which defines the generators of a monoidal category.

\begin{notation}
  For a set $O$, let $w(O)$ be the free monoid over $O$, i.e. the set of lists with elements taken from $O$ where multiplication is concatenation and the unit is the empty list. For a function $f : O \rightarrow O'$, let $w(f) : w(O) \rightarrow w(O')$ be the lifting of $f$ to lists:
\[ w(f)([A, B, C, \ldots]) = [f(A), f(B), f(C), \ldots] \]
\end{notation}

\begin{definition}\label{def:monoidal-signature}
  A (small, strict) monoidal signature $T = (O, M, \dom, \cod)$ consists of a set of objects $O$, a set of morphisms $M$, and a pair of functions $\dom : M \rightarrow w(O)$ and $\cod : M \rightarrow w(O)$.
\end{definition}

The maps $\dom$ and $\cod$ should be interpreted as giving input and output types to a morphism $m \in M$. For instance, if $\dom(m) = [A,B,C]$ and $\cod(m) = [D]$, then $m$ represents a morphism from $A \otimes B \otimes C$ to $D$. The empty list is interpreted as the tensor unit $I$.

\begin{example}
  Define a monoidal signature $T = (O, M, \dom, \cod)$ where
  \begin{align*}
    O    & = \{ A, B, C \} \\
    M    & = \{ f, g \} \\
    \dom & = \{ f \mapsto [A,B],\ g \mapsto [C] \} \\
    \cod & = \{ f \mapsto [C],\ g \mapsto [C] \}
  \end{align*}
  
  This signature defines three (primitive) objects $A, B, C$ and two morphisms $f : A \otimes B \rightarrow C$ and $g : C \rightarrow C$. We will often write a signature as a set of boxes, representing the diagram generators:
  \ctikzfig{box_set}
\end{example}

There is also a notion of a non-strict monoidal signature. In that case, $w(O)$ is replaced with the free $(\otimes, I)$-algebra over $O$. We treat the strict case for simplicity, but many of the results translate immediately, replacing equality with coherent natural isomorphism.

\begin{definition}\label{def:monsig}
  For monoidal signatures $S$, $T$, a monoidal signature homomorphism $f$ consists of functions $f_O : O_{S} \rightarrow O_{T}$ and $f_M : M_{S} \rightarrow M_{T}$ such that the following diagrams commute.
  \begin{center}
    \begin{tikzpicture}
      \matrix (m) [cdiag] {
        M_{S} & w(O_{S}) \\
        M_{T} & w(O_{T}) \\
      };
      \path [arrs]
        (m-1-1) edge node {$\dom_{S}$} (m-1-2)
        (m-2-1) edge node [swap] {$\dom_{T}$} (m-2-2)
        (m-1-1) edge node [swap] {$f_M$} (m-2-1)
        (m-1-2) edge node {$w(f_O)$} (m-2-2);
    \end{tikzpicture}
    \qquad\qquad
    \begin{tikzpicture}
      \matrix (m) [cdiag] {
        M_{\mathcal S} & w(O_{S}) \\
        M_{\mathcal T} & w(O_{T}) \\
      };
      \path [arrs]
        (m-1-1) edge node {$\cod_{S}$} (m-1-2)
        (m-2-1) edge node [swap] {$\cod_{T}$} (m-2-2)
        (m-1-1) edge node [swap] {$f_M$} (m-2-1)
        (m-1-2) edge node {$w(f_O)$} (m-2-2);
    \end{tikzpicture}
  \end{center}
  \noindent $\catMonSig$ is the category of monoidal signatures and monoidal signature homomorphisms.
\end{definition}

A monoidal signature is essentially a strict monoidal category without composition or identity maps. A monoidal signature homomorphism is thus a monoidal functor, minus the condition that it respect composition and identity maps.

\begin{definition}
  A monoidal signature is called \textit{simple} if the images of $\dom$ and $\cod$ are restricted to single-element lists.
\end{definition}

There are evident forgetful functors from $\catMonCat_s$, $\catSymMonCat_s$, $\catSymTraceCat_s$, and $\catCCCat_s$ into $\catMonSig$. If this forgetful functor has a left adjoint $F$, the image of a signature $T$ under $F$ is called the \textit{free monoidal category} over $T$.

To get a better feel for these objects, we unroll the universal property of the free category. Fix a monoidal category $\mathcal V$. Then, for a monoidal signature $T$, a monoidal signature homomorphism from $T$ to $U (\mathcal V)$ is called a \textit{valuation}. This homomorphism gives a value in $\mathcal V$ for each of the generators in $T$. The universal property of the adjunction then guarantees there is a unique strong monoidal functor $\tilde v : FT \rightarrow \mathcal V$ such that:

\begin{center}
  \begin{tikzpicture}
    \matrix (m) [cdiag] {
      T & UF(T)       \\
        & U\mathcal V \\
    };
    \path [arrs]
      (m-1-1) edge [right hook-latex] node {} (m-1-2)
      (m-1-1) edge node [swap] {$v$} (m-2-2)
      (m-1-2) edge [dashed] node {$U(\tilde v)$} (m-2-2);
  \end{tikzpicture}
\end{center}

In the non-strict case, $\tilde v$ is only unique up to unique, coherent natural isomorphism. In \cite{KellyLaplaza1980}, Kelly and Laplaza gave a prescription for constructing the free category on any ``algebraically-defined'' additional structure on a category. They went on to describe concretely the free compact closed category on a category (or equivalently, a simple signature). Shum proved a similar result in 1994~\cite{Shum1994} for tortile categories, i.e. braided monoidal categories with coherently-defined left and right duals. In the next section, we will discuss work to define the free symmetric, traced, and compact categories on an \textit{arbitrary} signature, using graphical language.

\section{Formalising Graphical Languages}\label{sec:graphical-language}

There are three ways in which one can formalise graphical languages for monoidal categories. The first formalisation is \textbf{algebraic}, where string diagrams are used as an equivalent representation of a tensor network defined using the \textit{abstract index} notation introduced by Penrose in his 1971 paper \cite{Penrose1971}. The second formalisation is \textbf{topological}. In this approach, topological graphs (i.e. realisations of 1D simplicial complexes) with added structure are used to represent morphisms \cite{JS}. The third formalisation is \textbf{combinatoric}. A special kind of typed graph called a \textit{string graph} is used to represent morphisms. This approach was developed by Dixon, Duncan, and Kissinger \cite{DixonDuncan2008,DixonKissinger2010}. In this section, we'll discuss the first two approaches. We'll discuss \textit{string graphs} in detail in Part \ref{part:rewriting}.

\subsection{Algebraic Approach: Abstract Tensor Systems}

In this section, we will look at the original formulation of string diagrams, due to Penrose in 1971 \cite{Penrose1971}.

Recall that a tensor is a set of real or complex numbers, indexed by one or more natural numbers. For example, the following is an $(n_1 \cdot n_2 \cdot n_3)$-dimensional tensor:
\[ \{ \chi_{i j}^{k} : i = 1..n_1; \ j = 1..n_2; \ k = 1..n_3 \} \]

Subscripts should be thought of as inputs and superscripts as outputs. Familiar examples of tensors are vectors, $v^i$ and matrices, $M_i^j$. We can compose tensors by \emph{contraction}, i.e. ``summing together'' a lower index and an upper index of the same dimension:
\[ \xi^i_j = \sum_{kl} \chi_{kl}^i \beta_j^{k} \rho^l \]

As expected, when we focus on vectors and matrices, we recover the usual notions of composition and application of linear maps. In order to simplify such expressions, we can use the Einstein summation convention, where any repeated indices are assumed to be summed over.
\[ \xi^i_j = \chi_{kl}^i \beta_j^{k} \rho^l \]

Penrose introduced \textit{abstract tensor systems} to express generalised tensors and contractions. In many ways, this formalism resembles that of monoidal categories. Natural number indices are replaced with formal \textit{labels}, which are simply names that can be used to identify inputs and ouputs. These are taken from a \textit{labelling set}.
\[ \mathcal L = \{ a, b, c, \ldots, a_0, b_0, \ldots, a_1, b_1, \ldots \} \]

Vector spaces are replaced by sets of formal tensors. For two lists of labels $U$, $L$ the set $\mathcal T_L^U$ has as elements \textit{formal tensors}. For $L = \{ a_0, \ldots, a_m \}$ and $U = \{ b_0, \ldots, b_n \}$, we write a formal tensor using the usual tensor notation:\footnote{Note that the data associated with the tensor includes a total ordering on the sets $L$ and $U$. This is implicit in the use of \textit{lists} of index names in the tensor.}
\[ \chi_{a_0,\ldots,a_m}^{b_0,\ldots,b_n} \in \mathcal T_L^U \]

It's useful to think of the sets $\mathcal T_L^U$ as something akin to a hom-set, whose elements $\xi$ are morphisms from $X^{\otimes |L|}$ (i.e. $|L|$ copies of $X$) to $X^{\otimes |U|}$. Penrose defines four operations over sets of abstract tensors.
\begin{align*}
  \textbf{Relabelling: \ \ } & R : \mathcal T_L^U \rightarrow
    \mathcal T_{L'}^{U'} \textrm{ for $L \cong L'$ and $U \cong U'$} \\
  \textbf{Addition: \ \ } & + : \mathcal T_L^U \times \mathcal T_L^U \rightarrow \mathcal T_L^U \\
  \textbf{Outer product: \ \ } & \otimes : \mathcal T_L^U \times \mathcal T_{L'}^{U'} \rightarrow \mathcal T_{L+L'}^{U+U'} \\
  \textbf{Contraction: \ \ } & C_p^q : \mathcal T_L^U \rightarrow \mathcal T_{L - \{p\}}^{U - \{q\}}
\end{align*}

These satisfy certain compatibility axioms (e.g. associativity and identity laws), mirroring those of normal tensor contraction.

\begin{remark}
  The usual, categorical composition can be defined in terms of outer product and contraction.
  \[ \beta_{y'_0,\ldots,y'_m}^{z_0,\ldots,z_n} \circ \alpha_{x_0,\ldots,x_l}^{y_0,\ldots,y_m} :=
      C_{y'_0}^{y_0} C_{y'_1}^{y_1} \cdots C_{y'_m}^{y_m}
      (\alpha_{x_0,\ldots,x_l}^{y_0,\ldots,y_m} \otimes
       \beta_{y'_0,\ldots,y'_m}^{z_0,\ldots,z_n}) \]

  If we then include for every pair of labels a Dirac delta tensor $\delta_a^b$ (i.e. an identity map), with suitable axioms, the data from an abstract tensor system defines a symmetric traced category. In chapter \ref{sec:sg-free-monoidal-categories} we make use of (essentially) the converse construction to prove the main theorem.
  
  With the inclusion of raising and lowering tensors $g_{a,b}$, $g^{a,b}$, it becomes a compact closed category with $X = X^*$. For this reason, abstract tensor systems are widely regarded as the prototype for compact closed categories, in their modern formulation.
\end{remark}

As in the concrete case, we represent outer (i.e. tensor) product as juxtaposition, contraction by repeating an index, and let relabeling be implicit.

\begin{example}
  The following is an abstract tensor contraction, followed by its explicit form, in terms of the functions above.
  \[ \alpha_{ab}^c \beta_c^d + \gamma_{ab}^d := C_c^{c'}(\alpha_{a,b}^c \otimes \beta_{c'}^d) + \gamma_{a,b}^d \]
\end{example}

However, even with this convention, contraction expressions can get quite complex. Consider this expression, involving six abstract tensors:
\begin{equation}\label{eqn:gross-tensor}
  \alpha_{abc}^{de} \beta_{f}^{bfg} \gamma_{dh}^{i}
  \rho_{i}^{h} \phi_{eg}^{jk} \delta_l^l  
\end{equation}

In order to work with this expression, one has to keep track of $11$ indices, which makes computations time-consuming and error-prone. To address this issue, Penrose introduced a second, graphical notation. Tensors are drawn as boxes, and contractions over pairs of indices as wires. The ``identity'' tensor (i.e. the Dirac delta $\delta_i^j$) is also drawn as a wire. The non-contracted, or ``free'' indices are left as dangling wires, and contractions $\delta_i^i$ are represented as circles. In the graphical notion, expression (\ref{eqn:gross-tensor}) becomes the following diagram:
\ctikzfig{tensor_diagram}

Such a diagram can always be interpreted, unambiguously as an abstract tensor network, up to a relabelling of indices. It is also clear how the data above forms a symmetric traced category.

\subsection{Topological Approach: Anchored Graphs}\label{sec:js-construction}


In 1991, Joyal and Street formalised the graphical Penrose notation as a generalised topological graph, with some added structure \cite{JS}. They went on to prove that variations of these graphs could be used to construct the free planar, symmetric, and braided categories on a monoidal signature. The usual notion of a finite topological graph is a Hausdorff space that forms the geometric realisation of a one-dimensional, finite simplicial (or equivalently, CW) complex. A generalised finite topological graph is a finite topological graph that is allowed to have some ``open ends''. That is, some edges in the graph look like the half-open interval or the open interval. We'll drop the adjective finite for the rest of this section, and assume that all of the graphs are finite.

\begin{definition}[Generalised Topological Graph] \label{def:genalised-top-graph}
  A \emph{generalised topological graph} is a pair $(G, G_0)$, where $G$ is a Hausdorff space and $G_0$ is a  finite subset of points in $G$ such that $G - G_0$ is isomorphic to a sum of open intervals $I_o := (0,1)\subseteq \mathbb R$ and copies of $S_1$. The points in $G_0$ are called \textit{vertices}. The compactification of an open interval $I_o \cong e \subseteq G - G_0$ is called an \emph{edge} $\hat e$. A copy of $S_1 \subseteq G - G_0$ is called a \emph{circle} $\hat c$. If a subgraph of $\hat G$ is an edge or a circle, we shall call it a \textit{wire} in $G$. Let $W(G)$ be the set of all wires in $G$.
\end{definition}

\begin{example}\label{ex:gen-top-graph}
  A generalised topological graph. The points in $G_0$ are marked as dots.
  \ctikzfig{generalised_top_graph}
\end{example}

This gives a bit more information than just the graph itself. The set $G_0$ is used to identify all of the ``logical'' vertices in the graph, not just the the ``topological'' vertices (i.e. those points lacking a neighbourhood of $[0,1]$). In particular, a vertex can occur along an edge. Also note that $G$ need not be compact. Let $\hat G$ be the compactification of $G$ obtained by freely adding endpoints to any open ends. Since $\hat{(-)}$ is left adjoint to the forgetful functor $U : \catCHaus \rightarrow \catHaus$, the embedding of $e$ in $G$ uniquely fixes an embedding from $\hat e$ into $\hat G$.

Note that all edges naturally embed in the compactification $\hat G \supseteq G$ obtained by adding endpoints to open edges.
\begin{center}
  \begin{tikzpicture}
    \matrix (m) [cdiag, row sep=0.5cm] {
      e &   & \hat e \\
        & G &        \\
        &   & \hat G \\
    };
    \path [arrs]
      (m-1-1) edge [right hook-latex] (m-1-3)
      (m-1-1) edge [right hook-latex] (m-2-2)
      (m-2-2) edge [right hook-latex] (m-3-3)
      (m-1-3) edge [dashed] (m-3-3);
  \end{tikzpicture}
\end{center}

String diagrams are not just graphs, but directed graphs. Thus, more data needs to be added to impose directedness to the wires in the graph. These wires can either be edges ($\cong [0,1]$) or circles ($\cong S^1$). In both cases, we have orientable manifolds, so we impose directedness by giving each of these manifolds an orientation.

\begin{definition}[Polarised Graph] \label{def:polarised-graph}
  A \emph{polarised graph} is a tuple $(G, G_0, (o_w), (p_{v,i}))$, where $o_w$ assigns each wire $w \in W(G)$ an orientation. We can therefore define an input $\hat e(0)$ and an output $\hat e(1)$ for each edge. For each vertex $v \in G_0$, $\textrm{in}(v)$ is the set of edges such that $\hat e(1) = v$ and $\textrm{out}(v)$ is the set of edges such that $\hat e(0) = v$. For all $v \in G_0$, $p_{v,0}$ is a total order on $\textrm{in}(v)$ and $p_{v,1}$ is a total order on $\textrm{out}(v)$, called a \emph{polarisation}. A polarised graph that contains no directed cycles is called \emph{progressive}.
\end{definition}

Polarised graphs come with a notion of boundary. We can furthermore put an ordering on this boundary.

\begin{definition}[Boundary of a polarised graph] \label{def:anchored-graph}
  For a polarised graph $(G, G_0, (o_w), (p_{v,i}))$, $\partial G := \hat G - G$ is a discrete space called the \emph{boundary} of $G$. Points in $\partial G$ that are the input of some edge are called \emph{inputs} of $G$, and outputs of edges in $\partial G$ are called \emph{outputs} of $G$. A polarised graph with a pair of total orders $\beta_0, \beta_1$ on its inputs and its outputs respectively is called an \emph{anchored graph}.
\end{definition}

\begin{example}
  We can make the topological graph from Example \ref{ex:gen-top-graph} into an anchored graph $\Gamma = (G, G_0, (o_w), (p_{v,i}), (\beta_i))$. We depict $\Gamma$ diagrammatically as follows. $G$ is drawn as before, but with the elements of $G_0$ shown as squares. The orientations provided by $\omega$ are depicted as arrow heads on the wires. The polarisation $p_{v,i}$ is shown by ordering from left to right the inputs and outputs to each square. The orderings $\beta_i$ are shown by ordering the inputs and outputs of the graph from left to right.
  \ctikzfig{anchored_graph}
\end{example}

This is very close to representing the diagrammatic language given at the beginning of this chapter. The only thing missing is labels on the boxes and wires.

\begin{definition}[Valuation] \label{def:valuation}
  For an anchored graph $\Gamma = (G, G_0, (o_w), (p_{v,i}), (\beta_i))$ and a monoidal signature $T = (O, M, \dom, \cod)$, a \emph{valuation} $\nu$ of $\Gamma$ is a function $\nu_0$ that assigns an element of $O$ to every edge or circle in $\Gamma$ and a function $\nu_1$ that assigns an element $m$ of $M$ to every point in $G_0$ in such a way that respects the domain on codomain of $m$.
\end{definition}

An isomorphism of anchored graphs with valuations $(\Gamma, \nu) \rightarrow (\Gamma', \nu')$ is an isomorphism of generalised topological graphs that respects orientations, polarisation, and input and output ordering, and is compatible with the valuations $\nu$ and $\nu'$.

Since an anchored graph gives a total order to inputs and outputs, we can associate input and output words to a pair $(\Gamma, \nu)$. Let $T = (O, M, \dom, \cod)$ be a monoidal signature. $\mathbb F_S(T)$ is the category whose objects are words in $w(O)$. For words $v$ and $w$, arrows are isomorphism classes of progressive anchored graphs with valuations into $T$ that have input word $v$ and output word $w$. Composition $g \circ f$ is defined by plugging the outputs of some representative of the isomorphism class of $f$ into the inputs of some representative of $g$, then taking the isomorphism class of the resultant graph.

\begin{theorem}\label{thm:js-free-smc}
  \cite{JS} $\mathbb F_S(T)$ is the free symmetric monoidal category over $T$.
\end{theorem}

In~\cite{JS}, Joyal and Street alluded to a sequel paper ``\textit{The Geometry of Tensor Calculus, II}'' (GTC-II) in which they would prove analogous results for traced and compact closed categories. For various reasons, this paper was never completed. In chapter \ref{sec:sg-free-monoidal-categories} we show that categories of string graphs can be used to construct the free symmetric traced and compact closed categories over a monoidal signature. This, along with the fact that the geometric realisation functor lifts to an equivalence of categories between the string graph-based free categories and the topologically-defined free categories~\cite{DixonKissinger2010} suffices to prove two of the missing GTC-II theorems: namely that suitably defined categories of topological graphs can be used to define (i) the free symmetric traced category and (ii) the free compact closed category over a monoidal signature.

\chapter{Algebraic Structures in Monoidal Categories}\label{ch:monoidal-algebra}

In a monoidal category $\mathcal V$, one can define various algebraic structures \textit{internal} to $\mathcal V$. A standard example is that of a monoid in $\mathcal V$.

\begin{definition}\label{def:monoid-in-category}
  For a monoidal category $\mathcal V$, a monoid in $\mathcal V$ is a triple $(X, \mu : X \otimes X \rightarrow X, \eta : I \rightarrow X)$ such that the following diagrams commute:
  \begin{equation}\label{eq:monoid-diagram}
    \begin{tikzpicture}
      \matrix (m) [cdiag] {
        X \otimes X \otimes X & X \otimes X \\
        X \otimes X & X \\
      };
      \path [arrs]
        (m-1-1) edge node {$X \otimes \mu$} (m-1-2)
        (m-2-1) edge node [swap] {$\mu$} (m-2-2)
        (m-1-1) edge node [swap] {$\mu \otimes X$} (m-2-1)
        (m-1-2) edge node {$\mu$} (m-2-2);
    \end{tikzpicture}
    \quad
    \begin{tikzpicture}
      \matrix (m) [cdiag] {
        X &  \\
        X \otimes X & X \\
      };
      \path [arrs]
        (m-2-1) edge node [swap] {$\mu$} (m-2-2)
        (m-1-1) edge node [swap] {$\eta \otimes X$} (m-2-1)
        (m-1-1) edge node {$1_X$} (m-2-2);
    \end{tikzpicture}
    \quad
    \begin{tikzpicture}
      \matrix (m) [cdiag] {
        X & X \otimes X \\
          & X \\
      };
      \path [arrs]
        (m-1-1) edge node {$X \otimes \eta$} (m-1-2)
        (m-1-2) edge node {$\mu$} (m-2-2)
        (m-1-1) edge node [swap] {$1_X$} (m-2-2);
    \end{tikzpicture}
  \end{equation}
\end{definition}

The diagrams in this definition establish associativity, left unit, and right unit.

\begin{example}
  A monoid in the category $(\catSet, \times, \{ * \})$ is the monoid in the usual sense. That is, a group without inverses. $X$ is the carrier set, $\mu : X \times X \rightarrow X$ is the associative multiplication operation, and $\nu : \{ * \} \rightarrow X$ is the map the picks out the unit, i.e. $\eta(*) = e$. If we write $\mu(x,y)$ as $x \cdot y$, then the commutative diagrams from Definition \ref{def:monoid-in-category} are equivalent to these equations:
  \begin{itemize}
    \item $(x \cdot y) \cdot z = x \cdot (y \cdot z)$,
    \item $e \cdot x = x$, and $x \cdot e = x$.
  \end{itemize}
\end{example}

\begin{examples}
  A monoid in $(\catVect_K, \otimes, K)$ is an associative, unital $K$-algebra. Since $\mu : V \otimes V \rightarrow V$ is a linear map, then due to the universal property of the tensor product, $\mu$ uniquely determines a \textit{bilinear} map $(-\cdot-) : V \times V \rightarrow V$. Again, associativity and unitality come from Diagrams (\ref{eq:monoid-diagram}). Similarly, a monoid in $(\catAb, \otimes, \mathbb N)$, the category of Abelian groups and group homomorphisms, is a ring.
\end{examples}

We can also express a monoid graphically, as a triple $(X,\, \mult : X \otimes X \rightarrow X,\, \unit : I \rightarrow X)$. Associativity, left unit, and right unit become the following equations.
\begin{equation}\label{eq:monoid-axioms}
\beginpgfgraphicnamed{monoid_axioms}
\InputIfFileExists{monoid_axioms.tikz}{}{\input{./figures/monoid_axioms.tikz}}
\endpgfgraphicnamed
\end{equation}

By simply turning everything upside-down, we define a comonoid in a monoidal category.

\begin{definition}\label{def:comonoid}
  For a monoidal category $\mathcal V$, a comonoid in $\mathcal V$ is a triple $(X, \delta : X \rightarrow X \otimes X, \epsilon : X \rightarrow I)$ such that the following diagrams commute:
  \begin{equation}\label{eq:comonoid-diagram}
    \begin{tikzpicture}
      \matrix (m) [cdiag] {
        X & X \otimes X \\
        X \otimes X & X \otimes X \otimes X \\
      };
      \path [arrs]
        (m-1-1) edge node {$\delta$} (m-1-2)
        (m-2-1) edge node [swap] {$X \otimes \delta$} (m-2-2)
        (m-1-1) edge node [swap] {$\delta$} (m-2-1)
        (m-1-2) edge node {$\delta \otimes X$} (m-2-2);
    \end{tikzpicture}
    \quad
    \begin{tikzpicture}
      \matrix (m) [cdiag] {
        X &  \\
        X \otimes X & X \\
      };
      \path [arrs]
        (m-2-1) edge node [swap] {$\epsilon \otimes X$} (m-2-2)
        (m-1-1) edge node [swap] {$\delta$} (m-2-1)
        (m-1-1) edge node {$1_X$} (m-2-2);
    \end{tikzpicture}
    \quad
    \begin{tikzpicture}
      \matrix (m) [cdiag] {
        X & X \otimes X \\
          & X \\
      };
      \path [arrs]
        (m-1-1) edge node {$\delta$} (m-1-2)
        (m-1-2) edge node {$X \otimes \epsilon$} (m-2-2)
        (m-1-1) edge node [swap] {$1_X$} (m-2-2);
    \end{tikzpicture}
  \end{equation}
\end{definition}

Graphically, these axioms are:
\ctikzfig{comonoid_axioms}

A \textit{commutative monoid} is a monoid satisfying the following equation:
\ctikzfig{commutative_monoid}

Similarly, a \textit{cocommutative comonoid} is a comonoid satisfying the same equation, but upside-down.
\ctikzfig{cocommutative_comonoid}

\section{Bi-algebras and Hopf Algebras}\label{sec:bialgebras}

Suppose we wished to define a \textit{group} in a monoidal category. We begin by looking at how groups are defined in $\catSet$. A group is a monoid, with an additional unary operation $(-)^{-1}$ such that $x \cdot x^{-1} = e = x^{-1} \cdot x$. How can we express these equations as commutative diagrams? This poses a challenge, because unlike in the monoid axioms, the free variable $x$ in the inverse axiom features \textit{twice}. Luckily, $\catSet$ is not just a monoidal category, but a cartesian monoidal category. A cartesian monoidal category is a category with finite products, where the monoidal product is defined using the categorical product and the monoidal unit is the terminal object. In the case of $\catSet$, these are the cartesian product on the one-element set.

In particular, cartesian monoidal categories always have a diagonal map for every object induced by the universal property of the product.

\begin{center}
  \begin{tikzpicture}
    \matrix (m) [cdiag] {
        &      X     &   \\
      X & X \times X & X \\
    };
    \path [arrs]
      (m-1-2) edge node [swap] {$1_X$} (m-2-1)
      (m-1-2) edge node {$1_X$} (m-2-3)
      (m-1-2) edge [dashed] node {$\Delta$} (m-2-2)
      (m-2-2) edge node {$\pi_1$} (m-2-1)
      (m-2-2) edge node [swap] {$\pi_2$} (m-2-3);
  \end{tikzpicture}
\end{center}

The diagonal map can be used to ``copy'' the free variable $x \in X$.

\begin{definition}
  A \textit{group object} in a cartesian monoidal category $(\mathcal C, \times, \top)$, is a tuple $(X, \mu, \eta, \iota)$ such that $(X, \mu, \eta)$ is a monoid and the following diagram commutes:
  \begin{center}
    \begin{tikzpicture}
      \matrix (m) [cdiag] {
         & X \times X  &      & X \times X &   \\
       X &             & \top &            & X \\
         & X \times X  &      & X \times X &   \\
      };
      \path [arrs]
        (m-2-1) edge [dashed] node {$!$} (m-2-3)
        (m-2-3) edge node {$\eta$} (m-2-5)
        (m-2-1) edge node {$\Delta$} (m-1-2)
        (m-1-2) edge node {$X \times \iota$} (m-1-4)
        (m-1-4) edge node {$\mu$} (m-2-5)
        (m-2-1) edge node [swap] {$\Delta$} (m-3-2)
        (m-3-2) edge node [swap] {$\iota \times X$} (m-3-4)
        (m-3-4) edge node [swap] {$\mu$} (m-2-5);
    \end{tikzpicture}
  \end{center}
  where $\Delta : X \rightarrow X \times X$ is the diagonal map and $! : X \rightarrow \top$ is the unique map from $X$ to $\top$, the terminal object.
\end{definition}

The picture becomes clearer by looking at the graphical versions of these identities. Note how the terminal map is used to ``delete'' the free variable, just as the $\Delta$ map is used to copy it.
\ctikzfig{group_antipode}

We have defined a group object in a cartesian monoidal category, but we have not quite answered the original question. Is there a notion of a group in a (non-cartesian) monoidal category? The maps $\Delta$ and $!$ in a monoidal category do not come for free. Instead, we require that they be part of the \textit{structure} of the the ``group'' object. To justify what this structure should be, we highlight some properties of $(\Delta, !)$.

\begin{proposition}\label{prop:cartesian-comonoid}
  For any object $X$ in a cartesian monoidal category, the triple $(X, \Delta, !)$ forms a comonoid. Furthermore, the following diagrams commute for all $X$, $Y$, and $f$.
  \begin{center}
    \begin{tikzpicture}
      \matrix (m) [cdiag] {
        X & Y \\
        X \times X & Y \times Y \\
      };
      \path [arrs]
        (m-1-1) edge node {$f$} (m-1-2)
        (m-2-1) edge node [swap] {$f \times f$} (m-2-2)
        (m-1-1) edge node [swap] {$\Delta_X$} (m-2-1)
        (m-1-2) edge node {$\Delta_Y$} (m-2-2);
    \end{tikzpicture}
    \qquad
    \begin{tikzpicture}
      \matrix (m) [cdiag] {
        X \times Y  & X \times Y \times X \times Y \\
        X \times X \times Y \times Y &       .     \\
      };
      \path [arrs]
        (m-1-1) edge node {$\Delta_{X,Y}$} (m-1-2)
        (m-2-1) edge node [swap] {$X \times \sigma_{X,Y} \times Y$} (m-1-2)
        (m-1-1) edge node [swap] {$\Delta_X \times \Delta_Y$} (m-2-1);
    \end{tikzpicture}
  \end{center}
\end{proposition}

\begin{proof}
  We write $\langle g,h \rangle : A \rightarrow B \times C$ for the unique induced map to the product of the pair of maps $g : A \rightarrow B$, $h : A \rightarrow C$. We can prove associativity using the definition of $\Delta$ and some well-known properties of $\langle -, -\rangle$.
  \begin{align*}
    (X \times \Delta) \circ \Delta
    & = \langle \pi_1, \Delta \circ \pi_2 \rangle \circ \Delta
      = \langle \pi_1 \circ \Delta, \Delta \circ \pi_2 \circ \Delta \rangle
      = \langle 1_X, \Delta \circ 1_X \rangle \\
    & = \langle 1_X, \langle 1_X, 1_X \rangle \rangle
      = \langle \langle 1_X, 1_X \rangle, 1_X \rangle
      = (\Delta \times X) \circ \Delta
  \end{align*}
  
  Right unit (and similarly, left unit) follows from terminality of $!$.
  \begin{align*}
    (X \times !_X) \circ \Delta
    & = \langle \pi_1, !_X \circ \pi_2 \rangle \circ \Delta
      = \langle \pi_1, !_{X\times X} \rangle \circ \Delta \\
    & = \langle \pi_1 \circ \Delta, !_{X\times X} \circ \Delta \rangle
      = \langle 1_X, !_X \rangle
      = 1_X
  \end{align*}
  
  The copying property is proved as follows:
  \begin{align*}
    \Delta \circ f
    & = \langle 1_X, 1_X \rangle \circ f
      = \langle f, f \rangle
      = \langle f \circ \pi_1 \circ \Delta, f \circ \pi_2 \circ \Delta \rangle \\
    & = \langle f \circ \pi_1, f \circ \pi_2 \rangle \circ \Delta
      = (f \times f) \circ \Delta
  \end{align*}
  
  For the final property, we establish that $(\Delta_X \times \Delta_Y) \circ (X \times \sigma_{X,Y} \times Y) = \langle 1_{X,Y}, 1_{X,Y} \rangle$. For simplicity, let $\pi_1, \pi_2$ be projections of the two-part products $(\cdot \times \cdot) \times (\cdot \times \cdot)$, and let $\pi'_1,\pi'_2,\pi'_3,\pi'_4$ be projections of the four-part products $(\cdot \times \cdot \times \cdot \times \cdot)$.
  \begin{center}
    \begin{tikzpicture}
      \matrix (m) [cdiag,column sep=2cm] {
        & X \times Y & X \times Y \\
        X \times Y &
        X \times X \times Y \times Y & 
        X \times Y \times X \times Y \\
        & X \times Y & X \times Y \\
      };
      \path [arrs]
        (m-2-1) edge node {$\Delta_X \times \Delta_Y$} (m-2-2)
        (m-2-2) edge node {$X \otimes \sigma \otimes Y$} (m-2-3)
        (m-2-1) edge [bend left=20] node {$1_{X,Y}$} (m-1-2)
        (m-1-2) edge node {$1_{X,Y}$} (m-1-3)
        (m-2-2) edge node [swap] {$\pi_1 \times \pi_1 = \langle \pi'_1, \pi'_3 \rangle$} (m-1-2)
        (m-2-3) edge node [swap] {$\langle \pi'_1, \pi'_2 \rangle = \pi_1$} (m-1-3)
        (m-2-1) edge [bend right=20] node [swap] {$1_{X,Y}$} (m-3-2)
        (m-3-2) edge node [swap] {$1_{X,Y}$} (m-3-3)
        (m-2-2) edge node {$\pi_2 \times \pi_2 = \langle \pi'_2, \pi'_4 \rangle$} (m-3-2)
        (m-2-3) edge node {$\langle \pi'_3, \pi'_4 \rangle = \pi_2$} (m-3-3);
    \end{tikzpicture}
  \end{center}
\end{proof}

We already saw graphical versions of the comonoid axioms after \ref{def:comonoid}. The remaining two diagrams from Proposition \ref{prop:cartesian-comonoid} are:
\ctikzfig{diagonal_axioms}

In particular, Proposition \ref{prop:cartesian-comonoid} implies that $\Delta$ copies $\mu$.
\begin{equation}\label{eq:cart-mu-copy}
\beginpgfgraphicnamed{cartesian_mu_copy}
\InputIfFileExists{cartesian_mu_copy.tikz}{}{\input{./figures/cartesian_mu_copy.tikz}}
\endpgfgraphicnamed
\end{equation}

$\Delta$ also copies the unit $\eta$. By terminality, $\mu$ also ``co-copies'' $!$.
\begin{equation}\label{eq:cart-unit-copy}
\beginpgfgraphicnamed{cartesian_unit_copy}
\InputIfFileExists{cartesian_unit_copy.tikz}{}{\input{./figures/cartesian_unit_copy.tikz}}
\endpgfgraphicnamed
\end{equation}

In a monoidal category, a monoid paired with a comonoid that behaves like copy and delete operations is called a \textit{bialgebra}.

\begin{definition}\label{def:bialg}
  In a monoidal category $\mathcal V$, a \textit{bialgebra} is a tuple $(X, \mu, \eta, \delta, \epsilon)$ where $(X, \mu, \eta)$ is a monoid, $(X, \delta, \epsilon)$ is a comonoid and the following diagrams commute.
  \begin{center}
    \begin{tikzpicture}
      \matrix (m) [cdiag] {
        X \otimes X                     & X & X \otimes X                     \\
        X \otimes X \otimes X \otimes X &   & X \otimes X \otimes X \otimes X \\
      };
      \path [arrs]
        (m-1-1) edge node {$\mu$} (m-1-2)
        (m-1-2) edge node {$\delta$} (m-1-3)
        (m-1-1) edge node [swap] {$\delta \otimes \delta$} (m-2-1)
        (m-2-1) edge node {$X \otimes \sigma \otimes X$} (m-2-3)
        (m-2-3) edge node [swap] {$\mu \otimes \mu$} (m-1-3);
    \end{tikzpicture}
    
    \begin{tikzpicture}
      \matrix (m) [cdiag] {
        I & X           \\
          & X \otimes X \\
      };
      \path [arrs]
        (m-1-1) edge node {$\eta$} (m-1-2)
        (m-1-2) edge node {$\delta$} (m-2-2)
        (m-1-1) edge node [swap] {$\eta \otimes \eta$} (m-2-2);
    \end{tikzpicture}
    \quad
    \begin{tikzpicture}
      \matrix (m) [cdiag] {
        X \otimes X &   \\
        X           & I \\
      };
      \path [arrs]
        (m-1-1) edge node [swap] {$\mu$} (m-2-1)
        (m-2-1) edge node [swap] {$\epsilon$} (m-2-2)
        (m-1-1) edge node {$\epsilon \otimes \epsilon$} (m-2-2);
    \end{tikzpicture}
  \end{center}
\end{definition}

These conditions are the same as equations (\ref{eq:cart-mu-copy}) and (\ref{eq:cart-unit-copy}), but with $(X,\Delta,!)$ replaced by an \textit{arbitrary} comonoid.
\ctikzfig{bialg_axioms}

In the general case, it is worth noting that the comonoid within a bialgebra acts like a copy and deletion operation \textit{relative to the associated monoid}, rather than globally. This distinction applied to the category of Hilbert spaces has deep connections to the \textit{no-cloning theorem} for quantum mechanics. See, for example \cite{AbramskyNoCloning}.
In a cartesian category, \textit{any} monoid can be made into a bialgebra using the comonoid $(X, \Delta, !)$.

In the general case of monoidal categories, the inverse map is replaced with an arbitrary map, called an \textit{antipode}. The analogue to a group object is a \textit{Hopf algebra}.

\begin{definition}\label{def:hopf-algebra}
  A \textit{Hopf algebra} is a tuple $(X, \mu, \eta, \delta, \epsilon, \iota)$ such that $(X, \mu, \eta, \delta, \epsilon)$ forms a bialgebra and the following diagram commutes:
  \begin{center}
    \begin{tikzpicture}
      \matrix (m) [cdiag,column sep=7mm] {
           & X \otimes X &   & X \otimes X &   \\
         X &             & I &             & X \\
           & X \otimes X &   & X \otimes X &   \\
      };
      \path [arrs]
        (m-2-1) edge node {$\epsilon$} (m-2-3)
        (m-2-3) edge node {$\eta$} (m-2-5)
        (m-2-1) edge node {$\delta$} (m-1-2)
        (m-1-2) edge node {$X \otimes \iota$} (m-1-4)
        (m-1-4) edge node {$\mu$} (m-2-5)
        (m-2-1) edge node [swap] {$\delta$} (m-3-2)
        (m-3-2) edge node [swap] {$\iota \otimes X$} (m-3-4)
        (m-3-4) edge node [swap] {$\mu$} (m-2-5);
    \end{tikzpicture}
  \end{center}
\end{definition}

\begin{example}
  For any group $G$, we can form the group algebra $(K[G],\mu,\eta)$ as follows. Let $K[G]$ be a vector space in $\catVect_K$ spanned by the elements $g \in G$. $\mu$ is then defined on basis vectors by the group multiplication: $\mu(g \otimes h) = gh$. $\eta(1) = e$, the group identity. This clearly forms a monoid in $\catVect_K$. We can turn this monoid into a bialgebra by adding maps that ``copy'' and ``delete'' the basis of group elements. We can make it into a Hopf algebra by adding a linear map that sends elements of the group basis to their inverse.
  \begin{align*}
    \delta   & :: g \mapsto g \otimes g \\
    \epsilon & :: g \mapsto 1           \\
    \iota    & :: g \mapsto g^{-1}     
  \end{align*}
  The Hopf algebra $(K[G], \mu, \eta, \delta, \epsilon, \iota)$ captures all of the structure of the group algebra $K[G]$ without relying on ``global'' copy an deletion operations, which are not present in the (non-cartesian) monoidal category $\catVect_K$.
\end{example}

\begin{examples}
  Here are a few more examples of Hopf algebras.
  \begin{itemize}
    \item Any group object in a cartesian monoidal category is automatically a Hopf algebra, with $\delta = \Delta$ and $\epsilon = !$.
    \item The universal enveloping algebra of a Lie algebra $U(\mathfrak g)$ has a natural Hopf algebra structure, given by the unique extension of the following maps on $\mathfrak g$ to $U(\mathfrak g)$:
    \begin{center}
      $\delta(x) = (x \otimes 1) + (1 \otimes x)$\quad
      $\epsilon(x) = 0$\quad
      $\iota(x) = -1$
    \end{center}
    \item More generally, quantum groups of the form $U_q(\mathfrak g)$ for $q \in \mathbb C$ carry a Hopf algebra structure.
  \end{itemize}
\end{examples}

\section{Frobenius Algebras}\label{sec:frobenius-algebras}

Associative algebras $(A,\mu,\eta)$ in a compact closed category of vector spaces always come with two canonical right representations: one over $A$ and one over the dual space $A^*$. Frobenius algebras are associative algebras that are \textit{self-dual}, i.e. these two representations are isomorphic. Frobenius algebras can be defined in any monoidal category. In order to talk about representations abstractly, it is convenient to use the (equivalent) language of modules.

\begin{definition}\label{def:module}
  For a monoid $\mathcal A = (A,\mu,\eta)$, a right $\mathcal A$-module $(M,m)$ is a morphism $m : A \otimes M \rightarrow M$ such that the following diagrams commute.
  \begin{center}
    \begin{tikzpicture}
      \matrix (m) [cdiag] {
        A \otimes A \otimes M & A \otimes M \\
        A \otimes M           & M           \\
      };
      \path [arrs]
        (m-1-1) edge node {$A \otimes m$} (m-1-2)
        (m-2-1) edge node [swap] {$m$} (m-2-2)
        (m-1-1) edge node [swap] {$\mu \otimes M$} (m-2-1)
        (m-1-2) edge node {$m$} (m-2-2);
    \end{tikzpicture}
    \qquad
    \begin{tikzpicture}
      \matrix (m) [cdiag] {
        M & A \otimes M \\
          & M           \\
      };
      \path [arrs]
        (m-1-1) edge node {$\eta \otimes M$} (m-1-2)
        (m-1-2) edge node {$m$} (m-2-2)
        (m-1-1) edge node [swap] {$1$} (m-2-2);
    \end{tikzpicture}
  \end{center}
\end{definition}

Left $\mathcal A$-modules are defined analogously, and for commutative monoids, the two concepts are equivalent. Graphically, the right module equations look like asymmetric versions of associativity and unit laws.
\ctikzfig{module_axioms}

\begin{definition}
  For a monoid $\mathcal A = (A,\mu,\eta)$, and right $\mathcal A$-modules $(M,m)$ and $(N,n)$, an $\mathcal A$-module homomorphism is a morphism $\phi : M \rightarrow N$ such that:
  \begin{center}
    \begin{tikzpicture}
      \matrix (m) [cdiag] {
        A \otimes M & A \otimes N \\
        M           & N           \\
      };
      \path [arrs]
        (m-1-1) edge node {$A \otimes \phi$} (m-1-2)
        (m-1-1) edge node [swap] {$m$} (m-2-1)
        (m-1-2) edge node {$n$} (m-2-2)
        (m-2-1) edge node [swap] {$\phi$} (m-2-2);
    \end{tikzpicture}
  \end{center}
\end{definition}

Module homomorphisms pass through the module map.
\begin{equation}\label{eq:module-hm}
\beginpgfgraphicnamed{module_hm}
\InputIfFileExists{module_hm.tikz}{}{\input{./figures/module_hm.tikz}}
\endpgfgraphicnamed
\end{equation}

Monoids in a monoidal category have a canonical right $\mathcal A$-module, namely the monoid itself.
\ctikzfig{regular_module}

This is called the \textit{regular right $\mathcal A$-module}. In a compact closed category, a monoid comes with a canonical right module over its dual object $\tilde\mu : A \otimes A^* \rightarrow A^*$.
\[ \tilde\mu := (A^* \otimes d) \circ (A^* \otimes \mu \otimes A^*) \circ (e \otimes A \otimes A^*) \]
We shall call this the \textit{dual right $\mathcal A$-module}. Graphically:
\ctikzfig{dual_module}

\begin{proposition}\label{prop:dual-module}
  $(A^*, \tilde\mu)$ is a right $\mathcal A$-module.
\end{proposition}

\begin{proof}
  We can show the two module conditions are satisfied graphically.
  \begin{align*}
    \textrm{Associativity:\quad} & %
\beginpgfgraphicnamed{dual_module_assoc}
\InputIfFileExists{dual_module_assoc.tikz}{}{\input{./figures/dual_module_assoc.tikz}}
\endpgfgraphicnamed \\
    \textrm{Unit:\quad}         & %
\beginpgfgraphicnamed{dual_module_unit}
\InputIfFileExists{dual_module_unit.tikz}{}{\input{./figures/dual_module_unit.tikz}}
\endpgfgraphicnamed 
  \end{align*}
\end{proof}

Frobenius is generally credited with being the first to study the class of algebras whose regular and dual modules are isomorphic. This is a particularly important aspect of many associative algebras, including algebras of the form $K[G]$ for some finite group $G$. The formal definition is due to Brauer and Nesbitt \cite{BrauerNesbitt1937}. An abstract version is provided here.

\begin{definition}
  A Frobenius algebra in a compact closed category is a monoid $(A, \mu, \eta)$ equipped with a module isomorphism $\chi : (A, \mu) \rightarrow (A^*, \tilde\mu)$.
\end{definition}

The existence of a Frobenius algebra on an object $A$ implies that $A$ is \textit{self-dual}, and furthermore, the monoid structure $\mu,\eta$ is compatible with that self-duality. Graphically, we represent $\chi$ and its inverse $\chi^{-1}$ as follows:
\ctikzfig{frobenius_module_iso}

Graphically, the axioms of a Frobenius algebra are the monoid axioms given by (\ref{eq:monoid-axioms}) and three additional identities:
\begin{equation}\label{eq:frobenius-module-ids}
\beginpgfgraphicnamed{frobenius_module_ids}
\InputIfFileExists{frobenius_module_ids.tikz}{}{\input{./figures/frobenius_module_ids.tikz}}
\endpgfgraphicnamed
\end{equation}

\begin{proposition}\label{prop:group-algebra-frobenius}
  Let $G$ be a finite group. Then, the associative algebra $K[G]$ is Frobenius.
\end{proposition}

\begin{proof}
  Let $(K[G],\mult,\unit)$ be the associative algebra of the finite group $G$. For $g \in G$, let $\hat g$ be the unique map defined on $h \in G$ as:
  \[ \hat g(h) =
     \begin{cases}
      1 & \textrm{ if } h = g \\
      0 & \textrm{ otherwise}
     \end{cases}
   \]
  For $\chi$, let $\chi(g) = \widehat{g^{-1}}$. Since $K[G]$ is finite-dimensional, this extends to an isomorphism $K[G] \cong K[G]^*$. It only remains to show that $\chi$ is a module homomorphism. It suffices to show that the two sides of equation (\ref{eq:module-hm}) agree on all basis elements $g_i \in G$.
  \ctikzfig{group_alg_frob_pf}
  The proof is completed by noting that $g_3 g_1 = g_2^{-1}$ if and only if $g_1 g_2 = g_3^{-1}$.
\end{proof}

There are actually many equivalent definitions for a Frobenius algebra. Two additional definitions replace the isomorphism condition given above with a non-degeneracy condition.

\begin{definition}
  A map of the form $\Phi : A \otimes A \rightarrow I$ is said to be non-degenerate if the following map is an isomorphism:
  \ctikzfig{left_transpose}
  
  Equivalently, there exists a unique map $\Phi' : I \rightarrow A \otimes A$ such that:
  \ctikzfig{non_degen}
\end{definition}

\begin{theorem}\label{thm:frobenius-defs}
  The following are equivalent.
  \begin{enumerate}
    \item $(A,\mu,\eta,\chi)$ is a Frobenius algebra,
    \item $(A, \mu, \eta)$ is a monoid with non-degenerate map $\Phi : A \otimes A \rightarrow I$ that is associative with respect to the monoid:
    \ctikzfig{frobenius_cup_def}
    \item $(A, \mu, \eta)$ is a monoid with a map $\epsilon : A \rightarrow I$ such that the following is non-degenerate:
    \ctikzfig{frobenius_counit_def}
    \item and $(A, \mu, \eta)$ is a monoid, $(A, \delta, \epsilon)$ is a comonoid satisfying the \textit{Frobenius identity}:
    \ctikzfig{frobenius_comonoid_def}
  \end{enumerate}
\end{theorem}

\begin{proof}
  ($1 \Rightarrow 2$) Define $\Phi$ in terms of $\chi$ as:
  \ctikzfig{frobenius_def_pf0}
  The associativity of $\Phi$ follows from the definition of dual module and the module homomorphism identity in equation (\ref{eq:frobenius-module-ids}).
  \ctikzfig{frobenius_def_pf1}
  
  ($2 \Rightarrow 3$) Let $\epsilon$ be defined as:
  \ctikzfig{frobenius_def_pf2}
  It can then easily be shown that $\epsilon \circ \mu = \Phi$, which is non-degenerate.
  
  ($3 \Rightarrow 4$) Let $\Phi' : I \rightarrow A \otimes A$ be the unique map such that:
  \ctikzfig{frobenius_def_pf4}
  
  We can show this induced cap commutes horizontally with the multiplication:
  \ctikzfig{frobenius_def_pf5}
  
  Define $\delta$ as follows:
  \begin{equation}\label{eq:delta-def}
\beginpgfgraphicnamed{frobenius_def_pf6}
\InputIfFileExists{frobenius_def_pf6.tikz}{}{\input{./figures/frobenius_def_pf6.tikz}}
\endpgfgraphicnamed
  \end{equation}
  
  Then it follows from the monoid identities on $(A, \mu, \eta)$ that $(A,\delta, \epsilon)$ is a comonoid and that the Frobenius identity holds.
  
  ($4 \Rightarrow 1$) Finally, let the isomorphism $\chi$ and its inverse be defined as follows:
  \ctikzfig{frobenius_def_pf7}
  
  It follows from the Frobenius identity and unit laws that this is indeed an isomorphism.
  \ctikzfig{frobenius_def_pf8}
  
  The only thing the remains to be checked is that this is an isomorphism of \textit{modules}:
  \ctikzfig{frobenius_def_pf9}
\end{proof}

The equation from definition $4$ first appeared in a paper by Carboni and Walters~\cite{CarboniWalters1987} in the context of categories of relations, though it was only later discovered to be a characterising equation for Frobenius algebras. This presentation is particularly interesting because, like in the case of bialgebras, it consists of (1) a monoid, (2) a comonoid, and (3) a rule for how they interact, called a \textit{distributive law}. In definition $4$, a Frobenius algebra consists of a monoid and a comonoid. Therefore a Frobenius algebra could be \textit{commutative} or \textit{cocommutative}. While these might seem like distinct conditions, it turns out they are the same.

\begin{theorem}\label{thm:comm-cocomm}
  For a Frobenius algebra $(A, \mu, \eta, \delta, \epsilon)$, $(A, \mu, \eta)$ is a commutative monoid if and only if $(A, \delta, \epsilon)$ is a cocommutative comonoid.
\end{theorem}

\begin{proof}
  First, we use Equation (\ref{eq:delta-def}) to derive another form for $\delta$ in terms of $\mu$.
  \ctikzfig{frobenius_delta_alt_form}
  
  Cocommutative then trivially follows from commutativity:
  \ctikzfig{frobenius_cocomm_pf}
  
  The opposite implication is the same proof, upside-down.
\end{proof}

It is worth noting that this is \textit{not} the case for bialgebras or Hopf algebras. As a simple counter-example, consider any non-commutative monoid in a cartesian category. The cocommutative monoid $(\Delta, !)$ automatically makes this into a bialgebra.

\begin{example}
  Revisiting the group algebra example, we can define the rest of the Frobenius algebra structure as follows.
  \begin{align*}
    \delta   & :: g \mapsto \sum_{g_1,g_2 \in G, g_1 g_2 = g} g_1 \otimes g_2 \\
    \epsilon & = \hat e                                                      
  \end{align*}
  The comultiplication can be though of an averaging operation. It takes a group element to a sum over all of its possible factorisations.
\end{example}

\begin{example}
  Any compact closed category automatically has a Frobenius algebra on $A^* \otimes A$ given by:
  \ctikzfig{compact_frob}
\end{example}

In a $\dagger$-monoidal category, we can introduce the notion of a $\dagger$-Frobenius algebra, where the monoid structure is just the dagger of the comonoid structure.

\begin{definition}
  A Frobenius algebra $(A, \mu, \eta, \delta, \epsilon)$ is called a $\dagger$-Frobenius algebra if $\mu = \delta^\dagger$ and $\eta = \epsilon^\dagger$.
\end{definition}

As we have already seen in the proof of Theorem \ref{thm:frobenius-defs}, the Frobenius identity implies that the Frobenius ``cap'' and ``cup'' maps form a compact structure.
\ctikzfig{frobenius_compact}

In fact, the Frobenius algebra induces a \textit{self-dual} compact structure, i.e. $A^* = A$. Using this compact structure, we can define a transposition operation $(-)^\blacktranspose$ relative to a particular Frobenius algebra.
\begin{equation}\label{eq:dot-transpose}
\beginpgfgraphicnamed{dot_transpose}
\InputIfFileExists{dot_transpose.tikz}{}{\input{./figures/dot_transpose.tikz}}
\endpgfgraphicnamed
\end{equation}

The situation here is a bit delicate. Whereas in most categories there is a canonical choice of a compact structure when $A$ and $A^*$ are distinct objects (e.g. a vector space and its dual space), defining a compact structure when $A$ and $A^*$ are the same object often involves a choice. Different Frobenius algebras defined on a single object $A$ will often define \textit{different} compact structures. It is a well-known fact that for finite-dimensional vector spaces, there is no canonical isomorphism connecting a space to its dual space. Picking a self-dual compact structure corresponds to choosing a \textit{particular} isomorphism $A \cong A^*$.

The distinction between compact structures, where duals are defined up to isomorphism, and (non-canonical) self-dual compact structures is worth bearing in mind, particularly with regards to the Frobenius algebras defined in Part \ref{part:entanglement}. However, the situation is simpler when it comes to the trace operation. In Equation \ref{eq:compact-partial-trace}, we showed that any compact structure can be used to construct a trace operation. It turns out that the compact structures generated by Frobenius algebras (and more generally, \textit{any} compact structures) always generate \textit{the same} trace.
\ctikzfig{frobenius_trace}

In the case of commutative Frobenius algebras, we usually represent the trace without a twist on the cap:
\ctikzfig{frobenius_trace_notwist}

Often it is useful to consider maps that can pass freely through the Frobenius algebra structure. That is, we consider maps that are module endomorphisms of the regular left and right $\mathcal A$ modules. By analogy to phase gates in quantum circuits, we call such maps \textit{phases}.

\begin{definition}
  For a Frobenius algebra $\mathcal A = (A,\mu,\eta,\delta,\epsilon)$, a map $\phi : A \rightarrow A$ is called a \textit{phase} for $\mathcal A$ if:
  \ctikzfig{frobenius_phase}
\end{definition}

We shall see the connection between abstract phases and quantum phase gates in sections \ref{sec:complementary-obs} and \ref{sec:zx-calculus}.

\begin{proposition}
  A phase for a Frobenius algebra $\mathcal A$ is also a comodule endomorphism for $\delta$, considered as a left and right comodule.
  \ctikzfig{frobenius_phase_comod}
\end{proposition}

\begin{proof}
  The proof follows straightforwardly from the Frobenius identities.
  \ctikzfig{frobenius_phase_comod_pf}
  
  The left comodule identity follows similarly.
\end{proof}

There is a canonical phase associated with any Frobenius algebra called its \textit{loop map}.
\ctikzfig{frobenius_loop_map}

\begin{proposition}\label{prop:loop-map-endo}
  For a Frobenius algebra $\mathcal A = (A, \mu, \eta, \delta, \epsilon)$, the loop map is a phase for $\mathcal A$. I.e.:
  \ctikzfig{frobenius_loop_module_endo}
\end{proposition}

\begin{proof}
  The proof follows from the Frobenius identity.
  \ctikzfig{frobenius_loop_module_endo_pf}
\end{proof}

For commutative Frobenius algebras, phases can always be expressed as a left (or right) multiplication by a point.

\begin{proposition}\label{prop:frobenius-phase-form}
  All phases for a given commutative Frobenius algebra are of the following form:
  \ctikzfig{phase_point_form}
\end{proposition}

\begin{proof}
  A map of the above form is clearly a phase. Let $\varphi : A \rightarrow A$ be an arbitrary phase map. Then:
  \ctikzfig{phase_point_form_pf}
\end{proof}

\subsection{Normal Form for Frobenius Algebras}\label{sec:nf-frobenius}

We shall primarily be interested in \textit{commutative Frobenius algebras} (CFAs). The primary purpose of this section is to show that CFAs have particularly nice normal forms. To do that, we will state and prove the so-called \textit{spider theorem}. First, we introduce the notion of a \textit{spider}.

\begin{definition}\label{def:spider}
  For a commutative Frobenius algebra $\mathcal A$, a \textit{spider} is a map $S_m^n$ defined as follows.
  \ctikzfig{spider_def}

  We represent the maps $S_m^n$ as single dots, with $m$ in-edges and $n$ out-edges.
  \ctikzfig{spider_def2}
\end{definition}

Since a CFA is a commutative, associative, cocommutative, and coassociative, we can freely interchange edges.

\begin{examples}
  Some spiders of various arities:
  \ctikzfig{spider_examples}
\end{examples}

As a minor technical point, it is occasionally necessary to distinguish morphisms in a monoidal category from their actual representations as string diagrams. For this section, we shall use $D, D', \ldots$ to represent formal string diagrams and $|D|, |D'|, \ldots$ to represent the associated morphisms in a monoidal category.

\begin{definition}\label{def:frobenius-map}
  For a commutative Frobenius algebra $\mathcal A = (A, \mu, \eta, \delta, \epsilon)$, an $\mathcal A$-diagram $D$ is a string diagram whose vertices are all labeled $\mu$, $\eta$, $\delta$, or $\epsilon$. An \textit{$\mathcal A$-tree} is an $\mathcal A$-diagram that contains no cycles. An $\mathcal A$-tree that is formally equal to the diagram of $S_m^n$ for some $m, n$ is said to be in \textit{spider-normal form}.
\end{definition}

\begin{lemma}\label{lem:cfa-tree-nf}
  Suppose $D$ is a connected $\mathcal A$-tree with $m$ in-edges and $n$ out-edges. Then $|D| = S_m^n$. I.e. the value of an $\mathcal A$-tree is uniquely determined by its input and output arities.
\end{lemma}

\begin{proof}
  We proceed by induction on the number of vertices in $D$. Note that any one-vertex $\mathcal A$-tree is trivially equal to a spider $S_1^0$, $S_0^1$, $S_2^1$, or $S_1^2$. Thus, assume for an $\mathcal A$-tree containing $N$ vertices, $|D| = S_m^n$. We show that for an $\mathcal A$-tree $D'$ containing $N+1$ vertices, $|D'| = S_{m'}^{n'}$. Since a spider is commutative and cocommutative, we can assume without loss of generality that the additional vertex is composed on the rightmost leg above or below the spider. These cases follow trivially from associativity, unitality, and the definition of spiders:
  \ctikzfig{spider_tree_easy_cases}
  
  The only cases remaining are post-composition by $\mu$ and pre-composition by $\delta$. First, consider pre-composition with $\delta$:
  \ctikzfig{spider_tree_hard_case1}
  
  If $m = 0$, this is already in spider-normal form. So, consider the case where $m \geq 1$. By definition of $S_m^n$, we can pull out a multiplication. Applying the Frobenius identity and the induction hypothesis:
  \ctikzfig{spider_tree_hard_case1_pf}
  
  We complete the proof by applying the same method upside-down for post-composition with $\mu$.
  \ctikzfig{spider_tree_hard_case2_pf}
\end{proof}

Because of this lemma, we could have equivalently defined $S_m^n$ in Definition \ref{def:spider} as $S_m^n = |D|$ for \textit{any} $\mathcal A$-tree $D$ with $m$ inputs and $n$ outputs. With the help of a normal form result for $\mathcal A$-trees and Proposition \ref{prop:loop-map-endo}, we are ready to state a general normal form result for commutative Frobenius algebras.

\begin{theorem}\label{thm:cfa-nf}
  For a commutative Frobenius algebra $\mathcal A$, any connected $\mathcal A$-diagram is equivalent to one of this form:
  \ctikzfig{frobenius_nf}
\end{theorem}

\begin{proof}
  We first prove a small identity relating to traces:
  \begin{equation}\label{eq:trace-to-loop}
\beginpgfgraphicnamed{frobenius_nf_pf1}
\InputIfFileExists{frobenius_nf_pf1.tikz}{}{\input{./figures/frobenius_nf_pf1.tikz}}
\endpgfgraphicnamed
  \end{equation}
  
  Using the axioms of a compact closed category, we can deform any $\mathcal A$-diagram into an $\mathcal A$-tree with traces.
  \ctikzfig{frobenius_nf_pf2}
  
  Applying Lemma \ref{lem:cfa-tree-nf}, we have:
  \ctikzfig{frobenius_nf_pf3}
  
  We can then turn each of the traces into a loop map, using equation (\ref{eq:trace-to-loop}), then push the loops up using Proposition \ref{prop:loop-map-endo}.
  \ctikzfig{frobenius_nf_pf4}
  
  Repeating the process for each of the traces, we obtain a diagram in normal form.
\end{proof}

\subsection{Special and Anti-special Commutative Frobenius Algebras}\label{sec:special-anti}

For any monoid in a traced category, we can define a map $\kappa : A \otimes A \rightarrow I$:
\begin{equation}\label{eq:killing-form}
\beginpgfgraphicnamed{killing_form}
\InputIfFileExists{killing_form.tikz}{}{\input{./figures/killing_form.tikz}}
\endpgfgraphicnamed
\end{equation}

In $\catFVect_K$, any map from $A \otimes A$ to the base field $K$ is called a \textit{bilinear form}. The map given by (\ref{eq:killing-form}) is a particularly important bilinear form called the \textit{Killing form}\footnote{The Killing form is much more commonly defined for Lie algebras than for associative algebras. The Killing form of an associative algebra is sometimes called simply its \textit{trace form}.}. This form plays a particularly important role in the representation theory of an algebra. For an algebraically closed field $K$, the Killing form of a finite-dimensional associative $K$-algebra is non-degenerate if and only if that algebra is semisimple. Since the Killing form automatically associates with the multiplication, any algebra with a non-degenerate Killing form (i.e. any finite-dimensional semisimple algebra) is automatically Frobenius. However, the converse is not true. There are many interesting Frobenius algebras that have degenerate Killing forms. First, note that we can relate the rank of the Killing form of a commutative Frobenius algebra with that of the loop map defined in the previous section.
\ctikzfig{killing_form_rank}

We shall primarily focus on the minimal and maximal cases of this rank: the cases where the loop map is full-rank or rank one. It is never rank zero, except in the case of the zero-dimensional space. We wish to abstract the notions of full-rank and rank one to an arbitrary category. Clearly full-rank maps are just isomorphisms. A linear map is rank-one if any only if it can be factored through the base field:
\begin{center}
  \begin{tikzpicture}[yshift=0.8mm]
    \matrix (m) [cdiag] {
      A & K \\
        & B \\
    };
    \path [arrs]
      (m-1-1) edge node [swap] {$f$} (m-2-2)
    (m-1-1) edge node {$\xi$} (m-1-2)
      (m-1-2) edge node {$v$} (m-2-2);
  \end{tikzpicture}
\end{center}

We therefore define disconnected morphisms as an abstraction of rank-one linear maps.

\begin{definition}
  A morphism $f : A \rightarrow B$ in a monoidal category is called \textit{disconnected} if it factors through the tensor unit.
  \ctikzfig{disconnected_map}
\end{definition}


For a Frobenius algebra, having an invertible loop is the abstract analogue to being semisimple. In particular, we shall consider commutative Frobenius algebras whose loop map is equal to the identity. 

\begin{definition}\label{def:special}
  A \textit{special commutative Frobenius algebra} (SCFA) is a commutative Frobenius algebra $\mathcal A = (A, \mu, \eta, \delta, \epsilon)$ such that $\mu \circ \delta = 1_A$. Graphically:
  \begin{equation}\label{eq:special}
\beginpgfgraphicnamed{special}
\InputIfFileExists{special.tikz}{}{\input{./figures/special.tikz}}
\endpgfgraphicnamed
  \end{equation}
\end{definition}

This is not a great loss of generality, given the following theorem.

\begin{theorem}\label{thm:semisimple-special}
  For any commutative Frobenius algebra $\mathcal A = (A, \mu, \eta, \delta, \epsilon)$ such that $\mu \circ \delta$ is invertible, there exists an invertible phase $L$ such that $(A, \mu, \eta, \delta \circ L^{-1}, \epsilon \circ L)$ is an SCFA.
\end{theorem}

\begin{proof}
  Since $\mathcal A$ is semisimple, the loop map is invertible. Therefore, let $L = \mu \circ \delta$. Module endormorphisms are closed under inversion, so since $L$ is a phase, $L^{-1}$ is a phase. For, we show $(A, \delta \circ L^{-1}, \epsilon \circ L)$ is a comonoid:
  \ctikzfig{frobenius_semi_pf1}
  
  The Frobenius identity:
  \ctikzfig{frobenius_semi_pf2}
  
  Specialness follows by definition of $L$.
  \ctikzfig{frobenius_semi_pf3}
\end{proof}

It is also worth noting that for special Frobenius algebras, defining just the monoid part (or just the comonoidal part) is enough to fix the entire structure. This is because the partial trace of $\mu$ is equal to $\epsilon$.
\ctikzfig{special_counit_trace}

By condition (3.) of Theorem \ref{thm:frobenius-defs}, $(A,\mu,\eta,\epsilon)$ is enough to define the entire Frobenius algebra.

\begin{example}\label{ex:basis-scfa}
  Fix a basis $e_i$ for a finite-dimensional vector space $V$. We can define the $e_i$ to be the vectors ``copied'' by $\delta$ and ``deleted'' by $\epsilon$.
  \[ \delta :: e_i \mapsto e_i \otimes e_i \qquad\qquad \epsilon :: e_i \mapsto 1 \]
  This clearly forms a cocommutative comonoid. We can complete the Frobenius algebra by letting the $e_i$ be a basis of idempotents of $\mu$.
  \[ \mu :: e_i \otimes e_i \mapsto e_i \qquad\qquad \eta :: \sum e_i \]
  This forms a commutative Frobenius algebra, and by definition, $\mu \delta = 1_V$. We can do this for any basis $e_i$, and its a well-known fact that \textit{any} semisimple commutative algebra $(A,\mu,\eta)$ over an algebraically closed field has a basis of idempotents, summing to $\eta$, so in particular, all SCFAs are of this form.
\end{example}

Now, we consider the other extreme: the case where the loop map is disconnected:
\begin{equation}\label{eq:disconnected-loop}
\beginpgfgraphicnamed{disconnected_loop}
\InputIfFileExists{disconnected_loop.tikz}{}{\input{./figures/disconnected_loop.tikz}}
\endpgfgraphicnamed
\end{equation}

Due to a theorem by Herrmann~\cite{Herrmann2010}, we can actually obtain an explicit form for equation (\ref{eq:disconnected-loop}).

\begin{theorem}\label{thm:disconnected-anitspecial}
  Let $\mathcal A$ be a commutative Frobenius algebra with a disconnected loop map. Then the following equation holds:
  \ctikzfig{antispecial}
\end{theorem}

\begin{proof}
  Assume the following equation holds, for any maps $x: A \rightarrow I$, $y: I \rightarrow A$:
  \ctikzfig{disconnected_loop}
  
  Then
  \ctikzfig{disconnected_antispecial_pf}
\end{proof}

A commutative Frobenius algebra with a disconnected loop is called \textit{anti-special}.

\begin{definition}
  An \textit{anti-special commutative Frobenius algebra} (ACFA) is a commutative Frobenius algebra such that:
  \ctikzfig{antispecial}
\end{definition}

In addition to having a unit and a counit, anti-special Frobenius algebras have canonical disconnecting points which we shall refer to as the \textit{anti-unit} and \textit{anti-counit}.

\begin{definition}
  For an ACFA, the \textit{anti-unit} $\widetilde\eta$ and the \textit{anti-counit} $\widetilde\epsilon$ are defined as follows:
  \ctikzfig{anti_unit}
\end{definition}

Special and anti-special Frobenius algebras have well-behaved normal forms.

\begin{theorem}\label{thm:scfa-spider}
  For an SCFA $\mathcal S$, any connected $\mathcal S$-diagram is equivalent to a spider.
\end{theorem}

\begin{proof}
  Any connected $\mathcal S$-diagram is equivalent to one in the form given in Theorem \ref{thm:cfa-nf}. We can then use equation (\ref{eq:special}) to remove all of the loops.
\end{proof}

\begin{lemma}\label{lem:lolli-copy}
  For an ACFA $\mathcal A$, the comultiplication map copies the anti-unit and the multiplication map copies the anti-counit.
  \ctikzfig{lolli_copy}
\end{lemma}

\begin{proof}
  Follows straightforwardly by an application of Theorem \ref{thm:cfa-nf} and the anti-specialness condition.
  \ctikzfig{lolli_copy_pf}
  
  The upside-down equation is proved similarly.
\end{proof}

\begin{theorem}\label{thm:acfa-spider}
  Suppose $\circl$ is an invertible scalar. For an ACFA $\mathcal A$, any connected $\mathcal A$-diagram is equivalent to one of the following, for scalar map $k : I \rightarrow I$.
\begin{center}
\beginpgfgraphicnamed{acfa_nf}
\InputIfFileExists{acfa_nf.tikz}{}{\input{./figures/acfa_nf.tikz}}
\endpgfgraphicnamed
\end{center}
\end{theorem}

\begin{proof}
  If a connected $\mathcal A$-diagram contains no loops, it is equivalent to a spider by Theorem \ref{thm:cfa-nf}. An $\mathcal A$-diagram containing more than one loop is equivalent to a scalar multiple of two disconnected diagrams containing one loop each.
  \begin{equation}\label{eq:acfa-nf-pf1}
\beginpgfgraphicnamed{acfa_nf_pf1}
\InputIfFileExists{acfa_nf_pf1.tikz}{}{\input{./figures/acfa_nf_pf1.tikz}}
\endpgfgraphicnamed
  \end{equation}
  
  Similarly, a diagram containing just a single loop can be made into two disconnected diagrams containing single loops.
  \begin{equation}\label{eq:acfa-nf-pf2}
\beginpgfgraphicnamed{acfa_nf_pf2}
\InputIfFileExists{acfa_nf_pf2.tikz}{}{\input{./figures/acfa_nf_pf2.tikz}}
\endpgfgraphicnamed
  \end{equation}
  
  If the diagram has no inputs or outputs, it is in the form of (i.), so assume it has at least one input or output. In that case, any diagram equivalent to the RHS of (\ref{eq:acfa-nf-pf1}) or (\ref{eq:acfa-nf-pf2}) can be put into the form of (iii.) using Lemma \ref{lem:lolli-copy}.
\end{proof}


	\part{Graphical Languages and Rewriting}\label{part:rewriting}
	



\chapter{Rewrite Systems}\label{ch:rewrite-systems}

Rewrite systems provide a model of computation that is particularly well suited to formalising dynamical systems, computing algebraic identities, and constructing proofs by automated or semi-automated means. The most well-studied type of rewrite systems are \textit{term rewrite systems}. A term rewrite system consists of a set of \textit{generators} (i.e. symbols with arities), \textit{variables}, and \textit{rewrite rules} between terms formed from generators and variables.

Term rewrite systems have applications in the study of programming languages, computer algebra systems, automated theorem proving, and many other areas of theoretical computer science. Rewriting terms is essentially the same as rewriting trees, and it was shown that many of these ``tree'' rewriting techniques could actually be applied to arbitrary graphs, or even objects of more general categories.

In this chapter, we shall review some of the basic principals of rewrite systems and the double-pushout approach to graph rewriting. We will then illustrate how the DPO approach is always well-defined in a particular class of graph-like categories called adhesive categories. However, to define a category suitable for rewriting string graphs, we shall need a weaker notion than adhesivity. For that reason, we introduce partial adhesive categories, and show how these categories inherit ``enough adhesivity'' from their ambient adhesive categories to do DPO rewriting.

\section{Term Rewriting}\label{sec:term-rewriting}

\begin{definition}
  A \textit{term signature} $\Sigma = (G, a : G \rightarrow \mathbb N)$ consists of a set $G$ of \textit{generators} and a function $a$ assigning each generator an arity.
\end{definition}

We define the set of \textit{terms} for a signature and a set of variables recursively.

\begin{definition}
  For a term signature $\Sigma = (G,a)$ and a set $X$ of variables, we can form the set $T(\Sigma, X)$ of terms as follows:
  \begin{itemize}
    \item For all $x \in X$, $x$ is a term.
    \item For $g \in G$ such that $a(g) = n$ and terms $t_i \in T(\Sigma, X)$, $g(t_1,\ldots,t_n)$ is a term.
  \end{itemize}
\end{definition}

Variables are place-holders for other terms, i.e. other elements of $T(\Sigma, X)$. The mechanism by which variables are assigned values is called \textit{substitution}.

\begin{definition}
  For a set of terms $T(\Sigma, X)$, a function $\sigma : X \rightarrow T(\Sigma, X)$ is called a \textit{substitution}. It can be lifted to a function $\hat\sigma : T(\Sigma, X) \rightarrow T(\Sigma, X)$ as follows. For a term $t \in T(\Sigma, X)$, replace every occurrence of any variable $x\in X$ in the term with $\sigma(x)$. The resulting term is $\hat\sigma(t)$.
\end{definition}

A rewrite rule is a pair $(l,r) \in T(\Sigma, X) \times T(\Sigma, X)$, usually written $l \rewritesto r$. A rewrite rule can be used to rewrite a term $t \in T(\Sigma, X)$ to a new term $t' \in T(\Sigma, X)$. This occurs in two stages.

\begin{enumerate}
  \item \textbf{Matching:} A substitution $\sigma$ is chosen such that $\hat\sigma(l)$ occurs as a sub-term of $t$.
  \item \textbf{Replacement:} The occurrence of $\hat\sigma(l)$ in $t$ is replaced by $\hat\sigma(r)$.
\end{enumerate}

We can elaborate on this process a bit. Terms are essentially just trees:
\ctikzfig{term_tree}

 For a given tree, a subtree can be uniquely identified by its \textit{lexicographic position}. The lexicographic position of a sub-tree $t'$ of $t$ is simply a list of natural numbers $p = [c_0, c_1, \ldots, c_{n-1}]$. The root of $t'$ can then be located by taking the $c_0$-th child of the root vertex of $t$, the $c_1$-th child of that vertex, and so on until $c_{n-1}$. In other words, if we represent a term $t$ as a list of lists,
\[ t' = t[c_0][c_1]\ldots[c_{n-1}] \]

Thus, in the case of term rewriting, finding a matching $l$ in $l \rewritesto r$ is simply a case of identifying a substitution and lexicographic position such that
\[ t[c_0][c_1]\ldots[c_{n-1}] = \hat\sigma(l) \]

Performing the rewrite is then just a case of replacing the subtree at that position:
\[ t[c_0][c_1]\ldots[c_{n-1}] := \hat\sigma(r) \]

A set $\mathcal R$ of rewrite rules is called a \textit{rewrite system}. We write $t \rewritesto_{\mathcal R} t'$ there exists a rule $l\rewritesto r \in \mathcal R$ that rewrites $t$ into $t'$ using the above procedure. This forms a relation $\rewritesto_{\mathcal R} \subseteq T(\Sigma, X) \times T(\Sigma, X)$ called a \textit{reduction relation} of $\mathcal R$. It is not always the case that the LHS of some rule in $\mathcal R$ will match a given term $t$. If there is no such matching, $t$ is called \textit{irreducible}. Otherwise, it is called \textit{reducible}.

\begin{example}
  Consider the algebraic theory of unital rings. This theory has two binary operations $(- + -)$ and $(- \cdot -)$ as well as two $0$-ary operations (i.e. constants) $1$ and $0$. The usual ring axioms can be turned into a rewrite system by directing each of the equations.
  \begin{enumerate}
    \item $((x + y) + z) \rewritesto (x + (y + z))$
    \item $((x \cdot y) \cdot z) \rewritesto (x \cdot (y \cdot z))$
    \item $(x + 0) \rewritesto x$
    \item $(0 + x) \rewritesto x$
    \item $(x \cdot 1) \rewritesto x$
    \item $(1 \cdot x) \rewritesto x$
    \item $(x \cdot (y + z)) \rewritesto ((x \cdot y) + (x \cdot z))$
    \item $((x + y) \cdot z) \rewritesto ((x \cdot z) + (y \cdot z))$
  \end{enumerate}
  
  We can define a term
  \[ t := (0 + (((1 + y) + y) \cdot (x + z))) \]
  
  Let $l \rewritesto r$ be rule $8$ from above. To find a matching of $l$ on $t$, first define a substitution:
  \[ \sigma :: \{ x \mapsto (1 + y),\ y \mapsto y,\ z \mapsto (x + z) \} \]
  
  Apply the substitution to $l$ and $r$:
  \begin{align*}
    \hat\sigma(l) & = (((1 + y) + y) \cdot (x + z))                 \\
    \hat\sigma(r) & = (((1 + y) \cdot (x + z)) + (y \cdot (x + z)))
  \end{align*}
  
  Note that $\hat\sigma(l)$ now occurs as a subterm of $t$.
  \[ t = (0 + \hat\sigma(l)) \]
  
  We form a new term $t'$ by replacing the occurrence of $\hat\sigma(l)$ with $\hat\sigma(r)$.
  \[ t' =  (0 + \hat\sigma(r)) = (0 + (((1 + y) \cdot (x + z)) + (y \cdot (x + z)))) \]
\end{example}

$\rewritetrans_{\mathcal R}$ is used to represent the reflexive, transitive closure of $\rewritesto_{\mathcal R}$, and $\rewriteequiv_{\mathcal R}$ the reflexive, transitive, symmetric closure. In the case where $\mathcal R$ represents the ``directed version'' of the axioms of an algebraic structure, it is often the goal to evaluate the truth of the following proposition:
\[ t \rewriteequiv_{\mathcal R} t' \]

Or, ``Is $t$ equivalent to $t'$ by the axioms of $\mathcal R$?'' As this is a statement of the word problem, it is not decidable for arbitrary rewrite systems $\mathcal R$. However, there are many ``good'' rewrite systems, where this proposition \textit{is} decidable. Two characteristics of a rewrite system that make it ``good'' are \textit{termination} and \textit{confluence}.

\begin{definition}\label{def:terminating}
  A rewrite system $\mathcal R$ is said to be \textit{terminating} if there exists no infinite chain of rewrites:
  \[ t_1 \rewritesto t_2 \rewritesto \ldots \rewritesto t_n \rewritesto \ldots \]
\end{definition}

In practice, one often proves termination by identifying a \textit{reduction order} on terms.

\begin{definition}
  A partially ordered set $(P, \leq)$ is called \textit{well-founded} (or Noetherian) if it has a smallest element and contains no infinite sequence of strictly decreasing elements:
  \[ p_1 > p_2 > \ldots > p_n > \ldots \]
\end{definition}

Well-foundedness is usually defined by every non-empty subset $P' \subseteq P$ having at least one minimal element (i.e. an element that is not strictly greater than any other). Up to the Axiom of Choice, these two definitions are equivalent. Intuitively, a well-founded poset is a poset over which one can perform (generalised) induction. A standard example is $(\mathbb N, \leq)$, which corresponds to the usual notion of induction over the natural numbers.

\begin{definition}\label{def:reduction-ordering}
  Let $(P, \leq)$ be a well-founded poset. A function $\omega : T(\Sigma, X) \rightarrow P$ is called a \textit{reduction ordering} for a rewrite system $\mathcal R$ if:
  \[ t_1 \rewritesto_{\mathcal R} t_2 \ \  \implies \ \ \omega(t_1) > \omega(t_2) \]
\end{definition}

Clearly any rewrite system that admits a reduction ordering is terminating. Termination guarantees that even a na\"ive rewriting strategy (choosing rules at random, applying until there are no more matchings) will terminate. If we can rewrite a term $t$ in any number of steps to an irreducible term $t'$, $t'$ is called a \textit{normal form} of $t$. However, with a rewrite system that is only terminating, there is no guarantee that two distinct sequences of rewrites will result in the same normal form. To get \textit{unique} normal forms, we need an additional property called \textit{confluence}.

\begin{definition}\label{def:confluent}
  A rewrite system $\mathcal R$ is said to be \textit{confluent} if sequences of rewrites starting with the same term can be rejoined. That is, for terms $t, t_1, t_2$ such that $t \rewritetrans t_1$ and $t \rewritetrans t_2$, there exists a term $t'$ such that  $t_1 \rewritetrans t'$ and $t_2 \rewritetrans t'$.
  \begin{center}
    \begin{tikzpicture}
      \matrix (m) [cdiag] {
        t   & t_1 \\
        t_2 & t'  \\
      };
      \path [arrs,-open triangle 45]
        (m-1-1) edge node {$*$} (m-1-2)
        (m-1-1) edge node [swap] {$*$} (m-2-1)
        (m-1-2) edge node {$*$} (m-2-2)
        (m-2-1) edge node {$*$} (m-2-2);
    \end{tikzpicture}
  \end{center}
\end{definition}

\begin{theorem}\label{thm:ct-unique-nf}
  If a rewrite system $\mathcal R$ is terminating and confluent, every term $t$ has a unique normal form $t\normalised$.
\end{theorem}

\begin{proof}
  The existence of at least one normal form is guaranteed by termination. Suppose $t_1$ and $t_2$ are normal forms for $t$. Then $t \rewritetrans t_1$ and $t \rewritetrans t_2$, so by confluence there exists $t'$ such that $t_1 \rewritetrans t'$ and $t_2 \rewritetrans t'$. However, since $t_1$ and $t_2$ are irreducible, the only possibility is that $t_1 = t' = t_2$.
\end{proof}

Terminating, confluent rewrite systems provide a decidable solution to the word problem.

\begin{theorem}\label{thm:ct-word-problem}
  For terms $t_1$ and $t_2$, $t_1 \rewriteequiv t_2$ iff $t_1 \normalised = t_2 \normalised$.
\end{theorem}

\begin{proof}
  ($\Leftarrow$) follows from the definition of $(-)\normalised$. For ($\Rightarrow$), assume $t_1 \rewriteequiv t_2$. Then, there exists a finite sequence of forward and backward rewrite steps between $t_1$ and $t_2$.
  \begin{equation}\label{eq:peak-sequence}
    t_1
    \leftrewritetrans p_1 \rewritetrans q_1 \leftrewritetrans \ldots
    \rewritetrans q_{n-1} \leftrewritetrans p_n \rewritetrans q_n \leftrewritetrans
    t_2
  \end{equation}
  Note that the rewrite sequence consists of ``peaks'' $p_i$ and ``troughs'' $q_i$. We proceed by induction on the number of peaks. An arbitrary rewrite sequence with $n$ peaks, as in equation (\ref{eq:peak-sequence}) can be reduced to a rewrite system with $n-1$ peaks by applying confluence to $q_{n-1} \leftrewritetrans p_n \rewritetrans q_n$.
  \begin{equation*}
    t_1
    \leftrewritetrans p_1 \rewritetrans q_1 \leftrewritetrans \ldots
    \rewritetrans q_{n-1} \rewritetrans p'_n \leftrewritetrans q_n \leftrewritetrans
    t_2
  \end{equation*}
  
  If there are zero peaks, then there exists at most one trough $q$, and by confluence $t_1 \normalised = q \normalised = t_2 \normalised$.
\end{proof}

\section{Graph Rewriting}\label{sec:graph-rewriting}

\begin{definition}
  Let $\catGraph$ be the category of (directed, multi-) graphs. It is defined as the functor category $[\mathbb G, \catSet]$, for $\mathbb G$ defined as:
  \begin{center}
    \cpair{E}{V}{s}{t}
  \end{center}
  $E$ identifies the edges of the graph, and $V$ the vertices. $s$ and $t$ are functions taking an edge to its source and target respectively. If $t(e) = v$ then $e$ is called an \emph{in-edge} of $v$ and if $s(e) = v$ then $e$ is called an \emph{out-edge} of $v$.
\end{definition}

For a graph $G : \mathbb G \rightarrow \catSet$, we shall write $V_G$ for $G(V)$, $E_G$ for $G(E)$, etc. and drop the subscripts when it is clear from context. Since a graph homomorphism $f : G \rightarrow H$ is just a natural transformation, it is a pair of functions $f_V, f_E$ such that:
\begin{center}
  \csquare{E_G}{V_G}{E_H}{V_H}{s}{s}{f_E}{f_V}
  \qquad\qquad
  \csquare{E_G}{V_G}{E_H}{V_H}{t}{t}{f_E}{f_V}
\end{center}

It is natural to ask if term rewriting (i.e. tree rewriting) techniques can be extended to arbitrary graphs. It turns out that most notions and techniques from term rewriting translate directly to graphs, but the concepts of matching and replacement become more complicated. This is because, unlike in the case of trees, there is no canonical ``root vertex'' of a graph, and hence no absolute notion of ``position''. Proceeding by analogy to the term rewriting case, suppose we wish to apply a graph rewrite rule $L \rewritesto R$ to a graph $G$. This poses two problems:
\begin{enumerate}
  \item How does one keep track of where a matching of the LHS of a graph rewrite rule has been made?
  \item How does one decide where to attach to the RHS of a graph rewrite rule?
\end{enumerate}

The first problem is solved by letting a matching be represented by an injective graph homomorphism $m : L \rightarrow G$, subject to certain conditions. The second problem is solved by requiring that all graph rewrite rules have a common subgraph embedded into the LHS and RHS called the \textit{invariant subgraph} of the rewrite rule.
\[ L \hookleftarrow I \hookrightarrow R \]

This invariant subgraph serves as ``glue'' to attach $R$ to $G$ after (the non-invariant part of) $L$ has been removed. We can sketch out diagrammatically how this works. Suppose we have matching $m : L \rightarrow G$. Since $m$ is an injection, we can think of $L$ as a subgraph of $G$. Before we insert $R$, we must remove $L$. To do this, we might consider using the graph theoretic subtraction.

\begin{definition}\label{def:graph-subtraction}
  For a subgraph $G' \subseteq G$, the \textit{graph theoretic subtraction} $G - G'$ is a new graph formed by removing $G'$, as well as any edges in to or out of $G'$, from $G$.
\end{definition}

So $G$ consists of a component $L$, a component $G-L$, and some edges between those components.
\medskip
\ctikzfig{graph_decompose}
\medskip

If we simply delete $L$ from $G$, there is no way to keep track of how $L$ was connected to $G$ in the first place. To get around this, we treat the invariant subgraph of the rewrite rule as an \textit{interface} for $L$, and further decompose $G$:
\begin{equation}
  \medskip
\beginpgfgraphicnamed{graph_decompose2}
\InputIfFileExists{graph_decompose2.tikz}{}{\input{./figures/graph_decompose2.tikz}}
\endpgfgraphicnamed
  \medskip
\end{equation}

The part of $L$ that is not contained in $I$ (i.e. $L - I$ and the edges between $L - I$ and $I$) is called the \textit{interior} of $L$. Note that there are edges between $G$ and $I$ and there are edges between $I$ and $L - I$, but there are no edges directly connecting $G$ to $L - I$. We require this to be the case for any valid matching $m$ of the rewrite rule. This is called the \textit{no-dangling-edges} condition. Now, when we replace $L$ with $R$, we know where the edges connected to $I$ should go, since \textit{$R$ also contains a copy of $I$}.
\medskip
\ctikzfig{rewrite_procedure}
\medskip

This procedure can be formalised elegantly using pushouts. Begin by noting that the graph $G$ with the interior of $I$ removed is the (unique) graph $G'$ such that the following square is a pushout:
\begin{center}
  \begin{tikzpicture}
    \matrix (m) [cdiag] {
      I  & L \\
      G' & G \\
    };
    \path [arrs]
      (m-1-1) edge [right hook-latex] (m-1-2)
      (m-1-1) edge node [swap] {$m'$} (m-2-1)
      (m-1-2) edge node {$m$} (m-2-2)
      (m-2-1) edge [right hook-latex] (m-2-2);
    \NWbracket{(m-2-2)}
  \end{tikzpicture}
\end{center}

In other words, $G$ is the result of gluing $L$ and $G'$ together along $I$. $G'$ is called the \textit{pushout complement} of $I \hookrightarrow L \overset{m}{\rightarrow} G$. Once we have $G'$ and $m' : I \rightarrow G'$, we can glue on $R$ by performing a second pushout:
\begin{center}
  \begin{tikzpicture}
    \matrix (m) [cdiag] {
      I  & R \\
      G' & H \\
    };
    \path [arrs]
      (m-1-1) edge [right hook-latex] (m-1-2)
      (m-1-1) edge node [swap] {$m'$} (m-2-1)
      (m-1-2) edge (m-2-2)
      (m-2-1) edge (m-2-2);
    \NWbracket{(m-2-2)}
  \end{tikzpicture}
\end{center}

So, to complete the rewrite $G \rewritesto H$, we first compute a pushout complement, then compute a pushout. This style of graph rewriting is called the \textit{double pushout} (DPO) approach. We typically express the whole rewrite as a single DPO diagram:
\begin{center}
  \begin{tikzpicture}
    \matrix (m) [cdiag] {
      L & I  & R \\
      G & G' & H \\
    };
    \path [arrs]
      (m-1-2) edge [left hook-latex] (m-1-1)
      (m-1-1) edge node [swap] {$m$} (m-2-1)
      (m-1-2) edge [right hook-latex] (m-1-3)
      (m-2-2) edge (m-2-1)
      (m-2-2) edge (m-2-3)
      (m-1-2) edge node {$m'$} (m-2-2)
      (m-1-3) edge (m-2-3);
    \NEbracket{(m-2-1)}
    \NWbracket{(m-2-3)}
  \end{tikzpicture}
\end{center}

Note that pushouts (and pushout complements) are defined up to isomorphism, so for $G \cong G'$, $H \cong H'$, $G \rewritesto H$ if and only if $G' \rewritesto H'$. As before, a set of graph rewrite rules $\mathcal R$ is called a \textit{graph rewrite system}, and we write $G \rewritesto_{\mathcal R} H$ if there exists a rule in $\mathcal R$ that rewrites $G$ to $H$, as above. The notions of confluence and termination are identical to those in Section \ref{sec:term-rewriting}, as well as the proofs of Theorems \ref{thm:ct-unique-nf} and \ref{thm:ct-word-problem} for graphs, replacing term equality with graph isomorphism.

\begin{example}
  Let the following be a rewrite rule $L \hookleftarrow I \hookrightarrow R$:
  \ctikzfig{rewrite_ex_rule}
  
  Then, we can find a matching of $L$ on a bigger graph $G$:
  \ctikzfig{rewrite_ex_matching}
  
  We perform the rewrite by first removing the interior of $L$, then gluing $R$ to remainder of the graph along the common subgraph $I$.
  \ctikzfig{rewrite_ex}
  
  The full double pushout diagram looks like this:
  \ctikzfig{rewrite_ex_dpo}
\end{example}

\section{Adhesive Categories}\label{sec:adhesive-categories}

In the previous section, we discussed double-pushout rewriting in the category $\catGraph$. However, we may also wish to do rewriting on the objects of many categories that look a bit like the category of graphs, such as typed graphs, graphs with extra data, Petri nets, or objects of an arbitrary topos. There have been two main approaches to carrying out this abstraction. The first was to define DPO rewriting in the context of high-level replacement (HLR) systems, introduced by Ehrig et al \cite{Ehrig1991}. The second approach, which we shall build on in this dissertation, relies on \textit{adhesive categories}, introduced by Lack and Soboci\'nski \cite{LackAdh2005}. These two notions are actually compatible, as was shown with the definition of \textit{adhesive HLR categories} in \cite{Ehrig2008}, which categories equipped with a class of monomorphisms that behave ``adhesively''. Our construction of \textit{partial adhesive categories} is somewhat in the spirit of adhesive HLR categories, but adhesive behaviour is localised to a certain class of well-behaved \textit{spans} within the category rather than a certain class morphisms.

Adhesive categories provide a general context in which rewriting on graph-like structures is well-defined. Their definition relies on the notion of a special kind of pushout called a van Kampen square.

\begin{definition}\label{def:van-kampen}
  A van Kampen square (VK-square) is a pushout
  \begin{center}
    \posquare{A}{B}{C}{D}{}{}{}{}
  \end{center}
  
  Such that for any commutative cube
  \begin{center}
    \vksquare
      {\bottomobjects {A}{B}{C}{D}}
      {\topobjects    {G}{H}{E}{F}}
      {\toparrows     {}{}{}{}}
      {\sidearrows    {}{}{}{}}
      {\bottomarrows  {}{}{}{}}
  \end{center}
  where the back and left faces ($ABEF$ and $ACEG$) are pullbacks, the following are equivalent:
  \begin{itemize}
    \item the front and right faces ($CDGH$ and $BDFH$) are pullbacks
    \item the top face ($EFGH$) is a pushout
  \end{itemize}
\end{definition}

A pushout of a span $\cspan{A}{f}{B}{g}{C}$ where either $f$ or $g$ is a monomorphism is called a \textit{pushout along a monomorphism}. An adhesive category is a category where pushouts along monomorphisms are van Kampen squares.

\begin{definition}\label{def:adhesive}
  A category $\mathcal A$ is said to be \emph{adhesive} if
  \begin{enumerate}
    \item $\mathcal A$ has pushouts along monomorphisms,
    \item $\mathcal A$ has pullbacks,
    \item and pushouts along monomorphisms in $\mathcal A$ are van
          Kampen squares.
  \end{enumerate}
\end{definition}

At first sight, the definition of van Kampen squares can seem rather opaque, so it is useful to consider a concrete example. Suppose we have a set $X = A \cup B$ which is composed to two (possibly overlapping) subsets $A$ and $B$. There are two equivalent ways to define a map into $X$:
\begin{enumerate}
  \item For some set $X'$, simply define a function $f : X' \rightarrow X$.
  \item For sets $A'$, $B'$, define functions $f_A : A' \rightarrow A$ and $f_B : B' \rightarrow B$ such that $f_A$ and $f_B$ agree on $A \cap B$. I.e. for the restrictions $f_A |(A \cap B) : I \rightarrow A \cap B$, $f_B |(A \cap B) : I' \rightarrow A \cap B$, $I = I'$ and $f_A |(A \cap B) = f_B |(A \cap B)$.
\end{enumerate}

First, we'll see how we can obtain (2.) from (1.). Starting with a map $f : X' \rightarrow X$, one can simply let $f_A$ and $f_B$ be the restrictions of $f$ to the subsets $A$ and $B$ respectively. The restriction of a function to a subset of its codomain is just a pullback. We can therefore define $f_A$ and $f_B$ by this diagram:
\begin{equation}\label{eq:adh-ex1}
  \begin{tikzpicture}
    \matrix (m) [cdiag] {
      A' & X' & B' \\
      A  & X  & B  \\
    };
    \path [arrs]
      (m-1-1) edge (m-1-2)
      (m-1-1) edge node [swap] {$f_A$} (m-2-1)
      (m-1-2) edge node {$f$} (m-2-2)
      (m-2-1) edge (m-2-2)
      (m-1-3) edge (m-1-2)
      (m-2-3) edge (m-2-2)
      (m-1-3) edge node {$f_B$} (m-2-3);
      \SEbracket{(m-1-1)}
      \SWbracket{(m-1-3)}
  \end{tikzpicture}
\end{equation}

To obtain (1.) from (2.), we start with maps $f_A, f_B$, such that $f_A |(A \cap B) = f_B |(A \cap B)$. Again using pullbacks to express restriction, this means the following diagram commutes (treating $I$ from above as $A' \cap B'$):
\begin{equation}\label{eq:adh-ex2}
  \begin{tikzpicture}
    \matrix (m) [cdiag] {
      A' & A' \cap B' & B' \\
      A  & A  \cap B  & B  \\
    };
    \path [arrs]
      (m-1-2) edge (m-1-1)
      (m-1-1) edge node [swap] {$f_A$} (m-2-1)
      (m-1-2) edge (m-2-2)
      (m-2-2) edge (m-2-1)
      (m-1-2) edge (m-1-3)
      (m-2-2) edge (m-2-3)
      (m-1-3) edge node {$f_B$} (m-2-3);
      \SWbracket{(m-1-2)}
      \SEbracket{(m-1-2)}
  \end{tikzpicture}
\end{equation}

Let $X' = A' \cup B'$ and define $f$ as:
\begin{equation*}
  f(x) =
  \begin{cases}
    f_A(x) & \textrm{ if } x \in A \\
    f_B(x) & \textrm{ if } x \in B
  \end{cases}
\end{equation*}

This function is defined unambiguously, because if $x$ is in both $A$ and $B$ then $f_A(x) = f_B(x)$. Since $X' = A' \cup B'$ is a pushout, we can equivalently define $f$ as the induced map in the following diagram:
\begin{equation}\label{eq:adh-ex3}
  \begin{tikzpicture}
    \matrix (m) [cdiag] {
      A' \cap B' & B' &   \\
      A'         & X' & B \\
                 & A  & X \\
    };
    \path [arrs]
      (m-1-1) edge (m-1-2)
      (m-1-1) edge (m-2-1)
      (m-1-2) edge (m-2-2)
      (m-2-1) edge (m-2-2)
      (m-2-1) edge node [swap] {$f_A$} (m-3-2)
      (m-1-2) edge node {$f_B$} (m-2-3)
      (m-2-2) edge [dashed] node {$f$} (m-3-3)
      (m-3-2) edge (m-3-3)
      (m-2-3) edge (m-3-3);
      \NWbracket{(m-2-2)}
  \end{tikzpicture}
\end{equation}

It is easy to show that these two procedures are the inverses of each other. If we combine diagrams (\ref{eq:adh-ex1}), (\ref{eq:adh-ex2}), and (\ref{eq:adh-ex3}) into a single diagram, we get the commutative cube from Definition \ref{def:van-kampen}.
\begin{center}
  \vksquare
    {\topobjects    {A'}{X'}{A' \cap B'}{B'}}
    {\bottomobjects {A \cap B}{B}{A}{X}}
    {\toparrows     {}{}{}{}}
    {\sidearrows    {f_A}{}{f_B}{f}}
    {\bottomarrows  {}{}{}{}}
\end{center}

Now, as in Definition \ref{def:van-kampen}, suppose we have a commutative cube where the bottom face is a pushout (i.e. $X = A \cup B$) and the left and back faces are pullbacks (i.e. $f_A |(A \cap B) = f_B |(A \cap B)$). We can now read the van Kampen square condition as follows: the front and right faces are pullbacks (i.e. $f$ restricts to $f_A$ and $f_B$ along $A$ and $B$ respectively) if and only if the top face is a pushout (i.e. $X' = A' \cup B'$).

\begin{examples}\label{ex:adh-cats}
  Some examples of adhesive categories:
  \begin{itemize}
    \item $\catSet$, $\catGraph$, and $\catSet_{*}$ are adhesive categories.
    \item For adhesive categories $\mathcal C$ and $\mathcal D$ and an object $X \in \mathcal C$ then $\mathcal C \times \mathcal D$, $X/\mathcal C$, and $\mathcal C/X$ are adhesive.
    \item For a small category $\mathbb X$ and an adhesive category $\mathcal C$, the functor category $[\mathbb X, \mathcal C]$ is adhesive.
    \item Any elementary topos is adhesive.
  \end{itemize}
  
  Note that unlike toposes, adhesive categories are stable under coslices. In particular, $\catSet_{*}$ is an adhesive category, but not a topos.
\end{examples}

Adhesive categories are useful for double-pushout rewriting because pushout complements for pushouts along monomorphisms are unique, when they exist.

\begin{definition}\label{def:pushout-complement}
  A \emph{pushout complement} for a pair of arrows $\crun{A}{m}{B}{g}{D}$, is an object $C$ and a pair of arrows $\crun{A}{f}{C}{n}{D}$ such that the following is a pushout:
  \begin{center}
    \posquare{A}{B}{C}{D}{m}{n}{f}{g}
  \end{center}
\end{definition}

The pushout complement above should be thought of as ``subtracting $B$ from $D$, modulo $A$''. If we think of $A$ as the interface of $B$, then another way to put it is ``removing the interior of $B$ from $D$''.

\begin{notation}
  We occasionally write $B +_{m,f} C$ for pushout of the span $\cspan{B}{m}{A}{f}{C}$ and $D -_{m,g} B$ for the pushout complement of $\crun{A}{m}{B}{g}{D}$.
\end{notation}

We first need a few basic lemmas before showing that pushout complements in an adhesive category are unique.

\begin{lemma}\label{lem:adh-po-pb}
  \cite{LackAdh2005} In an adhesive category:
\begin{itemize} 
  \item monomorphisms are stable under pushout, and
  \item pushouts along monomorphisms are also pullbacks.
\end{itemize}
\end{lemma}

\begin{proof}
  Let $m$ be a monomorphism, and let the following diagram be a pushout:
   \begin{center}
    \posquare{A}{B}{C}{D}{m}{n}{f}{g}
  \end{center}
We need to show that $n$ is mono and the above square is also a pullback. Construct a commutative cube, containing two copies of the given pushout square: one on the bottom face and one on the right face. Place a copy of $f$ in the upper-left corner, and fill in the rest with identities.
  \begin{center}
    \vksquare
      {\bottomobjects {A}{B}{C}{D}}
      {\topobjects    {C}{C}{A}{A}}
      {\toparrows     {1}{f}{f}{1}}
      {\sidearrows    {1}{1}{m}{n}}
      {\bottomarrows  {m}{f}{g}{n}}
  \end{center}
All of the faces of this cube commute trivially. Commutative squares involving identities often (trivially) form pullbacks or pushouts. In particular, the top face is pushout, the left face is a pullback, and the back face is a pullback iff $m$ is a monomorphism. Since $m$ is defined to be a monomorphism and the pushout we started with is a VK-square, adhesivity shows that the front and right faces must be pullbacks. Since the front face is a pullback, $n$ is a monomorphism and thus monomorphisms are stable under pushout. Furthermore, since the right face is the pushout we started with, we have also shown that pushouts along monomorphisms are pullbacks.
\end{proof}

We now provide a few lemmas about pullbacks and pushouts that hold in \emph{any} category.

\begin{lemma}\label{lem:double-pb-square}
  For the following commutative diagram, where the right square is a pullback:
  \begin{center}
    \begin{tikzpicture}
    \matrix (m) [cdiag] {
      A & B & C \\
      D & E & F \\
    };
    \path [arrs]
      (m-1-1) edge node {$f$} (m-1-2)
      (m-1-2) edge node {$g$} (m-1-3)
      (m-1-1) edge node [swap] {$h$} (m-2-1)
      (m-1-2) edge node {$i$} (m-2-2)
      (m-1-3) edge node {$j$} (m-2-3)
      (m-2-1) edge node [swap] {$k$} (m-2-2)
      (m-2-2) edge node [swap] {$l$} (m-2-3);
    \end{tikzpicture}
  \end{center}
  the outer square is a pullback iff the left square is a pullback.
\end{lemma}

\begin{proof}
  Left implies outer is trivial. For the other direction, assume the outer square is a pullback. For an object $A'$ let $f' : A' \rightarrow B$ and $h' : A' \rightarrow D$ be arrows such that $i f' = k h'$.
  \begin{center}
    \begin{tikzpicture}
      \matrix (m) [cdiag] {
        A & B & C \\
        D & E & F \\
      };
      \node (pb) [xshift=-7mm,yshift=7mm] at (m-1-1) {$A'$};
      \path [arrs]
        (m-1-1) edge node [swap] {$f$} (m-1-2)
        (m-1-2) edge node {$g$} (m-1-3)
        (m-1-1) edge node {$h$} (m-2-1)
        (m-1-2) edge node {$i$} (m-2-2)
        (m-1-3) edge node {$j$} (m-2-3)
        (m-2-1) edge node [swap] {$k$} (m-2-2)
        (m-2-2) edge node [swap] {$l$} (m-2-3)
        (pb) edge [bend left=20] node {$f'$} (m-1-2)
        (pb) edge [bend right=20] node [swap] {$h'$} (m-2-1);
      \end{tikzpicture}
  \end{center}
  
  Then $gf'$ and $h'$ form a cone under the outer pullback, so there exists a unique $u : A' \rightarrow A$ such that $gf' = gfu$ and $h' = hu$. From the universal property of the right pullback, it follows that $f' = fu$, so the left square is a pullback.
\end{proof}

\begin{lemma}\label{lem:three-faces}
  For a commutative cube, where the front, right, and back faces are pullbacks:
  \begin{center}
    \vksquare
      {\bottomobjects {A}{B}{C}{D}}
      {\topobjects    {E}{F}{G}{H}}
      {\toparrows     {}{}{}{}}
      {\sidearrows    {}{}{}{}}
      {\bottomarrows  {}{}{}{}}
  \end{center}
  the left face is also a pullback.
\end{lemma}

\begin{proof}
  Since the back and right faces are both pullback squares, we can apply Lemma \ref{lem:double-pb-square} to show that the back-right two faces form a larger pullback square. Then by commutativity of the cube, the square formed by the left-front two faces is also a pullback square. Applying Lemma \ref{lem:double-pb-square} from right to left concludes that the left face is a pullback.
\end{proof}

\begin{lemma}\label{lem:po-stable-iso}
  Isomorphisms are stable under pushout. For the following pushout, if $\phi$ is an isomorphism, then so too is $q$.
  \begin{center}
    \posquare{A}{B}{C}{D}{\phi}{q}{f}{p}
  \end{center}
\end{lemma}

\begin{proof}
  Let $f\phi^{-1}$ and $1_C$ form a cocone to $C$ over the pushout. It follows straightforwardly that the induced map $q' : D \rightarrow C$ is the inverse of $q$.
\end{proof}

We are now ready to prove the uniqueness theorem for pushout complements.

\begin{theorem}\label{lem:unique-pushout-complements}
  \cite{LackAdh2005}. If a pair of arrows $(m,g)$, where $m$ is mono, has a pushout complement, it is unique up to isomorphism. That is, for any two pushout complements, $(f,n)$ and $(f',n')$, there exists an isomorphism $\phi$ making the following diagram commute:
  \begin{equation}\label{dia:unique-compl}
    \csquareslant{A}{C}{C'}{D}{f}{n'}{f'}{n}{\phi}
  \end{equation}
\end{theorem}

\begin{proof}
  Suppose both of the following are pushout squares:
  \begin{center}
    \posquare{A}{B}{C}{D}{m}{n}{f}{g}
    \qquad
    \posquare{A}{B}{C'}{D}{m}{n'}{f'}{g}
  \end{center}
  
 Following a similar strategy to Lemma~\ref{lem:adh-po-pb}, we use these squares to build a commutative cube:
  \tikzstyle{vkarrow2}=[arrow dashed]
  \begin{center}
    \vksquare
      {\bottomobjects {A}{B}{C}{D}}
      {\topobjects    {C''}{C'}{A}{A}}
      {\toparrows     {k}{h}{f'}{1}}
      {\sidearrows    {l}{1}{m}{n'}}
      {\bottomarrows  {m}{f}{g}{n}}
  \end{center}
  \vkcleararrows
\noindent where the first pushout square forms the bottom face, and the second pushout square forms the right face.  The front face is the pullback of $n$ and $n'$, and the back face is the pullback of $m$ with itself. This pullback consists of identities because $m$ is mono. Now, $f$ and $f'$ form a cone under the pullback of $n$ and $n'$, so let $h$ be the induced map. 

The right face is a pushout along a monomorphism, so by Lemma \ref{lem:adh-po-pb}, it is also a pullback. We can therefore conclude from Lemma \ref{lem:three-faces} that the left face is also a pullback. From the VK-square property, we can then conclude that the top face is a pushout. By Lemma \ref{lem:po-stable-iso}, $k$ is an isomorphism. 
We can also form a similar cube, but with the positions of the two pushouts interchanged:
  \tikzstyle{vkarrow2}= [arrow dashed]
  \tikzstyle{vkarrow1}= [arrow bold]
  \tikzstyle{vkarrow3}= [arrow bold]
  \tikzstyle{vkarrow4}= [arrow bold]
  \tikzstyle{vkarrow5}= [arrow bold]
  \tikzstyle{vkarrow6}= [arrow bold]
  \tikzstyle{vkarrow8}= [arrow bold]
  \tikzstyle{vkarrow10}=[arrow bold]
  \tikzstyle{vkarrow12}=[arrow bold]
  \begin{center}
    \vksquare
      {\bottomobjects {A}{B}{C'}{D}}
      {\topobjects    {C''}{C}{A}{A}}
      {\toparrows     {l}{h}{f}{1}}
      {\sidearrows    {k}{1}{m}{n}}
      {\bottomarrows  {m}{f'}{g}{n'}}
  \end{center}
  \vkcleararrows
  and conclude similarly that $l$ is an isomorphism. By construction, all faces commute, and thus by letting $\phi := kl^{-1}$, we can read off the statement of the theorem from the commutative cube (with the relevant arrows shown in bold).
\end{proof}

As in the case of graphs, a \emph{rewrite rule} in an adhesive category $\mathcal A$ is a span of monomorphisms.
\[ L \rewritesto R := \cspan{L}{b_1}{I}{b_2}{R} \]

DPO rewriting consists of three steps.

\begin{enumerate}
  \item Identify a matching $m : L \rightarrow G$.
  \item Compute the pushout complement $G'$ of $L$ in $G$.
  \item Push out $G'$ and $R$ to obtain a rewritten graph $H$.
\end{enumerate}

First, we define \textit{matchings}. Adhesive categories only guarantee the \emph{uniqueness} of pushout complements, not the existence. In most categories, these will not exist for arbitrary pairs of morphisms. Thus, we build the existence condition into the definition of matching.

\begin{definition}\label{def:matching}
  For a rewrite rule
  \[ L \rewritesto R := \cspan{L}{b_1}{I}{b_2}{R} \]
  a \emph{matching} $m$ of $L \rewritesto R$ on $G$ is a monomorphism $m : L \rightarrow G$ such that the following pushout complement exists:
  \[ \posquare{I}{L}{G'}{G} {b_1}{}{}{m} \]
\end{definition}

This existential condition is not particularly useful in determining which morphisms are matchings. Luckily, in specific adhesive categories, we can often do better. For graph-like categories, we can usually ensure the existence of a pushout complement using some version of the \emph{no-dangling-edges condition}.

\begin{definition}\label{def:no-dangling-edges-condition}
  A monomorphism $m : L \rightarrow G$ of a rewrite rule $\cspan{L}{b_1}{I}{b_2}{R}$ in the category $\catGraph$ is said to satisfy the \emph{no-dangling-edges} condition if for any vertex $v \in V_L - V_I$, all edges incident to $m(v)$ must be in the image of $m$.
\end{definition}

\begin{theorem}
  In $\catGraph$, a monomorphism $m$ is a matching if and only if it satisfies the no-dangling-edges condition.
\end{theorem}

\begin{proof}
  Since $\crun{I}{b_1}{L}{m}{G}$ are both monos, we will assume without loss of generality that $I \subseteq L \subseteq G$. For ($\Leftarrow$), we form the pushout complement $G'$ as a graph with vertices $V_G - (V_L - V_I)$ and edges $E_G - (E_L - E_I)$. By the no-dangling-edges condition, if $s(e) \in V_L - V_I$ or $t(e) \in V_L - V_I$ then $e \in E_L$. If $e$ were to be in $E_I$ then $s(e) \in V_I$ and $t(e) \in V_I$ so $e \in V_L - V_I$. Therefore, the maps $s$ and $t$ have well-defined restrictions to $G'$. Furthermore, $G' \cup L = G$ and $G' \cap L = I$, so $G'$ is the pushout complement.
  
  For ($\Rightarrow$), suppose $m$ is a matching. Then there exists a pushout square:
  \begin{center}
    \posquare{I}{L}{G'}{G}{b_1}{b_1'}{m'}{m}
  \end{center}
  
  Since monos in an adhesive category are stable over pushout, $b_1'$ and $m'$ are also monos. As such, we identify $I$ with its image under $m'$ regard $G$ as $G' \cup L$. Let $v$ be a vertex in $V_L - V_I$, and $e$ an edge incident to $v$ in $G$. Suppose $e$ is not in $L$ (i.e. the image of $m : L \rightarrow G$), then $e$ must be in $G'$. Thus $v$ must also be in $G'$, so $v \in V_{G'} \cap V_L = V_I$, which is a contradiction. We therefore conclude that $e \in E_L$.
\end{proof}

We can now define the notion of a rewrite in an adhesive category.

\begin{definition}
  For a rewrite rule $L \rewritesto R$ and a matching $m : L \rightarrow G$, the \textit{rewrite} of $G$ into $H$ is a double pushout diagram:
  \begin{center}
    \dposquares
      {\dpoobjects{L}{I}{R}  {G}{G'}{H}}
      {\dpoarrows{}{}  {m}{}{}  {}{}}
  \end{center}
\end{definition}

\section{Partial Adhesive Categories}\label{sec:partial-adhesive}

Adamek sums up a procedure by which many categories are defined in \emph{Abstract and Concrete Categories} \cite{Adamek2009}:

\begin{quote}
  Many familiar constructs of an ``algebraic'' or ``topological'' nature have natural descriptions that can be accomplished in two steps. The first step [...] consists of defining algebraic (resp. topological) categories by means of certain functors. The second step consists of singling out full, concrete subcategories by imposing certain axioms on the objects.
\end{quote}

We shall follow this prescription to construct the category of string graphs in the next chapter. We already saw in Examples \ref{ex:adh-cats} how adhesivity is inherited by categories defined ``by means of certain functors'': namely slice, coslice, and functor categories. There is no reason for a full subcategory of an adhesive category to also be adhesive. However, we can still characterise certain subcategories of an adhesive category that inherit ``enough adhesiveness'' to do rewriting.

\begin{definition}[Partial Adhesive Category]\label{def:partial-adh-category}
  $\mathcal C$ is called a \textit{partial adhesive category} if it is a full subcategory of an adhesive category $\mathcal A$ and the embedding functor $S : \mathcal C \rightarrow \mathcal A$ preserves monomorphisms.
\end{definition}

The category $\mathcal C$ inherits the unique pushout complement property for a certain class of pushouts in $\mathcal C$, which we shall call $S$-pushouts.

\begin{definition}[$S$-spans and $S$-pushouts]\label{def:adhesive-span}
  Let $\mathcal C$ be a partial adhesive category and $S : \mathcal C \rightarrow \mathcal A$ the embedding functor. A span $\cspan{A}{f}{B}{g}{C}$ in $\mathcal{C}$ is called an \emph{$S$-span} if it has a pushout and that pushout is preserved by $S$. Such pushouts are called \emph{$S$-pushouts}.
\end{definition}

Recall that full and faithful functors reflect colimits. Since $S$ reflects \emph{all} pushouts, we could equivalently define $S$-spans in $\mathcal C$ as those spans which have a pushout \textit{reflected by $S$}.

\begin{definition}[$S$-pushout complement]\label{def:adhesive-complement}
  An \emph{$S$-pushout complement} for a pair of arrows $(b,f)$ is a pushout complement, where the following diagram is an $S$-pushout:
  \begin{center}
    \posquare{I}{L}{G'}{G}{b}{}{c}{f}
  \end{center}
  
  We call $b$ the \textit{boundary} of $L$ and $c$ the \textit{coboundary} of $L$ in $G$.
\end{definition}

\begin{lemma}\label{lem:adhesive-complement}
  If a pair of arrows $(b, f)$, where $b$ is mono, have an $S$-pushout complement, it is unique up to isomorphism.
\end{lemma}
\begin{proof}
  Let $(c,g)$ and $(c',g')$ be $S$-pushout complements. Then the following diagrams are pushouts in the adhesive category $\mathcal A$:
  \begin{center}
    \posquare{S B}{S L}{S G'}{S G}{S b}{S g}{S c}{S f}
    \qquad
    \posquare{S B}{S L}{S G''}{S G}{S b}{S g'}{S c'}{S f}
  \end{center}
  
  Since $S$ preserves monos, these are both pushout complements of $(S b, S f)$ for $S b$ mono. So the following diagram commutes in $\mathcal A$, for $\phi'$ an isomorphism:
  \begin{center}
    \csquareslant{S B}{S G'}{S G''}{S G}{S c}{S g'}{S c'}{S g}{\phi'}
  \end{center}
  
  Since $S$ is full and faithful, there exists an isomorphism $\phi : G' \rightarrow G''$ such that $S \phi = \phi'$. Replacing $\phi'$ in the above diagram yields:
  \begin{center}
    \csquareslant{S B}{S G'}{S G''}{S G}{S c}{S g'}{S c'}{S g}{S \phi}
  \end{center}
  
  Diagram (\ref{dia:unique-compl}) commutes by the faithfulness of $S$.
\end{proof}

\begin{definition}[$S$-matching]\label{def:adhesive-matching}
  For a rewrite rule $L \rewritesto R$, a monomorphism $m : L \rightarrow G$ is called an \emph{$S$-matching} if $\crun{B}{b_1}{L}{m}{G}$ has an $S$-pushout complement.
\end{definition}

Clearly if $m$ is an $S$-matching, then $S m$ is a matching. For the converse to be true, it suffices for the image of $S$ to be closed under subobjects.

\begin{definition}[$S$-rewrite]\label{def:adhesive-rewrite}
Let $L \rewritesto R := \cspan{L}{b_1}{I}{b_2}{R}$ be a rewrite rule and $m : L \rightarrow G$ be an $S$-matching. Then for $G'$ the $S$-pushout complement of $\crun{B}{b_1}{L}{m}{G}$, if the right pushout square in the following diagram exists and is an $S$-pushout:
  \begin{center}
    \dposquares
      {\dpoobjects{L}{B}{R}  {G}{G'}{H}}
      {\dpoarrows{b_1}{b_2}  {m}{}{}  {}{}}
  \end{center}
  
  Then $H$ is the $S$-rewrite of $L \rewritesto R$ at $m$.
\end{definition}

We write $H$ defined as in \ref{def:adhesive-rewrite} as $G[L \rewritesto R]_m$, or more explicitly $G[\cspan{L}{b_1}{I}{b_2}{R}]_m$.

\subsection{Example: The Category of Simple Graphs}

Partial adhesive categories should be thought of as adhesive categories, with some extra axioms imposed on the objects. In the presence of these axioms, one may need to verify by hand the relevant $S$-pushouts and $S$-pushout complements exist for a particular class of rewrite rules or matchings. In practice, this tends to be fairly straightforward. In this section, we give the derivation of these properties for the category of simple graphs.

Let $\catGr$ be the category of simple graphs, i.e. graphs where every pair of vertices is connected by at most one edge in either direction. Equivalently, a simple graph is just a binary relation from a set to itself. An object in $\catGr$ consists of a set $V$ of vertices, $E$ of edges and an injection $e : E \hookrightarrow V \times V$. A simple graph homomorphism is a pair $f_V, f_E$ such that
\begin{center}
  \begin{tikzpicture}
    \matrix (m) [cdiag] {
      E_G & V_G \times V_G \\
      E_H & V_H \times V_H \\
    };
    \path [arrs]
      (m-1-1) edge node {$e$} (m-1-2)
      (m-1-1) edge node [swap] {$f_E$} (m-2-1)
      (m-1-2) edge node {$f_V \times f_V$} (m-2-2)
      (m-2-1) edge node [swap] {$e'$} (m-2-2);
  \end{tikzpicture}
\end{center}

There is an evident embedding of $\catGr$ into $\catGraph$:
\begin{equation*}
  S : (V_G, E_G, e : E \rightarrow V) \mapsto (V_G, E_G, \pi_1 \circ e, \pi_2 \circ e)
\end{equation*}

Under the identifications $s = \pi_1 \circ e$, $t = \pi_2 \circ e$, the notions of graph homomorphism in $\catGr$ and $\catGraph$ are equivalent, so $S$ is a full subcategory embedding. $\catGr$ is a reflective subcategory of $\catGraph$, so $S$ has a left adjoint. As a right adjoint, $S$ preserves limits and, in particular, monomorphisms. Therefore $\catGr$ is a partial adhesive category.

\begin{lemma}\label{lem:gr-po-regular-mono}
  Let $\cspan{A}{m}{B}{n}{C}$ be a span in $\catGr$, where $m$ is a mono and $n$ is a regular mono. Then $m,n$ has an $S$-pushout.
\end{lemma}

\begin{proof}
  Since $S$ is full and faithful, it suffices to show for $m$ a mono, $n$ a regular mono, that $D$ defined by the following pushout in $\catGraph$ is a simple graph.
  \begin{center}
    \posquare{S A}{S B}{S C}{D}{S m}{}{S n}{}
  \end{center}
  
  We can consider $D$ to be a union of a simple graph $S B$ and another simple graph $S C$. Regular monos in $\catGr$ are precisely the full subgraph embeddings, so $S A \cong S B \cap S C$ is a full subgraph of $S C$. Consider two vertices $v,v'$ in $D$ and edges $e,e'$ such that $s(e)=s(e')=v$ and $t(e)=t(e')=v'$. The only way these edges can possibly be distinct is if $e$ is in $S B$ and $e'$ is in $S C$. Thus $v$ and $v'$ are in $S B \cap S C$. Since the intersection is a full subgraph of $S C$, $e'$ is in $S B \cap S C$, so $e = e'$.
\end{proof}

We can an $S$-matching that is a regular a monomorphism a \textit{regular $S$-matching}. We can show that when $m$ is a regular $S$-matching in $\catGr$ then the associated $S$-rewrite exists and is unique.

\begin{theorem}
  For a rewrite rule $L \rewritesto R := \cspan{L}{b_1}{I}{b_2}{R}$ and a regular $S$-matching $m : L \rightarrow G$, the associated $S$-rewrite is well-defined. That is, the following two $S$-pushout squares exist:
  \begin{equation}\label{eq:gr-dpo}
    \dposquares
      {\dpoobjects{L}{B}{R}  {G}{G'}{H}}
      {\dpoarrows{b_1}{b_2}  {m}{n}{}  {}{}}
  \end{equation}
\end{theorem}

\begin{proof}
  The existence of the left $S$-pushout square follows from the fact that any subgraph of a simple graph is also a simple graph. Since the following is a pushout along a monomorphism in $\catGraph$, it is also a pullback, by Lemma \ref{lem:adh-po-pb}.
  \begin{center}
    \posquare{S I}{S L}{S G'}{S G}{S b_1}{S f}{S n}{S m}
  \end{center}
  
  Since $S$ is full and faithful, it reflects pullbacks. So the following is a pullback in $\catGr$.
  \begin{center}
    \pbsquare{I}{L}{G'}{G}{b_1}{f}{n}{m}
  \end{center}
  
  Regular monos are stable under pullback, so $n$ is a regular mono. By Lemma \ref{lem:gr-po-regular-mono}, the right square in (\ref{eq:gr-dpo}) is an $S$-pushout.
\end{proof}

Thus, if we restrict to regular monomorphisms for matchings, DPO rewriting in the partial adhesive category $\catGr$ is well-defined, and the procedure is \textit{identical} to that in the category $\catGraph$.

\subsection{Commutation of \textit{S}-Pushouts and \textit{S}-Pushout Complements}

Often we wish to glue objects together using $S$-pushouts. If we perform an $S$-rewrite that is confined to a single component of this glued-together object, it should not matter if we first apply the rewrite then compose or if we first compose then apply the rewrite. In later sections, we shall define categories whose morphisms consist of graphs modulo a rewrite system. In such categories, this commutation of gluing and rewriting is crucial to ensuring that the composition operation well-defined. Therefore, we now prove two lemmas regarding the compatibility of $S$-pushouts, $S$-pushout complements, and $S$-rewrites.

\begin{lemma}\label{lem:pushout-pushout-compl-commute}
  $S$-pushout complements commute with $S$-pushouts. Let the following diagram be an $S$-pushout:
  \begin{center}
    \posquare{P}{H}{G}{G +_{p,q} H}{q}{i_1}{p}{i_2}
  \end{center}
  
  Assume $\crun{B}{b}{K}{m}{G}$ and $\crun{B}{b}{K}{i_1 m}{G +_{p,q} H}$ both have pushout complements,
  \begin{center}
    \posquare{B}{K}{G -_{b,m} K}{G}{b}{s}{c}{m}
    \qquad
    \posquare{B}{K}{(G +_{p,q} H) -_{b,i_1 m} K}{G +_{p,q} H}{b}{s'}{c'}{i_1 m}
  \end{center}
  
  \noindent and the $S$-span $(p,q)$ factors through $G -_{b,m} K$, i.e. there exists $p'$ such that $s p' = p$ and $(p',q)$ is an $S$-span.
  
  Then, for a second $S$-pushout:
  \begin{center}
    \posquare{P}{H}{G}{(G -_{b,m} K) +_{p',q} H}{q}{j_1}{p'}{j_2}
  \end{center}
  
  \noindent there is an isomorphism $\phi : (G +_{p,q} H) -_{b,i_1 m} K \overset{\sim}{\rightarrow} (G -_{b,m} K) +_{p',q} H$, commuting with the coboundaries $c$ and $c'$ of $K$ in $G$ and $G +_{p,q} H$ respectively.
  \begin{equation}\label{dia:compatible-with-coboundary}
    \csquare{B}{(G +_{p,q} H) -_{b,i_1 m} K}
            {G -_{b,m} K}{(G -_{b,m} K) +_{p',q} H}
            {c'}{j_1}{c}{\phi}
  \end{equation}
\end{lemma}

\begin{proof}
  The proof follows from the associativity of pushouts and the uniqueness of pushout complements. First, note that, in the following diagram, [1] commutes and is a pushout because $s p' = p$:
  \begin{center}
    \begin{tikzpicture}
    \node at (1,0) {\small [1]};
    \matrix (m) [cdiag] {
      & P           &  H            \\
    B & G -_{b,m} K &               \\
    K & G           &  G +_{p,q} H  \\
    };
    \path [arrs]
     (m-1-2) edge node {$q$} (m-1-3)
     (m-1-2) edge node {$p'$} (m-2-2)
     (m-2-1) edge node {$c$} (m-2-2)
     (m-2-1) edge node {$b$} (m-3-1)
     (m-3-1) edge node [swap] {$m$} (m-3-2)
     (m-2-2) edge node {$s$} (m-3-2)
     (m-1-3) edge node {$i_2$} (m-3-3)
     (m-3-2) edge node [swap] {$i_1$} (m-3-3);
    \NWbracket{(m-3-2)}
    \NWbracket{(m-3-3)}
    \end{tikzpicture}
  \end{center}
  Next, we make the two pushouts in the opposite order.
  \begin{center}
    \begin{tikzpicture}
    \node at (-0.75,-0.9) {\small [2]};
    \matrix (m) [cdiag] {
      & P           & H                        \\
    B & G -_{b,m} K & (G -_{b,m} K) +_{p',q} H \\
    K &             & Q                        \\
    };
    \path [arrs]
     (m-1-2) edge node {$q$} (m-1-3)
     (m-1-2) edge node {$p'$} (m-2-2)
     (m-2-1) edge node {$c$} (m-2-2)
     (m-2-1) edge node {$b$} (m-3-1)
     (m-3-1) edge node [swap] {$k_1$} (m-3-3)
     (m-2-3) edge node {$k_2$} (m-3-3)
     (m-2-2) edge node {$j_1$} (m-2-3)
     (m-1-3) edge node {$j_2$} (m-2-3);
    \NWbracket{(m-2-3)}
    \NWbracket{(m-3-3)}
    \end{tikzpicture}
  \end{center}
  
  By associativity of pushouts, there exists an isomorphism $\psi$ such that:
  \begin{center}
    \begin{tikzpicture}
      \matrix (m) [cdiag] {
        G +_{p,q} H & H \\
        K & Q \\
      };
      \path [arrs]
        (m-1-2) edge node [swap] {$i_2$} (m-1-1)
        (m-2-1) edge node {$i_1 m$} (m-1-1)
        (m-1-2) edge node {$k_2 j_2$} (m-2-2)
        (m-2-1) edge node [swap] {$k_1$} (m-2-2)
        (m-1-1) edge node {$\psi$} (m-2-2);
    \end{tikzpicture}
  \end{center}
  
  Since $Q$ is only defined up to isomorphism, we are free to take $\psi = 1_{G +_{p,q} H}$, in which case $k_1 = i_1 m$. Then, square [2] from above becomes:
  \begin{center}
    \posquare{B}{(G -_{b,m} K) +_{p',q} H}
             {K}{G +_{p,q} H}
             {j_1 c}{i_1 m}{b}{k_2}
  \end{center}
  
  Compare this to the definition of $(G +_{p,q} H) -_{b,i_1 m} K$ as a pushout complement:
  \begin{center}
    \posquare{B}{(G +_{p,q} H) -_{b,i_1 m} K}
             {K}{G +_{p,q} H}
             {c'}{i_1 m}{b}{s'}
  \end{center}
  
  The result then follows from uniqueness of pushout complements.
\end{proof}


\begin{theorem}\label{thm:pushout-and-rewrite}
  $S$-adhesive rewrites commute with $S$-adhesive pushouts. Let $m : L \rightarrow G$ be an $S$-matching of $L \rewritesto R := \cspan{L}{b_1}{I}{b_2}{R}$. The rewrite is computed as the double pushout:
  \begin{center}
    \begin{tikzpicture}
    \matrix (m) [cdiag] {
    L & B           & R                   \\
    G & G -_{b,m} L & G[L \rewritesto R]_m \\
    };
    \path [arrs]
     (m-1-2) edge node [swap] {$b_1$} (m-1-1)
     (m-1-2) edge node {$b_2$} (m-1-3)

     (m-2-2) edge node {$s$} (m-2-1)
     (m-2-2) edge node [swap] {$s'$} (m-2-3)

     (m-1-1) edge node [swap] {$m$} (m-2-1)
     (m-1-2) edge node {$c$} (m-2-2)
     (m-1-3) edge node {$m'$} (m-2-3);
    \NEbracket{(m-2-1)}
    \NWbracket{(m-2-3)}
    \end{tikzpicture}
  \end{center}
  Let $(p,q)$, $(\widehat p, q)$ and $(p', q)$ be three $S$-spans, such that:
  \begin{equation}\label{dia:compatible}
    \begin{tikzpicture}
    \matrix (m) [cdiag] {
    G                   &   &   \\
    G -_{b_1,m} L       & P & H \\
    G[L \rewritesto R]_m &   &   \\
    };
    \path [arrs]
      (m-2-2) edge [bend right] node [swap] {$p$} (m-1-1)
      (m-2-2) edge node [swap] {$\widehat p$} (m-2-1)
      (m-2-2) edge [bend left] node {$p'$} (m-3-1)
      (m-2-1) edge node {$s$} (m-1-1)
      (m-2-1) edge node [swap] {$s'$} (m-3-1)
      (m-2-2) edge node {$q$} (m-2-3);
    \end{tikzpicture}
  \end{equation}
  Then, for the pushout injection $i_1 : G \rightarrow G +_{p,q} H$, if $i_1 m$ is mono, the following is an isomorphism:
  \[ (G[L \rewritesto R]_m) +_{p',q} H \cong
     (G +_{p,q} H)[L \rewritesto R]_{i_1 m} \]
\end{theorem}

\begin{proof}
  $G[L\rewritesto R]_m$ is the pushout of $R$ and $G -_{b_1,m} L$ along $B$, so by uniqueness of pushout complements, we can choose $(G -_{b_1,m} L)$ to be equal to $((G[L \rewritesto R]_m) -_{b_2,m'} R)$, for the same coboundary $c$. Then, by two applications of Lemma \ref{lem:pushout-pushout-compl-commute}, we can choose $(G -_{b_1,m} L) +_{\widehat p,q} H = ((G[L \rewritesto R]_m) -_{b_2,m'} R) +_{\widehat p,q} H$ as the pushout complement of both of the following squares.
  \begin{center}
    \begin{tikzpicture}
    \matrix (m) [cdiag] {
    L           & B                        & R \\
    G +_{p,q} H & (G -_{b_1,m} L) +_{\widehat p,q} H &
         G[L \rewritesto R] +_{p', q} H \\
    };
    \path [arrs]
     (m-1-2) edge (m-1-1)
     (m-1-2) edge (m-1-3)

     (m-2-2) edge (m-2-1)
     (m-2-2) edge (m-2-3)

     (m-1-1) edge node [swap] {$i_1 m$} (m-2-1)
     (m-1-2) edge node {$c'$} (m-2-2)
     (m-1-3) edge (m-2-3);
    \NEbracket{(m-2-1)}
    \NWbracket{(m-2-3)}
    \end{tikzpicture}
  \end{center}
  
  Note that $c'$ becomes the coboundary for both squares because diagram (\ref{dia:compatible-with-coboundary}) commutes. This is then exactly the computation of the rewrite $(G +_{p,q} H)[L \rewritesto R]_{i_1 m}$. The theorem holds because $S$-rewrites are unique up to isomorphism.
\end{proof}

We shall use these two theorems throughout this dissertation to show that rewriting is compatible with several notions of composing graphs.

\chapter{String Graphs and Monoidal Theories}\label{ch:string-graphs}

In this chapter, we come to one of the primary contributions of this dissertation: the formalisation of the diagrammatic language of monoidal categories using certain typed graphs called \textit{string graphs}. In the previous chapter, we introduced the notion of a partial adhesive category. Our primary reason for doing so was to enable double-pushout graph rewriting in the category of string graphs.

Before passing to string graphs, we will look at the category of typed graphs, defined as a slice category $\catGraph / G_T$. Objects of the slice category are pairs $(G, \tau_G : G \rightarrow G_T)$. These should be thought of as a graph $G$ along with a map $\tau_G$ giving a \textit{type} in $G_T$ to every vertex and edge in $G$. Morphisms are simply graph homomorphisms respecting this type map: $f : G \rightarrow H$ such that $\tau_H \circ f = \tau_G$. It may seem odd at first that the type map is a graph homomorphism rather than simply a pair of functions defined on vertices and edges. However, the first can encode the latter, so these ``unresticted'' vertex and edge typings are merely a special case of homomorphic graph typings. As an example, suppose we fix a set $X$ of vertex types and have only one edge type. Then $G_T$ can be defined as the connected graph whose vertices are the elements of $X$. For instance, if we let $X = \{ \textbf{black}, \textbf{white}, \textbf{grey} \}$, we can form a connected graph:
\ctikzfig{typegraph_colours}

Since $G_T$ is isomorphic to the totally connected directed graph with $3$ vertices, any function $\tau_V : V_G \rightarrow V_{G_T}$ extends uniquely to a graph homomorphism $\tau : G \rightarrow G_T$. We can think of the fibres of $\tau_V$ (i.e. the inverse images $\tau_V^{-1}(\textbf{black})$, $\tau_V^{-1}(\textbf{white})$, and $\tau_V^{-1}(\textbf{grey})$) as sets of vertices in $G$ that are coloured ``black'', ``white'', or ``grey'' respectively. We can represent a graph with coloured vertices as a pair consisting of a graph $G$ and a typing function $\tau$.
\[
\beginpgfgraphicnamed{coloured_graph}
\InputIfFileExists{coloured_graph.tikz}{}{\input{./figures/coloured_graph.tikz}}
\endpgfgraphicnamed :=
\left(
  \ \ %
\beginpgfgraphicnamed{coloured_graph_G}
\InputIfFileExists{coloured_graph_G.tikz}{}{\input{./figures/coloured_graph_G.tikz}}
\endpgfgraphicnamed\ \ ,\ \ 
  \tau_V :: \left\{
  \begin{matrix}
    v_1 \mapsto \textbf{white}, & v_2 \mapsto \textbf{grey}, \\
    v_3 \mapsto \textbf{black}, & v_4 \mapsto \textbf{grey}  \\
  \end{matrix} \right\}
\right)
\]

Typegraphs can express edge typing as well. For instance, to add a set of edge types $Y = \{ +,\ -\}$, we can simply add a copy of $Y$ for every pair of vertices $(v_1, v_2)$ in $G_T$ connecting $v_1$ to $v_2$.
\ctikzfig{typegraph_colours_and_edges}

One can even restrict which vertex types can be connected by which edge types by taking subgraphs of $G_T'$. In the coming sections, we will use this restriction to make sure that our diagrams are well-typed, i.e. composition of ``boxes'' in the language of string diagrams should respect the types on ``wires''.

\section{String Graphs}\label{sec:string graphs}

String diagrams consist of \textit{boxes}, which represent morphisms in a monoidal category, and \textit{wires}, which are used to connect boxes together. We turn string diagrams into a \textit{string graphs}, which are typed graphs with two distinguished kinds of vertex: \textit{box-vertices} and \textit{wire-vertices}. As the name suggests, \textit{box-vertices} represent the boxes (i.e. morphisms/generators) in the diagrammatic language. These represent the ``logical'' or ``semantic'' vertices of the string graph. An important characteristic of wires, which distinguishes them from normal edges in a graph, is that they are not required to have boxes at both ends and they can be connected to themselves to form \textit{circles}. For that reason, we will represent wires as chains of special ``place-holder'' vertices called \textit{wire-vertices}.
\ctikzfig{wire}

Note how the wire-vertices carry the type of the wire. Representing boxes as \textit{box-vertices}, it is possible to translate any string diagram into a string graph.
\ctikzfig{string_diagram_to_graph}

Edge types (not shown) are used to keep track of the ordering of inputs and outputs to the boxes (see Definition \ref{def:derived-typegraph}).

The number of wire-vertices making up a wire is irrelevant, so for the purposes of representing a string diagram, the following two string graphs are equivalent:
\ctikzfig{string_graph_equiv}

If we were to treat the two graphs above as 1D simplicial complexes that define topological graphs, then the geometric realisation of the wires on the left are homeomorphic to those of the wires on the right. For that reason, the two graphs above are called \textit{wire-homeomorphic}. We can formalise the notion of wire-homeomorphism as a confluent, terminating graph rewrite system and prove that string graphs, up to wire-homeomorphism can be used to construct free monoidal categories. This construction is essentially the graph version of the topological construction described in Section \ref{sec:js-construction}.

For a monoidal signature $T$, we can define a category $\catSGraph_T$ of string graphs with generators taken from $T$. We do this by turning the monoidal signature $T$ into a typegraph $G_T$ and defining $\catSGraph_T$ is a full subcategory of $\catGraph/G_T$. We shall then prove that the embedding of $\catSGraph_T$ into $\catGraph/T$ preserves monos, so $\catSGraph_T$ is a partial adhesive category.

First, we show the construction of $G_T$. Recall that a monoidal signature $T = (O,M,\dom,\cod)$ consists of sets $O,M$ and functions $\dom : M \rightarrow w(O), \cod : M \rightarrow w(O)$ into the set $w(O)$ of lists of elements of $O$.

\begin{definition}\label{def:derived-typegraph}
  For a monoidal signature $T$, the \textit{derived typegraph} $G_T$ of $T$ has vertices $O + M$ and the following edges:
  \begin{itemize}
    \item a self-loop $\textrm{mid}_{X}$ for every $X \in O$,
    \item an edge $\textrm{in}_{f,i}$ for $f \in M$, $0 \leq i < \textrm{Length}(\dom(f))$, connecting $\dom(f)[i]$ to $f$, and
    \item an edge $\textrm{out}_{f,j}$ for $f \in M$, $0 \leq j < \textrm{Length}(\cod(f))$, connecting $f$ to $\cod(f)[j]$.
  \end{itemize}
\end{definition}

\begin{example}
  Let $T$ be the following monoidal signature:
  \ctikzfig{box_set}
  Then, the derived typegraph $G_T$ is:
  \ctikzfig{typegraph_gt}
\end{example}

\begin{definition}
  Let $(G, \tau : G \rightarrow G_T)$ be an object in $\catGraph/G_T$, for a monoidal signature $T$. By definition of $G_T$, the vertices of $G_T$ are $O + M$. A vertex $v \in V_G$ is called a \textit{box-vertex} of $\tau_V(v) \in M$. It is called a \textit{wire-vertex} if $\tau_V(v) \in O$. Let $B(G) \subseteq V_G$ be the set of all box-vertices, and $W(G)$ be the set of wire-vertices.
\end{definition}

Note that since we have omitted self-loops on box-vertices in $G_T$, any path between two box-vertices in a $G_T$-typed graph must pass \textit{through} a wire-vertex. This is important to the definition of string graphs, as the ``object'' types from the signature $T$ are carried by wire-vertices.

There are two restrictions that we place on typed graphs $(G,\tau)$ to make string graphs. The first is that wires should not split or merge. Namely, any wire-vertex in $G$ should have at most one in-edge and one out-edge. In other words, for a graph $G = (V_G,E_G,s,t)$, the restrictions of $s$ and $t$ to the wire vertices $W(G)$ are both monomorphisms, i.e. $s',t'$ defined by the pullbacks below are monomorphisms:
\begin{center}
  \begin{tikzpicture}
    \matrix (m) [cdiag] {
      E'  & W(G) & E'' \\
      E_G & V_G  & E_G \\
    };
    \path [arrs]
      (m-1-1) edge node {$s'$} (m-1-2)
      (m-1-1) edge (m-2-1)
      (m-1-2) edge [right hook-latex] (m-2-2)
      (m-2-1) edge node [swap] {$s$} (m-2-2)
      (m-1-3) edge node [swap] {$t'$} (m-1-2)
      (m-1-3) edge (m-2-3)
      (m-2-3) edge node {$t$} (m-2-2);
      \SEbracket{(m-1-1)}
      \SWbracket{(m-1-3)}
  \end{tikzpicture}
\end{center}

The second condition is that the incident edges of a box-vertex $b$ in $G$ should match those of its type, $\tau_V(b)$. In other words, there should be the same number of inputs and outputs to $b$ as there are to $\tau_V(b)$. We formalise this condition by introducing the notion of a \textit{local isomorphism}.

\begin{definition}
  For a vertex $v \in V_G$, the \textit{edge neighbourhood} of $v$ is the set of edges $N(v) := s^{-1}(v) \cup t^{-1}(v)$.
\end{definition}

Fix a graph homomorphism $f : G \rightarrow H$. Then for a vertex $v$ in $G$ and an adjacent edge $e$, the edge $f(e)$ is adjacent to $f(v)$. Thus $f_E(N(v)) \subseteq N(f(v))$. Let $f^v : N(v) \rightarrow N(f(v))$ be the function defined by $f^v(e) = f_E(e)$, for $e \in N(v)$.

\begin{definition}
  For $G_T$-typed graphs $(G,\tau_G)$, $(H, \tau_H)$, a typed graph homomorphism $f : G \rightarrow H$ is called a \textit{local isomorphism} if for every box $b \in B(G)$, $f^b : N(b) \rightarrow N(f(b))$ is a bijection.
\end{definition}

In particular for a $G_T$-typed graph $(G,\tau)$, the typing map $\tau$ can be considered as a typed graph homomorphism from $(G,\tau)$ to $(G_T,1_{G_T})$. Thus, we can require that it be a local isomorphism.

\begin{definition}
  A $G_T$-typed graph $(G, \tau)$ is called a \textit{string graph} if $\tau$ is a local isomorphism and every wire-vertex in $G$ has at most one in-edge and one out-edge. The category $\catSGraph_T$ is the full subcategory of $\catGraph/T$ whose objects are string graphs.
\end{definition}

Since the typing maps $\tau$ in $\catSGraph_T$ are local isomorphisms, we can show that \textit{every} arrow in $\catSGraph_T$ is a local isomorphism.

\begin{proposition}\label{prop:all-local-iso}
 Every arrow in $\catSGraph_T$ is a local isomorphism.
\end{proposition}

\begin{proof}
  Let $(G, \tau_G)$, $(H, \tau_H)$ be $G_T$-typed graphs. By definition $\tau_G$ and $\tau_H$ are both local isomorphisms. For any $f : (G,\tau_G) \rightarrow (H,\tau_H)$ in $\catGraph/G_T$, the following diagram commutes:
  \begin{center}
    \ctri{G}{G_T}{H}{\tau_G}{f}{\tau_H}
  \end{center}
  
  Thus, for any $v$ in $G$ we get this triangle in \catSet:
  \begin{center}
    \ctri{N(v)}{N(\tau_G(v))}{N(f(v))}{\tau_G^v}{f^v}{\tau_H^{f(v)}}
  \end{center}
  
  Since $\tau_G^v$ and $\tau_H^{f(v)}$ are both bijections, $f^v$ is a bijection, so $f$ is a local isomorphism.
\end{proof}

\begin{proposition}\label{prop:sgraph-monos}
  A morphism in $\catSGraph_T$ is a monomorphism iff it is injective.
\end{proposition}
\begin{proof}
  Suppose $m : G \rightarrow H$ in $\catSGraphTG$ is not injective. $m$ is a local isomorphism, so if $m$ takes two distinct edges $e_1$ and $e_2$ to the same edge, then the adjacent vertices of $e_1$ and $e_2$ must also be distinct.
  
  To show that $m$ cannot be mono, we will define a string graph $K$ and distinct maps $f,g : K \rightarrow G$ such that $m f = m g$. If $m$ takes two distinct box-vertices $v_1$, $v_2$ in $G$ to a single box-vertex in $H$, then let $K$ be the subgraph of $G$ consisting of just $v_1$ and its neighbourhood. If $m$ takes two distinct wire-vertices to a single wire-vertex in $H$, then let $K$ be a $T$-string graph consisting of a single, isolated wire-vertex. In either case, there are two distinct maps $f,g$ such that $m f = m g$.
\end{proof}

\begin{corollary}\label{cor:sgraph-partial-adhesive}
  $\catSGraph_T$ is a partial adhesive category.
\end{corollary}

\subsection{\textit{S}-pushouts of String Graphs}

For the constructions to come, it is useful to characterise the $S$-pushouts in $\catSGraph_T$. To simplify matters, it suffices to characterise the $S$-pushouts of spans $\cspan{G}{m}{K}{n}{H}$ where $m$ and $n$ are both monomorphisms. We show that a span of monos is an $S$-span if and only if it is \textit{boundary-coherent}. Before we define boundary-coherence, we need the notion of a boundary.

\begin{definitions}
  If a wire-vertex has no in-edges, it is called an \emph{input}. We write the set of inputs of a string graph $G$ as $\In(G)$. Similarly, a wire-vertex with no out-edges is called an \emph{output}, and the set of outputs is written $\Out(G)$. The inputs and outputs define a string graph's \emph{boundary}, written $\Bound(G)$. If a boundary vertex has no in-edges and no out-edges, (it is both and input and output) it is called an \emph{isolated wire-vertex}. An string graph consisting of only isolated wire-vertices is called a \emph{point graph}.
\end{definitions}

By abuse of notation, we may treat $\In(G)$, $\Out(G)$, and $\Bound(G)$ as sets or as point graphs. The intended usage will be clear from context.

\begin{definition}\label{def:boundary-coherent}
  A pair of morphisms $f : K \rightarrow G$, $g : K \rightarrow H$ in $\catSGraph_T$ is called \emph{boundary-coherent} if:
  \begin{enumerate}
    \item for all $v \in \In(K)$ at least one of $f(v)$ and $g(v)$ is an input, and
    \item for all $v \in \Out(K)$ at least one of $f(v)$ and $g(v)$ is an output. 
  \end{enumerate}
\end{definition}
 
\begin{theorem}
  A span of monomorphisms $\cspan{G}{m}{I}{n}{H}$ in $\catSGraph_T$ has an $S$-pushout if and only if $m$ and $n$ are boundary-coherent.
\end{theorem}
 
\begin{proof}
  For ($\Leftarrow$), it suffices to show that $K$ is a string graph, for the following pushout in $\catGraph/G_T$:
  \begin{equation}\label{eq:boundary-coherence-posquare}
    \posquare{S I}{S G}{S H}{K}{S m}{f}{S n}{g}
  \end{equation}
  
  Since $K$ is isomorphic to the union of $S G$ and $S H$, the images of $f$ and $g$ cover $K$. Let $f(b) \in B(K)$ be a box-vertex in the image of $f$ but not $g$. Then the neighbourhood of $b$ is identical to that of $f(b)$, so the inclusion of edges $f^b : N(b) \rightarrow N(f(b))$ is a bijection. Similarly, for $g(b) \in B(K)$ not in the image of $f$, $g^b : N(b) \rightarrow N(g(b))$ is a bijection. Finally, pick a box-vertex $b \in B(I)$, then since $m$ and $n$ are local isomorphisms, $S m$ and $S n$ restrict to bijections on $N(b) \subseteq E_I$. Thus, all of the edges in $N(Sm(b))$ are identified with edges in $N(Sn(b))$ in $K$, and the inclusions $f$ and $g$ restrict to bijections on $N(Sm(b))$ and $N(Sn(b))$ respectively. So $f$ and $g$ are local isomorphisms. Since the images of $f$ and $g$ cover $K$, it follows that the typing map $\tau_K$ is a local isomorphism (cf. the proof of Proposition \ref{prop:all-local-iso}).
  
  Now, suppose a wire-vertex $v \in W(K)$ has out-edges $e_1$, $e_2$. The only way these can possibly be distinct is if $v$ is in the image of both $f$ and $g$. Then, there must be a $v'$ such that $f \circ Sn(v') = g \circ Sm(v') = v$. If $v'$ is not an output in $I$, then the out-edges of $m(v') \in B(G)$ and $n(v') \in B(H)$ must be in the images of $m$ and $n$, respectively, so $e_1 = e_2$. If it is an output, then at least one of $m(v')$, $n(v')$ must be an output, so $v$ has at most one out-edge. We can show similarly that $v$ must have at most one in-edge. So $K$ is a string graph.
  
  For ($\Rightarrow$), suppose the span $\cspan{G}{m}{I}{n}{H}$ is not boundary-coherent. If the pushout $K$ given by (\ref{eq:boundary-coherence-posquare}) is not a string graph, then either (a) the span $m,n$ does not have a pushout or (b) $m,n$ does have a pushout, but it is not preserved by $S$. In either case, the span does not have an $S$-pushout, so it suffices to show that $K$ is not a string graph. If the span $m,n$ is not boundary-coherent, there exists wire-vertex $v$ in $I$ such that either (a) $v$ is an input and $f(v)$ and $g(v)$ both have in-edges or (b) $v$ is an output and $f(v)$ and $g(v)$ both have out-edges. If (a) is true, then the image of $v$ will have two distinct in-edges in $K$. If (b) is true, it will have two distinct out-edges. In either case, $K$ is not a string graph.
\end{proof}
 
\begin{example}
  Consider the following span of string graphs, which is not boundary-coherent:
  \ctikzfig{not_boundary_coherent}
  \noindent If we push out this span in $\catGraph/G_T$, we get:
  \ctikzfig{not_boundary_coherent2}
  \noindent which is clearly not a string graph. Therefore this span does not have an $S$-pushout.
\end{example}

Note that the empty graph is the initial object in $\catSGraph_T$, so we have an easy corollary.

\begin{corollary}\label{cor:spgraph-coproducts}
  $\catSGraph_T$ has finite coproducts and $S$ preserves them.
\end{corollary}

\begin{proof}
  The proof follows from the observation that for any string graphs $G,H$, the span of initial arrows $\cspan{G}{!}{\emptyset}{!}{H}$ is trivially boundary-coherent.
\end{proof}

One particularly important type of boundary-coherent span is a \textit{plugging}. These are used to ``plug'' the inputs of one graph into the outputs of another graph.

\begin{definition}\label{def:plugging}
  A boundary-coherent span $\cspan{G}{m}{P}{n}{H}$ where $P$ is a point graph and $G$ and $H$ contain no isolated vertices is called a \textit{plugging}.
\end{definition}

Recall that wire-vertices in a point graph are both inputs and outputs. Thus boundary-coherence forces the images $m(p)$ and $n(p)$ of a wire-vertex $p \in P$ to have opposite polarities. That is to say, exactly one of $m(p),n(p)$ is an input and exactly one is an output.

\begin{example}
  The following $S$-span defines a plugging:
  \ctikzfig{plugging_ex}
  
  Pushing out the span yields a string graph with the two smaller graphs plugged together:
  \ctikzfig{plugging_ex_result}
\end{example}

\section{Rewriting with String Graphs}\label{sec:rewriting-string-graphs}

If a string graph contains no isolated points, then the set of inputs and outputs is disjoint, i.e. $\Bound(G) \cong \In(G) + \Out(G)$. We will define string graph rewrite rules as pairs of string graphs with no isolated wire-vertices whose inputs and outputs are in bijection. Such a pair uniquely induces a span $\cspan{L}{b_1}{I}{b_2}{R}$ for $I \cong \In(L) + \Out(L) \cong \In(R) + \Out(R)$.

\begin{definition}
  A \textit{string graph rewrite rule} $L \rewritesto R$ is a span $\cspan{L}{b_1}{I}{b_2}{R}$ where:
  \begin{enumerate}
    \item $L$ and $R$ contain no isolated wire-vertices,
    \item $\In(L) \cong \In(R)$, $\Out(L) \cong \Out(R)$,
    \item $I \cong \In(L) + \Out(L) \cong \In(R) + \Out(R)$, and
    \item the following diagram commutes for $b_1$ and $b_2$ the induced maps of the coproduct inclusions $i,j$ and $i',j'$ respectively:
    \begin{center}
      \begin{tikzpicture}
        \matrix (m) [cdiag] {
          &  \In(L) &   & \In(R)  &   \\
        L &         & I &         & R \\
          & \Out(L) &   & \Out(R) &   \\
        };
        \path [arrs]
          (m-2-3) edge [dashed] node [swap] {$b_1$} (m-2-1)
          (m-2-3) edge [dashed] node {$b_2$} (m-2-5)
          (m-1-2) edge [left hook-latex] (m-2-1)
          (m-3-2) edge [left hook-latex] (m-2-1)
          (m-1-4) edge [right hook-latex] (m-2-5)
          (m-3-4) edge [right hook-latex] (m-2-5)
          (m-1-2) edge node {$i$} (m-2-3)
          (m-3-2) edge node [swap] {$j$} (m-2-3)
          (m-1-4) edge node [swap] {$i'$} (m-2-3)
          (m-3-4) edge node {$j'$} (m-2-3)
          
          (m-1-2) edge node {$\sim$} (m-1-4)
          (m-3-2) edge node [swap] {$\sim$} (m-3-4);
      \end{tikzpicture}
    \end{center}
  \end{enumerate}
\end{definition}

Since $L$ contains no isolated wire-vertices, the images of $\In(L)$ and $\Out(L)$ are disjoint, so $b_1$ is injective. For the same reason, $b_2$ is also injective.

\begin{example}
  Let $f : A \rightarrow A \otimes A$ and $g : A \otimes A \rightarrow A$ be morphisms in a monoidal category such that $f \circ g = 1_A$. This can be expressed as an equation between string diagrams:
  \ctikzfig{string_diagram_eq_ex}
  \noindent ...or as a rewrite rule between string graphs:
  \ctikzfig{string_rewrite_rule_ex}
\end{example}

\begin{theorem}\label{thm:all-monos-matchings}
  Let $L \rewritesto R := \cspan{L}{b_1}{I}{b_2}{R}$ be a string graph rewrite rule. \textit{Any} monomorphism $m : L \rightarrow G$ is an $S$-matching.
\end{theorem}

\begin{proof}
  Since $m$ is a local isomorphism and $b_1$ covers the boundary of $L$, $m$ satisfies the no-dangling-edges condition, so it has a pushout complement $G'$ in $\catGraph/G_T$. For $m$ to be an $S$-matching, it suffices to show that $G'$ is a string graph. $I$ consists only of wire-vertices, so by the no-dangling-edges condition, the adjacent edges of the box-vertices in $G'$ are the same as they were in $G$, so the typing map $\tau_{G'}$ is still a local isomorphism. The fact that $G'$ is a string graph then follows from it being a subgraph of $G$.
\end{proof}

\subsection{Wires and Wire-Homeomorphism}\label{sec:homeomorphism}

String graphs are meant to be the discrete version of (topological) string diagrams. In string diagrams, wires can be thought of as copies of the unit interval $[0,1] \subset \mathbb R$, considered as an oriented manifold. Boxes are distinguished points, to which we ascribe semantic meaning. Clearly if we replace a wire in a string diagram with a homeomorphic wire, we get the same string diagram. In other words, the meaning of a string diagram is unaffected by shortening or lengthening wires.

A \textit{simple chain} is a connected, acyclic graph where each vertex has at most one in-edge and one out-edge. A vertex in a chain with only an in-edge or only an out-edge is called an \textit{endpoint}. A \textit{simple cycle} is connected graph where each vertex has exactly one in-edge and one out-edge.

\begin{definition}\label{def:wire}
  For a string graph $G$, a \emph{closed wire} $\widehat W \subseteq G$ is a simple cycle of wire-points or a simple chain such that
  \begin{enumerate}
    \item the endpoints of $\widehat W$ are either box-vertices or in the boundary of $G$, and
    \item all other vertices in $\widehat W$ are wire-vertices.
  \end{enumerate}
\end{definition}

It is worth noting that wires are not necessarily string graphs, as their typing function need not be a local isomorphism at the endpoints.

\begin{example}
  A string graph and its $4$ closed wires:
  \ctikzfig{string_graph_chains}
\end{example}

Every wire-vertex in a wire will be of a single type in $G_T$. This is called the \textit{wire type}.

\begin{definition}
  Two wires $\widehat W$ and $\widehat W'$ are said to be \textit{homeomorphic} if they have the same wire type and (a) they are both simple cycles, or (b) they are simple chains where the endpoints and edges adjacent to the endpoints are of the same type.
\end{definition}

Informally, we justify this terminology by noting that, if we consider $\widehat W$ and $\widehat W'$ as simplicial complexes, then for homeomorphic wires, the geometric realisations of the two complexes as topological spaces (i.e. as a circle or the unit interval) are homeomorphic.

\begin{definition}\label{def:homeo-with-manifolds}
  Two string graphs $G$ and $G'$ are called \emph{wire-homeomorphic} if $G'$ can be obtained from $G$ by replacing any number of closed wires $\widehat W$ with homeomorphic wires $\widehat W'$.
\end{definition}

For any string graph $G$, there is a unique smallest wire-homeomorphic graph $G\normalised$ obtained by replacing every wire in $G$ with a homeomorpic wire containing a single wire-vertex. We justify the notation $G\normalised$ by formalising wire-homeomorphism using a rewrite system $\mathbb H$.

\begin{definition}\label{def:homeo-rewrite-system}
  For a monoidal signature $T = (O, M, \dom, \cod)$, the rewrite system $\mathbb H$ is defined as follows. For every $X \in O$, we define a loop contraction rule $h^L_X$ and a wire contraction rule $h^W_X$:
  \ctikzfig{wire_homeo1}
  
  For every $f \in M$ and $0 \leq i < \textrm{Length}(\dom(f))$, $0 \leq j < \textrm{Length}(\cod(f))$, we define an input contraction rule $h^I_{f,i}$ and an output contraction rule $h^O_{f,j}$:
  \ctikzfig{wire_homeo2}
\end{definition}

\begin{proposition}
  Two string graphs $G,H$ are wire-homeomorphic if and only if $G \rewriteequiv_{\mathbb H} H$.
\end{proposition}

\begin{proof}
  For ($\Leftarrow$), we observe the all of the rules in $\mathbb H$ leave the endpoints of a wire fixed, while decreasing the number of other wire-points. For ($\Rightarrow$), it suffices to show that $\mathbb H$ lets us increase or decrease the number of wire-vertices in any wire in $G$. Consider the types of wire $\widehat W$ that can occur in $G$:
  \begin{enumerate}
    \item a simple cycle,
    \item a chain starting with a boundary and ending with a boundary,
    \item a chain ending with a box-vertex, or
    \item a chain starting with a box-vertex.
  \end{enumerate}
  
  These alternatives are not mutually exclusive, but they are exhaustive. In the case of (1.) the rules $h^L_X$ or $h^W_X$ can always decrease the size of $\widehat W$ if applied forwards, and increase the size of $\widehat W$ if applied backwards. In the case of (2.) apply $h^W_X$, for (3.) apply $h^I_{f,i}$, and for (4.) apply $h^O_{f,j}$.
\end{proof}

\begin{lemma}
  The rewrite system $\mathbb H$ is confluent (up to graph isomorphism) and terminating.
\end{lemma}

\begin{proof}
  Termination comes from observing that each contraction rule strictly decreases the total number of wire-vertices. Confluence follows from noting that any forward-directed rewrite procedure starting with $G$ terminates at the unique minimal wire-homeomorphic graph $G\normalised$.
\end{proof}

\begin{example}
  Normalising a string graph with respect to $\mathbb H$:
  \ctikzfig{normalising_wrt_homeo}
  Note that there are multiple ways this string graph could be normalised, but since $\mathbb H$ is confluent, the end result will always be the same.
\end{example}

\section{Cospans over the Category of String Graphs}

For any category $\mathcal C$ with pushouts, we can form its \textit{cospan bicategory} $\Csp(\mathcal C)$. The $0$-cells of $\Csp(\mathcal C)$ are the objects from $\mathcal C$, the $1$-cells are cospans $\ccospan{X}{f}{F}{g}{Y}$, and the $2$-cells are cospan homomorphisms. A cospan homomorphism from $\ccospan{X}{f}{F}{g}{Y}$ to $\ccospan{X}{f'}{F'}{g'}{Y}$ is a morphism $\alpha : F \rightarrow F'$ in $\mathcal C$ that commutes with the cospan maps.
\begin{center}
  \begin{tikzpicture}
    \matrix (m) [cdiag, row sep=5mm] {
        & F  &   \\
      X &    & Y \\
        & F' &   \\
    };
    \path [arrs]
      (m-2-1) edge node {$f$} (m-1-2)
      (m-2-3) edge node [swap] {$g$} (m-1-2)
      (m-2-1) edge node [swap] {$f'$} (m-3-2)
      (m-2-3) edge node {$g'$} (m-3-2)
      (m-1-2) edge node {$\alpha$} (m-3-2);
  \end{tikzpicture}
\end{center}

$2$-cell composition is the usual composition of morphisms in $\mathcal C$. For two cospans $\ccospan{X}{f}{F}{g}{Y}$ and $\ccospan{Y}{h}{G}{i}{Z}$, the composition is formed by pushing out over $Y$.
\begin{center}
  \begin{tikzpicture}
    \matrix (m) [cdiag] {
        &   & Y         &   &   \\
      X & F &           & G & Z \\
        &   & G \circ F &   &   \\
    };
    \path [arrs]
      (m-2-1) edge node {$f$} (m-2-2)
      (m-1-3) edge node [swap] {$g$} (m-2-2)
      (m-1-3) edge node {$h$} (m-2-4)
      (m-2-5) edge node [swap] {$i$} (m-2-4)
      (m-2-2) edge node [swap] {$p_1$} (m-3-3)
      (m-2-4) edge node {$p_2$} (m-3-3);
    \THETAbracket{45}{(m-3-3)};
  \end{tikzpicture}
\end{center}

The composed cospan is then $\ccospan{X}{p_1 f}{G \circ F}{p_2 i}{Z}$. Defining $\textrm{Id}_X := X$, identity $1$-cells are cospans of identity maps: $\ccospan{X}{1}{\textrm{Id}_X}{1}{X}$. It follows from general properties of pushouts that the following are cospan isomorphisms:
\[ (H \circ G) \circ F \cong H \circ (G \circ F) \qquad\qquad
   \textrm{Id}_Y \circ G \cong G \cong G \circ \textrm{Id}_X \]

We can form the ordinary category $\csp(\mathcal C)$ by taking objects to be the $0$-cells from $\Csp(\mathcal C)$ and arrows to be isomorphism-classes of $1$-cells in $\Csp(\mathcal C)$.

Soboci\'{n}ski and Sassone~\cite{SobocinskiThesis,Sassone2005} have extensively studied rewrite systems in the context of bicategories, and in particular bicategories of cospans over adhesive categories. In this section, we will focus particular types of cospan constructions over string graphs that will be used in sections~\ref{sec:sg-cospan-rewriting} and ~\ref{sec:sg-free-monoidal-categories} to construct free monoidal categories and free monoidal categories containing algebraic structures.

Recall that for string graphs, certain $S$-pushouts called \textit{pluggings} perform the function of composition. Using cospans of string graphs $\ccospan{X}{d}{G}{c}{Y}$ where $X$ and $Y$ are point graphs, we can ``pin'' the inputs and outputs for a particular graph in place (i.e. distinguish domain from codomain and induce a total order) and allow us to define composition unambiguously. To ensure that the cospans we consider are meaningful in terms of morphisms in monoidal categories and every cospan composition is a plugging, we shall restrict our attention to \textit{framed point graphs} and \textit{framed cospans} in $\catSGraph_T$.

\begin{definition}
  A \textit{framed point graph} is a triple $(X, <, \sgn)$ where $X$ is a point graph, $<$ is a total order on $V_X$, and $\sgn : V_X \rightarrow \{ +,\ - \}$ is called a \textit{signing map}. A cospan $\ccospan{X}{d}{G}{c}{Y}$ is called a \textit{framed cospan} if:
  \begin{enumerate}
    \item $X$ and $Y$ are framed point graphs,
    \item $G$ contains no isolated wire-vertices,
    \item the induced map $[d,c] : X + Y \rightarrow G$ restricts to an isomorphism $[d,c]' : X + Y \overset{\sim}{\rightarrow} \Bound(G)$,
    \item for every $v \in V_X$, $d(v) \in \In(G) \Leftrightarrow \sgn(v) = +$, and
    \item for every $v \in V_Y$, $c(v) \in \Out(G) \Leftrightarrow \sgn(v) = +$.
  \end{enumerate}
\end{definition}

The function $\sgn$ assigns a polarity to each element of the boundary. A positive polarity marks a wire that runs in the usual (downward) direction, whereas a negative polarity marks a wire that runs in the dual (upward) direction.

\begin{notation}
  For a framed point graph $X$, let $X^*$ be the same framed point graph with the signs reversed.
\end{notation}

A totally downward-directed framed cospan is called \textit{positive}.

\begin{definition}
  A framed point graph $X$ is called \textit{positive} if $\sgn(v) = +$ for all $v \in V_X$. A framed cospan $\ccospan{X}{d}{G}{c}{Y}$ is called positive if both $X$ and $Y$ are positive.
\end{definition}

\begin{proposition}
  For two framed cospans:
  \[ \ccospan{X}{d}{G}{c}{Y} \qquad\qquad \ccospan{Y}{d'}{H}{c'}{Z} \]
  The span $\cspan{G}{c}{Y}{d'}{H}$ is a plugging.
\end{proposition}

\begin{proof}
  For each $v \in V_Y$, if $\sgn(v) = +$ then $c(v) \in \Out(G)$ and $d(v) \in \In(H)$. If $\sgn(v) = -$, then $c(v) \in \In(G)$ and $d(v) \in \Out(H)$. Thus $\cspan{G}{c}{Y}{d'}{H}$ is boundary-coherent. $Y$ is a point graph, so the span is a plugging.
\end{proof}

As a consequence, for framed cospans $G,H$, the composition $H \circ G$ exists and is computed by $S$-pushout.

\begin{example}
  Composing framed cospans by plugging:
  \ctikzfig{cospan_plugging}
\end{example}

\begin{definition}
  For a framed point graph $(X,\leq,\sgn)$, the \textit{pseudo-identity} cospan
  \[\ccospan{X}{d}{\mathbb{1}_X}{c}{X}\]
  \noindent is constructed as follows. $W_X$ has vertices $V_X + V_X$. $d : X \rightarrow \mathbb{1}_X$ maps the vertices of $X$ into the first copy and $c : X \rightarrow \mathbb{1}_X$ maps them into the second copy. $E_{\mathbb{1}_X} = \{ e_v : v \in V_X \}$ and:
  \begin{align*}
    s(e_v) & = \begin{cases}
      d(v) & \textrm{ if } \sgn(v) = + \\
      c(v) & \textrm{ if } \sgn(v) = -
    \end{cases} \\
    t(e_v) & = \begin{cases}
      c(v) & \textrm{ if } \sgn(v) = + \\
      d(v) & \textrm{ if } \sgn(v) = -
    \end{cases}
  \end{align*}
\end{definition}

If we try to form the category of framed cospans over $\catSGraph_T$, we run into a problem. For the identity cospans $\ccospan{X}{1}{\textrm{Id}_X}{1}{X}$, the string graph $\textrm{Id}_X$ contains isolated wire-vertices, so they are not framed cospans. In other words, the ``category'' of framed cospans has no identities! However, there are cospans that come quite close to identity maps.

\begin{example}
  Let $X$ be a framed point graph:
  \ctikzfig{fp_graph}
  
  \noindent The pseudo-identity $\mathbb 1_X$ is defined as:
  \ctikzfig{pseudo_id}
  
  These are called pseudo-identities because composing with them yields a string graph that is wire-homeomorphic to the original graph.
  \ctikzfig{compose_pseudo_id}
\end{example}

In order to obtain honest identities from pseudo-identities, we shall define a category of framed cospans modulo a rewrite system, called a \textit{rewrite category}.

\section{Rewriting on Cospans and Rewrite Categories}\label{sec:sg-cospan-rewriting}

Rewrite categories are categories of framed cospans, modulo a rewrite system. For these to be well-defined, we need to define a notion of rewriting on cospans and show that rewriting cospans is compatible with cospan composition. In particular, we show for $|G|$ an equivalence class of cospans over a given rewrite system, we can define $|H| \circ |G| := |H \circ G|$ in a way that does not depend on the choice of representatives $G$ and $H$.

\begin{definition}
  Let $\ccospan{X}{d}{G}{c}{Y}$ and $\ccospan{X}{d'}{H}{c'}{Y}$ be cospans in $\catSGraph_T$, and let $L \rewritesto R$ be a string graph rewrite rule. For a matching $m$, the following rewrite:
  \begin{center}
    \dposquares
      {\dpoobjects{L}{B}{R}  {G}{G'}{H}}
      {\dpoarrows{b_1}{b_2}  {m}{}{}  {i_1}{i_2}}
  \end{center}
  \noindent is called a \textit{cospan rewrite} if the maps $i_1$ and $i_2$ are cospan homomorphisms. That is, there exist maps $\widehat d$, $\widehat c$ such that the following diagram commutes:
  \begin{center}
    \begin{tikzpicture}
      \matrix (m) [cdiag] {
          & G  &   \\
        X & G' & Y \\
          & H  &   \\
      };
      \path [arrs]
        (m-2-1) edge node {$d$} (m-1-2)
        (m-2-3) edge node [swap] {$c$} (m-1-2)
        (m-2-2) edge node [swap] {$i_1$} (m-1-2)
        (m-2-1) edge node {$\widehat d$} (m-2-2)
        (m-2-3) edge node [swap] {$\widehat c$} (m-2-2)
        (m-2-1) edge node [swap] {$d'$} (m-3-2)
        (m-2-3) edge node {$c'$} (m-3-2)
        (m-2-2) edge node {$i_2$} (m-3-2);
    \end{tikzpicture}
  \end{center}
\end{definition}

We can show that any rewrite on string graphs lifts to a cospan rewrite, for unique maps $\widehat d, \widehat c$. This result relies on the fact that morphisms in a framed cospan factor uniquely through the $S$-pushout complements associated with a string graph rewrite.

\begin{lemma}\label{boundary-factors-complement}
  Let $\cspan{L}{b_1}{I}{b_2}{R}$ be a string graph rewrite rule, $\ccospan{X}{d}{G}{c}{Y}$ a framed cospan, and $m : L \rightarrow G$ an $S$-matching on $G$. Then, for the associated $S$-pushout complement:
  \begin{center}
    \posquare{I}{L}{G'}{G}{b_1}{i}{}{m}
  \end{center}
  \noindent the cospan maps $d$ and $c$ factor uniquely through $i$. That is, there exists unique $\widehat d, \widehat c$ such that:
  \begin{center}
    \begin{tikzpicture}
      \matrix (m) [cdiag] {
        X & G' & Y \\
          & G  &   \\
      };
      \path [arrs]
        (m-1-1) edge [dashed] node {$\widehat d$} (m-1-2)
        (m-1-3) edge [dashed] node [swap] {$\widehat c$} (m-1-2)
        (m-1-1) edge node [swap] {$d$} (m-2-2)
        (m-1-3) edge node {$c$} (m-2-2)
        (m-1-2) edge node {$i$} (m-2-2);
    \end{tikzpicture}
  \end{center}
\end{lemma}

\begin{proof}
  $G'$ is the subgraph of $G$ that has the interior of $L$ removed and $i$ is the inclusion of $G'$ in $G$. For any vertex $v$ in the interior of $L$, $v$ cannot be in $\Bound(L)$ by definition of string graph rewrite rule. Thus, there exists no graph homomorphism $m$ such that $m(v) \in \Bound(G)$. By the definition of framed cospan, the image of $d$ is contained in $\Bound(G)$, so for all $v' \in V_X$, $d(v') \in V_{G'}$. Letting $\widehat d(v') = d(v')$, we have $i \circ \widehat d = d$. The existence of $\widehat c$ follows similarly. Uniqueness is automatic, since $i$ is a monomorphism.
\end{proof}

\begin{corollary}
  A string graph rewrite $G \rewritesto H$ lifts uniquely to a rewrite of framed cospans:
  \[ \ccospan{X}{d}{G}{c}{Y} \ \ \rewritesto\ \ \ccospan{X}{d'}{H}{c'}{Y} \]
\end{corollary}

\begin{notation}
  For a string graph rewrite system $\mathcal R$ and a cospan $\ccospan{X}{d}{G}{c}{Y}$, we write $|G|_{\mathcal R}$ for the set of cospans $G'$ such that:
  \[ \ccospan{X}{d}{G}{c}{Y}\ \ \rewriteequiv_{\mathcal R}\ \ \ccospan{X}{d'}{G'}{c'}{Y} \]
\end{notation}

We drop the subscript $\mathcal R$ when it is clear from context.

\begin{theorem}\label{thm:rw-compose-commute}
  Rewriting commutes with composition. Let the following diagram define a composition of framed cospans:
  \begin{equation}\label{eq:composition-for-rw}
    \begin{tikzpicture}
      \matrix (m) [cdiag] {
          &   & Y         &   &   \\
        X & G &           & H & Z \\
          &   & H \circ G &   &   \\
      };
      \path [arrs]
        (m-2-1) edge node {$d_1$} (m-2-2)
        (m-1-3) edge node [swap] {$c_1$} (m-2-2)
        (m-1-3) edge node {$d_2$} (m-2-4)
        (m-2-5) edge node [swap] {$c_2$} (m-2-4)
        (m-2-2) edge node [swap] {$i_1$} (m-3-3)
        (m-2-4) edge node {$i_2$} (m-3-3);
      \THETAbracket{45}{(m-3-3)};
    \end{tikzpicture}
  \end{equation}
  
  For a string graph rewrite rule $L \rewritesto R$ and an $S$-matching $m : L \rightarrow G$, the composed morphism $i_1 m : L \rightarrow H \circ G$ is an $S$-matching and:
  \begin{equation}\label{eq:rw-commute-iso}
    H \circ (G[L \rewritesto R]_m) \cong (H \circ G)[L \rewritesto R]_{i_1 m}
  \end{equation}
  
  Similarly, for $n : L \rightarrow H$, $i_2 n : L \rightarrow H \circ G$ is an $S$-matching and:
  \begin{equation}\label{eq:rw-commute-iso2}
    (H[L \rewritesto R]_n) \circ G \cong (H \circ G)[L \rewritesto R]_{i_2 n}
  \end{equation}
\end{theorem}

\begin{proof}
  $i_1$ is a monomorphism because the pushout in (\ref{eq:composition-for-rw}) is a plugging, so $i_1 m$ is a monomorphism, hence an $S$-matching by Theorem \ref{thm:all-monos-matchings}. By Lemma \ref{boundary-factors-complement}, the map $c_1$ factors through the pushout complement $G -_{b_1,m} L$ as in diagram (\ref{dia:compatible}):
  \begin{center}
    \begin{tikzpicture}
    \matrix (m) [cdiag] {
      G                    &   &   \\
      G -_{b_1,m} L        & Y & H \\
      G[L \rewritesto R]_m &   &   \\
    };
    \path [arrs]
      (m-2-2) edge [bend right] node [swap] {$c_1$} (m-1-1)
      (m-2-2) edge node [swap] {$\widehat c_1$} (m-2-1)
      (m-2-2) edge [bend left] node {$c_1'$} (m-3-1)
      (m-2-1) edge (m-1-1)
      (m-2-1) edge (m-3-1)
      (m-2-2) edge node {$d_2$} (m-2-3);
    \end{tikzpicture}
  \end{center}
  All three pairs $(c_1,d)$, $(\widehat c_1,d)$, and $(c_1',d)$ are pluggings, so they are $S$-spans. Isomorphism (\ref{eq:rw-commute-iso}) follows from Theorem \ref{thm:pushout-and-rewrite}. Isomorphism (\ref{eq:rw-commute-iso2}) follows similarly.
\end{proof}

\begin{corollary}
  For a string graph rewrite system $\mathcal R$, composition given by $|H|_{\mathcal R} \circ |G|_{\mathcal R} := |H \circ G|_{\mathcal R}$ is well defined, and does not depend on the choices of $G$ and $H$.
\end{corollary}

\begin{proof}
  If $G' \in |G|_{\mathcal R}$ and $H' \in |G|_{\mathcal R}$, then by Theorem \ref{thm:rw-compose-commute}, we can always find matchings to rewrite $H' \circ G'$ into $H \circ G$, so $H' \circ G' \in |H \circ G|_{\mathcal R}$.
\end{proof}

\begin{theorem}\label{thm:rw-cat}
  For a rewrite system $\mathcal R$ containing the wire-homeomorphism rules $\mathbb H$, $\FCsp(\mathcal R, T)$ is a category where:
  \begin{itemize}
    \item objects are framed point graphs,
    \item arrows are equivalence classes $|G|_{\mathcal R}$ of framed cospans,
    \item identities are defined by pseudo-identity cospans: $|\mathbb 1_X|_{\mathcal R}$, and
    \item composition is given by: $|H|_{\mathcal R} \circ |G|_{\mathcal R} := |H \circ G|_{\mathcal R}$.
  \end{itemize}
\end{theorem}


\begin{proof}
  Composition is associative, because it is defined by pushouts as in cospan categories, so it remains to show that the pseudo-identity string graphs yield genuine identities in $\FCsp(\mathcal R, T)$. For a framed cospan $G$, consider the composition:
  \begin{center}
    \begin{tikzpicture}
      \matrix (m) [cdiag] {
          &   & Y                   &             &   \\
        X & G &                     & \mathbb 1_Y & Y \\
          &   & \mathbb 1_Y \circ G &             &   \\
      };
      \path [arrs]
        (m-2-1) edge node {$d_1$} (m-2-2)
        (m-1-3) edge node [swap] {$c_1$} (m-2-2)
        (m-1-3) edge node {$d_2$} (m-2-4)
        (m-2-5) edge node [swap] {$c_2$} (m-2-4)
        (m-2-2) edge node [swap] {$i_1$} (m-3-3)
        (m-2-4) edge node {$i_2$} (m-3-3);
      \THETAbracket{45}{(m-3-3)};
    \end{tikzpicture}
  \end{center}
  
  $\mathbb 1_Y \circ G$ contains a copy of $G$, as well as some extra edges and wire-vertices. For every wire-vertex $v \in W(\mathbb 1_Y \circ G) - W(G)$, there exists a unique edge $e \in E_{\mathbb 1_Y \circ G} - E_G$ adjacent to that vertex. We can always find a rewrite rule in $\mathbb H \subseteq \mathcal R$ that produces a string graph isomorphic to $\mathbb 1_Y \circ G$ with $v$ and $e$ removed. Since all edges $e \in E_{\mathbb 1_Y \circ G} - E_G$ arise in this way, we can repeat the process until the resultant string graph is isomorphic to $G$. So $|\mathbb 1_Y \circ G|_{\mathcal R} = |G|_{\mathcal R}$. Similarly, $|G \circ \mathbb 1_X|_{\mathcal R} = |G|_{\mathcal R}$.
\end{proof}

We often want to consider the full subcategory of $\FCsp(\mathcal R, T)$ containing only equivalence classes of ``downward-directed'' cospans. This category is called $\FCsp^+(\mathcal R, T)$.

\begin{definition}
  $\FCsp^+(\mathcal R, T)$ be the full subcategory of $\FCsp(\mathcal R, T)$ whose objects are positive framed point graphs.
\end{definition}

\begin{theorem}
  $\FCsp^+(\mathcal R, T)$ is a symmetric traced category and $\FCsp(\mathcal R, T)$ is a compact closed category.
\end{theorem}

\begin{proof}
  First, we show that $\FCsp(\mathcal R, T)$ is a compact closed category. For any framed graphs $A, B$, $A \otimes B$ has vertices $V_A + V_B$ which can be given a total order by placing all of the elements in $V_B$ above those in $V_A$ (i.e. the usual disjoint union of totally ordered sets). The monoid product of framed cospans is given coproducts in $\catSGraph_T$. For cospans:
  \[ \ccospan{A}{d_1}{G}{c_1}{C} \qquad\qquad \ccospan{B}{d_2}{H}{c_2}{C} \]
  \noindent the underlying string graphs of $A \otimes B$ and $C \otimes D$ are coproducts $A + B$ and $C + D$, so there is an induced cospan over $G \otimes H := G + H$:
  \begin{center}
    \begin{tikzpicture}
      \matrix (m) [cdiag] {
        A & A \otimes B & B \\
        G & G \otimes H & H \\
        C & C \otimes D & D \\
      };
      \path [arrs]
        (m-1-1) edge [right hook-latex] (m-1-2)
        (m-1-3) edge [left hook-latex] (m-1-2)
        (m-1-1) edge node [swap] {$d_1$} (m-2-1)
        (m-1-2) edge [dashed] node {$[d_1,d_2]$} (m-2-2)
        (m-1-3) edge node {$d_2$} (m-2-3)
        (m-2-1) edge [right hook-latex] (m-2-2)
        (m-2-3) edge [left hook-latex] (m-2-2)
        (m-3-1) edge node {$c_1$} (m-2-1)
        (m-3-2) edge [dashed] node [swap] {$[c_1,c_2]$} (m-2-2)
        (m-3-3) edge node [swap] {$c_2$} (m-2-3)
        (m-3-1) edge [right hook-latex] (m-3-2)
        (m-3-3) edge [left hook-latex] (m-3-2);
    \end{tikzpicture}
  \end{center}
  
  It can easily be verified that this induced cospan is framed. The rest of the structure maps are analogous to those from string diagrams:
  \ctikzfig{structure_string_graphs}
  
  All of the axioms of a compact closed category then follow from string graph isomorphism and edge homeomorphism. For example:
  \ctikzfig{sg_compact_pf}
  
  Since any full subcategory of a compact closed category is a symmetric traced category, $\FCsp^+(\mathcal R, T)$ is a symmetric traced category, with the trace operation defined using the compact structure from $\FCsp(\mathcal R, T)$.
\end{proof}

\section{Free Monoidal Categories}\label{sec:sg-free-monoidal-categories}


In this section, we will prove that the category $\FCsp^+(\mathbb H, T)$ is the free symmetric trace category on a monoidal signature $T$. Free categories are special in that they are sound and complete with respect to the axioms of that type of category. That is, we shall prove that two morphisms are equal by the axioms of a symmetric traced category \textit{if and only if} they are equal in $\FCsp^+(\mathbb H, T)$. Once we prove this, it is a simple matter to show that $\FCsp(\mathbb H, T)$ is the free \textit{compact closed} category on a monoidal signature, by proving that the category $\FCsp(\mathbb H, T)$ is equivalent to the result of performing the ``Int'' construction~\cite{JSV} on the free symmetric traced category $\FCsp^+(\mathbb H, T)$. For simplicity, we will focus on strict categories in this section.

Before we get to the bulk of the proof, we introduce some notation. The first thing we introduce is the notion of \textit{indexing} a morphism.

\begin{definition}\label{def:x-word}
  For a small, strict monoidal category $\mathcal V$, fix a set of \textit{atomic objects} $O$, such that any object in $\mathcal V$ is a monoidal product of elements of $O$. For an object $X \in \textrm{ob}\mathcal V$, an \textit{$X$-word} is a monoidal product $X_{i_1} \otimes \ldots \otimes X_{i_M} = X$ such that all $i_k$ are distinct and $X_{i_k} \in O$.
\end{definition}

We will assume that $O$ contains ``enough copies'' of every atomic object to find $X$-words for every object $X$. Replacing $X$ with an $X$-word is simply the act of binding each position in the monoidal product to a unique index that we can refer to later.

\begin{definition}\label{def:indexing}
  For a morphism $f : X \rightarrow Y$, an \textit{indexing} of $f$ is a choice of an $X$-word and a $Y$-word such that
  \[ f = f' : X_{i_1} \otimes \ldots \otimes X_{i_M} \rightarrow  Y_{j_1} \otimes \ldots \otimes Y_{j_N} \]
  
  A morphism from an $X$-word to a $Y$-word for any $X, Y$ is called \textit{indexed}.
\end{definition}

For an $X$-word $X_{i_1} \otimes \ldots \otimes X_{i_M}$, and an index $i \in \{ i_1,\ldots,i_M \}$, $\sigma_{X:i}$ is the (unique) symmetry map that permutes the object $X_i$ to the end of the list and leaves the other objects fixed.
\begin{align*}
  \sigma_{X:i} & = %
\beginpgfgraphicnamed{sigma_i}
\InputIfFileExists{sigma_i.tikz}{}{\input{./figures/sigma_i.tikz}}
\endpgfgraphicnamed
\end{align*}

In any strict symmetric traced category, we can define a contraction operator $C_i^j(-)$ which ``traces together'' the $i$-th input with the $j$-th output on an indexed morphism.

\begin{definition}
  Let $f : X_{i_1} \otimes \ldots \otimes X_{i_M} \rightarrow  Y_{j_1} \otimes \ldots \otimes Y_{j_N}$ be an indexed morphism in a symmetric traced category such that for indices $i \in \{ i_1,\ldots,i_M \}$ and $j \in \{ j_1,\ldots,j_N \}$, $X_i = Y_j$. Then we define the \textit{contraction} $C_i^j(f)$ as follows:
  \begin{align*}
    C_i^j(f) & := \Tr^{X_i = Y_j}(\sigma_{Y:j} \circ f \circ \sigma_{X:i}^{-1}) \\
             &  = \ \ %
\beginpgfgraphicnamed{contraction_def}
\InputIfFileExists{contraction_def.tikz}{}{\input{./figures/contraction_def.tikz}}
\endpgfgraphicnamed
  \end{align*}
\end{definition}

Note that a contraction of an indexed morphism yields an indexed morphism, so we can contract many times. Also, the resulting morphism does not depend on the order in which we perform contractions.

\begin{lemma}\label{lem:contractions-commutative}
  Contractions are commutative. For an indexed morphism $f$ distinct indices $i,i'$ and $j,j'$:
  \[ C_{i}^{j}(C_{i'}^{j'}(f)) = C_{i'}^{j'}(C_{i}^{j}(f)) \]
\end{lemma}

\begin{proof}
  Let $X=X_{i_{1}}\otimes\ldots\otimes X_{i_{M}}$ and $Y=Y_{j_{1}}\otimes\ldots\otimes Y_{j_{N}}$, let $i,i'\in\{i_{1},\ldots,i_{M}\}$ be distinct, and let $j,j'\in\{j_{1},\ldots,j_{N}\}$ be distinct.
  \begin{align*}
  C_{i}^{j}(C_{i'}^{j'}(f)) & =\Tr^{X_{i}=Y_{j}}(\sigma_{Y:j}\circ\Tr^{X_{i'}=Y_{j'}}(\sigma_{Y:j'}\circ f\circ\sigma_{X:i'}^{-1})\circ\sigma_{X:i}^{-1})\\
   & =\Tr^{X_{i}=Y_{j}}(\Tr^{X_{i'}=Y_{j'}}((\sigma_{Y:j}\otimes1_{Y_{j'}})\circ\sigma_{Y:j'}\circ f\circ\sigma_{X:i'}^{-1}\circ(\sigma_{X:i}^{-1}\otimes1_{X_{i'}})))\\
   & =\Tr^{X_{i}\otimes X_{i'}=Y_{j}\otimes Y_{j'}}((\sigma_{Y:j}\otimes1_{Y_{j'}})\circ\sigma_{Y:j'}\circ f\circ\sigma_{X:i'}^{-1}\circ(\sigma_{X:i}^{-1}\otimes1_{X_{i'}}))\\
   & =(*)
  \end{align*}

  Let $X'$ be a new $X$-word equal to $X$ with the factors $X_{i}$ and $X_{i'}$ deleted. Let $Y'$ be a $Y$-word equal to $Y$ with the factors at $Y_{j}$ and $Y_{j'}$ omitted.
  \begin{align*}
    (*)
     & =\Tr^{X_{i}\otimes X_{i'}=Y_{j}\otimes Y_{j'}}((\sigma_{Y:j}\otimes1_{Y_{j'}})\circ\sigma_{Y:j'}\circ f\circ\sigma_{X:i'}^{-1}\circ(\sigma_{X:i}^{-1}\otimes1_{X_{i'}})\circ(1_{X'}\otimes1_{X_{i}}\otimes1_{X_{i'}})\\
     & =\Tr^{X_{i}\otimes X_{i'}=Y_{j}\otimes Y_{j'}}((\sigma_{Y:j}\otimes1_{Y_{j'}})\circ\sigma_{Y:j'}\circ f\circ\sigma_{X:i'}^{-1}\circ(\sigma_{X:i}^{-1}\otimes1_{X_{i'}})\circ(1_{X'}\otimes(\sigma_{X_{i},X_{i'}}^{-1}\circ\sigma_{X_{i},X_{i'}})))\\
     & =\Tr^{X_{i'}\otimes X_{i}=Y_{j'}\otimes Y_{j}}((1_{Y'}\otimes\sigma_{X_{i},X_{i'}})\circ(\sigma_{Y:j}\otimes1_{Y_{j'}})\circ\sigma_{Y:j'}\circ f\circ\sigma_{X:i'}^{-1}\circ(\sigma_{X:i}^{-1}\otimes1_{X_{i'}})\circ(1_{X'}\otimes\sigma_{X_{i},X_{i'}}^{-1}))\\
     & =\Tr^{X_{i'}\otimes X_{i}=Y_{j'}\otimes Y_{j}}((1_{Y'}\otimes\sigma_{Y_{j},Y_{j'}})\circ(\sigma_{Y:j}\otimes1_{Y_{j'}})\circ\sigma_{Y:j'}\circ f\circ\sigma_{X:i'}^{-1}\circ(\sigma_{X:i}^{-1}\otimes1_{X_{i'}})\circ(1_{X'}\otimes\sigma_{X_{i},X_{i'}}^{-1}))\\
     & =(*)
  \end{align*}

  Now, we need to have a look at the symmetry maps. By the coherence theorem for symmetric monoidal categories, a composition of symmetry maps is uniquely defined by the permutation it performs. First note that $(1_{Y'}\otimes\sigma_{Y_{j},Y_{j'}})\circ(\sigma_{Y:j}\otimes1_{Y_{j'}})$ sends $Y_{j}$ to the end of the list, so:
  \[ (1_{Y'}\otimes\sigma_{Y_{j},Y_{j'}})\circ(\sigma_{Y:j}\otimes1_{Y_{j'}})=\sigma_{Y:j} \]

  Note that the $\sigma_{Y:j}$ on the LHS refers to a different map than the one on the RHS, as it has a different domain and codomain. By naturality:
  \[ \sigma_{Y:j}\circ\sigma_{Y:j'}=(\sigma_{Y:j'}\otimes1_{Y_{j}})\circ\sigma_{Y:j}\]

  It can be shown similarly that:
  \[ \sigma_{X:i'}^{-1}\circ(\sigma_{X:i}^{-1}\otimes1_{X_{i'}})\circ(1_{X'}\otimes\sigma_{X_{i},X_{i'}}^{-1})=\sigma_{X:i'}^{-1}\circ\sigma_{X:i}^{-1}=\sigma_{X:i}^{-1}\circ(\sigma_{X:i'}^{-1}\otimes1_{X_{i}}) \]

  Substituting in to the above expression, we complete the proof:
  \begin{align*}
    (*)
     & =\Tr^{X_{i'}\otimes X_{i}=Y_{j'}\otimes Y_{j}}((\sigma_{Y:j'}\otimes1_{Y_{j}})\circ\sigma_{Y:j}\circ f\circ\sigma_{X:i}^{-1}\circ(\sigma_{X:i'}^{-1}\otimes1_{X_{i}}))\\
     & =\Tr^{X_{i'}=Y_{j'}}(\Tr^{X_{i}=Y_{j}}((\sigma_{Y:j'}\otimes1_{Y_{j}})\circ\sigma_{Y:j}\circ f\circ\sigma_{X:i}^{-1}\circ(\sigma_{X:i'}^{-1}\otimes1_{X_{i}})))\\
     & =\Tr^{X_{i'}=Y_{j'}}(\sigma_{Y:j'}\circ\Tr^{X_{i}=Y_{j}}(\sigma_{Y:j}\circ f\circ\sigma_{X:i}^{-1}\circ)\circ\sigma_{X:i'}^{-1})\\
     & =C_{i'}^{j'}(C_{i}^{j}(f))
  \end{align*}
\end{proof}

\begin{definition}\label{def:disconnected-indexed}
  For a (small, strict) traced symmetric category $\mathcal V$, define a set $M$ of \textit{atomic} morphisms, such that any morphism in $\mathcal V$ can be obtained from those morphisms and the traced symmetric structure. An indexed morphism is called \textit{disconnected} if it is of the form $f = f_1 \otimes \ldots \otimes f_K$, where each $f_k$ is an indexed morphism in $M$:
  \[ f_k : X_{i_{k,1}} \otimes \ldots X_{i_{k,M_k}} \rightarrow Y_{j_{k,1}} \otimes \ldots \otimes Y_{j_{k,N_k}} \]
\end{definition}

\begin{definition}\label{def:cnf}
  Let $f = f_1 \otimes \ldots \otimes f_M$ be a disconnected indexed map. For distinct indices $\{i_1,\ldots,i_P\} \subseteq \{i_{1,1},\ldots i_{K,M_K} \}$ and $\{j_1,\ldots,j_P\} \subseteq \{j_{1,1},\ldots j_{K,N_K} \}$, a map $f'$ is said to be in \textit{contraction normal form} (CNF) if:
  \[ f' = C_{i_1}^{j_1}(C_{i_2}^{j_2}(\ldots(C_{i_P}^{j_P}(f))\ldots)) \]
\end{definition}

\begin{definition}\label{def:totally-contracted}
  Let $f$ and $f'$ be given as in Definition \ref{def:cnf}. A component $f_k$ of $f$ is said to be \textit{totally contracted} if the indices of all of its inputs occur in $\{i_1,\ldots,i_P\}$ and the indices of all of its outputs occur in $\{j_1,\ldots,j_P\}$.
\end{definition}

\begin{lemma}\label{lem:reorder-tc}
  Let $f$ and $f'$ be given as in Definition \ref{def:cnf}. By re-indexing the contraction, we can reorder the totally contracted components of $f$ arbitrarily.
\end{lemma}

\begin{proof}
  It suffices to show that we can send any totally-contracted $f_k$ to the far right side of $f$ by re-indexing the contraction.
  \[ C_{i_1}^{j_1}(C_{i_2}^{j_2}(\ldots(C_{i_P}^{j_P}(
        f_1 \otimes \ldots \otimes f_{k} \otimes \ldots \otimes f_K))\ldots)) =
     C_{i'_1}^{j'_1}(C_{i'_2}^{j'_2}(\ldots(C_{i'_P}^{j'_P}(
        f_1 \otimes \ldots \otimes f_K \otimes f_{k}))\ldots))
  \]
  
  We can show this by applying the previous lemma and naturality of symmetries. First, write $f$ as $f = f^L \otimes f_k \otimes f^R$ for $f^L : X^L \rightarrow Y^L$, $f^R : X^R \rightarrow Y^R$, and $f_k : X^k \rightarrow Y^k$.
  \[ C_{i_1}^{j_1}(\ldots(C_{i_P}^{j_P}(f^L \otimes f_{k} \otimes f^R)\ldots)) = (*) \]
  
  We then pre-compose with the indentity (i.e. a swap map and its inverse). Since $f_k$ and all the maps after it are totally-contracted, we can eliminate one of the swap maps by re-indexing the $i$'s. Similarly, we can introduce a swap after $f$ by re-indexing the $j$'s.
  \begin{align*}
    (*)
    &= C_{i_1}^{j_1}(\ldots(C_{i_P}^{j_P}(
        (f^L \otimes f_{k} \otimes f^R) \circ (X^L \otimes \sigma_{X^R,X^k}) \circ (X^L \otimes \sigma_{X^k,X^R})
        )\ldots)) \\
    &= C_{i'_1}^{j'_1}(\ldots(C_{i'_P}^{j'_P}(
        (f^L \otimes f_{k} \otimes f^R) \circ (X^L \otimes \sigma_{X^R,X^k})
        )\ldots)) \\
    &= C_{i'_1}^{j'_1}(\ldots(C_{i'_P}^{j'_P}(
        (Y^L \otimes \sigma_{Y^k,Y^R}) \circ (f^L \otimes f_{k} \otimes f^R) \circ (X^L \otimes \sigma_{X^R,X^k})
        )\ldots)) = (*)
  \end{align*}
  
  An application of naturality completes the proof:
  \begin{align*}
    (*)
    &= C_{i'_1}^{j'_1}(\ldots(C_{i'_P}^{j'_P}(
        (f^L \otimes f^R \otimes f_{k}) \circ (X^L \otimes \sigma_{X^k,X^R}) \circ (X^L \otimes \sigma_{X^R,X^k})
        )\ldots)) \\
    &= C_{i'_1}^{j'_1}(\ldots(C_{i'_P}^{j'_P}(f^L \otimes f^R \otimes f_{k})\ldots))
  \end{align*}
  
  We can therefore send $f_k$ to an arbitrary position by performing this procedure in reverse.
\end{proof}

A totally contracted identity map whose input is connected to its output is called a \textit{minimal circle}.

\begin{lemma}\label{cor:totally-contracted-ids}
  Let $f, f'$ be Defined as in \ref{def:cnf}. If $f_k = 1_{X_{i_{k,1}}} = 1_{Y_{j_{k,1}}}$ is a totally contracted identity map that is not a minimal circle, then it can be removed by re-indexing.
\end{lemma}

\begin{proof}
  From Lemma \ref{lem:reorder-tc}, we can assume the totally contracted identity map is on the far right. Let $M$ be the index of the input and $N$ be the index of the output. By \ref{lem:contractions-commutative}, we can move the two contractions involving the identity map all the way to the inside, so $f'$ is of the form:
  \[ f' = C(C(\ldots(C_M^j(C_i^N(f'' \otimes 1_{X_M})))\ldots)) \]
  
  We can reduce the inner map using the definition of $\sigma_{Y:j}$ and the trace axioms.
  \begin{align*}
    C_{M}^{j}(C_{i}^{N}(f''\otimes1_{X_{M}}))
     & =\Tr^{X_{M}=Y_{j}}(\sigma_{Y:j}\circ\Tr^{X_{i}=Y_{N}}((f''\otimes1_{X_{M}})\circ\sigma_{X:i}^{-1}))\\
     & =\Tr^{X_{M}=Y_{j}}(\sigma_{Y:j}\circ\Tr^{X_{i}=Y_{N}}((f''\otimes1_{X_{M}})\circ(\sigma_{X:i}^{-1}\otimes1_{X_{M}})\circ(1_{X'}\otimes\sigma_{X_{M},X_{i}}^{-1})))\\
     & =\Tr^{X_{M}=Y_{j}}(\sigma_{Y:j}\circ\Tr^{X_{i}=Y_{N}}(((f''\circ\sigma_{X:i}^{-1})\otimes1_{X_{M}})\circ(1_{X'}\otimes\sigma_{X_{M},X_{i}}^{-1})))\\
     & =\Tr^{X_{M}=Y_{j}}(\sigma_{Y:j}\circ f''\circ\sigma_{X:i}^{-1}\circ\Tr^{X_{i}=Y_{N}}((1_{X'}\otimes\sigma_{X_{M},X_{i}}^{-1})))\\
     & =\Tr^{X_{M}=Y_{j}}(\sigma_{Y:j}\circ f''\circ\sigma_{X:i}^{-1}\circ(1_{X'}\otimes\Tr^{X_{i}=Y_{N}}(\sigma_{X_{M},X_{i}}^{-1})))\\
     & =\Tr^{X_{M}=Y_{j}}(\sigma_{Y:j}\circ f''\circ\sigma_{X:i}^{-1}\circ(1_{X'}\otimes1_{X_{M}=X_{i}}))\\
     & =\Tr^{X_{i}=Y_{j}}(\sigma_{Y:j}\circ f''\circ\sigma_{X:i}^{-1})=C_{i}^{j}(f'')
  \end{align*}
\end{proof}

\begin{theorem}\label{thm:sg-free-traced}
  $\FCsp^+(\mathbb H, T)$ is the free strict symmetric traced category on the monoidal signature $T$.
\end{theorem}

\begin{proof}
  Let $T = (O, M, \dom, \cod)$ be a monoidal signature, $\mathcal V$ be a strict symmetric traced category, and $F : T \rightarrow \mathcal V$ be a monoidal signature homomorphism. There is an evident signature homomorphism $E : T \rightarrow \FCsp^+(\mathbb H, T)$, taking each object $Z \in O$ to the framed point graph consisting of a single wire-vertex of type $Z$ and taking each $g \in M$ to a homeomorphism-class of cospans consisting only of a box-vertex of type $g$ and its inputs and outputs.
  
  Let $X$ be a framed, positive point graph, where $V_X = \{ p_1 < p_2 < \ldots < p_N \}$. Then, for all $p_i$, $F(\tau_X(p_i))$ is an object in $\mathcal V$. Let:
  \[ \widehat F(X) = F(\tau_X(p_1)) \otimes \ldots \otimes F(\tau_X(p_N)) \]
  
  This map is uniquely specified, and it respects the monoidal product on objects in $\FCsp^+(\mathbb H, T)$. Let $|G| : X \rightarrow Y$ be an arrow in $\FCsp^+(\mathbb H, T)$, represented by a cospan $\ccospan{X}{d}{G}{c}{Y}$. Let $\{ x_1 < x_2 < \ldots < x_M \}$ be set of wire-vertices in $G$ that are in the image of $d$, inheriting the total ordering from $X$. Let $\{ y_1 < y_2 < \ldots < y_N \}$ be the same for $Y$. Define a function $F'$ from $V_G$ to the morphisms of $\mathcal V$ as follows:
  \[ 
  F'(v) =
  \begin{cases}
    F(\tau_G(v))     & \textrm{ if $v$ is a box-vertex}  \\
    1_{F(\tau_G(v))} & \textrm{ if $v$ is a wire-vertex}
  \end{cases}
  \]
  
  Let $\{ z_1, \ldots, z_K \} \subseteq V_G$ be the set of vertices not in the image of $d$ or $c$. Define a disconnected indexed morphism in $\mathcal V$.
  \[ f = F'(x_1) \otimes \ldots \otimes F'(x_M) \otimes 
         F'(y_1) \otimes \ldots \otimes F'(y_N) \otimes
         F'(z_1) \otimes \ldots \otimes F'(z_K) \]
         
  We can form a CNF term from $|G|$ by adding a contraction $C_i^j(-)$ for every edge $e$ in $G$. We define the indices for these contractions as follows. If $\tau_G(e) = \textrm{mid}_Z$, then $F'(s(e))$ and $F'(t(e))$ are both identity maps. Let $i$ be the unique input of $F'(t(e))$ and $j$ is the unique output of $F'(s(e))$. If $\tau_G(e) = \textrm{in}_{g,k}$, then $F'(s(e))$ is an identity map. Let $j$ be the unique output of $F'(s(e))$ and let $i$ be the $k$-th input of $F'(t(e))$. If $\tau_G(e) = \textrm{out}_{g,k}$, then $F'(t(e))$ is an identity map. Let $i$ be the unique input of $F'(t(e))$ and let $j$ be the $k$-th output of $F'(s(e))$. Define $\widehat F(|G|)$ as:
  \[ \widehat F(|G|) = C_{i_1}^{j_1}(C_{i_2}^{j_2}(\ldots(C_{i_P}^{j_P}(
        f_1 \otimes \ldots \otimes f_{k} \otimes \ldots \otimes f_K))\ldots)) \]
  For this to be well-defined, we need to show that this does not depend on the choice of $G$. Choosing an isomorphic string graph cospan amounts to picking a different order for the internal vertices $z_i$ and the contractions $C_i^j(-)$. Any internal vertices in $G$ will be totally contracted in $\widehat F(|G|)$, so by Lemmas \ref{lem:contractions-commutative} and \ref{lem:reorder-tc}, choosing an isomorphic cospan will not affect the value of $\widehat F(|G|)$. Choosing a homeomorphic cospan amounts to varying the number of totally contracted identities that are not minimal circles. By Lemma \ref{cor:totally-contracted-ids}, this does not affect the value of $\widehat F(|G|)$ either.
  
  It is straightforward to show that $\widehat F$ takes identities to identities and respects composition and traces, so $\widehat F$ is a symmetric traced functor from $\FCsp^+(\mathbb H, T)$ to $\mathcal V$. Thus $\widehat F$ satisfies the universal property of the free symmetric traced category.
  \begin{center}
    \begin{tikzpicture}
      \matrix (m) [cdiag] {
        T & \FCsp^+(\mathbb H, T) \\
          & \mathcal V            \\
      };
      \path [arrs]
        (m-1-1) edge node {$E$} (m-1-2)
        (m-1-2) edge [dashed] node {$\widehat F$} (m-2-2)
        (m-1-1) edge node [swap] {$F$} (m-2-2);
    \end{tikzpicture}
  \end{center}
  
  Since it is possible to build any positive framed cospan using the cospans in the image of $E$ and the traced symmetric structure, $\widehat F$ is the unique map making this diagram commute.
\end{proof}

\begin{theorem}
  There is a symmetric monoidal equivalence of categories $\Int(\FCsp^+(\mathbb H, T)) \cong \FCsp(\mathbb H, T)$, such that the following diagram commutes:
  \begin{equation}\label{eq:int-embed-commute}
    \begin{tikzpicture}
      \matrix (m) [cdiag] {
        \FCsp^+(\mathbb H, T) & \Int(\FCsp^+(\mathbb H, T)) \\
                              & \FCsp(\mathbb H, T)         \\
      };
      \path [arrs]
        (m-1-1) edge [right hook-latex] (m-1-2)
        (m-1-2) edge node {$\cong$} (m-2-2)
        (m-1-1) edge [right hook-latex] (m-2-2);
    \end{tikzpicture}
  \end{equation}
\end{theorem}

\begin{proof}
  The category $\Int(\FCsp^+(\mathbb H, T))$ has as objects pairs of positive framed point graphs $(A, X)$ and as arrows $\mathbb H$-equivalence classes of cospans. An arrow $|G| : (A, X) \rightarrow (B, Y)$ is represented by a cospan $\ccospan{A + Y}{d}{G}{c}{B + X}$. Composition is done by composing on the first object and ``composing backwards'' via the trace on the second object.
  \ctikzfig{int_compose}
  
  Since the domain and codomain are disjoint unions of framed point graphs, we can decompose the map $d$ into maps $d_1 : A \rightarrow G$ and $d_2 : Y \rightarrow G$ and the map $c$ into maps $c_1 : B \rightarrow G$ and $c_2 : X \rightarrow G$.
  \[ \ccospan{A + Y}{[d_1,d_2]}{G}{[c_1,c_2]}{B + X} \]
  
  We can then form a new cospan by interchanging $X$ and $Y$ and flipping their sign maps. Since $d_2$ and $c_2$ are just string graph homomorphisms, we can consider them to have domains $X^*$ and $Y^*$, respectively. Thus we obtain a cospan representing an arrow in $\FCsp(\mathbb H, T)$.
  \begin{equation}\label{eq:flipped-cospan}
    \ccospan{A + X^*}{[d_1,c_2]}{G}{[c_1,d_2]}{B + Y^*}
  \end{equation}
  
  Let $F((A,X)) = A + X^*$ and let $F(|G|)$ be the $\mathbb H$-equivalence class of cospans represented by (\ref{eq:flipped-cospan}). We can show that $F$ respects composition:
  \ctikzfig{int_equiv_pf1}
  
  It can also be verified that $F$ preserves all of the traced symmetric structure, up to isomorphism. Flipping $X$ and $Y$ induces a bijection of hom-sets:
  \[ \hom_{\FCsp^+(\mathbb H, T)}(A + Y, B + X) \cong
     \hom_{\FCsp(\mathbb H, T)}(A + X^*, B + Y^*) \]
  \noindent so $F$ is full and faithful. For all positive framed point graphs $A, X$, objects $A + X^*$ are in the image of $F$. For an arbitrary object $Z$ in $\FCsp(\mathbb H, T)$, we can define a symmetry isomorphism that sends all the positive points to the left and all the negative points to the right.
  \ctikzfig{int_equiv_pf2}
  
  So $F$ is essentially surjective. Therefore there is a traced symmetric equivalence of categories $\Int(\FCsp^+(\mathbb H, T)) \cong \FCsp(\mathbb H, T)$. $\FCsp^+(\mathbb H, T)$ fully embeds into $\Int(\FCsp^+(\mathbb H, T))$ as objects $(A,\emptyset)$ and cospans $\ccospan{A + \emptyset}{d}{G}{c}{B + \emptyset}$. The functor $F$ acts trivially on these cospans, so diagram (\ref{eq:int-embed-commute}) commutes.
\end{proof}

\begin{corollary}
  $\FCsp(\mathbb H, T)$ is the free compact closed category on a monoidal signature $T$.
\end{corollary}

\begin{proof}
  The embedding of $T$ in $\FCsp(\mathbb H, T)$ factors through $\FCsp^+(\mathbb H, T)$. The universal property follows from composing free constructions.
  \begin{center}
    \begin{tikzpicture}
      \matrix (m) [cdiag] {
        T & \FCsp^+(\mathbb H, T) & \FCsp(\mathbb H, T) \\
          &                       & \mathcal V          \\
      };
      \path [arrs]
        (m-1-1) edge node {$E$} (m-1-2)
        (m-1-1) edge [bend right=15] node [swap] {$F$} (m-2-3)
        (m-1-2) edge [dashed] node {$\widehat F$} (m-2-3)
        (m-1-2) edge [right hook-latex] (m-1-3)
        (m-1-3) edge [dashed] node {$\widetilde{F}$} (m-2-3);
    \end{tikzpicture}
  \end{center}
  $F$ is a monoidal signature homomorphism, $\widehat F$ is the unique traced symmetric functor induced by $F$, and $\widetilde F$ is the unique compact closed functor induced by $\widehat F$.
\end{proof}

A consequence of this construction, is we can quite easily build ``free categories containing an X'' (e.g. the traced symmetric and compact closed analogues of Lawvere theories, PROPs, etc.) as rewrite categories. For an algebraic theory $\mathcal E$, we can translate all of the equations in $\mathcal E$ into rewrite rules, forming a rewrite system $\mathcal R$. Applying the axioms of an algebraic structure to a morphism in a monoidal category corresponds precisely to rewriting with the rules in $\mathcal R$. Purely by construction, a category of the form $\FCsp(\mathcal R, T)$ satisfies the axioms of $\mathcal E$, \textit{and nothing else}. We can formalise this using another unique factorisation property.

For a rewrite rule $L \rewritesto R \in \mathcal R$, the string graphs $L$ and $R$ (regarded as cospans over their inputs and outputs) define morphisms $|L|_{\mathbb H}, |R|_{\mathbb H} \in \hom_{\FCsp(\mathbb H, T)}(X,Y)$ in the free compact closed category over $T$.

\begin{theorem}
  Let $\mathcal R$ be a string graph rewrite system containing $\mathbb H$. Let $F : T \rightarrow \mathcal V$ be a valuation of a monoidal signature $T$ such that, for the functor $\widetilde F : \FCsp(\mathbb H, T) \rightarrow \mathcal V$ and every rewrite rule $L \rewritesto R \in \mathcal R$, $\widetilde F(|L|_{\mathbb H}) = \widetilde F(|R|_{\mathbb H})$, the valuation $F$ factors uniquely as:
  \begin{center}
    \begin{tikzpicture}
      \matrix (m) [cdiag] {
        T & \FCsp(\mathcal R, T) \\
          & \mathcal V           \\
      };
      \path [arrs]
        (m-1-1) edge node {$E_{\mathcal R}$} (m-1-2)
        (m-1-2) edge [dashed] node {$\widetilde F_{\mathcal R}$} (m-2-2)
        (m-1-1) edge node [swap] {$F$} (m-2-2);
    \end{tikzpicture}
  \end{center}
\end{theorem}

\begin{proof}
  Note that for any rewrite system $\mathcal R \supseteq \mathbb H$, then $\rewriteequiv_{\mathcal R} \supseteq \rewriteequiv_{\mathbb H}$, so we can think of the morphisms of $\FCsp(\mathcal R, T)$ as equivalence classes of morphisms in $\FCsp(\mathbb H, T)$, where the quotient functor $Q$ is identity-on-objects, and takes an equivalence class $|G|_{\mathbb H}$ to the unique (possibly larger) equivalence class $|H|_{\mathcal R} \supseteq |G|_{\mathbb H}$. This is a well-defined functor by Theorem~\ref{thm:rw-compose-commute}, and it is possible to define $E_{\mathcal R} := Q \circ E$. Thus is suffices to show that $\widetilde F$ factors uniquely through $Q$:
  \begin{center}
    \begin{tikzpicture}
      \matrix (m) [cdiag] {
        T & \mathcal V \\
        \FCsp(\mathbb H, T) & \FCsp(\mathcal R, T) \\
      };
      \path [arrs]
        (m-1-1) edge node {$F$} (m-1-2)
        (m-1-1) edge node [swap] {$E$} (m-2-1)
        (m-2-2) edge [dashed] node [swap] {$\widetilde F_{\mathcal R}$} (m-1-2)
        (m-2-1) edge node [swap] {$Q$} (m-2-2)
        (m-2-1) edge node [yshift=-1mm] {$\widetilde F$} (m-1-2);
    \end{tikzpicture}  
  \end{center}

  For existence, let $\widetilde F_{\mathcal R}(X) := \widetilde F(X)$ on objects. On arrows, let $\widetilde F_{\mathcal R}(|G|_{\mathcal R}) := \widetilde F(|G|_{\mathbb H})$. For this to be well-defined, it must not depend on the choice of representative $G \in |G|_{\mathcal R}$. First, assume that $G \rewritesto_{\mathcal R} G'$ for some $G'$. Then, there exists a rule $L \rewritesto R \in \mathcal R$ and a matching $m : L \rightarrow G$ that rewrites $G$ into $G'$. By deforming $G$ and $G'$, we can always find edge-homeomorphic cospans $H \in |G|_{\mathbb H}$ and $H' \in |G'|_{\mathbb H}$ such that:
  \begin{center}
    $H = H_{\textrm{out}} \circ (L \otimes H_{\textrm{side}}) \circ H_{\textrm{in}}$\ \ \ and\ \ \ $H' = H_{\textrm{out}} \circ (R \otimes H_{\textrm{side}}) \circ H_{\textrm{in}}$
  \end{center}
  where $\otimes$ and $\circ$ are tensor and composition of framed cospans. Since $|\cdot|_{\mathbb H}$ respects tensor and composition:
  \begin{center}
    $|G|_{\mathbb H} = |H_{\textrm{out}}|_{\mathbb H} \circ (|L|_{\mathbb H} \otimes |H_{\textrm{side}}|_{\mathbb H}) \circ |H_{\textrm{in}}|_{\mathbb H}$\ \ \ and\ \ \ $|G'|_{\mathbb H} = |H_{\textrm{out}}|_{\mathbb H} \circ (|R|_{\mathbb H} \otimes |H_{\textrm{side}}|_{\mathbb H}) \circ |H_{\textrm{in}}|_{\mathbb H}$.
  \end{center}
  
  Then, since $\widetilde F$ is a strict monoidal functor, we have that:
  \begin{center}
    $\widetilde F(|G|_{\mathbb H}) =
     \widetilde F(|H_{\textrm{out}}|_{\mathbb H})
       \circ (\widetilde F(|L|_{\mathbb H})
         \otimes \widetilde F(|H_{\textrm{side}}|_{\mathbb H}))
       \circ \widetilde F(|H_{\textrm{in}}|_{\mathbb H}) =
     \widetilde F(|H_{\textrm{out}}|_{\mathbb H})
       \circ (\widetilde F(|R|_{\mathbb H})
         \otimes \widetilde F(|H_{\textrm{side}}|_{\mathbb H}))
       \circ \widetilde F(|H_{\textrm{in}}|_{\mathbb H}) = \widetilde F(|G'|_{\mathbb H})$
  \end{center}
  
  Since $|G|_{\mathcal R}$ represents the closure of $G$ with respect to $\rewritesto_{\mathcal R}$, this argument can be iterated to show that for any $G,G' \in |G|_{\mathcal R}$, $\widetilde F(|G|_{\mathbb H}) = \widetilde F(|G'|_{\mathbb H})$. Thus $\widetilde F_{\mathcal R}$ is well-defined. Regarding the hom-sets of $\FCsp(\mathcal R, T)$ is quotients of hom-sets in $\FCsp(\mathbb H, T)$, $\widetilde F_{\mathcal R}$ and is unique functor induced by the universal property of quotients in $\catSet$.
\end{proof}

The exact same argument carries through replacing $\FCsp$ with $\FCsp^+$ everywhere. Thus, rewrite categories form the most general compact closed (or symmetric traced) categories containing the equivalence closure of a rewrite system. They form the formal basis for reasoning about all of graphical theories we introduce in the next part.


	\part{Entanglement, Graphically}\label{part:entanglement}
	



\chapter{Quantum Information and Entanglement}\label{ch:quantum}

Quantum information theory is the study of how data can be encoded and manipulated using microscopic systems subject to quantum effects. Over the past two decades, it has grown into a large and diverse field, with applications in security, foundations of physics, and perhaps most notably quantum computing. In this chapter, we introduce the basics of quantum mechanics, quantum information theory, and models of quantum computing.

\section{Quantum Mechanics}\label{sec:quantum-mechanics}

This section is provided for the non-physicist to briefly introduce the basic concepts of quantum mechanics used in this dissertation. Those familiar with QM can safely skip it.

Quantum mechanics is a strange but very successful theory of the universe at small scales. Its Hilbert-space formulation has four key components, which we shall focus on in detail:
\begin{enumerate}
\item \textbf{States} encode all of the information about a quantum system. These are represented as normalised vectors in a (complex) Hilbert space.
\item \textbf{Observables} give us ``questions'' to ask about a quantum system, and provide the mathematical means to turn a state into a probability distribution over measurement outcomes. These are given as self-adjoint operators on a Hilbert space.
\item \textbf{Dynamics} describe how a state evolves over time. These are expressed as unitary operators.
\item \textbf{Compound systems} are expressed as tensor products of simpler systems.
\end{enumerate}

To describe these components, we use Dirac's \textit{bra-ket} notation. Recall that for any vector $v$ in a Hilbert space $\mathcal H$, there is a natural way to get a linear map $\phi_{v} : \mathcal H \rightarrow \mathbb{C}$ (i.e. a vector in the dual space $\mathcal H^*$), using the inner product defined on $\mathcal H$:
\[ \phi_{v}(u) = \braket{v}{u} \]

Vectors in the dual space of $\mathcal H$ are used so often that we employ the following notational trick. We write a vector $v \in \mathcal H$ as a \emph{ket} $\ket v\in \mathcal H$, and we write the associated linear map $\phi_{v}$ as a \emph{bra} $\bra v$. $\bra v$ is a function from $\mathcal H$ to $\mathbb{C}$, so we can apply it to a vector $\ket u\in \mathcal H$. Then, some notational magic happens:
\[ \bra v \ket u = \braket {v}{u} \]

For this reason, we refer to the inner product as a \emph{bra-ket}. We define an operation $(-)^\dagger$ taking bras to kets and vice-versa.
\[ \left(\ket v\right)^\dagger = \bra{v} \textrm{\ \ \ and\ \ \ } \left(\bra{v}\right)^{\dagger}=\ket{v}\]

This operation naturally extends to linear maps:
\[ L^\dagger \ket{u} = \left(\bra uL\right)^\dagger \]

Fixing an orthonormal basis, with respect to $\braket{-}{-}$, we can represent $\ket v$ as a column vector and $L$ as a matrix. Then, $(-)^\dagger$ just becomes the conjugate-transpose:
\[
\left(\begin{matrix} v_1 \\ v_2 \end{matrix}\right) \mapsto
\left(\begin{matrix} \overline{v_1} & \overline{v_2} \end{matrix}\right)
\qquad\qquad
\left(\begin{matrix}
  a_{11} & a_{12} \\
  a_{21} & a_{22}
\end{matrix}\right) \mapsto
\left(\begin{matrix}
  \overline{a_{11}} & \overline{a_{21}} \\
  \overline{a_{12}} & \overline{a_{22}}
\end{matrix}\right)
\]

$(-)^\dagger$ extends to a contravariant functor $\dagger : \catHilb^{\textrm{op}} \rightarrow \catHilb$, giving $(\catHilb, \otimes, \mathbb C)$ the structure of a $\dagger$-monoidal category.

A vector $\ket\psi \in \mathcal H$ is the same thing as a linear map $\mathbb C \rightarrow \mathcal H$ sending $1 \in \mathbb C$ to $\ket\psi$. We use these two notions interchangeably. As string diagrams, we represent kets as triangles with a single out-edge and bras as triangles with a single in-edge.
\ctikzfig{bra_ket_graphical}

In quantum mechanics, pure states are unit vectors in a Hilbert space $\mathcal H$. Crucially, the fact that $\mathcal H$ is a vector space gives us a way to super-impose several quantum states to form a new state.
\[ \ket\xi = \sum_i \alpha_i \ket{\psi_i} \]
In this case, we say that the state $\ket\xi$ is in a superposition of the states $\ket{\psi_i}$.

If $\ket{\psi}$ is a state, then $\bra{\psi}$ can be thought of as a function that measures the extent to which some given state is $\ket{\psi}$. This is the essential content of the Born rule, which provides a method for turning a state and an observable into a probability distribution on measurement outcomes. An observable $O$ is a self-adjoint ($O = O^\dagger$) operator from a Hilbert space $\mathcal H$ to itself. If $\mathcal H$ is finite-dimensional then all self-adjoint operators diagonalise:
\[ O = \sum_i \alpha_i \ketbra{v_i}{v_i} \]
We can choose $\alpha_i$ such that each of the eigenvectors $\ket{x_i}$ are normalised, in which case we call them the \textit{eigenstates} of $O$. We can therefore interpret $O$ as a set $\left\{ \alpha_{i}\right\} $ of \emph{measurement outcomes} and a set $\left\{ \ket{v_{i}}\right\} $ of possible outcome states. Picture an experimental setup, where we have some quantum state $\ket{\psi}$ in a box and a measuring apparatus hooked to it. We dial in $O$ as the thing we want to measure and push a button. Suppose for simplicity that $O$ only has two possible outcomes:
\[ O = 1\ket{v_1}\bra{v_1}+2\ket{v_2}\bra{v_2} \]
\noindent When we push the button, the screen says $2$. Thus, we know the second measurement outcome occurred. In that sense, the eigenvalue $\alpha_i$ should be thought of as a ``marker'' for the $i$-th measurement outcome. The most informative observables are the ones where all of these markers are distinct and non-zero. These are called \textit{non-degenerate} observables. Three non-degenerate observables that play a particularly important role in quantum information and quantum computing are the Pauli spin operators:
\begin{equation}\label{eq:paulis}
  X = \left(\begin{matrix} 0 & 1 \\  1 &  0 \end{matrix}\right) \qquad\qquad
  Y = \left(\begin{matrix} 0 & i \\ -i &  0 \end{matrix}\right) \qquad\qquad
  Z = \left(\begin{matrix} 1 & 0 \\  0 & -1 \end{matrix}\right)
\end{equation}

Furthermore, upon getting outcome $2$, we know that the state $\ket{\psi}$ must be $\ket{v_2}$. Suppose beforehand that we had prepared $\ket{\psi}$ in some superposition of $\ket{v_1}$ and $\ket{v_2}$. As soon as we measured $\ket{\psi}$, this superposition \emph{collapsed} to a single state $\ket{v_2}$. This phenomenon is known as the collapse of the quantum state, or ``collapse of the wavefunction''. What this means physically is a question of interpretation, but mathematically it means after the measurement occurs, we can treat the quantum system as if it is in the state $\ket{v_2}$.

We compute the probability of getting outcome $i$ using the \textit{Born rule}.
\[ \textrm{Prob}(i,\ket{\psi}) = \vert\braket{v_{i}}{\psi}\vert^2 = \braket{\psi}{v_{i}}\braket{v_{i}}{\psi} \]
The key point here is that the Born rule is a function of the inner product. Because of the role it plays in measurement, we sometimes refer to elements of the dual space $\bra\psi \in \mathcal H^*$ as \textit{effects}.

Since $\ket\psi$ is normalised, the sum of the probabilities of all outcomes is $1$, so $\textrm{Prob}(i,\ket{\psi})$ is a probability distribution. These probabilities provide our only access to the ``real'' quantum state, so we consider two states to be equal if they give the same probability distributions with respect to any observable. However, the Born rule yields the same probabilities for states $\ket{\psi}$ and $e^{i\theta}\ket{\psi}$.
\[ \overline{e^{i\theta}} e^{i\theta} \braket{\psi}{v_{i}}\braket{v_{i}}{\psi}
   = e^{-i\theta} e^{i\theta} \braket{\psi}{v_{i}}\braket{v_{i}}{\psi}
   = \braket{v_{i}}{\psi}\braket{\psi}{v_{i}} \]
The scalar factor $e^{i\theta}$ is called a \textit{global phase}. We always identify states (and hence operators) differing only by a global phase.

The main reason the Pauli matrices are so interesting is that every distinct pair of them is complementary. Two observables are called \textit{complementary} if their associated bases of eigenstates are equally-far apart. Let $O,O'$ be observables in a $D$-dimensional Hilbert space with eigenstates $\{ \ket{v_i} \}$ and $\{ \ket{v'_j} \}$. Then $O$ and $O'$ are called complementary if their bases are \textit{mutually unbiased}. That is, they satisfy the following equation for all $i,j$:
\begin{equation}\label{eq:mutually-unbiased}
  |\braket{v_i}{v'_j}|^2 = \frac{1}{D}
\end{equation}

We can interpret this definition using the Born rule. If a state $\ket\psi$ is in the $i$-th eigenstate of $O$, then measuring $O$ will obtain outcome $i$ with \textit{certainty}. So, if we know $\ket\psi$ is in an eigenstate of $O$, we have maximal knowledge about the $O$ observable. However, if we measure $\ket\psi$ with respect to the $O'$ observable, we are equally likely to get any outcome. So, maximal knowledge about $O$ implies minimal knowledge about $O'$.

Unlike measurements which cause a state to collapse, dynamic evolution of quantum states is always reversible. We can evolve a quantum state in time by applying a unitary ($U^\dagger = U^{-1}$) operator. One should interpret the dagger of a unitary operator as that some operator, ``done backwards''. To be a bit more explicit, suppose we represent the evolution of a state according to the Schr\"odinger equation, for a self-adjoint operator $H$ called a \textit{Hamiltonian}.
\begin{equation}\label{eq:tdse}
  i\hbar \frac{d}{dt} \ket{\psi(t)} = H \ket{\psi(t)}
\end{equation}

If $H$ does not depend on $t$ (e.g. if the forces acting on a particle are constant), then the value of a solution at time $t$ has a simple expression in terms of the value at some initial time $0$.
\[ \ket{\psi(t)} = e^{-(i/\hbar) t H}\ket{\psi(0)} \]

Letting $U(t) = e^{-(i/\hbar) t H}$, then since $H$ is self-adjoint, $U(t)$ is unitary for all $t$ and describes the evolution of $\ket{\psi(0)}$ under $ H$ for time $t$. $U(-t)$ then corresponds to the same evolution, but with time running backwards. Since $U(t)U(-t) = U(-t)U(t) = U(0) = 1_{\mathcal H}$, it must be the case that $U(t)^\dagger = U(-t)$.

\section{Compound Systems and Entanglement}

Suppose we have a particle in the state $\ket{\psi}$ sitting in some potential well and another particle $\ket{\phi}$ sitting in another one far away. We use the tensor product to ``pair up'' the states of the two particles. That is, the overall state of the system is $\ket{\psi}\otimes\ket{\phi}$. But, since this is a quantum state, it could be in some superposition:
\[ \frac{1}{\sqrt{2}} \left(\ket{\psi_1}\otimes\ket{\phi_1} + \ket{\psi_2}\otimes\ket{\phi_2}\right) \]

Suppose one state is in a superposition $\ket{\psi_1}+\ket{\psi_2}$. Then the combined state is also in a superposition:
\[ \left(\ket{\psi_1}+\ket{\psi_2}\right)\otimes\ket{\phi}=\left(\ket{\psi_1}\otimes\ket{\phi}\right)+\left(\ket{\psi_2}\otimes\ket{\phi}\right)\]

This bilinearity justifies our use of the tensor product, since $\otimes$ provides the most general bilinear pairing for two spaces. Some states in the tensor product space $\mathcal H_1\otimes \mathcal H_2$ can be written in the form $\ket{\psi}\otimes\ket{\phi}$ for $\ket{\psi}\in \mathcal H_1$ and $\ket{\phi}\in \mathcal H_2$. Such states are called \emph{separable} or \textit{product states}. However the vast majority of states $\ket{\Psi}\in \mathcal \mathcal H_1\otimes \mathcal \mathcal H_2$ cannot be written this way. These are called \emph{entangled states}.

A measurement on one subsystem of an entangled state collapses the entire state, thus affecting the other subsystem. For instance, suppose we have some observable
\[ O=1\ket{v_1}\bra{v_1}+2\ket{v_2}\bra{v_2} \]
as before, and an entangled state
\[ \ket{\Psi}=\frac{1}{\sqrt{2}}\left(\ket{v_1}\otimes\ket{v_1}+\ket{v_2}\otimes\ket{v_2}\right) \]
and we measure $O$ on the left subsystem, getting outcome $1$. The whole state is now:
\[ \ket{\Psi'}=\ket{v_1}\otimes\ket{v_1} \]

If $\ket{\Psi}$ were a product state, then the second system would be unaffected, so entanglement is the crucial property that allows such correlations at a distance.

Entanglement is also a source of computational complexity for many-body systems. For finite dimensional Hilbert spaces, the dimension of the tensor product of two spaces is the product of the dimensions of each space.
\[ \dim\left(\mathcal H_1\otimes \mathcal H_2\right)=\dim\left(\mathcal H_1\right)\dim\left(\mathcal H_2\right) \]
So, the dimension of a compound system increases \emph{exponentially} with the number of subsystems.
\[ \dim(\,\underbrace{\mathcal H\otimes\ldots\otimes \mathcal H}_{N}\,)=\dim\left(\mathcal H\right)^{N} \]
Computing with such states quickly becomes untenable, even for low-dimensional $\mathcal H$, which is the main reason for trying to understand quantum phenomena like entanglement from a more structural level.

\section{Mixed State Quantum Mechanics}

Often it is more convenient to work with probabilistic mixtures of quantum states, rather than states that are totally determined. This is because nearly all procedures for preparing a quantum state in a lab only succeed with some probability. A set of quantum pure states along with associated probabilities $\{ (\ket{\psi_i}, p_i) \}$ is called an \textit{ensemble}. We can compute the probability of getting a particular measurement outcome on an ensemble using the Born rule, adjusting for probabilities.
\begin{equation}\label{eq:mixed-state-born}
\textrm{Prob}(i,\{p_j, \ket{\psi_j}\})
      = \sum_j p_j \vert\braket{v_i}{\psi_j}\vert^2
      = \sum_j p_j \braket{v_{i}}{\psi_j}\braket{\psi_j}{v_{i}}
      = \bra{v_i}\left(\sum_j p_j \ketbra{\psi_j}{\psi_j}\right)\ket{v_i}
\end{equation}

As in the case with pure states, we consider two ensembles to represent the same ``state'' if they produce identical probability distributions for all observables under the Born rule. From (\ref{eq:mixed-state-born}), we can see that two ensembles produce the same probability distributions if and only if:
\[ \sum_i p_i \ketbra{\psi_i}{\psi_i} = \sum_i q_i \ketbra{\phi_i}{\phi_i} \]

For normalised vectors $\ket{\psi_i}$ and $\sum_i p_i = 1$, this is the general form for a trace-1 positive operator. For that reason, we call trace-1 positive operators \textit{mixed states}. Because they encode probability densities, they are also sometimes called \textit{density matrices}. Pure states can be represented as density matrices of the form $\ketbra{\psi}{\psi}$.

Just as states cannot in practice be prepared with certainty, so too is the case for quantum evolutions and measurements. The mixed version of a unitary evolution is a completely positive map (CPM). For a finite-dimensional Hilbert space $\mathcal H$, let $\mathcal L(\mathcal H)$ be the vector space of linear maps $\mathcal H \rightarrow \mathcal H$. CPMs are just linear maps $\Phi : \mathcal L(H) \rightarrow \mathcal L(\mathcal H')$ that take positive operators to positive operators.

Outcomes for pure measurements span an orthonormal basis. A particular observable $O$ corresponds to the decomposition of the identity into 1-dimensional projectors corresponding to each outcome for $O$.
\[ 1_{\mathcal H} = \sum_i \ketbra{v_i}{v_i} \]

The probability of getting a \textit{particular} measurement outcome on a mixed state can be computed as in \ref{eq:mixed-state-born} by tracing the composition of the projector and the state's density matrix.
\[ \textrm{Prob}(i,\rho) = \bra{v_i}\rho\ket{v_i} = \Tr(\bra{v_i}\rho\ket{v_i}) = \Tr(\rho\ketbra{v_i}{v_i}) \]

For that reason, pure measurements are often referred to as \textit{projective measurements}. The mixed version of a projective measurement is a \textit{positive operator-valued map} (POVM). In finite dimensions, this is just a set of positive operators $P_i \in \mathcal L(\mathcal H)$ that sum to the identity. As in the projective case, probabilities are computed by tracing the composition of the positive operator and the state's density matrix.
\[ \textrm{Prob}(i,\rho) = \Tr(\rho P_i) \]

For the majority of this dissertation, we will only need concepts from pure-state quantum mechanics. However, when we look at multipartite entanglement in chapter \ref{ch:monoidal-entanglement}, it will occasionally be useful to ignore a subsystem of an entangled pure state. This can be done probabilistically by tracing out that subsystem and renormalising (if necessary). This is called the \textit{reduced density matrix} of an entangled state. For a state $\ket\Psi \in \mathcal H_1 \otimes \mathcal H_2$ and $\rho_{12} := \ketbra{\Psi}{\Psi}$, we can ignore the system $\mathcal H_2$ by tracing it out.
\[ \rho_1 = \Tr^{\mathcal H_2}(\rho_{12}) \]

\section{Quantum Computation}\label{sec:quantum-computation}

Quantum computation refers to the encoding of data into quantum states and the use of evolution and measurements to perform computations on that data. As in programming language design, quantum computation can be carried out using one of a variety of paradigms. We shall introduce two of them here: the \textit{circuit model} and \textit{measurement-based quantum computation} (MBQC).

For a $D$-dimensional Hilbert space, we often fix an orthonormal basis $\ket{0}, \ket{1}, \ldots, \ket{D-1}$ called the \textit{computational basis}. The two-dimensional Hilbert space $\mathbb C^2$ plays a special role in quantum computation, and is called the space of quantum bits, or \textit{qubits}. The basis vectors $\ket{0},\ket{1} \in \mathbb C^2$ can be thought of as classical bits, embedded in the bigger space of qubits. We also introduce a special notation for tensor products of basis vectors using bit strings:
\[ \ket{0010110} = \ket 0 \otimes \ket 0 \otimes \ket 1 \otimes 
                   \ket 0 \otimes \ket 1 \otimes \ket 1 \otimes \ket 0 \]

In $\mathbb C^2$ the space of qubits, the basis $\ket 0, \ket 1$ corresponds to the eigenbasis of the $Z$ observable. We use the term ``measuring in the computational basis'' to mean performing a measurement with respect to $Z$. We also define the $X$ basis $\ket +, \ket -$ and the $Y$ basis $\ket{i}, \ket{-i}$ as follows.
\begin{align*}
  \ket{+} & := \frac{1}{\sqrt{2}} \left(\ket 0 + \ket 1\right) & \ket{i}  & := \frac{1}{\sqrt{2}} \left(\ket 0 + i \ket 1\right) \\
  \ket{-} & := \frac{1}{\sqrt{2}} \left(\ket 0 - \ket 1\right) & \ket{-i} & := \frac{1}{\sqrt{2}} \left(\ket 0 - i \ket 1\right)
\end{align*}

\subsection{The Circuit Model}\label{sec:circuit-model}

In the circuit model, quantum computation proceeds in three steps:
\begin{enumerate}
  \item Prepare an $N$-qubit quantum state (usually a product state). Some of these qubits are treated as inputs and others simply as ``helper'' qubits called ancillas (which are usually initialised to $\ket 0$).
  \item Evolve the prepared state using small (usually 1- or 2-qubit), fixed-time evolutions called quantum gates.
  \item Measure some or all of the qubits, yielding the result of the computation. Unmeasured qubits are sometimes treated as outputs.
\end{enumerate}

Graphically, we represent a circuit evaluation as a string diagram.
\ctikzfig{quantum_circuit_ex}

Usually measurements are performed over the computational basis. That is, we measure the Pauli $Z$ observable given in equation (\ref{eq:paulis}). Using the Pauli operators as the Hamiltonian in the Schr\"odinger equation (\ref{eq:tdse}), we can produce unitary evolutions we call the \textit{phase gates}.
\begin{equation*}
  X_{\theta} = e^{-(i \theta/2) X} \qquad\qquad
  Y_{\theta} = e^{-(i \theta/2) Y} \qquad\qquad
  Z_{\theta} = e^{-(i \theta/2) Z}
\end{equation*}

Up to a global phase, we can recover identities and the Pauli gates themselves as phase gates.
\[ X_0 = Y_0 = Z_0 = 1_{\mathbb C^2} \qquad
   X_{\pi} = X \qquad Y_{\pi} = Y \qquad Z_{\pi} = Z \]

It is a well known fact that \textit{any} single-qubit unitary can be constructed from phase gates as $U = Z_\alpha X_\beta Z_\gamma$. This is called the \textit{Euler decomposition} of the unitary. An important single-qubit gate that is not a phase gate is the Hadamard gate $H$.
\[ H = Z_{-\pi/2} X_{-\pi/2} Z_{-\pi/2} = \frac{1}{\sqrt{2}}
   \left( \begin{matrix} 1 & 1 \\ 1 & -1 \end{matrix}\right) \]

$H$ interchanges the eigenbasis $\ket 0$, $\ket 1$ of $Z$ with the eigenbasis $\ket +, \ket -$ of $X$. A simple consequence is that $H Z_\theta H = X_\theta$. Perhaps the most common 2-qubit gate is the \textit{controlled-NOT} or CNOT gate.
\[
\beginpgfgraphicnamed{cnot}
\InputIfFileExists{cnot.tikz}{}{\input{./figures/cnot.tikz}}
\endpgfgraphicnamed\ \  ::\ 
\begin{cases}
  \ \ket{00} \mapsto \ket{00} \\
  \ \ket{01} \mapsto \ket{01} \\
  \ \ket{10} \mapsto \ket{11} \\
  \ \ket{11} \mapsto \ket{10} \\
\end{cases}
\]

It is called the controlled-NOT gate because the first qubit controls whether the second qubit has a NOT (a.k.a. $X$) gate applied to it. A generalisation of CNOT gates are controlled-unitary gates, which conditionally apply a unitary to the second qubit.
\[
\beginpgfgraphicnamed{cunitary}
\InputIfFileExists{cunitary.tikz}{}{\input{./figures/cunitary.tikz}}
\endpgfgraphicnamed\ \  ::\ 
\begin{cases}
  \ \ket{0} \otimes \ket\psi \mapsto \ket{0} \otimes \ket\psi \\
  \ \ket{1} \otimes \ket\psi \mapsto \ket{1} \otimes (U \ket\psi) \\
\end{cases}
\]

These are examples of gates that can take product states to entangled states. As such, they are sometimes referred to as \textit{entangling gates}. For example, applying a CNOT gate to the state $\ket + \otimes \ket 0$ yields an entangled state $\frac{1}{\sqrt{2}} (\ket{00} + \ket{11})$ called the \textit{Bell state}. Examples of gates that are not entangling gates are tensor products of 1-qubit gates, such as $X \otimes Z$. These will \textit{always} take product states to product states.

Just as AND, OR, and NOT gates can be used to construct arbitrary classical circuits, we have a notion for a set of quantum gates being able to construct arbitrary quantum circuits. We say a set of gates is \textit{universal for quantum computation} if any $N$-qubit unitary map can be constructed from compositions of those gates.

\begin{theorem}[\cite{BBC1995}]
  The gates $Z_\theta$, $H$, and CNOT are universal for quantum computation.
\end{theorem}

Many important quantum algorithms, such as the quantum Fourier transform, Shor's factoring algorithm, and Grover's search algorithm, can be presented in the circuit model.

\subsection{Measurement-based Quantum Computation}\label{sec:mbqc}

Measurement based quantum computation (MBQC), which is sometimes called one-way quantum computation, provides a different, equally powerful paradigm for quantum computation. For comparison to the previous section, we can organise the MBQC procedure into three steps.
\begin{enumerate}
  \item Prepare a known, highly-entangled state called a \textit{graph state}. This graph state may be entangled to some qubits in an unknown state called input qubits.
  \item Perform measurements at arbitrary angles, where the choice of angles can depend on previous measurement outcomes.
  \item Optionally, perform single-qubit corrections on unmeasured, output qubits.
\end{enumerate}

Graph states are constructed by preparing a collection of (non-input) qubits in the $\ket +$ state, then applying controlled-$Z$ gates to pairs of qubits to introduce entanglement. We represent such a state by drawing a vertex for every qubit and an edge whenever a controlled-$Z$ gate is applied.
\ctikzfig{graph_state}

This representation is unambiguous because controlled-$Z$ gates are symmetric and commute past each other. We define generalised measurements, or \textit{measurements with angles} as follows. An $X_\theta$ measurement on single qubit consists of first applying the unitary gate $X_\theta$ then measuring the $X$ observable. A $Z_\theta$ measurement consists of first applying the unitary gate $Z_\theta$ then measuring the $Z$ observable. Since they are measurements, neither of these operations is deterministic.

Each measurement has two possible outcomes: a ``desired'' outcome, and an ``erroneous'' outcome. The latter can be thought of as the quantum version of an ``accidental bit flip'' during the course of the computation. The key to MBQC is that if we choose our measurement angles wisely, we can \textit{correct} these errors as we go. If there are no measurements left to perform, we complete the calculation by applying any remaining corrections as single-qubit unitaries on the output qubits. This is known as \emph{feeding forward} corrections. There are several techniques for identifying and using graph states for deterministic MBQC, such as identifying a \textit{generalised flow} for the graph~\cite{BrowneGFlow2007}. To give a feel for how feed-forward works, we provide a simple example.

\begin{example}
Prepare the following graph state, where $q_1$ and $q_3$ are in an unknown state (i.e. they are inputs) and $q_2$ is prepared in the $\ket +$ state.
\ctikzfig{cnot_graph_state}

First measure $q_3$ in $Z$, getting outcome $i \in \{ 0, 1 \}$. Letting $Z_0 = 1_{\mathbb C^2}$ and $Z_1 = Z$, apply $Z_i$ to $q_1$ and apply $Z_i H$ to $q_2$. Treating $q_1$ and $q_2$ as outputs, this procedure computes the CNOT of $q_1$ and $q_3$, regardless of the outcome of the measurement of $q_3$. As $Z_i$ depends on the outcome of measuring $q_3$, this is an example of feeding-forward a measurement outcome to a correction.
\end{example}

This may seem like an excessively roundabout way to apply a CNOT gate, especially since a na\"ive implementation involves applying 2 controlled-$Z$ gates to prepare the graph state. However, if one assumes we have a stock of suitably nice graph states, we can actually perform \textit{arbitrary} quantum computations using (comparatively easy) single qubit unitaries and measurements.

\chapter{Categorical Quantum Mechanics}\label{ch:categorical-quantum-mechanics}

Categorical quantum mechanics (CQM) refers to a broad program initiated by Abramsky and Coecke in 2004 \cite{AC2004} that emphasises the abstract, categorical, and compositional aspects of quantum mechanics. The core thesis is that the structure-rich setting of Hilbert spaces and linear operators obscures the causes of many quantum phenomena. Therefore reasoning in the comparatively sparse context of an arbitrary monoidal category yields practical and foundational insights that were previously hidden. Depending on the aims of a particular project under the umbrella of CQM, this proscription of Hilbert spaces can be taken literally, as in Paquette's PhD thesis~\cite{PaquetteThesis2008} for example, in order to obtain structural results about quantum theory that are independent of Hilbert space formulation. This route could prove fruitful in light of doubts that the Hilbert space formalism is the ``correct'' way to think about quantum mechanics, as expressed in these now-infamous words of its progenitor.

\begin{quote}
  \textit{I would like to make a confession which may seem immoral: I do not believe absolutely in Hilbert space any more.}
  \cite[John von Neumann (1935)]{VNLetter}
\end{quote}

Alternatively, one can take a less ``hard-line'' approach by using categorical techniques to complement and expand upon concrete results based on Hilbert spaces. We adopt this approach in the sections to come.

This chapter offers an introduction to categorical quantum mechanics and a handful of illustrative examples. We employ notions from CQM to show how complementary observables can be studied as interacting Frobenius algebras and offer several new results about special types of complementary observables called strongly complementary observables. Most notably, we give a classification theorem for strongly complementary pairs of observables and show that a set of pairwise strongly complementary observables must contain no more than 2 distinct observables.

\section{Compact Categories and Teleportation}

Quantum protocols exploit the unique features of quantum mechanics (typically entanglement and superposition) to perform a task that would be difficult or impossible classically. The canonical example of a quantum protocol is quantum teleportation, whereby one party (called Alice) can transmit an arbitrary quantum state to a second party (called Bob) using only a shared entangled state and a classical data channel. To start, Alice has a quantum state $\ket\psi$ and Alice and Bob share a Bell pair. That is, a pair of qubits in the state $\ket{00} + \ket{11}$ (ignoring normalisation factors).

Alice performs an entangled, 2-qubit measurement called a \textit{Bell measurement} on $\ket\psi$ and her half of the Bell pair. A Bell measurement consists of measuring the two qubits in the Bell basis:
\[ \ket{\Psi_0} = \ket{00} + \ket{11} \quad\ \ 
   \ket{\Psi_1} = \ket{00} - \ket{11} \quad\ \ 
   \ket{\Psi_2} = \ket{01} + \ket{10} \quad\ \ 
   \ket{\Psi_3} = \ket{01} - \ket{10} \]

She gets an outcome $i \in \{ 0, 1, 2, 3 \}$, which she then sends to Bob. Bob then applies a unitary correction to his half of the Bell pair, based on $i$:
\[ U_0 = 1_{\mathbb C^2} \qquad\qquad U_1 = Z \qquad\qquad U_2 = X \qquad\qquad U_3 = XZ \]

Once this is done, Bob's qubit will be in the state $\ket\psi$, i.e. the state of $\ket\psi$ has been teleported to Bob. We can represent this protocol in circuit language:
\ctikzfig{teleportation}

This protocol works because, by performing a Bell basis measurement, Alice projects out her two qubits using the associated bra $\bra{\Psi_i}$.
\ctikzfig{teleportation2}

We can express all four elements of the Bell basis in terms of the Bell state and the corrections we defined before.\footnote{Note that applying $U_i$ to the right qubit of $\ket{\Psi}$ has the same affect as applying $U_i^\dagger$ on the left qubit because all of the maps $U_i$, written as matrices over the computational basis, have real entries. A more generalised scheme is provided in~\cite{AC2004}, replacing $(-)^\dagger$ with $(-)^\whitetranspose$.}
\[ \ket{\Psi_i} = (1 \otimes U_i)\ket{\Psi_0} = (U_i^\dagger \otimes 1)\ket{\Psi_0} \qquad
   \bra{\Psi_i} = \bra{\Psi_0}(1 \otimes U_i^\dagger) = \bra{\Psi_0}(U_i \otimes 1) \]

Thus, we can pull the $U_i$ out of the measurement all the way to the end, and prove the teleportation protocol works for all $i$.
\ctikzfig{teleportation3}

The crucial step is ($*$). This identity should look familiar. Teleportation, like many quantum protocols exploits the fact that $\catFHilb$ is compact-closed and finite-dimensional Hilbert spaces are all self-dual $(X \cong X^*)$.
\ctikzfig{teleportation4}

Using this insight, one can perform teleportation in any self-dual compact-closed category, including $\catRel$, $\catMat(\mathbb R)$, and (perhaps surprisingly) $\catSpek$, the category Rob Spekkens' toy theory. The last example is surprising, because $\catSpek$ can be defined using a local hidden variable model. Thus teleportation succeeds even in the absence of non-locality for a physical theory. For more details, see \cite{EdwardsPhase2011,EdwardsSpek2011}.

\section{Complementary Observables as Frobenius Algebras}\label{sec:complementary-obs}

The eigenstates of an observable play a key role in quantum mechanics. They form the set of possible outcome states one obtains by performing a measurement. Classical data is obtained from a quantum system via measurements, so an orthonormal basis of measurement outcomes can be thought of as a particular \textit{classical context} embedded in the overall quantum state space. In studying the interaction of multiple classical contexts (especially complementary ones), we can see the unique features of quantum mechanics. The question is, can we study the concept of a ``basis'' over an object in an arbitrary $\dagger$-compact category?

Recall in example \ref{ex:basis-scfa}, we used a basis of a vector space to construct a special commutative Frobenius algebra.
\[ \delta :: e_i \mapsto e_i \otimes e_i \qquad\qquad \epsilon :: e_i \mapsto 1 \qquad\qquad
   \mu :: e_i \otimes e_i \mapsto e_i \qquad\qquad \eta :: \sum e_i \]

It was also noted that \textit{all} SCFAs over an algebraically closed field are of this form. So, there is a one-to-one correlation between SCFAs and arbitrary bases. However, projective measurements like the ones we have described have outcomes in an \textit{orthonormal} basis. Orthonormal bases can be captured in a $\dagger$-compact closed category using $\dagger$-special commutative Frobenius algebras.

\begin{definition}\label{def:dag_scfa}
  A \textit{$\dagger$-special commutative Frobenius algebra}, or $\dagger$-SCFA, $(A,\delta^\dagger,\epsilon^\dagger,\delta,\epsilon)$ is a $\dagger$-Frobenius algebra such that $\delta^\dagger \delta = 1_A$.
\end{definition}

In \cite{CPV2008} Coecke, Pavlovic, and Vicary showed that $\dagger$-SCFAs in $\catFHilb$ are in one-to-one correspondence with orthonormal bases. So, for any orthonormal basis $\{ \ket{i} \}$ in a finite-dimensional Hilbert space, there exists a unique $\dagger$-SCFA whose comultiplication copies the basis vectors and whose counit deletes them.
\[ \delta :: \ket{i} \mapsto \ket{ii} \qquad\qquad \epsilon :: \ket{i} \mapsto 1 \]

The basis vectors $\ket{i}$ are called the \textit{classical points} of $\delta$. This respects the no-cloning principal in quantum mechanics, because $\delta$ cannot copy any arbitrary state, only those in $\{ \ket{i} \}$. In fact, one can prove that these are the \textit{only} vectors copied by $\delta$ and deleted by $\epsilon$, so a basis can always be recovered from a $\dagger$-SCFA by taking the set of classical points. For an observables $O$ and $O'$, let the associated $\dagger$-SCFAs be:
\ctikzfig{dag_scfa_O}

Then, we write classical points as triangles of the same colour, and their associated bras as upside-down triangles of that colour.
\ctikzfig{classical_points}

Two observables are complementary if their bases of eigenstates are mutually unbiased. That is, for any $i,j$, $|\braket{v_i}{v_j'}|^2 = 1/D$. In the graphical notation:
\ctikzfig{mub}

A question posed by Coecke and Duncan \cite{Coecke2008} was, ``Can we represent complementarity purely in terms of interacting Frobenius algebras?'' It turns out that complementarity is equivalent to a simple diagrammatic identity between two $\dagger$-SCFAs. First, we can move $1/D$ in the above equation to the other side and express it as a circle, as the trace of the identity always equals $D$. Then, replace $1$ on the RHS with ``deleted points''.
\begin{equation}\label{eq:mub2}
\beginpgfgraphicnamed{mub2}
\InputIfFileExists{mub2.tikz}{}{\input{./figures/mub2.tikz}}
\endpgfgraphicnamed
\end{equation}

Frobenius algebras fix an isomorphism of a space with its dual space. In the case of $\dagger$-SCFAs, this isomorphism takes a classical point to its adjoint: $\ket{v_i}^\whitetranspose = \bra{v_i}$ and $\ket{v_j'}^\graytranspose = \bra{v_j'}$. Graphically:
\ctikzfig{dag_frob_transpose}

We can simplify the LHS of equation (\ref{eq:mub2}) using this fact.
\ctikzfig{mub3}

The equation $(*)$ is due to two applications of the spider theorem to merge the grey and white vertices, leaving a single edge connecting grey to white. Plugging this back into equation (\ref{eq:mub2}), we get:
\ctikzfig{op_dir_hopf_plugged}

Since this equation holds for all $i,j$ and the classical points span the entire space, we can conclude that a more general identity holds:
\begin{equation}\label{eq:op-dir-hopf}
\beginpgfgraphicnamed{op_dir_hopf}
\InputIfFileExists{op_dir_hopf.tikz}{}{\input{./figures/op_dir_hopf.tikz}}
\endpgfgraphicnamed
\end{equation}

Suppose we define a map $S$, serving as an antipode:
\begin{equation}\label{eq:mub-antipode}
\beginpgfgraphicnamed{mub_antipode}
\InputIfFileExists{mub_antipode.tikz}{}{\input{./figures/mub_antipode.tikz}}
\endpgfgraphicnamed
\end{equation}

Then, equation (\ref{eq:op-dir-hopf}) resembles the equation found in Definition \ref{def:hopf-algebra} of a Hopf algebra.
\begin{equation*}
\beginpgfgraphicnamed{antipode_hopf}
\InputIfFileExists{antipode_hopf.tikz}{}{\input{./figures/antipode_hopf.tikz}}
\endpgfgraphicnamed
\end{equation*}

For that reason, we refer to (\ref{eq:op-dir-hopf}) as the \textit{Hopf law}.

\begin{theorem}\label{thm:hopf-complementary}
  Two $\dagger$-SCFAs correspond to complementary observables if and only if they satisfy the Hopf law.
\end{theorem}

All complementary observables satisfy the Hopf law. When a complementary pair of observables actually extends to a (scaled) Hopf algebra, we call them \textit{strongly complementary}.

\begin{definition}\label{def:strongly-complementary}
  Two observables $O$ and $O'$ are called \textit{strongly complementary} if their associated $\dagger$-SCFAs satisfy the following equations, called the \textit{scaled bialgebra equations}.
  \ctikzfig{zx_bialg}
\end{definition}

Note that we have only required that $(A, \mu_O, \eta_O, \delta_{O'}, \epsilon_{O'})$ be a \textit{bialgebra}, up to scalar factors. However, we can show that bialgebras consisting of $\dagger$-SCFAs automatically satisfy equation (\ref{eq:op-dir-hopf}), so they are Hopf algebras. Before we can show this, we need a couple of lemmas. For the remainder of the section, let $O = (\whitemult, \whiteunit, \whitecomult, \whitecounit)$ and $O' = (\graymult, \grayunit, \graycomult, \graycounit)$ be strongly complementary $\dagger$-SCFAs in $\catFHilb$.

\begin{lemma}\label{lem:monoid-cp}
  Up to a scalar, $\whitemult$ is a monoid over the classical points of $\graycomult$. For all $i,j$, the following are classical points for $O'$:
  \begin{equation}\label{eq:cp_mult}
\beginpgfgraphicnamed{cp_monoid}
\InputIfFileExists{cp_monoid.tikz}{}{\input{./figures/cp_monoid.tikz}}
\endpgfgraphicnamed
  \end{equation}
\end{lemma}

\begin{proof}
  We can show that the point labelled $i \cdot j$ is copied using the first bialgebra rule.
  \ctikzfig{cp_mult_copy}
  
  Deletion follows from the dagger of the second bialgabra rule.
  \ctikzfig{cp_mult_delete}
  
  We can apply the bialgebra:
  \ctikzfig{cp_unit}
  
  Since the scalar $\graycounit \circ \whiteunit$ is non-zero, $e$ is a classical point.
\end{proof}

This is the property that Coecke and Duncan refer to as \textit{closure}.

\begin{lemma}\label{lem:sc-transpose}
  For a strongly complementary pair of observables, $(\whitemult)^\graytranspose = \whitecomult$ and $(\whiteunit)^\graytranspose = \whitecounit$.
\end{lemma}

\begin{proof}
  We can use the previous lemma to evaluate over classical points for $\graycomult$. For multiplication:
  \ctikzfig{strong_transpose_pf}
  
  And for unit:
  \ctikzfig{strong_transpose_pf2}
\end{proof}

A simple consequence of this lemma is:
\ctikzfig{scalar_circ_pf}.

\begin{lemma}\label{lem:sc-antipode}
  The antipode map $S$ defined in figure (\ref{eq:mub-antipode}) is self-adjoint and is an automorphism for both Frobenius algebras.
\end{lemma}

\begin{proof}
  For $S$ to be self-adjoint it suffices to show we can interchange the caps and cups.
  \ctikzfig{mub_antipode_sa_pf}
  
  To show $S$ is a Frobenius algebra automorphism, we can use the previous identity and the fact that it is copied by $\whitecomult$.
  \ctikzfig{mub_antipode_automorphism_pf}
\end{proof}

\begin{theorem}\label{thm:strong-compl-hopf}
  A strongly complementary pair of observables forms a scaled Hopf algebra with antipode $S$.
  \ctikzfig{antipode_hopf}
\end{theorem}

\begin{proof}
  The proof follows straightforwardly from the bialgebra identities and Lemmas \ref{lem:sc-transpose} and \ref{lem:sc-antipode}.
  \ctikzfig{sc_hopf_pf}
\end{proof}

Using the results from this section, we can prove a stronger classification result than the ones given in \cite{CoeckeDuncan2009} for strongly complementary observables in $\catFHilb$.

\begin{theorem}\label{thm:sc-classification}
  Let $(G,\cdot,e)$ be a finite Abelian group of order $D$, and $\{ \ket{g} : g \in G \}$ be a $D$-dimensional orthonormal basis. Every strongly complementary pair $O,O'$ of $\dagger$-SCFAs is of the following form.
  \begin{align*}
    \delta_{O}   & :: \ket{g} \mapsto \ket{g} \otimes \ket{g}                            & 
    \epsilon_{O} & :: \ket{g} \mapsto 1                                                  \\
    \mu_{O'}     & :: \ket{g} \otimes \ket{h} \mapsto \frac{1}{\sqrt{D}} \ket{g \cdot h} & 
    \eta_{O'}    & ::  1 \mapsto \sqrt{D} \ket{e}                                                   
  \end{align*}
\end{theorem}

\begin{proof}
  First, we show this is indeed a strongly complementary pair. $(\delta_O, \epsilon_O)$ copies and deletes and orthonormal basis, so it extends to a $\dagger$-SCFA. Also, up to a scalar, $(\mu_{O'},\eta_{O'},\mu_{O'}^\dagger,\eta_{O'}^\dagger)$ is the induced Frobenius algebra of the group algebra $\mathbb C[G]$. It is a routine calculation to show that the factors of $1/\sqrt{D}$ and $\sqrt{D}$ cancel out where necessary in the monoid, comonoid, and Frobenius identities. We can give explicit forms for $\delta_{O'} = \mu_{O'}^\dagger$ and $\epsilon_{O'} = \eta_{O'}^\dagger$.
  \begin{align*}
    \delta_{O'}   & :: \ket{g} \mapsto \frac{1}{\sqrt{D}} \sum_{g_1 \cdot g_2 = g} \ket{g_1} \otimes \ket{g_2} &
    \epsilon_{O'} & =  \sqrt{D} \bra{e}
  \end{align*}
  
  To show specialness, evaluate $\mu \circ \delta$ for any $g \in G$:
  \[ \mu_{O'} \delta_{O'} \ket{g}
      = \mu_{O'} \left(\frac{1}{\sqrt{D}} \sum_{g_1 \cdot g_2 = g} \ket{g_1} \otimes \ket{g_2} \right)
      = \frac{1}{D} \sum_{g_1 \cdot g_2 = g} \ket{g} \]
  Every element in $G$ has exactly $|G| = D$ distinct factorisations (i.e. pairs $(gh, h^{-1})$ for all $h \in G$), so $\mu_{O'}\delta_{O'}\ket{g} = \ket{g}$. We can also compute the explicit form for the ``cap'' $\delta_{O'} \circ \eta_{O'}$.
  \[ \delta_{O'} \eta_{O'} = \delta_{O'} \left( \sqrt{D} \ket{e} \right)
                           = \sqrt{D} \sum_{g_1 \cdot g_2 = e} \ket{g_1} \otimes \ket{g_2}
                           = \sqrt{D} \sum_{g \in G} \ket{g} \otimes \ket{g^{-1}} \]
  From this it follows that $S\ket{g} = \ket{g^{-1}}$ and $(\mu_{O'}, \eta_{O'}, \delta_{O}, \epsilon_{O}, S)$ is (up to a scalar) the induced Hopf algebra of the group algebra $\mathbb C[G]$. Therefore $O$ and $O'$ are strongly complementary.
  
  Conversely, let $(\whitemult, \whiteunit, \whitecomult, \whitecounit)$ and $(\graymult, \grayunit, \graycomult, \graycounit)$ be strongly complementary $\dagger$-SCFAs. From Lemma \ref{lem:monoid-cp}, $\graymult$ is a monoid over the classical points of $\whitecomult$. We can then evaluate both sides of the equation from Theorem \ref{thm:strong-compl-hopf} over an arbitrary classical point.
  \ctikzfig{sc_hopf_inverse}
  
  Since $S$ is a Frobenius algebra automorphism, it is a permutation of the classical points of $\whitecomult$. Thus the previous equation implies that for all $i$, there exists $i'$ such that:
  \ctikzfig{sc_hopf_inverse2}
  
  In other words, all of the classical points of $\whitecomult$ have inverses, so $\graymult$ is isomorphic to the group algebra $\mathbb C[G]$ for some Abelian group $G$.
\end{proof}

For a $D$-dimensional space in $\catFHilb$, fixing $\whitecomult$ as the $\dagger$-SCFA corresponding to the computational basis $\mathcal B = \{ \ket{i} \}$, we can find the basis inducing the group algebra $\mathbb C[\mathbb Z_D]$ corresponding to the order-$D$ cyclic group by applying the $D$-dimensional Fourier transform $F_D$ to the elements of the computation basis: $\mathcal B' = \{ F_D \ket{i} \}$. Using the fact that $\mathbb C[G \times G'] \cong \mathbb C[G] \otimes \mathbb C[G']$ for Abelian groups $G, G'$, it is possible to compute the strongly complementary pair corresponding to \textit{any} finite Abelian group $G$ by decomposing it into its cyclic components.

\begin{example}
  In $\mathbb C^4$, let $\mathcal B$ be the computational basis. Then, we can compute the strongly complementary basis $\mathcal B' = \{ \ket{e_i} \}$ corresponding to $\mathbb Z_2 \times \mathbb Z_2$ as follows. The 2D Fourier transform is just the Hadamard matrix $H$. So $\ket{e_i} = (H \otimes H)\ket{i}$. Writing these as column vectors in the compuational basis, we have:
  \[ \mathcal B' = \left\{
    \frac{1}{\sqrt{4}} \left(\begin{matrix}  1 \\  1 \\  1 \\  1 \end{matrix} \right),\ 
    \frac{1}{\sqrt{4}} \left(\begin{matrix}  1 \\ -1 \\  1 \\ -1 \end{matrix} \right),\ 
    \frac{1}{\sqrt{4}} \left(\begin{matrix}  1 \\  1 \\ -1 \\ -1 \end{matrix} \right),\ 
    \frac{1}{\sqrt{4}} \left(\begin{matrix}  1 \\ -1 \\ -1 \\  1 \end{matrix} \right)
  \right\} \]
\end{example}

In the case of complementary observables, it is often useful to know how big a \textit{complete set} of mutually unbiased bases is for a given dimension. That is, a maximal set of bases such that is pairwise mutually unbiased. In the case of strongly complementary observables, there can only be two.

\begin{theorem}
  Let $\graycomult$ be a $\dagger$-SCFA of dimension $D \geq 2$ and let $(\graycomult, \comult)$ and $(\graycomult,\whitecomult)$ be strongly complementary pairs. Then $\comult$ and $\whitecomult$ cannot be strongly complementary.
\end{theorem}

\begin{proof}
  By contradiction. Suppose $(\comult,\whitecomult)$ is a strongly complementary pair. The units $\unit$ and $\whiteunit$ must both be proportional to classical points of $\graycomult$. We already showed that for any strongly complementary pair, $\counit \circ \whiteunit = \whitecounit \circ \unit = \sqrt{D} \neq 0$, so $\unit$ and $\whiteunit$ must be proportional to the \textit{same} classical point. Then:
  \ctikzfig{two_strong_compl_pf}
  
  This is a contradiction because the LHS is invertible, while the RHS is rank $1 < D$.
\end{proof}

The classification of strongly complementary observables is much simpler than the general case. Whereas the maximum number of mutually unbiased bases of dimension $6$ is still unknown, there is (up to isomorphism) one strongly complementary pair, corresponding to the cyclic group $\mathbb Z_6 \cong \mathbb Z_2 \times \mathbb Z_3$.

Phases have a special status for $\dagger$-SCFAs in $\catFHilb$. Recall that phases are maps such that:
\ctikzfig{frobenius_phase_white}

From Proposition \ref{prop:frobenius-phase-form}, we can put any phase in the form of the right multiplication by an arbitrary vector. Suppose a $\dagger$-SCFA corresponds to a basis $\{ \ket{i} \}$, then $\mu = \sum \ket{i}\bra{ii}$. For an arbitrary vector $\ket{\psi} = \sum \alpha_i \ket{i}$, this is:
\[ \mu(1 \otimes \ket{\psi}) = \sum \alpha_i \ketbra{i}{i} \]

So, phases are precisely the maps that are diagonal in the basis defined by a $\dagger$-SCFA. \textit{Unitary} phases are precisely the phase gates familiar from quantum computing.

For strongly complementary observables, the phases for $\graycomult$ associated with classical points of $\whitecomult$ are Frobenius algebra automorphisms of $\whitecomult$, up to a scalar. This follows from the bialgebra law.
\begin{equation}\label{eq:classical-phase-copy}
\beginpgfgraphicnamed{classical_phase_copy}
\InputIfFileExists{classical_phase_copy.tikz}{}{\input{./figures/classical_phase_copy.tikz}}
\endpgfgraphicnamed
\end{equation}

\section{The Z/X Calculus and Quantum Computation}\label{sec:zx-calculus}

We have already seen that strongly complementary observables satisfy many graphical identities. This collection of identities is sometimes referred to as the \textit{calculus of complementary observables} to emphasise that it can be used as a computational tool. We will now restrict our attention to the complementary pair $Z$ and $X$, and show how we can use the $Z/X$-calculus to perform calculations on quantum circuits.

\begin{definitions}
  Let $\mathcal Z = (\mathbb C^2, \delta_Z, \epsilon_Z)$ be the $\dagger$-SCFA corresponding to the $Z$-observable and let $\mathcal X = (\mathbb C^2, \delta_X, \epsilon_X)$ be the $\dagger$-SCFA corresponding to the $X$-observable.
  \begin{align*}
    \delta_Z   & :: \ket{0} \mapsto \ket{00}, \ket{1} \mapsto \ket{11} & 
    \epsilon_Z & :: \ket{0} \mapsto 1,        \ket{1} \mapsto 1        \\
    \delta_X   & :: \ket{+} \mapsto \ket{++}, \ket{-} \mapsto \ket{--} & 
    \epsilon_X & :: \ket{+} \mapsto 1,        \ket{-} \mapsto 1        
  \end{align*}
\end{definitions}

For the remainder of this section, we will define:
\ctikzfig{dag_scfa_ZX}

Up a global $e^{i\theta}$ factor, the unitary phases for $\mathcal Z$ are the phase gates $Z_\theta$ and the unitary phases for $\mathcal X$ are the phase gates $X_\theta$. We represent these as dots with a phase angle.
\ctikzfig{ZX_phases}

More generally, we can write arbitrary spiders with phases.
\ctikzfig{phase_spider}

Since the phase commutes with all of the Frobenius structure, it does not matter which leg of the spider we place the phase gate on. For this section, we will ignore (non-zero) scalars, as they will not be important for the calculations. Up to scalar factors, the following equations hold.
\ctikzfig{phase_classical_points}

Since $\mathcal Z$ and $\mathcal X$ are strongly complementary, the phase gates corresponding to classical points are Frobenius algebra automorphisms. Using equation (\ref{eq:classical-phase-copy}) and the fact that $X_\alpha Z_\pi \propto Z_\pi X_{-\alpha}$ for all $\alpha$, we have:
\ctikzfig{zx_phase_commute}

There is only one Abelian group of order $2$, so by Theorem \ref{thm:sc-classification}, $\mathcal X$ is the group algebra $\mathbb C[\mathbb Z_2]$ defined over the basis given by $\mathcal Z$. Both elements of $\mathbb Z_2$ are self-inverse, so the antipode of the strongly complementary pair is trivial.
\ctikzfig{zx_trivial_antipode}

As a consequence, we can freely change the direction of any edge between dots of different colours, and we can delete any two parallel edges between dots of different colours. Finally, we introduce the Hadamard gate, which exchanges the colours of dots.
\begin{center}
  $H := \frac{1}{\sqrt{2}} \left(\begin{matrix} 1 & 1 \\ 1 & -1 \end{matrix}\right)$
  \qquad\qquad
\beginpgfgraphicnamed{zx_hadimard}
\InputIfFileExists{zx_hadimard.tikz}{}{\input{./figures/zx_hadimard.tikz}}
\endpgfgraphicnamed
\end{center}

We refer to the bialgebra identities along with these additional rules as the \textit{Z/X calculus}.

\subsection{Example: Building and Rewriting Circuits}\label{sec:zx-building-circuits}

Consider the following map from $\mathbb C^2 \rightarrow \mathbb C^2$.
\ctikzfig{zx_cnot}

By evaluating the first qubit at $\ket{0}$ and $\ket{1}$, we can see that this map selectively applies $X_\pi$.
\ctikzfig{zx_cnot_on_points}

Therefore, it is a CNOT gate. From this, we have a universality theorem.

\begin{theorem}
  The generators of the $Z/X$ calculus are universal for quantum computation.
\end{theorem}

\begin{proof}
  We have already constructed a CNOT gate, so it suffices to show we can construct an arbitrary 1-qubit unitary. This is possible because every unitary map $\mathbb C^2 \rightarrow \mathbb C^2$ admits an Euler decomposition $U = Z_\gamma X_\beta Z_\alpha$. That is:
  \begin{equation}\label{eq:euler-decomp}
\beginpgfgraphicnamed{euler_decomp}
\InputIfFileExists{euler_decomp.tikz}{}{\input{./figures/euler_decomp.tikz}}
\endpgfgraphicnamed
  \end{equation}
\end{proof}

\begin{example}
  It is a basic property of CNOT gates that three alternating applications yields a qubit swap:
  \ctikzfig{cnot_swap}
  
  We prove this using the $Z/X$ calculus.
  \ctikzfig{cnot_swap_pf}
\end{example}

More examples like this can be found in \cite{CoeckeDuncan2009}. Hillebrand applied to $Z/X$ calculus to a wide variety of security protocols in~\cite{Hillebrand2011}, and Duncan and Perdrix applied it to MBQC in~\cite{DuncanPerdrix2010}.

\chapter{Monoidal Algebra in Quantum Entanglement Theory}\label{ch:monoidal-entanglement}

We now turn out attention to a different topic in quantum information theory: multipartite entanglement. In this chapter, we review several major results from the study of multipartite quantum entanglement. We then give an algebraic (i.e. diagrammatic) characterisation of a special class of highly entangled, symmetric states called Frobenius states. Frobenius states always induce commutative Frobenius algebras, and it can be shown that the two canonical maximally-entangled states on qubits, GHZ and W, can be distinguished by a simple property of this induced Frobenius algebra: specialness or anti-specialness.

In studying GHZ and W states abstractly, we introduce the notion of a GW-pair. A GW-pair consists of a special commutative Frobenius algebra and an anti-special commutative Frobenius algebra, and it exhibits an interaction theory characteristic of the algebras induced by the GHZ and W states. We provide a behavioural intuition for the generators of a GW-pair, show that they are universal for quantum computing, and use them to encode arithmetic on the complex projective line.

\section{Classifying Entanglement}\label{sec:classifying-entanglement}

Characterising general $N$-system entangled states is a very hard open problem in quantum information theory. Before talking about applications of categorical diagrams to the study of entangled states, we will briefly give some background and major results from the field.

\textit{Bipartite states}, i.e. quantum states consisting of two entangled systems, are fairly well understood. We can characterise bipartite states by ``how much'' entanglement they have: ranging from product states (which have no entanglement), to perfectly correlated states (which have maximal entanglement).
\bigskip
\begin{equation}\label{eq:entanglement-chart}
\beginpgfgraphicnamed{entanglement_chart}
\InputIfFileExists{entanglement_chart.tikz}{}{\input{./figures/entanglement_chart.tikz}}
\endpgfgraphicnamed
\end{equation}
\bigskip

This characterisation can be formalised using the \textit{majorisation order} on bipartite states. This is done using the \textit{Schmidt decomposition}. For any bipartite state $\ket\Psi \in \mathcal H \otimes \mathcal H$, there exist orthonormal bases $\{ \ket{u_i} \}$ and $\{ \ket{v_i} \}$ and non-negative real numbers $\alpha_i$ such that:
\[ \ket{\Psi} = \sum_i \alpha_i \ket{u_i} \otimes \ket{v_i} \]

The numbers $\alpha_i$ are called the \textit{Schmidt coefficients} of $\ket\Psi$ and are uniquely determined, up to permutation, by $\ket\Psi$. The number of non-zero Schmidt coefficients is called the \textit{Schmidt rank}. By reordering the associated basis vectors, we can always assume these coefficients are in decreasing order $\alpha_0 \geq \alpha_1 \geq \ldots \geq \alpha_{D-1}$. We can define the majorisation ordering on states $\ket{\Psi}, \ket{\Phi} \in \mathcal H \otimes \mathcal H$ using their associated with Schmidt coefficients $\{ \alpha_i \}$ and $\{ \beta_i \}$.
\[ \ket{\Psi} \leq_M \ket{\Phi} \ \ \Leftrightarrow\ \ 
   \forall k\ .\ \left(\, \sum_{i = 0}^k \alpha_i^2 \geq \sum_{i = 0}^k \beta_i^2 \right) \]

Intuitively, states whose Schmidt coefficients are more evenly spread are higher in the majorisation order. For instance, the product state $\ket{00}$ has Schmidt coefficients $(1,0)$, whereas the Bell state $\ketBell = \frac{1}{\sqrt{2}} \left(\ket{00} + \ket{11} \right)$ has Schmidt coefficients $(\frac{1}{\sqrt{2}}, \frac{1}{\sqrt{2}})$. Since $\frac{1}{2} \leq 1$ and $\frac{1}{2} + \frac{1}{2} \leq 1 + 0$, $\ket{00} \leq_M \ket{\textrm{Bell}}$.

This relation is transitive, reflexive, and anti-symmetric up to a change of orthonormal bases $\ket{u_i}$ and $\ket{v_i}$. We call two states that are equal up to a change of orthonormal basis on each subsystem \textit{LU-equivalent}, or equivalent up to \textit{local unitaries}

\begin{definition}
  Two states $\Psi, \Phi \in \mathcal H \otimes \ldots \otimes \mathcal H$ are said to be \textit{LU-equivalent} if there exist unitary maps $U_i : \mathcal H \rightarrow \mathcal H$ such that:
  \ctikzfig{lu_equiv}
\end{definition}

So, $\leq_M$ forms a partial order on LU-classes of bipartite states. For qubits, $\leq_M$ is a total order, but it is not total in general. Consider two states in $\mathbb C^3 \otimes \mathbb C^3$:
\[
\ket\Psi = \frac{1}{\sqrt{4}} \left( \sqrt{2} \ket{00} + \ket{11} + \ket{22} \right)
\qquad\qquad
\ket\Phi = \frac{1}{\sqrt{5}} \left( \sqrt{2} \ket{00} + \sqrt{2} \ket{11} + \ket{22} \right)
\]

Then, it is neither the case that $\ket\Psi \leq_M \ket\Phi$ nor that $\ket\Phi \leq_M \ket\Psi$. However the minimal and maximal elements of $\leq_M$ are always unique, up to LU-equivalence.

We can get a better feel for what characterises these states by thinking of bipartite states as processes, or \textit{quantum channels}, over which information can flow. We do this by employing the principal of \textit{channel-state duality}.

Fixing an orthonormal basis $\{ \ket{i} \} \in \mathcal H$ fixes a unitary isomorphism $\chi : \mathcal H^* \rightarrow \mathcal H$ sending $\bra{i}$ to $\ket{i}$. As we have seen in the previous chapter, this is the same as fixing a $\dagger$-special commutative Frobenius algebra. The map $\chi$ is then the following induced isomorphism $\mathcal H^* \cong \mathcal H$:
\ctikzfig{white_dualiser}

Up to normalisation factors, we can consider linear maps $L : \mathcal H \rightarrow \mathcal H$ as bipartite states $\ket{\Psi_L} \in L \otimes L$ and bipartite states $\ket{\Psi} \in \mathcal H \otimes \mathcal H$ as linear maps $L_{\Psi} : \mathcal H \rightarrow \mathcal H$.
\begin{equation}\label{eq:channel-state-duality}
\beginpgfgraphicnamed{map_as_state}
\InputIfFileExists{map_as_state.tikz}{}{\input{./figures/map_as_state.tikz}}
\endpgfgraphicnamed
\end{equation}

This is the ``pure'' version of the \textit{Choi-Jamio\l{}kowski isomorphism}. The general statement, which includes mixed states, is as follows. For a finite-dimensional Hilbert space $\mathcal H$, let $\mathcal L(\mathcal H)$ be the vector space of linear maps $\mathcal H \rightarrow \mathcal H$. Positive operators correspond to mixed states and completely positive maps correspond to (mixed) quantum operations.

\begin{theorem}\label{thm:choi-jam}
  Positive operators $L \in \mathcal L(\mathcal H \otimes \mathcal H)$ are in one-to-one correspondence with completely positive maps $\Phi : \mathcal L(\mathcal H) \rightarrow \mathcal L(\mathcal H)$.
\end{theorem}

Since pure states (up to a global phase) are the same thing as positive operators of the form $\rho := \ketbra{\psi}{\psi}$ and pure maps $L : \mathcal H \rightarrow \mathcal H$ are the same thing as CPMs of the form $\Phi(\rho) = L \rho L^\dagger$, equation (\ref{eq:channel-state-duality}) can be thought of as the pure fragment of channel-state duality. Under this correspondence, the Schmidt decomposition of a state is essentially the same as the singular value decomposition of its associated map. Let $\bra{v'_i} = \chi^\dagger(\ket{v_i})$, i.e. the transposition of $\ket{v_i}$ in the basis $\{ \ket{i} \}$. This clearly forms an orthonormal basis for $\mathcal H^*$. Using this basis, we can decompose $L_\Psi$ such that its singular values are the same as the Schmidt coefficients of $\ket\Psi$.
\[ \ket\Psi = \sum_i \alpha_i \ket{u_i} \otimes \ket{v_i}
   \qquad \leftrightarrow \qquad
   L = \sum_i \alpha_i \ketbra{u_i}{v'_i}
\]

In particular, the Schmidt rank of a bipartite state is the same as the rank of the associated map.

Quantum teleportation is an archetypal example of regarding a bipartite state as a channel. The entangled state that Alice and Bob share provides the ``medium'' over which the unknown state is teleported from Alice to Bob. Suppose we considered variants of the teleportation protocol over an arbitrary finite-dimensional Hilbert space $\mathcal H$, replacing the (perfectly correlated) Bell state with other kinds of bipartite states from chart (\ref{eq:entanglement-chart}). The states at the left extreme are the worst for teleportation. Regarded as channels, product states correspond to rank-one maps. No matter what we put in to such a channel, we always get the same output, up to a scalar. Such channels cannot be used to send \textit{any} quantum data. At the other extreme are the perfectly correlated states, which correspond to unitary maps under channel-state duality. Any such state can be used to construct a teleportation protocol which will always succeed. If this unitary $U_\Psi$ is not the identity (as in the case of the Bell state), we simply need to undo it by applying $U_\Psi^\dagger$ to one of the sub-systems.

In between these two extremes are maps $L$ of rank $2 \leq r \leq D$, which can be thought of as noisy channels. Suppose $r = D$, then $\ket{\Psi_L}$ can at least in principal teleport a state, but it might be impossible for Bob to recover Alice's state deterministically. However, in most cases, Bob can at least recover the state with non-zero probability. This is because an arbitrary linear map can be ``applied'' to a quantum system by first applying a big unitary to the system and an ancilla state, then measuring the ancilla.
\ctikzfig{stochastic_map}

If Bob gets outcome $0$, then he has successfully applied $L$, otherwise he has applied some other (unwanted) map $L'$. Maps $L$ that can be applied with non-zero probability are known as \textit{stochastic maps}. In the case where $r < D$, states can only be teleported in a lossy sense, i.e. they are (non-deterministically) projected on to a subspace of $\mathcal H$ before being sent. Under channel-state duality, chart (\ref{eq:entanglement-chart}) becomes:
\bigskip
\ctikzfig{entanglement_map_chart}
\bigskip

When we move from \textit{bipartite} entanglement to \textit{multipartite} entanglement, the picture becomes less clear. For one thing, this is no canonical analogue to the majorisation order involving three or more systems. However, we can define an \textit{operational} ordering on states that is equivalent to the majorisation ordering in the bipartite case. Operational orderings relate states by the existence of certain types of quantum protocols that can convert one state into another. The most well-known type of protocol used for this purpose is a LOCC protocol.

\begin{definition}\label{def:locc}
  An $N$-partite state $\ket\Psi$ can be converted into $\ket\Phi$ by \textit{Local Operations and Classical Communication} (LOCC) if there exists an $N$-party protocol that can \textit{deterministically} convert $\ket\Psi$ into $\ket\Phi$ consisting of any number of the following operations:
  \begin{enumerate}
    \item Party $p_i$ performs a local measurement of the $i$-th subsystem (possibly with an ancillary system which only $p_i$ can access).
    \item Party $p_j$ performs a local unitary, which can be conditioned on any previous measurement outcome in the protocol.
  \end{enumerate}
  
  In such a case, we write $\ket\Phi \leq_{\textrm{LOCC}} \ket\Psi$. If two states are LOCC-interconvertible (i.e. $\ket\Phi \leq_{\textrm{LOCC}} \ket\Psi$ and $\ket\Psi \leq_{\textrm{LOCC}} \ket\Phi$), we say they are \textit{LOCC-equivalent}, written $\ket\Psi \sim_{\textrm{LOCC}} \ket\Phi$.
\end{definition}

Nielsen and Vidal showed that the majorisation order is the same as the LOCC order on bipartite states~\cite{NielsenVidal2001}. As a consequence, two states are LOCC-equivalent if an only if they are LU-equivalent. In fact, this relationship between LU-equivalence and LOCC-equivalence is true for any number of systems.

\begin{theorem}[\cite{Bennett1999}]
  Two $N$-partite quantum states are LOCC-equivalent if and only if they are LU-equivalent.
\end{theorem}

\begin{example}
  The following two bipartite states are LOCC-equivalent:
  \begin{center}
    $\ketBell := \frac{1}{\sqrt{2}} \left( \ket{00} + \ket{11} \right)$
    \qquad
    $\ketEPR := \frac{1}{\sqrt{2}} \left( \ket{01} + \ket{10} \right)$
  \end{center}
\end{example}

Even though there are infinitely many LOCC classes in the bipartite case, the majorisation order gives us a straightforward way to characterise LOCC-equivalence, as well as those states that are minimal and maximal with respect to $\leq_{\textrm{LOCC}}$. However, for three or more systems, the picture becomes much more complicated. For that reason, it is convenient to introduce a course-graining of LOCC called Stochastic LOCC, or SLOCC. For SLOCC, we relax the requirement that the LOCC protocol succeeds deterministically and merely require that it succeed with some non-zero probability.

\begin{definition}\label{def:slocc}
  An $N$-partite state $\ket\Psi$ can be converted into $\ket\Phi$ by \textit{Stochastic Local Operations and Classical Communication} (SLOCC) if there exists an $N$-party LOCC protocol that converts $\ket\Psi$ into $\ket\Phi$ with a \textit{non-zero probability}.
\end{definition}

Again, we use $\leq_{\textrm{SLOCC}}$ and $\sim_{\textrm{SLOCC}}$ to represent SLOCC-convertibility and SLOCC-equivalence, respectively. It was shown by D\"ur, Vidal, and Cirac that SLOCC-equivalence can be characterised in a similar manner to LOCC-equivalence, but replacing \textit{unitary} maps with \textit{invertible} maps.

\begin{definition}\label{def:ilo-equiv}
  Two states $\ket\Psi, \ket\Phi \in \mathcal H \otimes \ldots \otimes \mathcal H$ are said to be ILO-equivalent if there exist invertible local operations $L_i : \mathcal H \rightarrow \mathcal H$ such that:
  \ctikzfig{ilo_equiv}
\end{definition}

\begin{theorem}[\cite{DVC}]\label{thm:slocc-ilo}
  Two states $\ket\Psi, \ket\Phi \in \mathcal H \otimes \ldots \otimes \mathcal H$ are SLOCC-equivalent iff they are ILO-equivalent.
\end{theorem}

We say a state $\ket\Psi$ is \emph{SLOCC-maximal} if its SLOCC-equivalence class is maximal with respect to $\leq_{\textrm{SLOCC}}$, i.e. $\ket\Psi \leq_{\textrm{SLOCC}} \ket\Phi \implies \ket\Phi \sim_{\textrm{SLOCC}} \ket\Psi$. We define SLOCC-minimal states similarly. For any number of systems, the unique minimal SLOCC-equivalence class is the class of product states $\ket{00\ldots{}0}$. For $\mathbb C^2 \otimes \mathbb C^2$, the only other SLOCC class is generated by the perfectly correlated states. So, the Hasse diagram of SLOCC classes has only two elements.
\ctikzfig{hasse_bipartite}

For three qubits, there are still only six SLOCC classes (or four classes, up to permutations of qubits), there is no longer a unique maximal SLOCC-class. The Hasse diagram of the tripartite SLOCC classes is:
\ctikzfig{hasse_tripartite}

Where $\ketGHZ = \frac{1}{\sqrt{2}} \left( \ket{000} + \ket{111} \right)$, $\ketW = \frac{1}{\sqrt{3}} \left( \ket{100} + \ket{010} + \ket{001} \right)$, and the $\ket{D_i}$ states represent the three separable configurations of a Bell state with $\ket{0}$.
\ctikzfig{tripartite_degen}

The complete classification for three qubits was given by D\"ur, Vidal, and Cirac~\cite{DVC}. They distinguished GHZ and W by using an entanglement measure called the \textit{3-tangle}, or \textit{residual tangle}. Intuitively, this is the entanglement ``left over'' after bipartite correlations are subtracted out. States that are SLOCC-equivalent to $\ketGHZ$ have a non-zero 3-tangle. These states have some true 3-body entanglement, even accounting for bipartite correlations. However, states that are SLOCC-equivalent to $\ketW$ always have a vanishing 3-tangle, which means that all the entanglement present can be accounted for by bipartite correlations. Informally, we might depict the correlations in GHZ and W as:
\ctikzfig{ghz_w_informal}

One could say that W has no ``true'' tripartite entanglement. One could also state this property positively, by stating that W-like entanglement, unlike GHZ-like entanglement, is (partially) robust to the loss of a one system. That is, if we trace out one of the three parties, the reduced density matrix of W becomes (ignoring normalisation):
\[ \rho'_{\textit{W}} = (\ket{10} + \ket{01})(\bra{10} + \bra{01}) + \ketbra{00}{00}
                    = \sqrt{2} \ketbra{\textit{EPR}}{\textit{EPR}} + \ketbra{00}{00} \]

So, $\rho'_{\textit{W}}$ corresponds to the probabilistic mixture of the (entangled) Einstein-Podolsky-Rosen state $\ketEPR$ and a separable pure state $\ket{00}$. However, if we do the same to the GHZ state, we get:
\[ \rho'_{\textit{GHZ}} = \ketbra{00}{00} + \ketbra{11}{11} \]

This is a mixture of two product states. Therefore it contains no entanglement if the third system is disregarded.

There are various schemes for classifying multipartite entanglement which have had some success for small (typically $< 6$) numbers of systems. These approaches have had mixed results, as the general problem becomes very difficult for four or more systems. For four or more systems, there is necessarily an infinite number of SLOCC classes~\cite{DVC}. This is because the continuous degrees of freedom in the state space of $N$ systems ($2^N - 2$) quickly overcomes the degrees of freedom available in a tensor product $L_1 \otimes \ldots \otimes L_N$ of local invertible maps ($6 N + 2$). Nevertheless, one sometimes give a finite set of SLOCC \textit{super-classes}, i.e. SLOCC classes with some free parameters, which span the state spaces of these larger systems. However, there is no one right choice for a parametrisation, so these super-classes are not uniquely determined and often reflect the methods that were used to obtain them.

Lamata et al introduced an inductive scheme for classifying multipartite states \cite{Lamata2006} for $N$ partite states in terms of the classification of $N-1$ partite states. It works by treating $N$ partite states as maps from $\mathcal H$ to $\mathcal H^{\otimes(N-1)}$.
\ctikzfig{n_party_state_map}

They then look at the image $V \subseteq \mathcal H^{\otimes(N-1)}$ of $L_\Psi$ and classify the vectors that span $V$ (i.e. $N-1$ partite states). Using their classification scheme, $\ketGHZ$ is the unique, tripartite qubit state whose image is spanned by two product states $\{ \ket{\psi_1} \otimes \ket{\phi_1}, \ket{\psi_2} \otimes \ket{\phi_2} \}$ where $\ket{\psi_1} \neq \lambda \ket{\psi_2})$ and $\ket{\phi_1} \neq \lambda' \ket{\phi_2}$. $\ketW$ is the unique state whose image is two dimensional, but only contains one product state. The whole classification for tripartite qubit states is as follows.

\begin{center}
  \bigskip
  \begin{tabular}{c|l}
    State     & $V$ spanned by                                                               \\
    \hline
    \ketGHZ   & $\{ \ket{\psi_1} \otimes \ket{\phi_1}, \ket{\psi_2} \otimes \ket{\phi_2} \}$ \\
    \ketW     & $\{ \ket{\psi} \otimes \ket{\phi}, \ket{\Psi} \}$                            \\
    \ket{D_1} & $\{ \ket{\psi} \otimes \ket{\phi_1}, \ket{\psi} \otimes \ket{\phi_2} \}$     \\
    \ket{D_2} & $\{ \ket{\psi_1} \otimes \ket{\phi}, \ket{\psi_2} \otimes \ket{\phi} \}$     \\
    \ket{D_3} & $\{ \ket{\Psi} \}$                                                           \\
    \ket{000} & $\{ \ket{\psi} \otimes \ket{\phi} \}$                                        \\
  \end{tabular}
  \bigskip
\end{center}

Since it can be proved that any two-dimensional subspace of $\mathbb C^2 \otimes \mathbb C^2$ contains at least one product state, this list is exhaustive, and matches the classification given by \cite{DVC}. It also reflects the fact that tracing out a qubit from $\ketGHZ$ yields a mixture of product states, while tracing out a qubit of $\ketW$ yields a mixture of a product state and an entangled state.

In \cite{Lamata2006four}, the authors went on to use this technique to enumerate a finite set of $4$-partite SLOCC super-classes spanning $(\mathbb C^2)^{\otimes 4}$. However, the amount of calculations to ensure these super-classes were non-overlapping and spanned the whole space was much greater than in the tripartite case, and it seems likely that enumeration of super-classes with more than four subsystems would be significantly more difficult than the $4$-partite case.

\subsection{Symmetric States}

The study of symmetric multipartite states is quite a bit simpler, because they have significantly fewer degrees of freedom. For instance, the space of $N$-partite symmetric states over qubits is spanned by the $N$-partite Dicke states:
\[ \ket{D^{(k)}_N} = \mathcal S_N(|\underbrace{0...0}_{N-k} \underbrace{1...1}_{k} \rangle) \]

Where $\mathcal S_N$ is the $N$-system symmetrisation map.
\[ \mathcal S_N :: \ket{i_1,i_2,\ldots,i_N} \mapsto
   \frac{1}{N!} \sum_{\pi \in \textrm{Perms}(N)}
     \ket{i_{\pi(1)},i_{\pi(2)},\ldots,i_{\pi(N)}} \]

So, for qubits, the space of $N$-system symmetric states is $(N+1)$-dimensional. It can also be shown that an arbitrary $N$-partite state is the symmetrisation of a product state.
\[ \ket\Psi = \mathcal S_N(\ket{\psi_1,\psi_2,\ldots,\psi_N}) \]

The factors $\ket{\psi_i}$ need not be distinct. Since the value of $\ket\Psi$ does not depend on the order of factors, $\ket\Psi$ is totally defined by the multiset of $M \leq N$ distinct factors $\ket{\phi_i}$:
\[ \ket\Psi = \mathcal S_N(
   \underbrace{\ket{\phi_1,\ldots,\phi_1}}_{d_1} \otimes
   \underbrace{\ket{\phi_2,\ldots,\phi_2}}_{d_2} \otimes
   \ldots \otimes
   \underbrace{\ket{\phi_M,\ldots,\phi_M}}_{d_M}) \]

We can also assume that the $d_i$ are in decreasing order. The number $M$ of distinct factors is called the \textit{diversity degree} of a symmetric state, and the list $\mathcal D_\Psi = [d_1,d_2,\ldots,d_M]$ is called its \textit{degeneracy configuration}.

\begin{definition}\label{def:symmetric-slocc-equiv}
  Two $N$-partite states $\ket\Psi, \ket\Phi \in \mathcal H^{\otimes N}$ are called \textit{symmetrically SLOCC-equivalent} if there exists a \textit{single} invertible map $L : \mathcal H \rightarrow \mathcal H$ such that:
  \ctikzfig{symmetric_slocc_equiv}
\end{definition}

Bastin et al provided a classification result for symmetric qubit states of any size with a small diversity degree~\cite{Bastin2009}. It relies crucially on the following fact about qubit states.

\begin{proposition}\label{prop:qubit-symmetric-slocc-equiv}
  Two $N$-partite symmetric states are SLOCC-equivalent if and only if they are symmetrically SLOCC-equivalent.
\end{proposition}

Using this fact, they showed that for small diversity degree, SLOCC-classes are uniquely fixed by their degeneracy configuration.

\begin{theorem}[\cite{Bastin2009}]
  Two  $N$-partite symmetric qubit states with diversity degree $M \leq 3$ are SLOCC-equivalent if and if only they have the same degeneracy configuration.
\end{theorem}

They proved that any state with diversity degree $M = N$ is SLOCC-equivalent to $\ket{\textit{GHZ}_N} = \ket{0\ldots{}0} + \ket{1\ldots{}1}$. As a consequence, the SLOCC classes of $N$-partite symmetric states are completely determined for $N \leq 4$. They also showed that for $4 \leq M < N$, there are necessarily infinitely-many SLOCC classes of symmetric states.

\section{Strong SLOCC-maximality and strong symmetry}\label{sec:strong-symmetry-strong-max}

In an effort to push state classification farther, we will identify the properties of that make GHZ and W states unique. In order to leverage techniques from categorical quantum mechanics and diagrammatic languages, we shall focus on properties that are \textit{compositional} in nature. That is, we shall look at how GHZ and W states, treated as maps via channel-state duality, behave when they are composed with other states or themselves.

Whereas bipartite states can be thought of as quantum channels, tripartite states can be thought of as \textit{algebraic operations} which have two channels of input and one channel of output (or, equivalently, as coalgebraic operations from 1 input to 2 outputs).
\ctikzfig{tripartite_map_state}

The GHZ and W states are characterised as the unique tripartite qubit states that are SLOCC-maximal, so we shall look at how to characterise SLOCC-maximality abstractly. For bipartite states, this condition is equivalent to the a bipartite state forming the ``cap'' of a self-dual compact structure.

\begin{proposition}\label{pro:bipartite-maximal}
  A bipartite state $\ket\Psi \in \mathcal H \otimes \mathcal H$ is SLOCC-maximal iff there exists an effect $\bra\Phi : \mathcal H \otimes \mathcal H \rightarrow \mathbb C$ such that
  \begin{equation}\label{eq:bipartite-maximal}
\beginpgfgraphicnamed{maximal_bipartite}
\InputIfFileExists{maximal_bipartite.tikz}{}{\input{./figures/maximal_bipartite.tikz}}
\endpgfgraphicnamed
  \end{equation}
\end{proposition}

\begin{proof}
  We first show that SLOCC-maximal maps have full Schmidt rank. Suppose $\ket\Psi$ has Schmidt decomposition:
  \[ \ket{\Psi} = \sum_i \alpha_i \ket{u_i} \otimes \ket{v_i} \]
  If for any $i$, $\alpha_i = 0$, then there exists a map $S$ with $\ket{u_i}$ in its null space such that $\ket\Psi = (S \otimes 1)\ket{\Psi'}$. Since $S$ is singular, it corresponds to a non-reversible stochastic local operation, so $\ket\Psi$ is not SLOCC-maximal.
  
  Since any SLOCC-maximal state must have full Schmidt rank, we can construct $\bra\Phi$ as follows.
  \[ \bra{\Phi} = \sum_i \frac{1}{\alpha_i} \bra{v_i} \otimes \bra{u_i} \]
  
  Then equation (\ref{eq:bipartite-maximal}) is satisfied. Conversely, if $\ket\Psi$ were not SLOCC maximal, then $\ket\Psi = (S \otimes 1)\ket{\Psi'}$ for some $\ket{\Psi'}$ and some singular map $S$. Clearly no bipartite state of this form could satisfy equation (\ref{eq:bipartite-maximal}).
\end{proof}

Characterising SLOCC-maximal tripartite states is trickier. However, GHZ and W satisfy a stronger version of SLOCC-maximality. As in \cite{Lamata2006}, we can study tripartite entangled states by studying the states that span the image of the associated map $\delta_\Psi$.
\ctikzfig{image_of_state}

For convenience, we define the 3 bipartite \textit{image spaces} of $\ket\Psi$ as follows.
\begin{align*}
  \textrm{Im}_1(\ket{\Psi}) & = \textrm{Span}\left\{ \, (\bra{i} \otimes 1 \otimes 1) \ket{\Psi} \, \right\} \\
  \textrm{Im}_2(\ket{\Psi}) & = \textrm{Span}\left\{ \, (1 \otimes \bra{i} \otimes 1) \ket{\Psi} \, \right\} \\
  \textrm{Im}_3(\ket{\Psi}) & = \textrm{Span}\left\{ \, (1 \otimes 1 \otimes \bra{i}) \ket{\Psi} \, \right\}
\end{align*}

For states in $(\mathbb C^D)^{\otimes 3}$ to be SLOCC-maximal, these spaces must all be $D$-dimensional. In fact, this condition is equivalent to SLOCC-maximality.

\begin{proposition}
  A tripartite state $\ket\Psi$ is SLOCC-maximal if and only if $\textrm{Im}_i(\ket\Psi)$ is $D$-dimensional for $i=1,2,3$.
\end{proposition}

\begin{proof}
  Suppose that $\textrm{Im}_1(\ket\Psi)$ is not $D$-dimensional. Then, there must exist $\bra{\psi} \in \mathbb C^D$ such that $(\bra\psi \otimes 1 \otimes 1)\ket\Psi = 0$. Fix an orthonormal basis $\bra{u_i}$ spanning $\bra{\psi}^\perp$. Then define a singular map $S = \sum_i \ketbra{u_i}{u_i}$, and note that there exists $\ket{\Psi'}$ such that $(S \otimes 1 \otimes 1)\ket{\Psi'} = \ket\Psi$, so $\ket\Psi$ cannot be SLOCC-maximal. The argument follows similarly for $\textrm{Im}_2(\ket\Psi)$ and $\textrm{Im}_3(\ket\Psi)$.
  
  Conversely, let $\ket\Psi$ not be SLOCC-maximal. Then it must be in one of the following forms:
  \[ (S \otimes 1 \otimes 1)\ket{\Psi'} \qquad
     (1 \otimes S \otimes 1)\ket{\Psi'} \qquad
     (1 \otimes 1 \otimes S)\ket{\Psi'} \]
  In each of these cases, one of the spaces $\textrm{Im}_i(\ket\Psi)$ must be less than $D$-dimensional.
\end{proof}

It is always the case for SLOCC-maximal tripartite states that each of the spaces $\textrm{Im}_i(\ket\Psi)$ contain entangled states. However, for $D > 2$, it need not contain a SLOCC-maximal bipartite state. The case of $D = 2$ is degenerate in the sense that \textit{any} entangled state is SLOCC-maximal.

\begin{definition}\label{def:strong-slocc-max}
  A tripartite state $\ket\Psi$ is called \textit{strongly SLOCC-maximal} if each of its associated image spaces $\textrm{Im}_i(\ket\Psi)$ contain a SLOCC-maximal bipartite state.
\end{definition}

\begin{theorem}\label{thm:string-slocc-slocc}
  Strong SLOCC-maximality implies SLOCC-maximality.
\end{theorem}

\begin{proof}
  It suffices to show that the image spaces $\textrm{Im}_i(\ket\Psi)$ are all $D$-dimensional. Let $\bra\psi$ be a state such that $\ket\Phi = (1 \otimes \bra\psi \otimes 1)\ket\Psi$ is a SLOCC-maximal bipartite state. Then, fixing an orthonormal basis $\bra{i}$ for $(\mathbb C^D)^*$, all states of the form 
  \ctikzfig{strong_slocc_pf}
  \noindent must be linearly independent for distinct values of $i$, so they span $\mathbb C^D$. Thus, non-projected states:
  \ctikzfig{strong_slocc_pf2}
  \noindent must span a $D$-dimensional subspace of $\mathbb C^D \otimes \mathbb C^D$. The result for $\textrm{Im}_2(\ket\Psi)$ and $\textrm{Im}_3(\ket\Psi)$ follows similarly.
\end{proof}

\begin{remark}
  In the case of $D=2$, any entangled bipartite state is SLOCC-maximal, so strong SLOCC-maximality and ``weak'' SLOCC-maximality are equivalent. However, for $D > 2$, the implication in Theorem \ref{thm:string-slocc-slocc} is strict. To see this, consider the following state in $\mathbb C^3 \otimes \mathbb C^3 \otimes \mathbb C^3$.
  \[ \ket\Psi = \ket{000} + \ket{101} + \ket{110} + \ket{202} + \ket{220} \]
  
  The spaces $\textrm{Im}_1(\ket\Psi)$, $\textrm{Im}_2(\ket\Psi)$, and $\textrm{Im}_3(\ket\Psi)$ are all $3$-dimensional, so $\ket\Psi$ is SLOCC-maximal. Every state in $\textrm{Im}_i(\ket\Psi)$ is of the form:
  \[ \ket\Phi = a\ket{00} + b(\ket{01} + \ket{10}) + c(\ket{02} + \ket{20}) \]
  
  But $\ket\Psi$ cannot be SLOCC-maximal. If $b=c=0$, then it is a product state, otherwise there exists a non-zero bra $\bra\xi := c\bra{1} - b\bra{2}$ such that $(\bra\xi \otimes 1)\ket\Phi = 0$.
\end{remark}

When a state is symmetric we can simplify the strong maximality condition. A state $\ket\Psi$ is strongly SLOCC-maximal if there exist effects $\bra\xi$, $\bra\Phi$ such that:
\begin{equation}\label{eq:symmetric-strong-max}
\beginpgfgraphicnamed{maximal_symmetric}
\InputIfFileExists{maximal_symmetric.tikz}{}{\input{./figures/maximal_symmetric.tikz}}
\endpgfgraphicnamed
\end{equation}

The GHZ and W states are strongly SLOCC-maximal because they are SLOCC-maximal qubit states. In particular, ignoring scalar factors, the following equations hold:
\begin{equation}\label{eq:ghz-w-strong-max}
\beginpgfgraphicnamed{ghz_w_bipartite_states}
\InputIfFileExists{ghz_w_bipartite_states.tikz}{}{\input{./figures/ghz_w_bipartite_states.tikz}}
\endpgfgraphicnamed
\end{equation}

The GHZ and W states have natural $N$-qubit symmetric analogues.
\begin{align*}
  \ket{\textit{GHZ}_N} & := \ket{00\ldots0} + \ket{11\ldots1}                            \\
  \ket{\textit{W}_N}   & := \ket{10\ldots0} + \ket{010\ldots0} + \ldots + \ket{0\ldots01}
\end{align*}

Not only do they have such $N$-partite versions, they come with a recipe for inductively constructing them. That is, for both of these states, there is a bipartite effect $\bra\Phi$ that can be used to ``glue'' a tripartite state on to an $N$-partite state by projecting out a pair of qubits to make an $(N+1)$-partite symmetric state.
\begin{align*}
  \ket{\textit{GHZ}_{N+1}} & = (1 \otimes \braBell \otimes 1)(\ket{\textit{GHZ}_N} \otimes \ketGHZ) \\
  \ket{\textit{W}_{N+1}}   & = (1 \otimes \braEPR \otimes 1)(\ket{\textit{W}_N} \otimes \ketW)     
\end{align*}

To inductively build a symmetric $N$-partite state, it suffices that the following condition hold.

\begin{definition}\label{def:strongly-symmetric}
  A symmetric state is said to be \emph{strongly symmetric} if there exists some bipartite effect $\bra\Phi$ such that
  \begin{equation}\label{eq:strongly-symmetric}
\beginpgfgraphicnamed{strongly_symmetric}
\InputIfFileExists{strongly_symmetric.tikz}{}{\input{./figures/strongly_symmetric.tikz}}
\endpgfgraphicnamed
  \end{equation}
\end{definition}

\section{Frobenius States and their Induced Frobenius Algebras}\label{sec:frobenius-states}

A state that is strongly SLOCC-maximal and strongly symmetric is called a Frobenius state.

\begin{definition}\label{def:Frobenius-state}\em 
  A symmetric tripartite state $\ket\Psi \in \mathcal H \otimes \mathcal H \otimes \mathcal H$ is said to be a \emph{Frobenius state} if there exist effects $\bra\Phi$, $\bra\xi$ such that:
  \begin{equation}\label{eq:frobenius-state-axioms}
\beginpgfgraphicnamed{maximal_symmetric}
\InputIfFileExists{maximal_symmetric.tikz}{}{\input{./figures/maximal_symmetric.tikz}}
\endpgfgraphicnamed \qquad\qquad %
\beginpgfgraphicnamed{strongly_symmetric}
\InputIfFileExists{strongly_symmetric.tikz}{}{\input{./figures/strongly_symmetric.tikz}}
\endpgfgraphicnamed
  \end{equation}
\end{definition}

Note that $\bra\Phi$ satisfying two equations in the definition must be the \emph{same effect}. This is a stronger condition than stating that a Frobenius state is both strongly SLOCC-maximal and strongly symmetric (i.e. these equations respectively hold for some possibly distinct effects $\bra\Phi$ and $\bra{\Phi'}$).

\begin{theorem}[Algebras as states]\label{thm:frob-state-from-cfa}
  For any commutative Frobenius algebra $(\mult, \unit, \comult, \counit)$, the following is a Frobenius state, with its two associated effects:
  \ctikzfig{frob_state_from_cfa}
\end{theorem}

\begin{proof}
  The induced triparite state is symmetric because $\comult$ is co-commutative. The two Frobenius state equations follow from the spider theorem of commutative Frobenius algebras.
\end{proof}

The Frobenius state conditions from Definition \ref{def:Frobenius-state} hold as a consequence of Theorem \ref{thm:cfa-nf}. Also, from any Frobenius state we can construct the associated commutative Frobenius algebra.

\begin{theorem}[States as algebras]\label{thm:frob-state-cfa}
  For any Frobenius state $\ket\Psi$, there exist effects $\bra\Phi$, $\bra\xi$ such that the following is a commutative Frobenius algebra:
  \ctikzfig{cfa_from_frob_state} 
\end{theorem}

\begin{proof} 
  Cocommutativity of $\comult$ follows from $\ket\Psi$ being symmetric. Coassociativity follows from the symmetry of $\ket\Psi$ and the strong symmetry equation.
  \ctikzfig{cfa_from_frob_state_pf1}
  
  Associativity and the Frobenius law follow similarly. Strong SLOCC-maximality implies that $\bra\xi$ is a counit for $\comult$ and that $(\bra\xi \otimes \bra\xi \otimes 1)\ket\Psi$ is a unit for $\mult$.
\end{proof}

We can therefore define a commutative Frobenius algebra using either the usual maps $(\mu, \eta, \delta, \epsilon)$, or a triple $(\ket\Psi, \bra\Phi, \bra\xi)$ consisting of a Frobenius state and its two associated effects. Also, note that for a given state $\ket\Psi$, there could be multiple induced commutative Frobenius algebras based upon the choice of $\bra\xi$. However, once $\bra\xi$ is fixed, $\bra\Phi$ is completely determined. This is analogous to the situation with Frobenius algebras where the maps $\mu$ and $\epsilon$ completely determine the other two.

\begin{theorem}\label{thm:fs-symmetric-slocc-iso}
  If $(\ket\Psi, \bra\Phi, \bra\xi)$ is a Frobenius state and $\ket{\Psi'}$ is symmetrically SLOCC-equivalent to $\ket\Psi$, then $\ket{\Psi'}$ extends to a Frobenius state in such a way that the induced commutative Frobenius algebras of $\ket\Psi$ and $\ket{\Psi'}$ are isomorphic.
\end{theorem}

\begin{proof}
  Let $\ket{\Psi'} = (L \otimes L \otimes L)\ket\Psi$ for some invertible map $L$. Then let $\bra{\xi'} = \bra\xi L^{-1}$ and $\bra{\Phi'} = \bra{\Phi}(L^{-1} \otimes L^{-1})$. It is straightforward to verify the Frobenius state axioms, and constructing the CFA as in Theorem \ref{thm:frob-state-cfa}, we have:
  \ctikzfig{iso_cfa}
\end{proof}

\subsection{Classification of Qubit Frobenius States}\label{sec:classification}

In section \ref{sec:special-anti}, we defined special and anti-special commutative Frobenius algebras. These are both defined as CFAs with one additional axiom.
\ctikzfig{scfa_acfa_recap}

We show in this section that GHZ and W, the two canonical SLOCC-maximal tripartite qubit states, are both Frobenius states. Furthermore, the specialness and anti-specialness condition serve to distinguish the (symmetric) SLOCC-classes of these two states.

It will be convenient to work with the unnormalised versions of $\ketGHZ$ and $\ketW$.
\[ \ketGHZ = \ket{000} + \ket{111} \qquad\qquad \ketW = \ket{100} + \ket{010} + \ket{001} \]

For the Frobenius state $\ketGHZ$, fixing $\bra\xi := \bra{0} + \bra{1}$ induces the following CFA, which we shall refer to as $\mathcal G$.
\begin{equation}\label{GHZ-SCFA}
  \begin{split}
    \whitemult & = \ket{0}\bra{00} + \ket{1}\bra{11} \qquad\qquad
    \whiteunit = \sqrt{2}\, \ket{+} = \ket{0}+\ket{1} \\
    \whitecomult & = \ket{00}\bra{0} + \ket{11}\bra{1} \qquad\qquad
    \whitecounit = \sqrt{2} \bra{+} = \bra{0}+\bra{1}
  \end{split}\vspace{-1.5mm}
\end{equation}

We can verify that $\mathcal G$ is \textit{special}.
\[ \whitemult \circ \whitecomult = \ketbra{0}{0} + \ketbra{1}{1} = 1_{\mathbb C^2} \]

For the Frobenius state $\ketW$, fixing $\bra\xi := \bra{0}$ induces the following CFA, called $\mathcal W$.
\begin{equation}\label{W-ACFA}
  \begin{split}
    \mult & = \ket{1}\bra{11} + \ket{0}\bra{01} + \ket{0}\bra{10} 
    \qquad\qquad\qquad
    \unit = \ket 1\qquad\qquad \\
    \comult & = \ket{00}\bra{0} + \ket{01}\bra{1} + \ket{10}\bra{1} 
    \qquad\qquad\qquad
    \counit = \bra 0\qquad\qquad
  \end{split}
\end{equation}    

Computing the partial traces of multiplication and comultiplication, we get:
\begin{align*}
  \lolli   & = \Tr^{\mathbb C^2}(\comult)
             = \sum_i (1 \otimes \bra{i})(\ket{00}\bra{0} + \ket{01}\bra{1} + \ket{10}\bra{1})\ket{i}
             = 2 \ket{0} \\
  \cololli & = \Tr^{\mathbb C^2}(\mult)
             = \sum_i \bra{i}(\ket{1}\bra{11} + \ket{0}\bra{01} + \ket{0}\bra{10})(1 \otimes \ket{i})
             = 2 \bra{1}
\end{align*}

We can then verify that $\mathcal W$ is \textit{anti-special}.
\[ \circl \otimes (\,\mult \circ \comult\,) = 2(\ket{0}\bra{1} + \ket{0}\bra{1}) = (2\ket{0})(2\bra{1}) = \lolli \circ \cololli \]

The next theorem is a straightforward consequence of the classification of tripartite qubit states and Theorem \ref{thm:fs-symmetric-slocc-iso}.

\begin{theorem}\label{thm:qubit-cfa-classification}
  Any SLOCC-maximal, symmetric state $\ket\Psi \in \mathbb C^2 \otimes \mathbb C^2 \otimes \mathbb C^2$ is a Frobenius state. It is SLOCC-equivalent to $\ketGHZ$ (resp. $\ketW$) if and only if the associated commutative Frobenius algebra is special (resp. anti-special).
\end{theorem}

\begin{proof}
  We first show $\ket\Psi$ is a Frobenius state. Since it is SLOCC-maximal, it must be SLOCC-equivalent to GHZ or W. Since it is a symmetric qubit state, it must be \textit{symmetrically} SLOCC-equivalent to GHZ or W by Proposition \ref{prop:qubit-symmetric-slocc-equiv}. Therefore it is a Frobenius state.
  
  If $\ket\Psi$ is SLOCC-equivalent to GHZ, its induced CFA is isomorphic to $\mathcal G$ by Theorem \ref{thm:fs-symmetric-slocc-iso}, so it must be special. If $\ket\Psi$ is SLOCC-equivalent to GHZ, its induced CFA is isomorphic to $\mathcal W$, so it must be anti-special.
  
  Conversely, let $\ket\Psi$ be the induced Frobenius state of a special commutative Frobenius algebra $\mathcal S$. It must be SLOCC-equivalent to GHZ or W. Suppose it were to be SLOCC-equivalent to W, then $\mathcal S$ is isomorphic to an anti-special commutative Frobenius algebra, which is a contradiction for dimensions $>1$. So, $\ket\Psi$ must be SLOCC-equivalent to GHZ. The result follows similarly for ACFAs.
\end{proof}

\subsection{Classification of Frobenius States for Higher Dimensions}\label{sec:higher-dim-classification}

We can push this classification a bit farther, but less is known in higher dimensions.

\begin{example}\label{ex:d-dimensional-frob-states}
  We can produce $D$-dimensional analogues of the GHZ state and the W state.
  \begin{align*}
    \ketGHZD & = \sum_{i} \ket{iii} \\
    \ketWD   & = \ket{100} + \ket{010} + \ket{001} +
                 \sum_{i} \left(\ket{0ii} + \ket{i0i} + \ket{ii0} \right)
  \end{align*}
\end{example}

\begin{theorem}
  $\ketGHZD$ and $\ketWD$ are both Frobenius states for any dimension. The induced CFA for $\ketGHZD$ is always special and the induced CFA for $\ketWD$ is always anti-special.
\end{theorem}

\begin{proof}
  For $\ketGHZD$, the associated effect in $(\mathbb C^D)^*$ is $\sum_i \bra{i}$ and for $\ketWD$, it is $\bra{0}$. The rest of the structure is uniquely determined, and the verification of SCFA and ACFA axioms is straightforward.
\end{proof}

\begin{theorem}
  For any $D$-dimensional SCFA, the induced Frobenius state is SLOCC-equivalent to $\ketGHZD$.
\end{theorem}

\begin{proof}
  Let $(\mathcal H, \mu, \eta, \delta, \epsilon)$ be an SCFA. Then so too is $(\mathcal H^*, \delta^*, \epsilon^*, \mu^*, \eta^*)$. Any special Frobenius algebra over $\mathbb C$ is semisimple, and any commutative semisimple algebra over $\mathbb C$ is isomorphic to the direct sum of $D$ copies of $(\mathbb C, \cdot)$. So, for a (not necessarily orthonormal) basis $\{ \ket{u_i} \}$, $\mu$ has the form:
  \[ \delta :: \ket{u_i} \mapsto \ket{u_i,u_i} \]
  For some arbitrary $\eta = \sum_i \alpha_i \ket{u_i}$, the tripartite state is:
  \[ \ket\Psi = (1 \otimes \delta) \circ \delta \circ \eta = \sum_i \alpha_i \ket{u_i,u_i,u_i} \]
  For $\ket\Psi$ to be SLOCC maximal, all of the scalars $\alpha_i$ must be non-zero. Let $L$ be defined as:
  \[ L :: \ket{u_i} \mapsto \frac{1}{\alpha_i} \ket{u_i} \]
  Then, applying $L$ to any of the subsystems yields
  \[ \ket{\Psi'} = \sum_i \ket{u_i,u_i,u_i} \]
  \noindent which is clearly SLOCC-equivalent to $\ketGHZD$.
\end{proof}

A complete classification of ACFAs is not yet complete. For dimensions $D < 6$, there are relatively few commutative algebras, up to isomorphism, so it is feasible to enumerate them and check which extend to ACFAs. For $D < 4$, there is only one ACFA, and its Frobenius state is SLOCC-equivalent to $\ketWD$. However, for $D = 4$, there are already two non-isomorphic ACFAs, so the classification of these types of Frobenius states may be more difficult in general.\footnote{Thanks to Alex Merry for performing these calculations.} A classification of \textit{all} Frobenius states for $D < 6$ is in progress, and will be included in a forthcoming sequel to \cite{CoeckeKissinger2010}.

\section{A Graphical Theory for Entanglement}\label{sec:ghzw-calculus}

We now look at some of the attributes of the pair of commutative Frobenius algebras corresponding to the GHZ and W states. $\mathcal G = (\whitemult, \whiteunit, \whitecomult, \whitecounit)$ is a $\dagger$-special commutative Frobenius algebra, so it has an orthonormal basis of classical points. For $\mathcal W = (\mult, \unit, \comult, \counit)$, note that $\mult$ defines a partial monoid over the classical points of $\mathcal G$, and $\comult$ a partial monoid over their adjoints. In other words, $\ket{i\updot j}$ and $\bra{i \downdot j}$ defined as follows must either be classical points for $\mathcal G$ or $0$.
\ctikzfig{ghzw_partial_monoid}

To handle the case where $i \updot j = \bot$ or $i \downdot j = \bot$, we use the following notation for ``undefined points'':
\ctikzfig{undef_points}

Since $\mathcal W$ is \textit{not} a $\dagger$-CFA, it is not necessarily true that $\left(\ket{i \updot j}\right)^\dagger = \bra{i \downdot j}$. We will shortly see that these two partial monoids are isomorphic, but they are never identical for non-trivial ACFAs. The other thing to note is that both the unit and the anti-unit are proportional to classical points of $\mathcal G$.
\ctikzfig{ghzw_unit_antiunit_classical}

We can build up these features abstractly. First, we generalise several results from section \ref{sec:complementary-obs} to the case where one of the Frobenius algebras is not $\dagger$-special. The next definition provides a condition for CFA to form a partial monoid over an orthonormal basis.

\begin{definition}\label{def:cfa-closed}
  An arbitrary commutative Frobenius algebra $(\mult, \unit, \comult, \counit)$ is said to be \textit{closed} on a $\dagger$-special commutative Frobenius algebra $(\whitemult, \whiteunit, \whitecomult, \whitecounit)$ if the following equations hold:
  \ctikzfig{ghzw_bialg}
\end{definition}

For the rest of this section, let  $(\mult, \unit, \comult, \counit)$ be a commutative Frobenius algebra closed on a $\dagger$-SCFA $(\whitemult, \whiteunit, \whitecomult, \whitecounit)$.

\begin{theorem}\label{thm:closed-dualiser}
  For a commutative Frobenius algebra $(\mult, \unit, \comult, \counit)$, closed on $(\whitemult, \whiteunit, \whitecomult, \whitecounit)$, the dualiser is unitary and self-adjoint. Furthermore, it is a permutation of the classical points of $\whitecomult$.
\end{theorem}

\begin{proof}
  We can compute caps and cups in terms of $\updot$ and $\downdot$. It follows from Definition \ref{def:cfa-closed} that $\unit = \ket{u}$ is a classical point for $\whitecomult$ and $\counit = \bra{c}$ is the adjoint of a classical point. So:
  \[ %
\beginpgfgraphicnamed{black_cap}
\InputIfFileExists{black_cap.tikz}{}{\input{./figures/black_cap.tikz}}
\endpgfgraphicnamed = \sum_{i\downdot j = u} \ket{i,j} \qquad \qquad \qquad
\beginpgfgraphicnamed{black_cup}
\InputIfFileExists{black_cup.tikz}{}{\input{./figures/black_cup.tikz}}
\endpgfgraphicnamed = \sum_{i\updot j = c} \bra{i,j} \]
     
  From these, we can compute the dualiser and its inverse.
  \[ d      = %
\beginpgfgraphicnamed{gw_dualiser}
\InputIfFileExists{gw_dualiser.tikz}{}{\input{./figures/gw_dualiser.tikz}}
\endpgfgraphicnamed = \sum_{i\downdot j = u} \ketbra{i}{j} \qquad \qquad \qquad 
     d^{-1} = %
\beginpgfgraphicnamed{gw_dualiser2}
\InputIfFileExists{gw_dualiser2.tikz}{}{\input{./figures/gw_dualiser2.tikz}}
\endpgfgraphicnamed = \sum_{i\updot j = c} \ketbra{i}{j} \]
  
  Written in the basis given by $\whitecomult$, these are both binary matrices (i.e. matrices whose entries are all either $0$ or $1$). But then, the only binary matrices whose inverses are also binary matrices are the permutations, so $d$ must be a permutation, and $d^{-1} = d^\dagger$ its inverse permutation. The fact that this map is self-adjoint follows from the symmetry of caps and cups. Evaluating the cup at classical points, we have:
  \ctikzfig{gw_self_adj_pf}
\end{proof}

For the remainder of this chapter, we will represent the dualiser by placing a tick on an edge.
\ctikzfig{dualiser_tick}

The dualiser generates a monoid isomorphism from $\updot$ to $\downdot$:
\begin{equation}\label{eq:tick-monoid-iso}
\beginpgfgraphicnamed{tick_monoid_iso}
\InputIfFileExists{tick_monoid_iso.tikz}{}{\input{./figures/tick_monoid_iso.tikz}}
\endpgfgraphicnamed
\end{equation}

\begin{definition}\label{def:ghz-w-pair}
  Let $\mathcal S = (\whitemult, \whiteunit, \whitecomult, \whitecounit)$ be a $\dagger$-SCFA and let $\mathcal A = (\mult, \unit, \comult, \counit)$ be an ACFA. Then $(\mathcal S, \mathcal A)$ is called a \textit{GW-pair} if $\mathcal A$ is closed on $\mathcal S$ and the following equations hold:
  \ctikzfig{antiunit_cp}
\end{definition}

Technically, we only have to require that $\lolli$ be \textit{proportional} to a classical point. Then, the fact that $\counit \circ \lolli = \circl = D$ uniquely fixes the scalar factor. Furthermore, we can prove that the dualiser interchanges the unit and anti-unit.

\begin{lemma}\label{lem:tick-lolli}
  The following equations hold for any GW-pair:
  \[ \circl\ \left(\unit\right)^\dagger = \cololli \qquad \qquad
     \circl\ \left(\counit\right)^\dagger = \lolli \qquad \qquad
     \circl\ \tickunit = \lolli \]
\end{lemma}

\begin{proof}
  By Definition \ref{def:ghz-w-pair}, $\lolli$ and $\cololli$ correspond to classical points for $\whitecomult$, up to scalar. By closure, $\unit$ and $\counit$ are also classical points. Since $\counit \circ \lolli = \circl \neq 0$, then $\left(\counit\right)^\dagger$ and $\lolli$ must be proportional. By closure, $\counit$ is \textit{equal} to a classical point, hence it is normalised. So $\circl\ \left(\unit\right)^\dagger = \cololli$. The second equation holds similarly. We can also show that $\tickunit = \left(\counit\right)^\whitetranspose$:
  \ctikzfig{tick_is_transpose}
  
  Since $\counit$ is the adjoint of a classical point, $\left(\counit\right)^\dagger = \left(\counit\right)^\whitetranspose = \tickunit$. The final equation then follows.
\end{proof}

A consequence of this lemma is that the monoids $\updot$ and $\downdot$ cannot be equal for dimensions $D \geq 2$. By anti-specialness and equation (\ref{eq:tick-monoid-iso}), that would imply the identity was rank $1$. In \cite{CoeckeKissinger2010}, Coecke and Kissinger identified four axioms for interacting GHZ- and W-like states.
\ctikzfig{ghz_w_axioms}

\begin{theorem}\label{thm:gw-pair-subsumes-axioms}
  Any \textit{GW-pair} satisfies axioms (i.)-(iv.).
\end{theorem}

\begin{proof}
  Axioms (i.) and (ii.) are consequences of Theorem \ref{thm:closed-dualiser}. Axiom (iii.) is by definition, and axiom (iv.) is part of Lemma \ref{lem:tick-lolli}.
\end{proof}

It is known that not all SCFA/ACFA pairs in $\catRel$ satisfy the closure identities, so it is likely that the GW-pair conditions are strictly stronger than axioms (i.)-(iv.). However, the Frobenius states $\ketGHZD$ and $\ketWD$ provide at least one example of a GW-pair for any dimension in $\catFHilb$.

\subsection{Symmetric Modules of an SCFA and Distributivity}

For any Hilbert space $\mathcal H$, we can form the space $\mathcal H \otimes_S \mathcal H$ of \textit{symmetric} vectors in $\mathcal H \otimes \mathcal H$. There is a canonical projection from an arbitrary vector in $\mathcal H \otimes \mathcal H$ onto a symmetric vector in $\mathcal H \otimes_S \mathcal H$, called the \textit{symmetriser map} $\mathcal S_2$:
\ctikzfig{symmetriser}

Any monoid $(\mathcal H, \mu, \eta)$ can be extended to a monoid $(\mathcal H \otimes_S \mathcal H, \mu_S, \eta_S)$ by using the symmetriser.
\ctikzfig{symm_monoid}

There is a canonical $\mu_S$-module over the whole space $\mathcal H \otimes \mathcal H$. It is essentially the regular module of $\mu_S$, but without the symmetriser maps on the right.
\ctikzfig{symm_adjoint_module}

We shall refer to this as the \textit{extended regular module} $x_\mu$. For any SCFA, there is also a $\mu_S$-module $k_{\bra{u_i}}$ over $\mathcal H$ for every classical point $\ket{u_i}$.
\ctikzfig{symm_cp_module}

For a GW-pair, $\bra{c} := \tickcounit$ is the adjoint of a classical point, so it induces a $\mu_S$-module on $\mathcal H$. We call a GW-pair \textit{distributive} if the multiplication $\mult : \mathcal H \otimes \mathcal H \rightarrow \mathcal H$ is a $\mu_S$-module homomorphism from $(\mathcal H\otimes \mathcal H, x_\mu)$ to $(\mathcal H, k_{\bra{c}})$.

\begin{definition}\label{def:dist-gw-pair}
  A special commutative Frobenius algebra $(\whitemult, \whiteunit, \whitecomult, \whitecounit)$ and an anti-special commutative Frobenius algebra $(\mult, \unit, \comult, \counit)$ form a \textit{distributive GW-pair} if they are a GW-pair and:
  \begin{equation}\label{eq:dist-gw-pair}
\beginpgfgraphicnamed{symm_module_hm}
\InputIfFileExists{symm_module_hm.tikz}{}{\input{./figures/symm_module_hm.tikz}}
\endpgfgraphicnamed
  \end{equation}
\end{definition}

We refer to this condition as distributivity, because it resembles the distributive law for rings. To see this most clearly, consider arbitrary vectors $\ket{a}, \ket{b}, \ket{c} \in \mathcal H$. Noting that $\ket{a} \otimes \ket{a}$ is symmetric, we can prove the following equation:
\begin{equation}\label{eq:distrib-ex}
\beginpgfgraphicnamed{distributive_ex}
\InputIfFileExists{distributive_ex.tikz}{}{\input{./figures/distributive_ex.tikz}}
\endpgfgraphicnamed
\end{equation}

Ignoring scalars, and writing ``$+$'' for $\mult$ and ``$\cdot$'' for $\whitemult$, this equation becomes:
\[ a \cdot (b + c) = (a \cdot b) + (a \cdot c) \]

It reduces to the distributive law familiar from arithmetic. A slightly different, but equivalent way to see this is $\mult$ copies phases for $\whitemult$, up to a scalar.
\ctikzfig{distributive_phase_copy}

This is similar to the case for strongly complementary observables, as in equation (\ref{eq:classical-phase-copy}) from section~\ref{sec:complementary-obs}. However, the condition here is much stronger, because $\mult$ copies \textit{arbitrary} phases for $\whitemult$, rather than just those corresponding to classical points.

\begin{example}
  The Frobenius algebras $\mathcal G, \mathcal W$ defined in section~\ref{sec:classification} form a distributive GW-pair.
\end{example}

The distributive law implies many identities between $\mult$ and $\whitemult$. In the next lemma are two equations that we shall find useful in the next section.

\begin{lemma}\label{lem:gw-hopf-like}
  For any distributive GW-pair, the following equations hold:
  \begin{center}
\beginpgfgraphicnamed{gw_hopf_like}
\InputIfFileExists{gw_hopf_like.tikz}{}{\input{./figures/gw_hopf_like.tikz}}
\endpgfgraphicnamed\qquad\qquad
\beginpgfgraphicnamed{tick_cancel}
\InputIfFileExists{tick_cancel.tikz}{}{\input{./figures/tick_cancel.tikz}}
\endpgfgraphicnamed
  \end{center}
\end{lemma}

\begin{proof}
  The first equation follows from distributivity, noting that $\whitecap$ is symmetric.
  \ctikzfig{gw_hopf_like_pf}
  
  The second equation then follows from the first one.
  \ctikzfig{tick_cancel_pf}
\end{proof}

\subsection{Universality}\label{sec:ghzw-universality}

Returning to the example of the GHZ and W Frobenius algebras, we treat $(\whitemult, \whiteunit, \whitecomult, \whitecounit)$ and $(\mult, \unit, \comult, \counit)$ as generators for quantum circuits and look at their computational power. While they are not gates themselves, as they are not unitary, we shall soon see that we can think of certain gates as begin composed from these generators, as we did for the $Z/X$ calculus in section \ref{sec:zx-building-circuits}. Alternatively, one can think of the GW-pair maps as stochastic gates, prepared using post-selection or some more sophisticated measurement-based scheme.

First, note that the dualiser interchanges $\unit = \ket{1}$ and $\tickunit = \ket{0}$, so it serves as a NOT (i.e. Pauli $X$) gate.
\[ \tick \circ \unit = \tickunit \qquad\qquad \tick \circ \tickunit = \unit \]

While the GHZ dot copies both $\unit$ and $\tickunit$, the W dot acts like a ``controlled'' copy.
\ctikzfig{gw_copy}

With this behaviour in mind, we can use these generators to build a CNOT gate.
\[ %
\beginpgfgraphicnamed{cnot}
\InputIfFileExists{cnot.tikz}{}{\input{./figures/cnot.tikz}}
\endpgfgraphicnamed\quad:=\quad%
\beginpgfgraphicnamed{gw_cnot}
\InputIfFileExists{gw_cnot.tikz}{}{\input{./figures/gw_cnot.tikz}}
\endpgfgraphicnamed \]

We can verify that this is indeed a CNOT gate by rewriting.
\ctikzfig{gw_cnot_pf}

\begin{theorem}
  The generators of $\mathcal G$ and $\mathcal W$, along with single-qubit states, are universal for quantum computation.
\end{theorem}

\begin{proof}
  We have already illustrated the construction of a CNOT gate. To complete the proof, it suffices to show that we can apply arbitrary single-qubit unitaries. We can actually do better than this by showing we can apply arbitrary single-qubit linear maps. We can write a general single-qubit diagonal matrix as a GHZ phase. Let $\ket{\psi_1} = a \ket{0} + b \ket{1}$.
  \begin{equation}\label{eq:diag-matrix}
\beginpgfgraphicnamed{ghz_phase}
\InputIfFileExists{ghz_phase.tikz}{}{\input{./figures/ghz_phase.tikz}}
\endpgfgraphicnamed = a \ketbra{0}{0} + b \ketbra{1}{1} =
    \left( \begin{matrix} a & 0 \\ 0 & b \end{matrix} \right)
  \end{equation}
  
  For $\ket{\psi_2} = c \ket{0} + \ket{1}$, we can construct an arbitrary unit-diagonal upper triangular matrix as a W phase.
  \begin{equation}\label{eq:up-tri-matrix}
\beginpgfgraphicnamed{w_phase}
\InputIfFileExists{w_phase.tikz}{}{\input{./figures/w_phase.tikz}}
\endpgfgraphicnamed = \ketbra{0}{0} + c \ketbra{0}{1} + \ketbra{1}{1} =
    \left( \begin{matrix} 1 & c \\ 0 & 1 \end{matrix} \right)
  \end{equation}
  
  For $\ket{\psi_3} = d \ket{0} + \ket{1}$, we can construct an arbitrary unit-diagonal lower triangular matrix by applying the dualiser.
  \begin{equation}\label{eq:down-tri-matrix}
\beginpgfgraphicnamed{w_tick_phase}
\InputIfFileExists{w_tick_phase.tikz}{}{\input{./figures/w_tick_phase.tikz}}
\endpgfgraphicnamed = \ketbra{0}{0} + d \ketbra{1}{0} + \ketbra{1}{1} =
    \left( \begin{matrix} 1 & 0 \\ d & 1 \end{matrix} \right)
  \end{equation}
  
  Any linear map decomposes as $M = PLDU$, where $P$ is a permutation, $L$ is a unit-diagonal lower triangular matrix, $D$ is diagonal, and $U$ is a unit-diagonal upper triangular matrix. Since the only permutations on $\mathbb C^2$ are the identity and $\tick$, we can construct any single-qubit map using $\tick$ and the maps above.
\end{proof}

\subsection{Arithmetic on the Complex Projective Line}\label{sec:arithmetic}

It is a well-known fact that points on the Bloch sphere correspond to points on the complex projective line. Any state $a \ket{0} + b \ket{1}$ where $b \neq 0$ can be represented, up to a scalar by the quotient $\frac{b}{a}$. Defining $\numket{b/a} := \ket{0} + \frac{b}{a} \ket{1}$ and $\numket\infty := \ket{1}$, we can cover the entire Bloch sphere. Then, the usual projection of a sphere on to $\cpline$ takes these states to their corresponding points in $\cpline$.
\ctikzfig{projection}

Under this correspondence, we can see that the algebra induced by the W state is addition on $\cpline$ and the GHZ algebra is multiplication. Before we illustrate this, we define addition and multiplication on $\cpline$ as commutative partial monoids. For $k_1, k_2 \in \mathbb C$, addition and multiplication are defined as usual. For the rest, let $k \in \mathbb C$ be a non-zero complex number, and let $\bot$ represent undefined.
\begin{equation}
\begin{aligned}
  k \cdot \infty      & = \infty & \qquad\qquad\qquad & k + \infty      & = \infty \\
  0 \cdot \infty      & = \bot   &                    & 0 + \infty      & = \infty \\
  \infty \cdot \infty & = \infty &                    & \infty + \infty & = \bot  
\end{aligned}
\end{equation}

Intuitively, these are addition and multiplication operations for ``formal fractions'' over $\mathbb C$. That is, equivalence classes of pairs of complex numbers:
\[ |(d,n)| = \{ (d,n) \sim (\lambda d, \lambda n) : \lambda \in \mathbb C - \{0\}\  \} \]

Letting $\infty := (0,1)$ and $\bot := (0,0)$, we can reproduce the above multiplication tables with:
\[ |(d_1,n_1)| + |(d_2,n_2)| = |(d_1 d_2, n_1 d_2 + n_2 d_1)|
\qquad\qquad
|(d_1,n_1)| \cdot |(d_2,n_2)| = |(d_1 d_2, n_1 n_2)| \]

Using the convention that $\numket{\bot} = 0$, it is straightforward to verify the following equations.
\ctikzfig{gw_plus_times}

It follows from equation \ref{eq:dist-gw-pair} that the finitary points in $\cpline$ distribute over addition. That is, for $k \in \mathbb C$:
\ctikzfig{cp1_distrib}

Distributivity fails for $\infty$, as is usually the case when formally introducing points at infinity.
\ctikzfig{cp1_distrib_fail}

The relationship between GHZ states, W states and the arithmetic of fractions is explored in detail by Coecke, Kissinger, Merry, Roy in \cite{CKMR2011}.


	\part{Automation}\label{part:automation}
	



\chapter{Automating String Graph Rewriting: Quantomatic}\label{ch:automation}

The Quantomatic Project~\cite{Quantomatic} provides a set of tools for working with string graphs and string graph rewrite systems. It is divided into three parts: QuantoCore, QuantoGUI, and QuantoCoSy. The first is called QuantoCore, which is an ML library that does most of the work in representing, manipulating, and rewriting string graphs. QuantoCore uses three basic kinds of files.
\begin{enumerate}
  \item \textbf{*.graph} files store string graphs.
  \item \textbf{*.theory} files store graphical theories. A graphical theory contains information about what kinds of vertices can occur in a string graph, what kinds of data can occur on vertices, and how that data should be matched.
  \item \textbf{*.rules} files contain sets of string graph rewrite rules.
\end{enumerate}

\begin{figure}
  \centering
  \includegraphics[width=14cm]{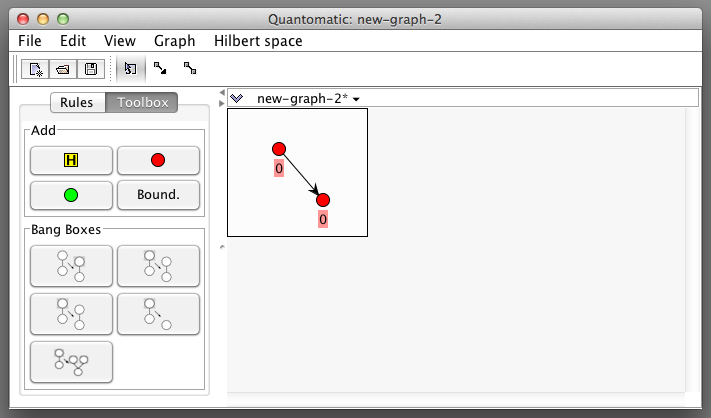}
  \caption{\label{fig:quanto-screen1} A string graph in Quantomatic}
\end{figure}

\begin{figure}
  \centering
  \includegraphics[width=14cm]{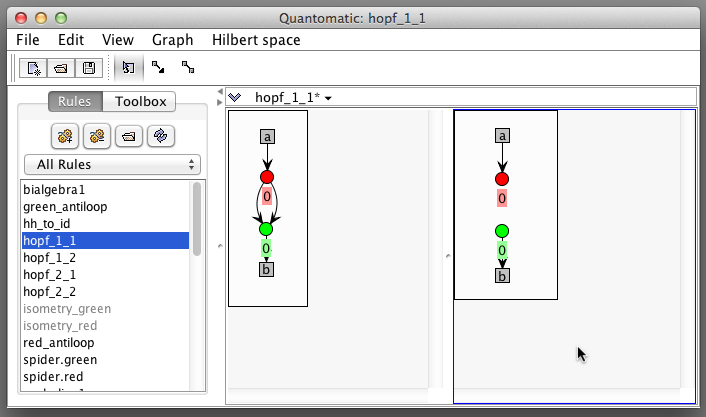}
  \caption{\label{fig:quanto-screen2} Editing a rule}
\end{figure}

\begin{figure}
  \centering
  \includegraphics[width=14cm]{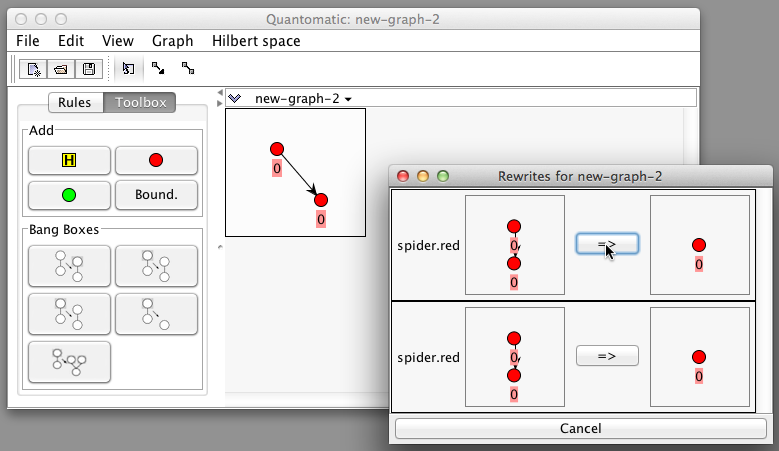}
  \caption{\label{fig:quanto-screen3} Applying a rewrite}
\end{figure}

\begin{figure}
  \centering
  \includegraphics[width=14cm]{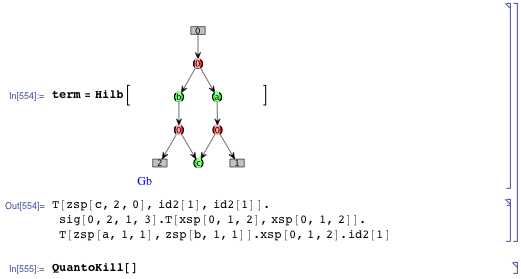}
  \caption{\label{fig:quanto-cas} Invoking Quantomatic from a Mathematica~\cite{Mathematica} notebook}
\end{figure}

On top of QuantoCore, we have developed a graphical user interface called QuantoGUI. Currently, QuantoGUI can:
\begin{itemize}
  \item create and edit string graphs (Figure \ref{fig:quanto-screen1}) and string graph rewrite systems (Figure \ref{fig:quanto-screen2}),
  \item search for rewrites in a selected subgraph and apply them manually (Figure \ref{fig:quanto-screen3}),
  \item display animated normalisations of string graphs with respect to a rewrite system,
  \item do ``fast-normalisation'' of string graphs and only display the output, and
  \item interact with computer algebra systems to perform concrete calculations of string graphs as linear maps (Figure \ref{fig:quanto-cas}, see~\cite{Calculemus2009} for details).
\end{itemize}

While it is currently quite minimal, we intend to make QuantoGUI into the graphical analogue of a proof assistant. Just as Isabelle~\cite{Isabelle} or Coq~\cite{Coq} exposes a variety of techniques for constructing formal proofs with respect to a term rewrite system (i.e. an algebraic theory), Quantomatic aims to do the same for string graph rewrite systems. The real power of systems like Isabelle and Coq comes not just from rewriting, but from the application of inductive reasoning techniques to prove theorems in first- or higher-order logic. These techniques do not always extend straightforwardly from terms to graphs. However graphs with quantifiers (defined by Rensink et al in ~\cite{Rensink2006}), pattern graphs, and graphs defined by grammars can give us the ability to reason about an infinite set of graphical equations simultaneously. We discuss how some of these more advanced techniques might work in section~\ref{sec:fw-pattern-graphs}.

For more details about the Quantomatic project and to download the software itself, visit the project's web page at
\href{http://sites.google.com/site/quantomatic}{\color{blue} \underline{http://sites.google.com/site/quantomatic}}.

\section{Conjecture Synthesis and QuantoCoSy}

In this section, we will discuss QuantoCoSy, the third component to the Quantomatic project. It performs automated theory creation using a technique called \textit{conjecture synthesis}.

One of the main goals of automated reasoning is to reproduce as much as possible on a machine the way a human mathematician thinks and a works. Consider a situation where a mathematician has the following:
\begin{enumerate}
  \item A set of generators for a new algebraic object $X$.
  \item A concrete model or set of models for $X$.
\end{enumerate}

Though $X$ is not defined yet, the models that the mathematician has in hand are things that ``morally'' should be $X$'s. For instance, if the mathematician were trying to develop an algebra for studying maximally entangled states, as in chapter \ref{ch:monoidal-entanglement}, these would be a set of known maximally entangled states. From this data, the mathematician now seeks to \textit{axiomatise} $X$.

One way he or she could start this process is to ``plug-and-chug''. That is, the mathematician could plug these generators together randomly and see which compositions equal other compositions. In actuality, this process is not totally random, as the mathematician calls upon experience and a handful of helpful heuristics for seeking out likely equations.

\textbf{Heuristic 1: seek familiarity.} the mathematician may discover that these generators are actually satisfying some properties of a known algebraic object, such as a Hopf algebra. From this, the mathematician deduces that the generators are more likely to satisfy \emph{the rest} of the identities of a Hopf algebra than they are to satisfy some other, randomly chosen identities.

\textbf{Heuristic 2: avoid redundancy.} Once the mathematician has a handful of identities, then while searching for new identities, he or she will avoid those which are trivially derivable from those already known. For instance, if the mathematician already knows $a \cdot (b + c) = (a \cdot b) + (a \cdot c)$, it is redundant to consider the terms $a \cdot (b + c + d)$ and $(a \cdot b) + (a \cdot c) + (a \cdot d)$.

\textbf{Heuristic 3: elegant identities are essential.} Though the generators in question may exhibit large, complex, and asymmetric identities, the mathematician is treating the generators as merely one example of a more abstract mathematical object. His experience tells him that simpler identities tend to be the most crucial in characterising abstract mathematical objects.

Though there is theoretically no end to this procedure, the mathematician may hold a conviction, as per heuristic 3, that he or she will eventually find no more ``interesting'' identities above a certain size. At this point, the procedure is effectively complete.

This section is about reproducing this process of (graphical) theory generation with a program called QuantoCoSy.

\subsection{Conjecture Synthesis for Terms}

Conjecture synthesis is a technique that automatically generates ``reasonable'' conjectures to test for an algebraic theory. This procedure for term-based theories was introduced by Johansson, Dixon, and Bundy in 2010~\cite{Johansson:2010tk}. The tool that implements their technique is called IsaCoSy (for Isabelle COnjecture SYnthesis). A single round of their algorithm proceeds as follows:
\begin{enumerate}
  \item Initialise a conjecture as an expression with \textit{holes}, e.g. ``$(*) = (*)$''. Holes mark places where there is more expansion to do. They also come with certain constraints on what terms can be instantiated in them. Initially these are maximum size or depth constraints to guarantee termination, but they will get updated later.
  \item Substitute holes with all possible terms-with-holes. This is done by a depth-first enumeration of possible substitutions, respecting the constraints on holes.
  \item Once there are no holes, save the expression as a possible conjecture.
  \item Perform (fast) post-filtering of conjectures that are obviously not true. IsaCoSy does this by using Isabelle's fast counter-example search.
  \item Try to prove the remaining conjectures using an automated proof search routine. When a proof is found, save the conjecture as a new rewrite rule.
\end{enumerate}

When viewed as a single round, this looks like the na\"ive technique that one would use to go about searching for conjectures. However, the interesting part is how constraints are updated between rounds. A given round of the synthesis procedure will produce a set of true equations, $E$. We can turn these equations into rewrite rules by putting a reduction ordering on terms (see Definition \ref{def:reduction-ordering} in section \ref{sec:term-rewriting}).
\[ S = \{ t_1 \rewritesto t_2 : \left((t_1, t_2) \in E \vee (t_2, t_1) \in E\right) \wedge \omega(t_1) > \omega(t_2)  \} \]

Let $R$ be the set of all terms occurring on the LHS of rules in $S$. These are called \textit{reducible expressions}, or \textit{redexes}. For a rewrite rule $\textbf{r} = (t_1 \rewritesto t_2)$, any equation containing $t_1$ can already be proved using $\textbf{r}$ and an equation containing $t_2$. As conjectures involving redexes are redundant, they are never considered. We can cut them out of the search space by updating the constraints on holes, such that the search performed in step (2.) above never generates terms containing redexes.

While this is a fairly simple procedure, in practice, this can exponentially reduce the search space.

\subsection{Adapting Conjecture Synthesis to String Graphs}

\begin{figure}
  \centering
  \includegraphics[width=14cm]{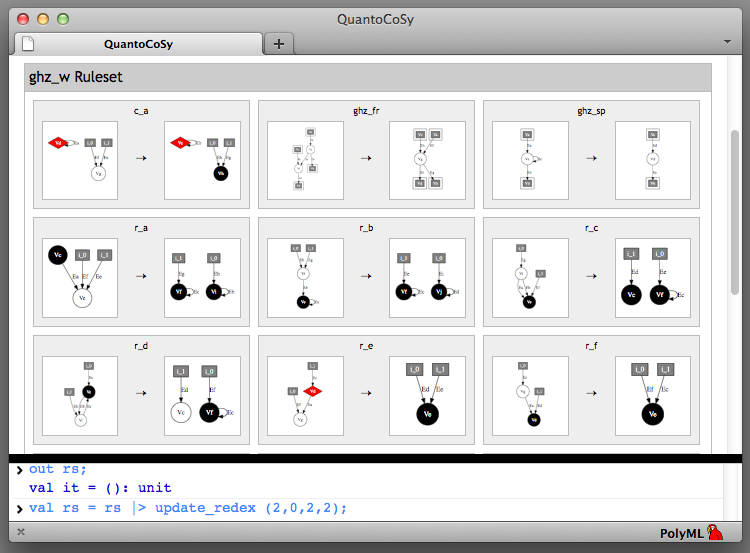}
  \caption{\label{fig:quantocosy} QuantoCoSy runs in Firefox using PolyChome~\cite{PolyChrome}, an extension for running Poly/ML in a web browser.}
\end{figure}

The synthesis procedure for string graph identities is similar to that for terms. The procedure described in this section is implemented on top of QuantoCore with a program called QuantoCoSy (Figure \ref{fig:quantocosy}).

In the term case, a single round of synthesis is parametrised by two natural numbers: the maximum term size (or term depth) and the maximum number of free variables occurring in the term. For string graphs, we parametrise a run with four natural numbers: the number of inputs $M$, the number of outputs $N$, the maximum number of box-vertices $B$, and the maximum number of pluggings $P$. We enumerate string graphs by starting with disconnected string graphs, i.e. graphs only containing box-vertices and their adjacent edges and wire-vertices. For instance:
\ctikzfig{sg_disconnected_graph}

The synthesis procedure for string graphs takes as input:
\begin{enumerate}
  \item A string-graph signature $T$. For simplicity, we will focus on string graph signatures with only a single wire type.
  \item An $(m,n)$-tensor for every box in $T$ with $m$ inputs and $n$ outputs.
  \item A function $\omega$ from string graphs to a well-ordered poset $(P,\leq)$. This will serve as a (candidate) reduction ordering.
\end{enumerate}

Throughout the synthesis, we maintain a rewrite system $S$, and a set $R$ of reducible string graphs. A single run given by natural numbers $(M,N,B,P)$ consists of the following steps:
\begin{enumerate}
  \item For all $p$ such that $0 \leq p \leq P$, generate all disconnected string graphs with $M + p$ inputs and $N + p$ outputs, up to isomorphism.
  \item For a disconnected string graph with $M + p$ inputs and $N + p$ outputs, there are $(M+p)\cdot(N+p)$ input/output pairs. Choose $p$ of them to plug together. These are chosen so that the enumeration is exhaustive and minimises the occurrences of isomorphic string graphs (any remaining isomorphic graphs will be filtered out later). After each plugging, if a string graph contains an element of $R$ as a subgraph, terminate that branch of the enumeration.
  \item Evaluate the string graphs as tensors, performing a tensor contraction for every edge (c.f. the construction of $\widehat F$ in Theorem \ref{thm:sg-free-traced}). Organise them into equivalence classes, up to scalar factors and permutations of inputs and outputs (which are stored with the string graphs). Filter out any remaining isomorphic graphs.
  \item For each equivalence class $C$, identify a set $C' \subseteq C$ of minimal elements with respect to $\omega$. Add any string graph in $C - C'$ to the set of reducible graphs $R$. Choose a string graph $s \in C'$ at random and add rules $t \rewritesto s$ to the rewrite system $S$ for all $t \in C - C'$. Add rules in both directions ($s \rewritesto s'$, $s' \rewritesto s$) for the other minimal graphs $s' \in C' - \{ s \}$.
\end{enumerate}

We postpone filtering out isomorphic graphs until step (3.) because tensor contraction is fast, and two graphs will not be isomorphic unless they are in the same equivalence class. We choose $\omega$ such that step (4.) picks out as few graphs in $C' \subseteq C$ as possible. If $C'$ is a singleton, all of the rewrites respect the reduction order. Rewrites that do not strictly decrease $\omega(G)$ (i.e. rewrites from an element of $C'$ to another element of $C'$) are called \textit{congruences}. We can retain a terminating rewrite system if we throw out all of the congruences, but not without losing some information about the model. Therefore, a large portion of the refinement of this technique has to do with eliminating congruences or handling them in smarter ways (i.e. building them in to the graph representation).

We applied QuantoCoSy to generators of the GW-pair given in section~\ref{sec:ghzw-calculus}. We preloaded the synthesis procedure with rewrites that merge two vertices of the same colour (i.e. the spider laws), and synthesised graphical identities for $B=3$, $P=3$, and $M + N \leq 3$. This yielded $223$ rewrite rules, most of which were versions of the four axioms given in Theorem~\ref{thm:gw-pair-subsumes-axioms}.

\begin{figure}\label{fig:redex-plot}
  \centering
\beginpgfgraphicnamed{naive_vs_redex_plot}
\InputIfFileExists{naive_vs_redex_plot.tikz}{}{\input{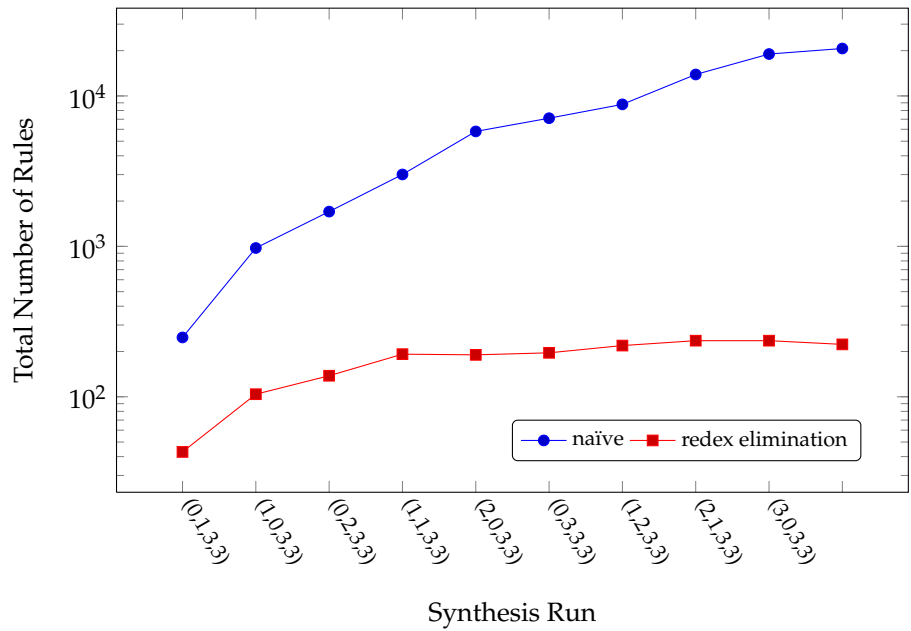}}
\endpgfgraphicnamed
  \caption{Rewrite rules synthesised from the generators of the GW-pair defined in section \ref{sec:ghzw-calculus}.}
\end{figure}

As was the case for terms, filtering out redexes has a huge impact on the number of string graphs that need to be checked. The na\"ive synthesis procedure with the same parameters yielded over $20,000$ rewrite rules. In Figure \ref{fig:redex-plot}, we plot the number of rewrite rules generated using a na\"ive graph enumeration algorithm against the number generated using the procedure above.

While the redex-elimination procedure yielded a much more manageable number of rewrite rules, a quick look at the rules by a human will show that many are still trivially consequences of each other. Therefore, there is still much work to be done in eliminating redundancy in rules and building more of the symmetries of a rewrite theory into the graphical representation itself. We discuss some of the ways in which we hope to accomplish this in sections~\ref{sec:graphical-knuth-bendix} and~\ref{sec:fw-pattern-graphs}.





\chapter{Conclusion}\label{ch:conclusion}

The main contributions of this dissertation fall under three categories: (1) the definition and properties of string graphs, (2) the application of string diagrams/string graphs to quantum computation, and (3) the automated generation and manipulation of string graph rewrite systems.

First, we defined string graphs and string graph rewrite systems using double-pushout graph rewriting. We also introduced the notion of composition for string graphs via certain pushouts called \textit{pluggings} and defined a category whose morphisms are string graphs modulo a rewrite system, where categorical composition is defined using plugging. Using these rewrite categories, we constructed the free traced symmetric and the free compact closed categories over a monoidal signature. These results allow us to prove identities in an arbitrary symmetric traced category using graph rewriting.

In~\cite{DixonKissinger2010}, Dixon and Kissinger proved that rewrite categories are equivalent to their topological analogues, as defined by Joyal and Street in~\cite{JS}. Since free categories are defined up to equivalence, the results in this dissertation imply two of the missing ``GTC-II'' theorems. Namely, a category defined using string diagrams forms the free symmetric traced (or compact closed) category over a monoidal signature. 

Next, we illustrated the application of diagrammatic languages to the study of quantum phenomena: namely complementary observables and multipartite entanglement. After reviewing Coecke and Duncan's characterisation of complementarity using interacting Frobenius algebras, we proved several new results, including a complete classification of pairs of strongly complementary observables in an finite dimensional Hilbert space. This classification theorem showed a 1-to-1 correspondence between $D$-dimensional strongly complementary pairs of observables and Abelian groups of order $D$.

After this, we showed how certain kinds of highly entangled states, called Frobenius states, induce Frobenius algebras. Furthermore, the two canonical tripartite qubit states---GHZ and W---are distinguished by a simple condition on these Frobenius algebras: specialness or anti-specialness.
\ctikzfig{scfa_acfa_recap}

A Frobenius state in $\mathbb C^2 \otimes \mathbb C^2 \otimes \mathbb C^2$ is SLOCC-equivalent to GHZ if and only if its induced commutative Frobenius algebra is special, and a Frobenius state is SLOCC-equivalent to W if and only if its induced CFA is anti-special.

Drawing on the interaction properties of these canonical qubit states, we introduced the theory of \textit{GW-pairs}. These abstract the interaction properties of the Frobenius algebras associated with GHZ and W, and can be used to study arbitrary multipartite entangled states. Focusing on the specific case of GHZ and W, we showed that the generators of the GW-pair are universal for quantum computing and characterised the behaviour of the GHZ- and W-algebras in terms of (partial) arithmetic operations defined on the complex projective line $\cpline$.

Finally, in part \ref{part:automation}, we introduced the Quantomatic project, which consists of tools for the automatic construction and manipulation of string graph rewrite theories. We illustrated how the process of conjecture synthesis introduced by Johansson, Dixon, and Bundy~\cite{Johansson:2010tk} can be adapted to the setting of string graphs, where each graphical generator is given a concrete valuation as a linear operator. Using this technique, it becomes practical to enumerate all of the graphical identities exhibited by a set of generators under composition for small- to medium-sized string graphs. These tools, along with the methods they employ, show great potential for changing the way we formulate and interact with a wide variety of theories involving interacting components.

\section{Future Work}\label{sec:future}

\subsection{Classifying Frobenius states}

In section \ref{sec:higher-dim-classification}, we gave a classification result for anti-special commutative Frobenius algebras of dimension $2$ and for special Frobenius algebras of \textit{any} dimension. The natural next step is to ask if a reasonable classification result can be constructed for ACFAs of dimension $D \geq 3$. This classification is likely to be more difficult in the case of SCFAs, where classification follows from the fact that there are not very many semisimple, commutative $K$-algebras. When $K$ is an algebraically closed field, the only semisimple algebras are direct sums of $K$ itself. In the case of ACFAs of dimension $D > 1$, the vector $\lolli$ always generates a non-trivial nilpotent ideal, so non-trivial ACFAs are never semisimple. However, it is our hope that anti-specialness will prove a strong enough condition on a Frobenius algebra to yield a straightforward classification.

We also intend to expand the classification results for Frobenius states on a different front: the classification of \textit{arbitrary} Frobenius states for low dimensions (e.g. $D \leq 5$). This problem is tractable because there are relatively few commutative, unital algebras of dimension $5$ or less. There is one such algebra for dimension $1$, two for $D=2$, four for $D=3$, nine for $D=4$, and $20$ for $D=5$. For dimensions $3$ and above, only some of these algebras extend to Frobenius algebras, and it becomes a straightforward task to enumerate those that are. The classification of Frobenius states for dimension $3$ was completed this year by Honda~\cite{Honda2011}. As in the case for two dimensions, a Frobenius state up to SLOCC is uniquely determined by the rank of its loop map. Thus there is only Frobenius state corresponding to an SCFA (full rank), one corresponding to an ACFA (rank 1), and one corresponding to what Honda calls an \textit{intermediate special commutative Frobenius algebra} (ICFA), which has a loop map of rank 2. The author, along with Coecke and Merry intend to incorporate this result into a complete classification of Frobenius states up to dimension $5$ in the sequel to~\cite{CoeckeKissinger2010}.

\subsection{Super-qubits and the W state bialgebra}

The W state exhibits some interesting properties that we have not yet fully explored. In section \ref{sec:frobenius-algebras}, we introduced the operation $(-)^\whitetranspose$ of transposition relative to a particular Frobenius algebra. In the case of the GW-pair corresponding to the GHZ and W states, let $\altwhitemult = \left(\comult\right)^\whitetranspose$ and let $\altwhiteunit = \left(\counit\right)^\whitetranspose$. The following then forms a bialgebra:
\[ \widehat{\mathcal W} := (\mathbb C \oplus \mathbb C, \altwhitemult, \altwhiteunit, \comult, \counit ) \]

But there is a catch: it does not form a bialgebra in $\catFHilb$, but rather the category $\catSuperHilb$ of \textit{super-Hilbert spaces and even maps}.\footnote{Thanks to Jamie Vicary for pointing this out.} This is the category where Hilbert spaces are graded into a ``bosonic'' and a ``fermionic'' part. Objects are $\mathbb Z_2$-graded Hilbert spaces $H_0 \oplus H_1$ and morphisms are linear maps $f : H_0 \oplus H_1 \rightarrow H'_0 \oplus H'_1$ that respect the grading. This category is monoidally equivalent to the category of unitary representations of $\mathbb Z_2$, but we use a different symmetry map that introduces a $-1$ factor when a ``fermionic'' element crosses another one:
\[ 
\sigma_{H_0 \oplus H_1, H'_0 \oplus H'_1}(\ket{i} \otimes \ket{j}) =
\begin{cases}
  -\ket{j} \otimes \ket{i} & \textrm{ if } \ket{i} \in H_1, \ket{j} \in H'_1 \\
  \ket{j} \otimes \ket{i}    & \textrm{ otherwise}
\end{cases}
\]

For instance, over the graded space $\mathbb C \oplus \mathbb C$ of \textit{super-qubits}, the swap map is defined by the following matrix:
\[
\sigma_{\mathbb C \oplus \mathbb C,\mathbb C \oplus \mathbb C} = 
\left(
\begin{matrix}
  1 & 0 & 0 & 0  \\
  0 & 0 & 1 & 0  \\
  0 & 1 & 0 & 0  \\
  0 & 0 & 0 & -1 \\
\end{matrix}
\right)
\]

Why is it interesting that $\widehat{\mathcal W}$ forms a bialgebra in $\catSuperHilb$, especially if $\catHilb$ is our primary category of interest? There is a faithful forgetful functor $U : \catSuperHilb \rightarrow \catHilb$ that is strongly monoidal, but does not preserve symmetries. As a result, any \textit{planar} diagrammatic identity we can prove in $\catSuperHilb$ also holds in $\catHilb$. For instance, we can prove the following equation for any commutative bialgebra.
\ctikzfig{bialg_path_count}

This is a bit of a contrived example, but its a special case of a general result for commutative bialgebras. A commutative bialgebra diagram is uniquely determined by the number of forward-directed paths from each input to each output. In the above equation, there are exactly 4 distinct paths connecting the input to the output on both the LHS and the RHS. When morphisms exhibit certain properties in one ``categorical context'' but not in another, we express this graphically, using the functorial boxes defined by Melli\`es in~\cite{Mellies2006}. Starting with a diagram in $\catHilb$, we can draw a box around a sub-diagram as long as the diagram (1) is planar and (2) contains only morphisms in the image of $U$. We can then rewrite the elements inside the box as if they are in $\catSuperHilb$, possibly breaking planarity along the way. Then, as long as the diagram is ultimately planar, we can erase the box.
\ctikzfig{func_boxes}

\subsection{GW-pairs and strongly complementary observables}

Aside from contrived examples we have given, one might wonder if there are \textit{useful} planar equations satisfied by bialgebras. In 2008, Melli\`{e}s showed that any commutative bialgebra whose comonoid is the transpose of its monoid, the following maps form a Frobenius algebra~\cite{MelliesTalk}:
\ctikzfig{frob_from_bialg}

Recall that the monoid, comonoid, and Frobenius identites are all planar, so the above yield a Frobenius algebra in both $\catSuperHilb$ and $\catHilb$. Furthermore, in the case of the W state algebra, this Frobenius algebra is precisely the $\dagger$-SCFA for the $X$ observable, defined in section \ref{sec:zx-calculus}. Expanding $\altwhitecomult$ and $\altwhitecounit$ using the dualiser, we obtain the following expression for the Frobenius algebra $(\mathbb C^2, \graymult, \grayunit, \graycomult, \graycounit)$.
\begin{equation}\label{eq:zx-encode}
\beginpgfgraphicnamed{zx_encode}
\InputIfFileExists{zx_encode.tikz}{}{\input{./figures/zx_encode.tikz}}
\endpgfgraphicnamed
\end{equation}

Under this encoding, we can relate the two constructions of the CNOT gate from sections \ref{sec:zx-building-circuits} and \ref{sec:ghzw-universality}.
\ctikzfig{cnot_comparison}

So, we know that in the particular case of the GHZ/W Frobenius algebra pair, we can construct the Frobenius algebras for the strongly complementary observables $Z$ and $X$. However the general relationship between GW-pairs and strongly complementary pairs is still unknown. However, we conjecture that the axioms of a distributive GW-pair subsume those of a strongly complementary pair.

\begin{conjecture}\label{conj:zx-encode}
  Let $(\whitecomult, \comult)$ be a distributive GW-pair. Then, for a third Frobenius algebra $(\graymult, \grayunit, \graycomult, \graycounit)$ defined as in (\ref{eq:zx-encode}) above, $(\whitecomult, \graycomult)$ forms a strongly complementary pair of $\dagger$-SCFAs.
\end{conjecture}

Another area for future work is the conceptual understanding of strong complementarity. In section \ref{sec:complementary-obs}, we provided a complete classification for pairs of strongly complementary observables. While we have a clear idea of what strong complementary means mathematically, a physical interpretation of strong complementary is still missing. The fact that a particular pair of strongly complementary observables (Pauli $Z$ and $X$) play such a central role in the study of complementarity in finite dimensions suggests that such an interpretation exists. As a first step toward finding this interpretation, we are looking for quantum protocols and theorems that rely crucially on certain forms of the bialgebra equations given in Definition \ref{def:strongly-complementary}.

\subsection{Knuth-Bendix completion for string graphs}\label{sec:graphical-knuth-bendix}

Knuth-Bendix completion is a procedure for turning terminating, non-confluent rewrite systems into terminating, confluent rewrite systems. It works by identifying \textit{critical pairs} for a rewrite system. For a finite set of term rewrite rules, we can always identify a finite set of terms $s \in S$ that represent all of the ``possible ways'' in which the left-hand sides of two rewrite rules can overlap. A critical pair is then a pair of distinct, normalised terms $t_1, t_2$ that were both rewritten from some such term $s$. A rewrite system is confluent precisely when it has no critical pairs. Knuth-Bendix completion takes a rewrite system $R$ and a strict, total reduction ordering $\omega$ on terms, and operates as follows:
\begin{enumerate}
  \item Compute all of the critical pairs for $R$.
  \item For each critical pair $(t_1, t_2)$, add $t_1 \rewritesto t_2$ to $R$ if $\omega(t_1) > \omega(t_2)$. Add $t_2 \rewritesto t_1$ otherwise.
  \item Repeat until there are no critical pairs.
\end{enumerate}

There are two possible outcomes for this procedure: (1) it halts and produces a confluent, terminating rewrite system, or (2) it keeps producing more an more rewrite rules forever. Since arbitrary word problems can be encoded as terminating term rewrite systems, there must exist some rewrite systems for which the Knuth-Bendix procedure does not halt. However, for many useful classes of rewrite systems, this procedure always halts, yielding a terminating, confluent rewrite system.

As an example, let $J \trianglelefteq K[X_1,\ldots,X_n]$ be some ideal of a polynomial ring. It is a well-known fact that the ideal membership problem is decidable precisely when one can find a special set of polynomials generating $J$ called a \textit{Gr\"obner basis}. There is a natural way to consider a particular polynomial $J_i$ as a rewrite rule on other polynomials.
\[ X_1^2 X_2 + 4 X_2 + 2 \qquad\longrightarrow\qquad
   \left( X_1^2 X_2 \rewritesto -4 X_2 - 2 \right) \]

This rule then rewrites certain polynomials $P$ that are ``matched'' by $J_i$ to $P - J_i$. Gr\"obner bases are then exactly those sets of polynomials generating $J$ that, considered as rewrite systems, are terminating and confluent. A crucial tool for computer algebra systems is Buchberger's algorithm, which derives Gr\"obner bases from arbitrary finite sets of polynomials. This algorithm is precisely Knuth-Bendix completion applied to polynomial rewrite systems~\cite{Baader1998}.

Confluence for diagrammatic rewrite systems can be more subtle than term rewrite systems. For example, it was shown by Lafont that for general diagrammatic rewriting, finite diagrammatic rewrite systems could lead to infinite families of critical pairs exhibiting what he calls \textit{global conflicts}~\cite{Lafont08}. However, this problem can be overcome by performing critical pair analysis on diagrams with ``gaps'', as Mimram did in his thesis~\cite{MimramThesis} to demonstrate a locally confluent presentation of $\catMat(\mathbb N)$.

The source of the subtlety here lies in the fact that the diagrams considered by Lafont et al exist in arbitrary monoidal categories. Any additional categorical structure (like symmetries and duals) are treated ``opaquely'' as morphisms in a 2D grid. String graph rewrite systems rely crucially on the fact that the traced symmetric structure of the category is \textit{built in} to the graphical representation of the string graph. Therefore much of the subtlety of the Lafont-style approach falls away (or more precisely, is absorbed in graph isomorphism) at the expense of generality.

In 2008, Kissinger defined a Knuth-Bendix procedure for diagrams of interacting commutative Frobenius algebras and applied it by hand to derive a confluent fragment of the Z/X-calculus~\cite{Kissinger2008}. Extending this to general string graph theories is straightforward, and its implementation in Quantomatic will be a useful tool both on its own and as a component of more sophisticated procedures.

\subsection{Pattern graphs and graphical inductive reasoning}\label{sec:fw-pattern-graphs}

We often wish to work with graphical generators that have commutative inputs and outputs with variable arities. We can encode this into the usual string graph formalism by adding a box type for every possible arity and adding rewrite rules for commutativity. However, experience has showed us that there is much to be gained by encoding as much symmetry into the \textit{representation} of an algebraic system as possible. For instance, in term rewrite theories, once a function is assumed to be commutative, its arguments are treated as a multiset, rather than an ordered list. We incorporate this into the theory of string graphs by introducing \textit{string graphs with arities}.

Monoidal signatures are replaced by \textit{signatures with arities} $\widehat T = (O, M, \dom, \cod)$. For $P\mathbb N$ the powerset of $\mathbb N$, the functions:
\[ \dom : M \rightarrow w(O \times P\mathbb N) \qquad\qquad
   \cod : M \rightarrow w(O \times P\mathbb N) \]
\noindent assign a morphism to a list of pairs $(o \in O, A \subseteq \mathbb N)$. The element $o$ is the type of the input or output as before, and the set $A$ is the set of allowed arities. The requirement that the typing maps $\tau_G : G \rightarrow G_T$ for string graphs must be local isomorphisms around box vertices is replaced with the requirement that these maps respect arities. For an edge $e$ in the typegraph, we require that $e$ occurs with an allowed arity around every box-vertex $v \in B(G)$. Formally, for all box-vertices $v \in B(V)$ where $\tau_G(v) = f$ and all edge types $\textrm{in}_{f,i},\textrm{out}_{f,j}$ such that $\dom(f)[i] = (o, A)$ and $\cod(f)[j] = (o', A')$ the following equations must hold:
\begin{equation}\label{eq:arity-restriction}
  \begin{split}
  |(\tau^v_G)^{-1}(\textrm{in}_{f,i})|  & \in A \\
  |(\tau^v_G)^{-1}(\textrm{out}_{f,j})| & \in A'
  \end{split}
\end{equation}

Using arities admits a great deal of flexibility in string graph signatures. The usual notion of string graphs without arities is recovered by letting all of the sets of allowed arities by $\{ 1 \}$, because $\tau_G$ is a local isomorphism iff the inverse images defined in (\ref{eq:arity-restriction}) are always of cardinality $1$. The other extreme is where all of the sets $A = \mathbb N$, where any arity is allowed. This is incredibly useful when working with spiders, as defined in section \ref{sec:nf-frobenius}. In fact, this is default mode for Quantomatic, as it was originally designed to work with ``spider-based'' graphical languages like those defined in part \ref{part:rewriting}. In between these two extremes, one could define optional (non-commutative) inputs and outputs by setting $A = \{ 0, 1 \}$ and (fixed arity) commutative inputs and outputs by setting $A = \{ k \}$.

One of the most useful things about having generators with variable arities (as opposed to just introducing new generators for every arity) is that we can define \textit{pattern graphs}. These are graphs where certain portions of a graph (and their incident edges) can be duplicated any number of times. We define pattern graphs by introducing \textit{$!$-boxes} (pronounced ``bang boxes'') around string subgraphs. Intuitively, pattern graphs represent a set of concrete string graphs where each of the $!$-boxes can occur $0$ or more times.

\begin{example}
  A pattern graph with two $!$-boxes, and the set of concrete graphs it represents.
  \ctikzfig{bang_graph_ex}
\end{example}

The real power of pattern graphs comes from the ability to define \textit{pattern graph rewrite rules}. These consist of a pair of pattern graphs, and a suitable bijection between the $!$-boxes on the LHS and RHS such that we can represent an infinite set of valid string graph rewrite rules. For instance, a rewrite rule to merge two vertices of any arity can be expressed as:
\begin{equation}\label{eq:bang-graph-sp}
\beginpgfgraphicnamed{bang_graph_rewrite_rule}
\InputIfFileExists{bang_graph_rewrite_rule.tikz}{}{\input{./figures/bang_graph_rewrite_rule.tikz}}
\endpgfgraphicnamed
\end{equation}

A pattern graph rule can be \textit{instantiated} into a concrete rule by replacing a single $!$-box with $N$ copies of that $!$-box on the LHS and replacing its corresponding $!$-box with $N$ copies on the RHS. For example, the following is an instance of the rule above:
\ctikzfig{bang_graph_rewrite_rule_instance}

Although a preliminary implementation of $!$-boxes already exists in Quantomatic, to assure its validity, we need to formalise string graphs with arities and pattern graph matching and rewriting within the framework of partial adhesive categories. In many ways this is a straight-forward task, but care must be taken when defining the correct notions of matching and rewriting.

We can also do rewriting on pattern graphs themselves. That is, we can use infinite sets of rules (i.e. pattern graph rewrite rules) to reason about infinite sets of graphs (i.e. pattern graphs).

One of the most interesting applications for pattern graph rewriting is combining the conjecture synthesis procedure with Knuth-Bendix completion to automatically generate new pattern graph rewrite rules. IsaScheme, a tool for ``scheme-based'' conjecture synthesis~\cite{MontanoRivas2010}, has already had some success in this area for term-based theories.

As it stands, the synthesis procedure can only discover new concrete rewrite rules, however, using Knuth-Bendix, it could automatically \textit{combine} its pre-existing knowledge (in the form of an initial set of pattern graph rewrite rules) with its findings.

Consider a case where we initiate the synthesis procedure with the rewrite rule given by (\ref{eq:bang-graph-sp}). At some point, it discovers a new identity:
\begin{equation}\label{eq:sg-black-dot-copy}
\beginpgfgraphicnamed{sg_black_dot_copy}
\InputIfFileExists{sg_black_dot_copy.tikz}{}{\input{./figures/sg_black_dot_copy.tikz}}
\endpgfgraphicnamed
\end{equation}

If we perform Knuth-Bendix on a rewrite system, we obtain a critical pair.
\ctikzfig{bang_graph_cp}

Then, under a suitable ordering on string graphs, we can consider the string graph on the right to be more reduced than the string graph on the left, so we introduce a new pattern-graph rewrite rule.
\ctikzfig{bang_graph_cp_new_rule}

This rule is stronger than (\ref{eq:sg-black-dot-copy}), and can be thought of as the ``spiderised'' version of that rule. Since this rule is stronger than the previous one, there are more reducible expressions in the rewrite system being synthesised, and hence a smaller search space for string graph enumeration. This suggests that incorporating a Knuth-Bendix step into the conjecture synthesis procedure could vastly improve its performance as well as generate fewer, more powerful graphical identities.

In addition to pure equational reasoning (i.e. rewriting proofs) with $!$-boxes, we can do some inductive reasoning as well. Suppose we extend the language of $!$-boxes, allowing them to be bound to expressions over natural numbers, possibly containing free variables.
\ctikzfig{bang_graph_rewrite_rule_bound}

We interpret this rule as the set of all concrete rules where the $i$-th $!$-box is duplicated $n_i$ times. For non-atomic expressions (i.e. expressions that are not just a single free variable or constant), it could become tricky to prove that any substitution of free variables yields valid concrete rewrite rule without ambiguity. However, if this can be done even in limited cases, such a language allows one to define a notion of induction, in the form of an inference rule on graphical identities.
\ctikzfig{induction_principal}

This rule says, ``If one can convert a single copy of a $!$-box to a single copy of another $!$-box, then one can convert $k$ copies of that $!$-box for any $k$.'' Note that the base case is trivially satisfied: if we kill both $!$-boxes, we are left with $G = G$.

\begin{example}\label{ex:inductive-pf}
  Using the \textbf{ind} inference rule to prove a new pattern graph identity. Take the following rules as given:
  \begin{center}
\beginpgfgraphicnamed{sg_black_dot_copy}
\InputIfFileExists{sg_black_dot_copy.tikz}{}{\input{./figures/sg_black_dot_copy.tikz}}
\endpgfgraphicnamed\qquad\qquad
\beginpgfgraphicnamed{bang_graph_rewrite_rule}
\InputIfFileExists{bang_graph_rewrite_rule.tikz}{}{\input{./figures/bang_graph_rewrite_rule.tikz}}
\endpgfgraphicnamed
  \end{center}
  
  Then, we will use the induction principal to prove:
  \ctikzfig{sg_many_copies}
  
  First, we derive the hypothesis for \textbf{ind} using rewriting.
  \ctikzfig{sg_induction_pf}
  
  The general rule is then constructed with an application of \textbf{ind}.
  \ctikzfig{sg_induction_pf2}
\end{example}

Note that this identity could not be produced using purely-equational means (e.g. using Knuth-Bendix completion), as none of the equations that we assumed contain the vertex $\blackdot$ inside of a $!$-box.

The proof above is quite similar in form to the types of inductive proofs checked by automated proof assistants. In principal, it would be straightforward to verify this proof in Quantomatic. Perhaps more advanced techniques like those employed by the proof-planning tool IsaPlanner~\cite{IsaPlanner} (e.g. rippling) could be adapted to string graphs to automatically \textit{search} for such proofs as well.

Of course, $!$-boxes are not the only way to compactly represent infinite sets of string graphs. We could also describe sets of graphs by introducing a ``meta'' rewrite system, where certain types of rewrite rules are treated as \textit{productions} in a graph grammar. In fact, this was the terminology originally used by the graph rewriting literature in the 1970s~\cite{Ehrig1973}. These rules would not need to respect inputs and outputs to a string graph, but some provision would need to be made when applying a meta-rule to a normal string graph rewrite rule to ensure that both the LHS and RHS are expanded in the same way. One could think of pattern graphs, as defined in this section, as something akin to regular languages, whereas sets of graphs described by graph grammars are richer (e.g. context-free or context-sensitive) languages. It is the hope that increasingly sophisticated graphical languages will lead to increasingly elegant and powerful graphical theories with applications in physics, linguistics, logic, and beyond.


	\bibliographystyle{akbib}
	\bibliography{bibfile}
\end{document}